%% file: ratell.tex
\documentclass[12pt]{amsart}
\usepackage{epic,eepic}
\usepackage{amscd,amssymb,xypic}
 \newlength{\baseunit}               
\newcommand{\bpf}{\noindent {\em Proof. }}
\newcommand{\epf}{\qed \vspace{+10pt}}

\newtheorem{tm}{Theorem}[section]
\newtheorem{pr}[tm]{Proposition}
\newtheorem{lm}[tm]{Lemma}

\newtheorem{co}[tm]{Corollary}

\newtheorem{defn}[tm]{Definition}
\newcommand{\proj}{\mathbb P}

\newcommand{\oh}{{\mathcal{O}}}

\newcommand{\ca}{{\mathcal{A}}}
\newcommand{\cb}{{\mathcal{B}}}
\newcommand{\cd}{{\mathcal{D}}}
\newcommand{\ce}{{\mathcal{E}}}
\newcommand{\cf}{{\mathcal{F}}}

\newcommand{\cll}{{\mathcal{L}}}
\newcommand{\cm}{{\mathcal{M}}}

\newcommand{\cu}{{\mathcal{U}}}
\newcommand{\cv}{{\mathcal{V}}}
\newcommand{\cw}{{\mathcal{W}}}
\newcommand{\cx}{{\mathcal{X}}}
\newcommand{\cy}{{\mathcal{Y}}}
\newcommand{\cz}{{\mathcal{Z}}}
\newcommand{\mbar}{{\overline{M}}}
\newcommand{\cmbar}{{\overline{\cm}}}
\newcommand{\mb}{{\overline{M}}_{0,\sum h_{m,n-1}+1}}
\newcommand{\cmb}{{\overline{\cm}}_{0,\sum h_{m,n-1}+1}}
\newcommand{\al}{\alpha}
\newcommand{\be}{\beta}
\newcommand{\ga}{\gamma}
\newcommand{\Ga}{\Gamma}

\newcommand{\de}{\delta}
\newcommand{\De}{\Delta}
\newcommand{\lcm}{\operatorname{lcm}}
\newcommand{\Sym}{\operatorname{Sym}}
\newcommand{\Aut}{\operatorname{Aut}}
\newcommand{\Bl}{\operatorname{Bl}}

\newcommand{\Spec}{\operatorname{Spec}}
\newcommand{\cehat}{\hat{\ce}}
\newcommand{\clehat}{\hat{{\mathcal{E}}'}}

\newcommand{\com}{\mathbb{C}}
\newcommand{\zed}{\mathbb{Z}}
\newcommand{\mbd}{\mbar_{g,m}(\proj^1,d)}
\newcommand{\mbgd}{\mbar_g(\proj^1,d)}
\newcommand{\Def}{\operatorname{Def}}
\newcommand{\Ext}{\operatorname{Ext}}
\newcommand{\Hom}{\operatorname{Hom}}
\newcommand{\hyperExt}{\operatorname{{\mathbb E}xt}}
\newcommand{\hyperHom}{\operatorname{{\mathbb H}om}}

\newcommand{\vh}{\vec{h}}
\newcommand{\vi}{\vec{i}}

\newcommand{\vep}{\vec{\epsilon}}

\newcommand{\cc}{{\mathcal{C}}}

\newcommand{\Pic}{\operatorname{Pic}}
\newcommand{\Hilb}{\operatorname{Hilb}}
\newcommand{\Ob}{\operatorname{Ob}}

\newcommand{\tmsoc}{in {\em The Moduli Space of Curves}, R. Dijkgraaf,
C. Faber and G. van der Geer eds., Birkhauser, 1995, pp }

\begin{document}
\pagestyle{plain}
\title{The Enumerative Geometry of Rational and Elliptic Curves in Projective 
Space} 
\author{Ravi Vakil}
\date{\today}
\begin{abstract}
We study the geometry of varieties parametrizing degree $d$ rational
and elliptic curves in $\proj^n$ intersecting fixed general linear
spaces and tangent to a fixed hyperplane $H$ with fixed multiplicities
along fixed general linear subspaces of $H$.  As an application, we
derive recursive formulas for the number of such curves when the
number is finite.  These recursive formulas require as ``seed data''
only one input: there is one line in $\proj^1$ through two points.
These numbers can be seen as top intersection products of various
cycles on the Hilbert scheme of degree $d$ rational or elliptic curves
in $\proj^n$, or on certain components of $\mbar_0(\proj^n,d)$ or
$\mbar_1(\proj^n,d)$, and as such give information about the Chow ring
(and hence the topology) of these objects.  The formula can also be
interpreted as an equality in the Chow ring (not necessarily at the
top level) of the appropriate Hilbert scheme or space of stable maps.
In particular, this gives an algorithm for counting rational and
elliptic curves in $\proj^n$ intersecting various fixed general linear
spaces.  (The genus 0 numbers were found earlier by Kontsevich-Manin,
and the genus 1 numbers were found for $n=2$ by Ran and
Caporaso-Harris, and independently by Getzler for $n=3$.)
\end{abstract}
\maketitle

\tableofcontents

\section{Introduction}
\label{intro}

In this article, we study the geometry of varieties (over $\com$)
parametrizing degree $d$ rational and elliptic curves in $\proj^n$
intersecting fixed general linear spaces and tangent to a fixed hyperplane
$H$ with fixed multiplicities along fixed general linear subspaces of $H$.
As an application, we derive recursive formulas (Theorem \ref{rrecursiveX2}
and Theorem 
\ref{erecursiveW}) for the number of such curves when the number is finite.
As with M. Kontsevich's marvelous formula of [KM], these recursive formulas
require as ``seed data'' only one input: there is one line in $\proj^1$
through two points.  These numbers can be seen as top intersection products
of various cycles on the Hilbert scheme of degree $d$ rational or elliptic
curves in $\proj^n$, or on $\mbar_0(\proj^n,d)$ or certain components of
$\mbar_1(\proj^n,d)$, and as such give information about the Chow ring (and
hence the topology) of these objects.  The formula can also be interpreted
as an equality in the Chow ring (not necessarily at the top level) of the
appropriate Hilbert scheme or space of stable maps.

\begin{figure}
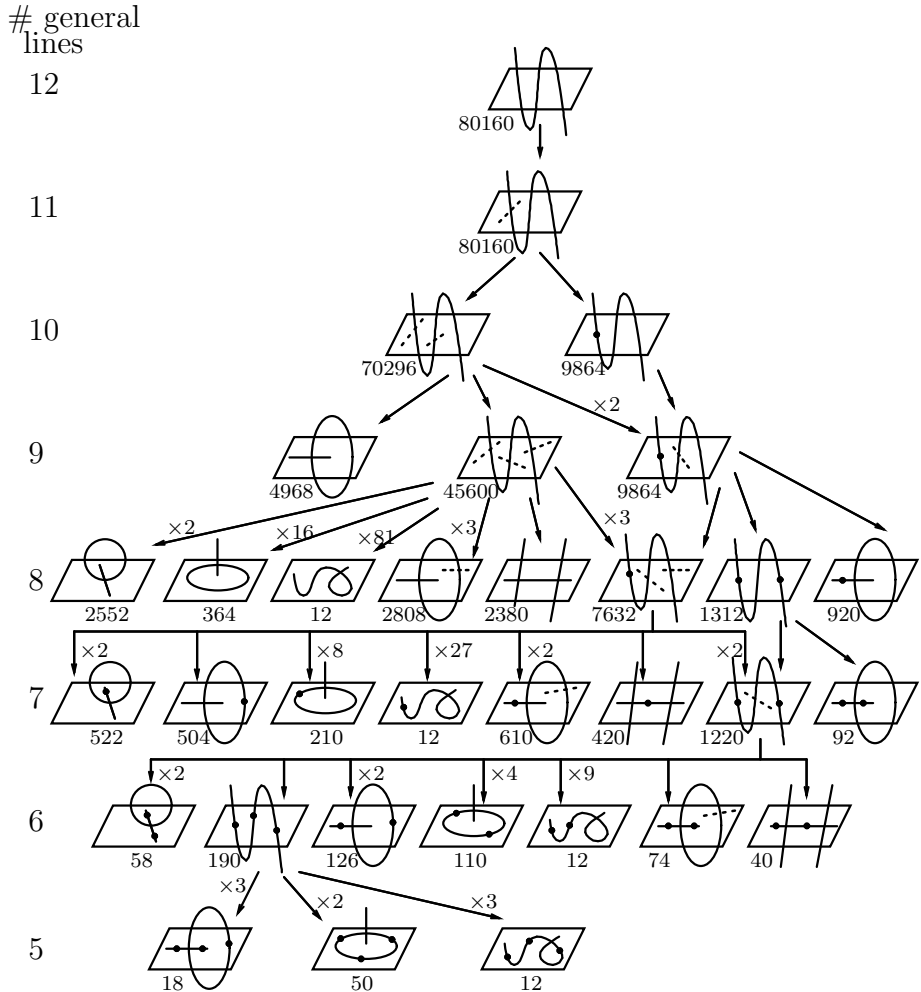

\begin{center}

	   \setlength{\unitlength}{.1\baseunit}
	    \input rcubics.tex 
\end{center}
\caption{Counting $80,160$ cubics in $\proj^3$ through 12 general lines}
\label{rcubics}
\end{figure}

The recursive formulas are convenient to program or use by hand, and provide
quick diagramatic enumerative calculations.  For example, Figure
\ref{rcubics} shows a diagramatic enumeration of the 80,160 twisted cubics
in $\proj^3$ intersecting 12 fixed general lines.  The methods used
are surprisingly elementary.  Little is assumed about Kontsevich's
space of stable maps.  The low genus of the curves under consideration
makes dimension counting straightforward.

Interest in such classical enumerative problems has been revived by recent
ideas from quantum field theory leading to the definition of quantum
cohomology and Gromov-Witten invariants, and by the subsequent discovery by
Kontsevich in 1993 of an elegant recursion solving the problem in genus 0
when no tangencies are involved (cf. [KM]; another proof, using different
techniques, was given independently by Y. Ruan and G. Tian in [RT]).  The
enumerative results here are a generalization of [KM], and the methods seem
more likely to generalize further than those of [KM].  In particular, such
ideas could apply to certain nonconvex rational spaces, in the same way
that the ideas of Caporaso and Harris in $\proj^2$ extend to rational ruled
surfaces (cf. [CH3] and [V1]).

In Section \ref{elliptic}, the same ideas are applied to the genus 1
case.  One of the motivations for this study was to gain more
information about higher genus Gromov-Witten invariants.  The
enumerative results for elliptic curves in $\proj^3$ have been
independently derived by E. Getzler (cf. [G3]), without tangency
conditions, by determining the genus 1 Gromov-Witten invariants of
$\proj^3$ and relating them to the enumerative problem.  (T. Graber
and R. Pandharipande ([GrP]) have also proposed a programme to
determine these numbers for all $\proj^n$ as well.)  There is some
hope of unifying the two methods, and getting some information about
Gromov-Witten invariants through degeneration methods.

There is also some hope that this method will apply to curves of higher genus.
The enumerative geometry of curves (of any genus) in the plane as described
in [CH3] can be seen as a variant of this perspective.  In [V1], the
corresponding problem (for curves of any genus in any divisor class) on any
rational ruled surface is solved by the same method.  In Subsection \ref{ehighgenus}, we will briefly discuss
the possibilities and potential obstructions to generalization to curves of
higher genus.  The results should also carry over to other highly symmetric
rational varieties, especially Flag manifolds and towers of
$\proj^1$-bundles, in the same way as the enumerative geometry of curves in
$\proj^2$ as described in [CH3] was generalized to rational ruled surfaces
in [V1].

This article contains the majority of the author's 1997 Harvard
Ph.D. thesis, and was partially supported by an NSERC 1967 Fellowship and a
Sloan Fellowship.  Tha author is extremely grateful to his advisor, Joe
Harris, and also to Dan Abramovich, Michael Thaddeus, Tony Pantev, Angelo
Vistoli, Rahul Pandharipande, Michael Roth, Lucia Caporaso, and Ezra
Getzler for useful discussions.

\subsection{Notation and basic results.}

The base field will always be assumed to be $\com$.  If $S$ is a set, define
$\Aut(S)$ to be the automorphisms (or symmetries) of $S$.

The main results of this article will be about varieties, but it will be
convenient to occasionally use the language of algebraic stacks (in the
sense of Deligne and Mumford, cf. [DM]).  Stacks have several advantages:
calculating Zariski-tangent spaces to moduli stacks is simpler than
calculating tangent spaces to moduli spaces, and moduli stacks (and
morphisms between them) are smooth ``more often'' than the corresponding moduli
spaces (and morphisms between them).  In general, for all definitions of
varieties made, there is a corresponding definition in the language of
stacks, and the corresponding stack will be written in a calligraphic font.
For example, $M_g$ is the moduli space of smooth genus $g$ curves, and
$\cm_g$ is the corresponding moduli stack.  Stacks are invoked as rarely as
possible, and the reader unfamiliar with stacks should have no problem
following the arguments.  An introduction to the theory of algebraic stacks
is given in the appendix to [Vi1].

\subsubsection{The moduli space of stable maps}
\label{itmsosm}
For the convenience of the reader we recall certain facts about moduli
spaces (and stacks) of stable maps to $\proj^n$, without proofs.  
A {\em stable marked curve}, denoted $(C, \{ p_1, \dots, p_m \} )$, of
genus $g$ with $m$ marked points $p_1$, \dots, $p_m$ is a connected nodal
complete marked curve with finite automorphism group; the marked points are
required to be smooth points of $C$.  The set of {\em special
points} on (the normalization of) a component of $C$ is the set of marked
points union the set of branches of nodes.  A marked curve is stable if
each rational component has at least three special points and each elliptic
component has at least one special point.

A degree $d$ {\em stable map} $(C, \{ p_1, \dots, p_m \}, \pi)$ to $\proj^n$
consists of a connected nodal complete marked curve $(C, \{ p_i \} )$
and a morphism $\pi:  C \rightarrow \proj^n$ such that $\pi_*[C] = d[L]$
(where $[L]$ is the class of a line in the first Chow group $A_1(\proj^n)$
of $\proj^n$, or equivalently in $H_2(\proj^n, \zed)$), such that the map $\pi$
has finite automorphism group.  This last condition is equivalent to
requiring each collapsed rational component to have at least three special
points and each contracted elliptic component to have at least one special
point.  

There is a coarse projective moduli space $\mbar_{g,m}(\proj^n,d)$ for such
stable maps with $p_a(C) = g$.  There is an algebraic stack
$\cmbar_{g,m}(\proj^1,d)$ that is a fine moduli space for stable maps.  The
stack $\cmbar_{g,m+1}(\proj^1,d)$ can be considered as the universal curve
over $\cmbar_{g,m}(\proj^1,d)$.
There is an open subvariety $M_{g,m}(\proj^n,d)$ of
$\mbar_{g,m}(\proj^n,d)$ that is a coarse moduli space of stable maps from
smooth curves, and an open substack $\cm_{g,m}(\proj^n,d)$ of
$\cmbar_{g,m}(\proj^n,d)$ that is a fine moduli space of stable maps from
smooth curves.  
The stack $\cmbar_{0,m}(\proj^n,d)$ is smooth of dimension $(n+1)(d+1)+m
-4$.

The versal deformation space to the stable map $(C, \{ p_i \}, \pi)$ in
$\cmbar_{g,m}(\proj^n,d)$ is obtained from the complex 
$$
\underline{\Omega}_\pi = \left( \pi^* \Omega_{\proj^n} \rightarrow
\Omega_C(p_1 + \dots + p_m)\right).
$$
(The versal deformation space depends only on the image of
$\underline{\Omega}_{\pi}$ in the derived category of coherent sheaves on
$C$.)  The vector space $\hyperHom (\underline{\Omega}_\pi, \oh_C)$
parametrizes infinitesimal automorphisms of the map $\pi$; as $(C, \{ p_i
\}, \pi)$ is a stable map, $\hyperHom (\underline{\Omega}_\pi, \oh_C) = 0$.
The space of 
infinitesimal deformations to the map $(C, \{ p_i \}, \pi)$ (i.e. the
Zariski tangent space to $\cmbar_{g,m}(\proj^n,d)$ at the point
representing this stable map), denoted $\Def_{(C, \{ p_i \}, \pi )}$, is
given by $\hyperExt^1(\underline{\Omega}_\pi, \oh_C)$ and the obstruction
space, denoted $\Ob_{(C, \{ p_i \},\pi)}$, is given by
$\hyperExt^2(\underline{\Omega}_\pi, \oh_C)$.

By applying the functor $\Hom( \cdot, \oh_C)$ to the exact sequence of complexes
$$
0 \rightarrow \Omega_C(p_1 + \dots + p_m)[-1] \rightarrow
\underline{\Omega}_\pi \rightarrow \pi^* \Omega_{\proj^n} \rightarrow 0
$$
we obtain the long exact sequence
\begin{eqnarray}
\nonumber
0 \rightarrow \hyperHom( \underline{\Omega}_\pi, \oh_C)
&\rightarrow& \Hom( \Omega_C(p_1 + \dots + p_m), \oh_C)  \rightarrow H^0(C,\pi^*
T_{\proj^n}) \\
\nonumber
\rightarrow \hyperExt^1( \underline{\Omega}_\pi, \oh_C)
&\rightarrow& \Ext^1( \Omega_C(p_1 + \dots + p_m), \oh_C)  \rightarrow H^0(C,\pi^*
T_{\proj^n}) \\
\nonumber
\rightarrow \hyperExt^2( \underline{\Omega}_\pi, \oh_C)
&\rightarrow& 0.
\end{eqnarray}

By the identifications given in the previous paragraph, and using
$\Hom(\Omega_C(p_1 + \dots + p_m), \oh_C) = \Aut(C, \{ p_i \})$ and
$\Ext^1(\Omega_C(p_1 + \dots + p_m), \oh_C) = \Def(C, \{ p_i \})$,
this long exact sequence can be rewritten as
\begin{eqnarray}
\nonumber
0 &\longrightarrow& \Aut (C, \{ p_i \}) \longrightarrow H^0(C,\pi^*
T_{\proj^n}) \\
\nonumber
 \longrightarrow  \Def (C, \{ p_i \}, \pi ) & \longrightarrow &
\Def (C, \{ p_i \}) \longrightarrow H^1(C,\pi^*
T_{\proj^n} ) \\
 \longrightarrow  \Ob (C, \{ p_i \}, \pi ) &\longrightarrow & 0.
\nonumber
\end{eqnarray}

The construction of the versal deformation space from
$\underline{\Omega}_\pi$ is discussed in [R3], [Vi2], and [LT2].  All other
facts described here appear in the comprehensive introduction [FP]. 

\subsection{Divisors on subvarieties of $\mbd$}
\label{ikey}
The results proved here will be invoked repeatedly in Sections
\ref{rational} and \ref{elliptic}.

\subsubsection{A property of stable maps from curves to curves}

We will make repeated use of a special property of maps from curves to
curves.  Let $(C, \{ p_i\}, \pi)$ be a stable map
of a complete marked curve to $\proj^1$.  The scheme $\pi^{-1}(p)$ consists
of reduced unmarked points for almost all $p$.  Let $A = A_1 \coprod \dots
\coprod 
A_l$ be the union of the connected components of fibers of $\pi$ that are
{\em not} reduced unmarked points.  Call the $A_j$ the {\em special loci}
of the map $\pi$.  Then each special locus $A_j$ is either a ramification
of $\pi$, a labeled point (that may also be a 
ramification), or a union of contracted components (possibly containing
labeled points, and possibly attached to other components at their
ramification points).

Let $\Def_{(C,\{ p_i \}, \pi)}$ be the versal deformation space to the
stable map
$(C, \{ p_i \} , \pi)$.  Let $\Def_{A_j}$ be the versal deformation space of
the map 
$$
(C \setminus \cup_{i \neq j} A_i, \{ p_i \} \cap A_j, \pi );
$$
the space $\Def_{A_j}$ parametrizes formal deformations of the map $\pi$ that
are trivial away from $A_j$.

\begin{pr}
The versal deformation space $\Def_{(C,\{ p_i \}, \pi)}$ is naturally
$\prod_j \Def_{A_j}$.
\label{ilocal}
\end{pr}

An informal argument in the analytic category is instructive.  Let $U_j
\subset C$ (for $1 \leq j \leq l$) be an open (analytic) neighborhood of the special
locus $A_j$ whose closure does not intersect the other special loci
$\cup_{i \neq j} A_i$.  Let $U$ be an open subset of $C$ whose closure does
not intersect the special loci $\cup A_i$, and such that 
$$
U \cup \left( \bigcup_{j=1}^l U_j \right) = C.
$$
Then ``small'' deformations of $(C, \{ p_i \}, \pi)$ are trivial on $U$, and thus
the deformations of $\pi$ on $U_1$, \dots, $U_l$ are mutually independent:
$$
\Def_{C,\{ p_i \}, \pi} = \prod_j \Def_{(U_j, \{ p_i \} \cap U_j, \pi)}.
$$

This argument can be carried out algebraically on the level of
formal schemes, which will give a rigorous proof.

{\noindent {\em Proof of the proposition. }}
The versal deformation to the stable map $(C, \{ p_i \}, \pi)$ is
constructed using the complex 
$$
\underline{\Omega}_{(C, \{ p_i \}, \pi)} = \left( \pi^* \Omega_{\proj^1} \rightarrow
\Omega_C \left( \sum p_i \right) \right)
$$
(see Subsubsection \ref{itmsosm}).
The complex $\underline{\Omega}_{(C, \{ p_i \}, \pi)}$ is exact on $C
\setminus \cup  A_j$, so it splits in the derived category as a direct sum
of objects 
$$
\underline{\Omega}_{  (C, \{ p_i \}, \pi)} = \oplus \underline{\Omega}_{A_j}
$$
where the cohomology sheaves of $\underline{\Omega}_{A_j}$ are supported on $A_j$.  
Then the entire construction of the versal deformation space naturally
factors through this direct sum, and the result follows.
\epf

Thus to understand the versal deformation space of a map to $\proj^1$ we
need only understand the versal deformations of the special loci
$\Def_{A_j}$.  

Let $(C, \{ p_i \}, \pi)$ (resp. $(C', \{ p'_i \}, \pi')$) be a stable map
to $\proj^1$ and let $A_1$ (resp. $A'_1$) be one of its special loci
consisting of a connected union of components contracted by $\pi$
(resp. $\pi'$).  Let $\tilde{C}$ be the closure of $C
\setminus A_1$ in $C'$, and $\tilde{C'}$ the closure of $C' \setminus A'_1$
in $C'$.  Consider $A_1$
(resp. $A_1'$) as a marked curve, where the markings are the intersection
of $\tilde{C}$ with $A_1$ (resp. $\tilde{C'}$ with $A'_1$) and the labeled
points of $C$ (resp. $C'$).  Assume that $A_1$ and $A_1'$ are isomorphic as
marked curves, and that the ramification orders of the points
of $\tilde{C} \cap A_1$  on $\tilde{C}$ are the same as those of
the corresponding points of $\tilde{C'} \cap A_1'$ on $\tilde{C'}$.

\begin{pr}
\label{ilocal2}
There is an isomorphism of versal deformation spaces $\Def_{A_1} \cong
\Def_{A_1'}$.
\end{pr}

As before, an analytic perspective is instructive.  The hypotheses of the
theorem imply that there is an analytic neighborhood  $U_{\text{an}}$ of $A_1$ that is
isomorphic to an analytic neighborhood $U'_{\text{an}}$ of $A_1'$, and
$\pi|_{U_{\text{an}}} = 
\pi'|_{U'_{\text{an}}}$.

\bpf
The hypotheses of the proposition imply that a formal neighborhood $U$ of
$A_1$ is isomorphic to a formal neighborhood $U'$ of $A_1'$, and $\pi|_U$
agrees with $\pi'|_{U'}$ (via this isomorphism).  As the cohomology sheaves
of $\underline{\Omega}_{A_1}$, $\underline{\Omega}_{A'_1}$ are supported in
$U$ and $U'$ respectively, the entire construction of $\Def_{A_1}$ and
$\Def_{A'_1}$ depends only on $(U, \pi|_U)$ and $(U', \pi'|_{U'})$.
\epf

\subsubsection{Subvarieties of $\mbar_{g,m}(\proj^1,d)$}

Fix a positive integer $d$, and let $\vec{h} = (h_1,h_2,\dots)$ represent a partition of
$d$ with $h_1$ 1's, $h_2$ 2's, etc., so $\sum_m m h_m = d$.  Fix a point
$z$ on $\proj^1$.  Let $X = X^{d,g}(\vh)$ be the closure in $\mbar_{g,\sum h_m + 1}(\proj^1,d)$ of
points representing stable maps $(C, \{ p^j_m \}, q, \pi)$ where $C$ is an
irreducible curve with $\sum h_m + 1$ (distinct) marked points $\{
p^j_m \}_{1 \leq j \leq h_m}$ and $q$, and $\pi^*(z) =
\sum_{j,m} m p^j_m$.  (In short, we have marked all the pre-images of $z$
and one other point.  The map $\pi$ has $h_m$ ramifications of order $m$
above the point $z$.)  Then by Riemann-Hurwitz,
\begin{equation}
\label{idimX}
\dim X = d+2g-1+\sum h_m,
\end{equation}
which is $d+2g-1$ plus the number of pre-images of $z$.  (We are implicitly
invoking the Riemann existence theorem here.)  Notice that for the map
corresponding to a general point in $X$, each special locus $A_j$ is either
a marked ramification above the point $z$, a simple unmarked ramification,
or the point $q$ (at which $\pi$ is smooth).  In these three cases, the
formal deformation space $\Def A_j$ is 0, $\Spec \com [[t]]$, and $\Spec
\com [[t]]$ respectively.

There is a corresponding stack $\cx = \cx^{d,g}(\vh) \subset \cmbar_{g,\sum h_m +
1}(\proj^1,d)$ as well.

Let $D$ be the
divisor $\{ \pi(q) = z \}$.  There are three natural questions to ask.
\begin{enumerate}
\item What are the components of the divisor $D$?
\item With what multiplicity do they appear?
\item What is the local structure of $X$ near these components?
\end{enumerate}

We will partially answer these three questions.

\begin{figure}
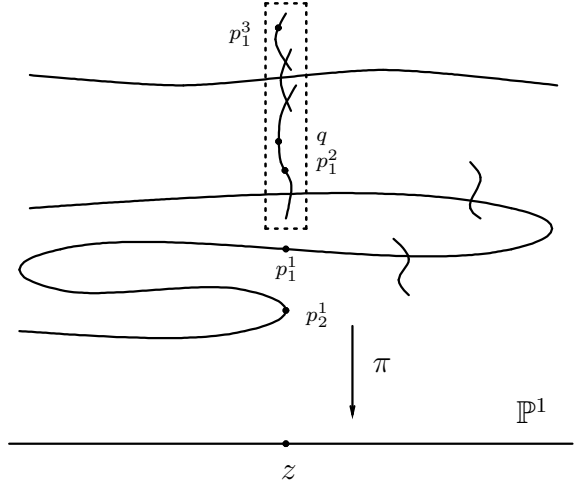

\begin{center}

	   \setlength{\unitlength}{.1\baseunit}
	    \input iegC.tex 
\end{center}
\caption{The map $(C, \{p^j_m \}, q, \pi \} ) \in Y$}
\label{iegC}
\end{figure}

Fix a component $Y$ of the divisor $D$ and a map $(C, \{ p^j_m
\}, q, \pi)$ corresponding to the general element of $Y$.
Notice that $\pi$ collapses a component of $C$ to $z$, as otherwise
$\pi^{-1}(z)$ is a union of points, and 
$$
\deg \pi^* (z) \geq \sum_{j,m} \deg \pi^* (z) |_{p^j_m} + \deg \pi^* (z) |_q \geq
\sum_{j,m} m + 1 > d.
$$
Let $C(0)$ be the connected component of $\pi^{-1} (z)$ containing $q$, and
let $\tilde{C}$ be the closure of $C \setminus C(0)$ in $C$ (see
Figure \ref{iegC}; $C(0)$ are those curves contained in the dotted rectangle, and
$\tilde{C}$ is the rest of $C$).

Let $h_m(0)$ be the number of $\{ p^j_m \}_j$ in $C(0)$, and $\tilde{h}_m =
h_m - h_m(0)$ the number in $\tilde{C}$.  
Let $\{ p^j_m(0) \}_{m,1 \leq j \leq h_m(0)}$ and 
$\{ \tilde{p}^j_m \}_{m,1 \leq j \leq \tilde{h}_m}$ be the partition of $\{
p^j_m \}_{1 \leq j \leq h_m}$ into those marked points lying on $C(0)$ and
those lying on $\tilde{C}$.
Let $s$ be the number of
intersections of $C(0)$ and $\tilde{C}$, and label these points $r^1$,
\dots, $r^s$.  Thus $g = p_a(C(0)) + p_a(\tilde{C})
+ s - 1$.   Let $m^k$ be the multiplicity of $(\pi|_{\tilde{C}})^*(z)$ at
$r^k$.  The data $(m^1, \dots, m^s)$ must be constant for any choice of $(C,
\{p^j_m \}, q, \pi)$ in an open subset of $Y$.  

\begin{pr}
The stable map $(\tilde{C}, \{ p^j_m(0) \}, \{ r^k \}, \pi)$ has no
collapsed components, and only simple ramification away from $\pi^{-1}(z)$.
The curve $\tilde{C}$ is smooth.
\end{pr}

The map $(\tilde{C}, \{ p^j_m(0) \}, \{ r^k \}, \pi)$ will
turn out to correspond to a general element in $X^{d,g'}(\vec{h'})$ for some
$\vec{h'}$, $g'$.  Note that $\tilde{C}$ may be reducible.

\bpf
Let $A_1$, \dots, $A_l$ be the special loci of $\pi$, and say $q \in A_1$.

The map $(C, \{ p^j_m \}, q, \pi)$ lies in $X$ and hence can be deformed to
a curve where each special locus is either a marked ramification above $z$,
a simple unmarked ramification, or an unramified marked point.  If $A_k$
($k>1$) is not one of these three forms then by Proposition \ref{ilocal}
there is a deformation of the map  $(C, \{ p^j_m \}, q, \pi)$ preserving
$\pi$ at $A_i$ ($i \neq k$) but changing $A_k$ into a combination of
special loci of these three forms.  Such a deformation (in which $A_1$
is preserved and thus still smoothable) is actually a
deformation in the divisor $D = \{ \pi(q) = z \}$, contradicting the
generality of $(C, \{ p^j_m \}, q, \pi)$ in $Y$.
\epf

The map $(\tilde{C}, \{ p^j_m(0) \}, \{ r^k \}, \pi)$ must lie in $X^{d, p_a(\tilde{C})}(
\vec{h'})$ where $\vec{h'}$ is the partition corresponding
to $(\pi |_{\tilde{C}})^*(z)$.  By Riemann-Hurwitz, $\tilde{C}$ moves in a family of
dimension at most
\begin{equation*}
d + 2 p_a(\tilde{C}) - 2 + \left( \sum \tilde{h}_m + s \right) ,
\end{equation*}
and the curve $C(0)$ (as a nodal curve with marked points $\{ r^k \}_{1
\leq k \leq s}$, $\{ p^j_m(0)
\}_{m, 1 \leq j \leq h_m(0)}$, and $q$) moves in a
family of dimension at most
$$
3 p_a (C(0)) - 3 + \sum h_m(0) + s +1
$$
so $Y$ is contained in a family of dimension
\begin{eqnarray}
\nonumber
d + 2 p_a(\tilde{C}) - 2 + \sum \tilde{h}_m + s &+& 3 p_a(C(0)) - 3 +
\sum h_m(0) + s +1\\
\nonumber
&=& d + 2g-1 + \sum h_m - 1 + p_a(C(0))\\
&=& \dim X - 1 + p_a(C(0)) 
\label{inaive}
\end{eqnarray}
by (\ref{idimX}).

We can now determine all components $Y$ of $D$ satisfying
$p_a(C(0)) = 0$.  For each choice of a partition of $\{ p^j_m \} $ into $\{
p^j_m(0) \} \cup \{ \tilde{p}^j_m \}$ (inducing a partition of $h_m$ into
$h_m(0) + \tilde{h}_m$ for all $m$), a positive integer $s$, and $(m^1,
\dots, m^s)$ satisfying $\sum m^s + \sum m \tilde{h}^j_m = d$, there is a
variety (possibly reducible) of dimension $\dim X - 1$ that is the closure
in $\mbar_{g, \sum h_m + 1}(\proj^1,d)$ of points corresponding to maps
$$
(C(0) \cup \tilde{C}, \{ p^j_m \}, q, \pi)
$$
where 
\begin{enumerate}
\item[A1.] The curve $C(0)$ is isomorphic to $\proj^1$ and has labeled points 
$$
\{ p^j_m(0) \} \cup \{ r^k \} \cup \{ q \},
$$
and $\pi(C(0)) = z$.
\item[A2.] The curve $\tilde{C}$ is smooth of arithmetic genus $g-s+1$ with
marked points $\{ \tilde{p}^j_m \} \cup \{ r^k \}$.  The map $\pi$ is
degree $d$ on $\tilde{C}$, and
$$
( \pi|_{\tilde{C}})^*(z) = \sum m \tilde{p}^j_m + \sum m^k r^k.
$$
\item[A3.] The curve $C(0) \cup \tilde{C}$ is nodal, and the curves $C(0)$ and $\tilde{C}$ intersect
at the points $\{ r^k \}$.
\end{enumerate}
Let $U$ be the union of these varieties.

An irreducible component $Y$ of the divisor $D$ satisfying $p_a(C(0)) = 0$
has dimension $\dim X - 1$ and is a subvariety of $U$, which also has
dimension $\dim X - 1$.  Hence $Y$ must be a component of $U$ and the
stable map corresponding to a general point of $Y$ satisfies properties
A1--A3 above.   (We don't yet know that all such $Y$ are subsets of $X$, but
this will follow from Proposition \ref{imultg0} below.)

\begin{figure}
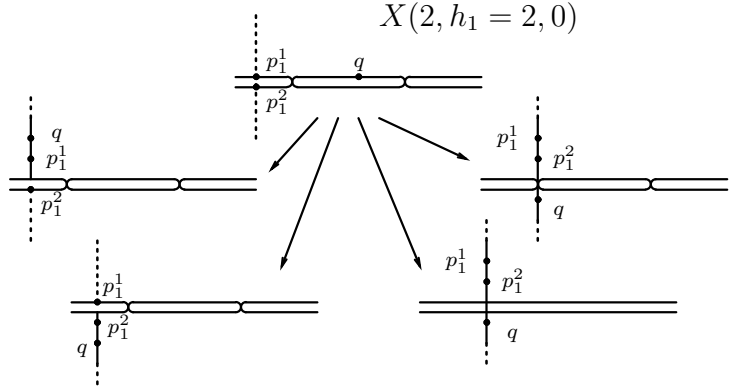

\begin{center}

	   \setlength{\unitlength}{.1\baseunit}
	    \input ideg0.tex 
\end{center}
\caption{The possible components of $D$ on $X^{2,0}(h_1=2)$}
\label{ideg0}
\end{figure}

For example, if $d=2$, $g=0$, and $h_1 = 2$, then there are four
components (see Figure \ref{ideg0}; $\pi^{-1}(z)$ is indicated by a dotted
line).  The components (from left to right) are a subset of the following.
\begin{enumerate}
\item
The curve $\tilde{C}$ is irreducible and maps with degree 2 to
$\proj^1$, ramifying over general points of $\proj^1$.  The marked points
$q$ and $p^1_1$ lie on $C(0)$, and $p^2_1$ lies on $\tilde{C}$.  The curve
$C(0)$ is attached to $\tilde{C}$ at the point 
$$
( \pi|_{\tilde{C}} )^{-1}(z) \setminus \{ p^2_1 \}.
$$
\item This case is the same as the previous one with $p^1_1$ and $p^2_1$
switched.
\item The curve $\tilde{C}$ is the disjoint union of two $\proj^1$'s, each
mapping to $\proj^1$ with degree 1.  Both
intersect $C(0)$, which contains all the marked points.  
\item The curve $\tilde{C}$ is irreducible and maps with degree 2 to
$\proj^1$, and one of its branch points is $z$.  All of the marked points
lie on $C(0)$. 
\end{enumerate}

Given a component $Y$ of $U$, we can determine the multiplicity of the
divisor $D = \{ \pi(q) = z \}$ along $Y$.  As this multiplicity will turn
out to be positive, we will have the corollary that, as sets, $U \subset D$.

For technical reasons, we use the language of stacks.

\begin{pr}
\label{imultg0}
Fix such a component $\cy$ with $p_a(C(0)) = 0$.  The multiplicity of $D$
along $\cy$ is $\prod_{k=1}^s m^k$.
\end{pr}

In particular, $Y$ is a subset of $X$.

\bpf
We make a series of reductions to simplify the proof.

{\em Step 1:  The deformations of $A_1$.}  Let $\Def_{(C, \{ p^j_m \}, q,
\pi), \cx}$ be the versal deformation space of $(C, \{ p^j_m \}, q,
\pi)$ in $\cx$; it is a subspace of $\Def_{(C, \{ p^j_m \}, q,
\pi)}$.  
We can compute the multiplicity on the versal deformation space of a stable
map $\Def_{(C, \{ p^j_m \}, q, \pi), \cx}$ corresponding to a general point in $\cy$.
If $A_1$ is the special locus of $\pi$ containing $q$, then the divisor
corresponding to $D$ on $\Def_{(C, \{ p^j_m \}, q, \pi), \cx}$ is the
pullback of a divisor $D_{\Def_{A_1,\cx}}$ on $\Def_{A_1,\cx}$.  We will study
$\Def_{A_1}$ by analyzing a simpler map.

{\em Step 2:  A simpler map.}
Consider the stable map 
$$
(C'(0) \cup \tilde{C}', \{ p^j_m(0) \}, q, \{r^k \}, \pi)
\in \cmbar_{0, \sum h_m(0)+1+s}(\proj^1,\sum m h_m(0))
$$
where:
\begin{enumerate}
\item[B1.] The marked curve $(C'(0),  \{ p^j_m(0) \}, q, \{r^k \})$ is
isomorphic to the marked curve $(C(0),  \{ p^j_m(0) \}, q, \{r^k \})$, and is collapsed to
$z$ by $\pi$.
\item[B2.] The stable map $(\tilde{C}', \pi)$ consists of $s$ rational curves
$C'(1)$, \dots, $C'(s)$ of degrees $m^1$, \dots, $m^s$ respectively, each
ramifying completely over $z$.
\item[B3.] The point of ramification of $C'(k)$ over $z$ is glued to $r^k$ on
$C'(0)$. 
\end{enumerate}
Let $A'_1$ be the special locus $C'(0)$ of this stable map.  By Proposition
\ref{ilocal2}, $\Def_{A_1} \cong \Def_{A'_1}$.  Thus without loss of
generality we may assume that $g=0$, $\vh = \vh(0)$, and  the point in $\cy$ is of the form
$$
(C = C(0) \cup \tilde{C}, \{ p^j_m \}, q, \{r^k \}, \pi)
\in \cmbar_{0, \sum h_m+1+s}(\proj^1,d)
$$
with properties B1--B3.  Note that in this case $d = \sum m h_m = \sum m^k r^k$.

{\em Step 3: Fixing the other special loci.} 
Since we are interested
in $\Def_{A_1}$, it will be to our advantage to hold the other special
loci constant.  Fix a point $y \neq z$ on the target $\proj^1$.  Let $\cx'$ be the
closed substack that is the stack-theoretic closure (in $\cmbar_{0,\sum
h_m + 1 + s}(\proj^1,d)$) of the points representing stable maps
$$
(C, \{ p^j_m \}, q, \{ y^k \}, \pi )
$$
where $\pi^*(z) = \sum m p^j_m(0)$ and $\pi^*(y) = \sum m^k y^k$. Let
$\cy'$ be the closure of points representing maps
$$
(C = C(0) \cup C(1) \cup \dots \cup C(s), \{ p^j_m \}, q, \{ y^k \},\pi )
$$
where:
\begin{enumerate}
\item[(i)] The curve $C(0)$ is rational, contains $\{ p^j_m \}$, intersects
all of the other $C(k)$, and is collapsed to $z$ by $\pi$.
\item[(ii)] The curve $C(k)$ ($k>0$) is rational.  The map $\pi$ is degree $m^k$
on $C(k)$, and $C(k)$ ramifies totally above $z$ (where it intersects
$C(0)$) and $y$ (where it is labeled $y^k$).
\end{enumerate}
If $(C, \{ p^j_m \}, q, \{y^k \}, \pi)$ is a map corresponding to a
general point in $\cy'$,
$$
\Def_{(C, \{ p^j_m \}, q, \{y^k \}, \pi), \cx'} = \Def_{A_1, \cx}:
$$
the only deformations of such a
map preserving the ramifications above $y$ are deformations of $A_1$.

{\em Step 4: Fixing the marked curve.}  We next reduce to the case where
$(C, \{ p^j_m \}, q, \{y^k \})$ is a fixed stable marked curve.  There is a
morphism of stacks $\al: \cx \rightarrow \cmbar_{0, \sum h_m + 1 + s}$ that
sends each map to stable model of the underlying pointed nodal curve.
Given any smooth marked curve $(C(0), \{ p^j_m \}, q, \{ y^k \} )$ in
$\cmbar_{0, \sum h_m + 1 + s}$, the stable map $(C, \{ p^j_m \},
q, \{ y^k
\},\pi )$ defined in Step 3 (where $C$ is a union of irreducible curves $C(0) \cup \dots \cup C(s)$)
corresponds to a point in $\al^{-1}(C(0), \{ p^j_m \}, q, \{ y^k \} )$.  (The
stable curve $\al(C, \{ p^j_m \}, q, \{ y^k
\},\pi )$ is constructed by forgetting $\pi$ and contracting the rational
tails $C(1)$, \dots, $C(s)$.)
Hence $\al|_{\cy'}$ is surjective.  Let $\cf_\al$ be a general fiber of
$\al$.  By Sard's theorem, $\al |_{\cy'}$ is regular in a Zariski-open
subset of $\cy'$, so $[\cy'] \cap [\cf_\al] = [\cy' \cap \cf_\al]$ in the
Chow group of $[\cx']$.

In order to determine the multiplicity of $D|_{\cx'}$ along $\cy'$, it suffices to
determine the multiplicity of the Cartier divisor $D|_{\cf_\al}$ along
$\cy' \cap \cf_\al$ (in the Chow group of $\cf_\al$).  ({\em Proof:} As $D$
is a Cartier divisor, $[D|_{\cf_\al}] = [D|_{\cx'}] \cdot [\cf_\al]$.  Thus
if $[D|_{\cx'}] = m [ \cy']$ in $A^1 \cx'$ then, intersecting with
$[\cf_\al]$, $[D|_{\cf_\al}] = [D|_{\cx'}] \cdot [\cf_\al] = m [\cy'] \cdot
[\cf_\al] = m [\cy' \cap \cf_\al ]$ in $A^1 \cf_\al$.)

With this in mind, fix a general $(C, \{ p^j_m \}, q, \{ y^k \} )$ in 
$\cmbar_{0, \sum h_m + 1 + s}$ and let $\cx''_o$ be the points of
$\cmbar_{0, \sum h_m + 1 + s}(\proj^1,d)$ 
representing stable maps $(C, \{ p^j_m \}, q, \{ y^k \}, \pi )$
where $\pi^*(z)= \sum m p^j_m$ and $\pi^*(y) = \sum m^k y^k$.  
Let $\cx'' = \cx' \cap \cf_\al$ be the closure  of $\cx''_o$, 
let $X''_o$, $X''$ be the corresponding
varieties, and define $\cy'' = \cy \cap \cf_\al$ and $Y''$ similarly.

{\em Step 5:  The variety $X''$ is actually $\proj^1$!}
Let $f$ and $g$ be sections of $\oh_C(d)$ with 
$$
(f=0) = \sum m p^j_m, \quad (g=0) = \sum m^k y^k;
$$
the maps in $\cx''_o$ are those of the form $[\be f, \ga g ]$ where
$$
[\be, \ga] \in \proj^1 \setminus \{ [0,1], [1,0] \}
$$
where $z = [0,1]$ and $y = [1,0]$.
The variety $X''$ is proper, so the normalization of the variety $X''$ is $\proj^1$.

The curve $X''$ has a rational map to $\proj^1$ given by
$$
(C, \{ p^j_m \}, q, \{ y^k \}, \pi ) \rightarrow \pi(q)
$$
and this map is an isomorphism from $X''_o$ to $\proj^1 \setminus \{ [0,1],
[1,0] \}$, so it must be an isomorphism from $X''$ to $\proj^1$.

{\em Step 6:  Calculating the multiplicity.}
Let $z'$ be a general point of the target $\proj^1$.  Then the divisor $\{
\pi(q) = z' 
\}$ is linearly equivalent to $D|_{X''} = \{ \pi(q) = z \}|_{X''}$, and
is $\oh_{X''}(1)$.

Thus, {\em as varieties}, $D|_{X''} = [1,0] = Y''$.  But the limit map has
automorphism group
$$
\zed / {m^1 \zed} \oplus \dots \oplus \zed / m^s \zed
$$
(as $\Aut(C(k), \pi|_{C(k)}) = m^k$) so {\em as stacks} $[D |_{\cx''} ] =
\left( \prod m^k \right) [\cy'']$.  Therefore $[D] = \prod m^k 
[\cy]$.
\epf

The above argument can be refined to determine the local structure of
$\cx$ near $\cy$ (and thus $X$ near $Y$ if the map corresponding to a
general point of $Y$ has no automorphisms):

\begin{co}
\label{ilocalst}
Let $\cy$ be the same component as in Proposition \ref{imultg0}.  
an \'{e}tale neighborhood of a general point of $\cy$, the
stack $\cx$ is isomorphic to
$$
\Spec \com [[a, b_1,\dots,b_s,c_1,\dots,c_{\dim X -1 }]] / (a = b_1^{m^1} = \dots =
b_s^{m^s})
$$
with $D$ given by $\{ a = 0 \}$, and $\cy$ given set-theoretically by the
same equation.
\end{co}
In particular, if $\gcd(m^i,m^j)>1$ for some $i$ and $j$, $\cx$ fails to be
unibranch at a general point of $\cy$.

\bpf
In the proof of Proposition \ref{imultg0} above, at the end of Step 4 we had reduced to considering a fixed
marked curve 
$(C, \{ p^j_m \}, q, \{ y^k \} )$ in 
$\cmbar_{0, \sum h_m + 1 + s}$ and maps $\pi$ from this marked curve to
$\proj^1$ where $\pi^*(z) = \sum m p^j_m$ and $\pi^*(y) = \sum m^k y^k$.
These maps are parametrized by the stack $\cx''$.

{\em Step 5${}'$: Rigidifying the moduli problem.} It will be more
convenient to work with varieties than stacks, so we rigidify the moduli
problem to eliminate nontrivial automorphisms.  Fix a point $x
\in \proj^1$ distinct from $y$ and $z$.  We will mark the $d$ pre-images of
$x$ with the labels $\{ x^1, \dots, x^d \}$.  Let $\cx''_x$ be the moduli
stack parametrizing maps $(C, \{p^j_m \}, q, \{ y^k \}, \{ x^l \}, \pi)$
where
\begin{enumerate}
\item[C1.]  The pointed curve $(C, \{p^j_m \}, q, \{ y^k \})$ is
fixed. 
\item[C2.]  $(C, \{p^j_m \}, q, \{ y^k \}, \pi)$ is a stable map,
\item[C3.] $\pi(x^i) = x$ for all $i$, and $x^i \neq x^j$ for $i \neq j$.
\end{enumerate}
The moduli stack $\cx''_x$
is actually a variety (as none of the maps
para\-metrized by this stack have automorphisms), and over an open
neighborhood of $\cy''$ in $\cx''$ the natural morphism
$$
\eta:  \cx''_x \rightarrow \cx''
$$
is an \'{e}tale (degree $d!$) morphism of proper stacks at a point of
$\cy''$.  (The variety $\cx''_x$ is an atlas for the stack $\cx''$.)
Define the Weil divisor $\cy''_x$ on $\cx''_x$ similarly; it is a union of
points.

We can now reprove Proposition \ref{imultg0} using the variety $\cx''_x$:  if
$\cy''$ is the point in $\cx''$ corresponding to the point in $\cy$
(i.e. the curve $C(0)$ with $l$ rational tails $C(1)$, \dots, $C(l)$
ramifying completely over $y$ and $z$) then $\eta^{-1}(\cy'')$ is
set-theoretically
$$
{\binom d {m^1, \dots, m^s} } \prod_{k=1}^s (m^k - 1)!
$$
points: there are $\binom d {m^1, \dots, m^s}$ ways to divide the $d$
points $\{ x^1, \dots, x^d \}$ above $x$ among $C(1)$, \dots, $C(s)$ and
$(m^k - 1)!$ possible choices of the markings above $x$ on $C(k)$ up to
automorphisms of $\pi |_{C(k)}$.  This is the number of partitions of $\{
x^1, \dots, x^d \}$ into cyclically-ordered subsets of sizes $m^1$, \dots,
$m^s$.  Hence the multiplicity at
each one of these points must be
$$
{\frac {d!} {
{\binom {d} {m^1, \dots, m^s} } \prod_{k=1}^s (m^k - 1)! } } =
\prod_{k=1}^s m^k.
$$

{\em Step 6${}'$:  The calculation.} We now determine the local structure at one of these points.  Fix sections
$f$ and $g$ of $\oh_C(d)$ with
$$
(f=0) = \sum m p^j_m, \quad (g=0) = \sum m^k y^k
$$
and let $\pi$ be the morphism to $\proj^1$ given by $\pi = [f,g]$.

Rather than considering elements of $\cx''_x$ as maps $[\be f, \ga g]$
(with $x$ fixed, and $\pi(q)$ varying), we now consider them as maps
$[f,g]$ with $x$ moving (and $\pi(q)$ fixed).  (We are degenerating the
point $(\proj^1, x, y, z, \pi(q)) \in \mbar_{0,4}$ in two ways.  Originally
we fixed $x$, $y$, $z$ and let $\pi(q)$ degenerate to $z$.  Now we fix
$y$, $z$, $\pi(q)$ and let $x$ tend to $y$.  They are equivalent as they
represent the same point in the curve $\mbar_{0,4} = \proj^1$.)
The Weil divisor is now
defined (set-theoretically) by $\{ x = y \}$, not $\{ \pi(q) = z \}$.  Then
$$
\cx''_x = \overline{\underbrace{C \times_{\pi} \dots \times_{\pi} C}_{d} \setminus \De}
$$
where $\De$ is the big diagonal (where any two of the factors are the same),
and the closure is in
$$
\underbrace{C \times_{\pi} \dots \times_{\pi} C}_{d}.
$$

Fix one of the points of $\cy''_x$, which corresponds to a partition of $\{
x^1, \dots, x^d \}$ into subsets of sizes $m^1$, \dots, $m^s$ and cyclic
orderings of these subsets.  Consider a neighborhood of this point in
$\cx''_x$.  By relabeling if necessary, we may assume that $x^k$ is in the
$k^{\text{th}}$ subset for $1 \leq k \leq s$ (so, informally, $x^k$ is close
to $t_k$; see Figure \ref{inap}).

\begin{figure}
\begin{center}

	   \setlength{\unitlength}{.1\baseunit}
	    \input inap.tex 
\end{center}
\caption{Near a point of $\cy''_x$}
\label{inap}
\end{figure}

In an \'{e}tale neighborhood of point in $\cy''_x$, 
$$
\cx''_x = \overline{\underbrace{C \times_{\pi} \dots \times_{\pi} C}_{d}
\setminus \De} 
= \overline{\underbrace{C \times_{\pi} \dots \times_{\pi} C}_s
\setminus \De} 
$$
where the second product consists of the first $s$ factors of the first.

Let $a$ be a local co-ordinate for $x$ near $y$.  As the ramification
of $\pi$ at $y^k$ is $m^k - 1$, there is an \'{e}tale-local co-ordinate
$b_k$ for $x^k$ near $y^k$ where $\pi$ is given by $a = b_k^{m^k}$.
Therefore in an \'{e}tale neighborhood of our point of $\cy''_x$, 
\begin{eqnarray*}
\cx''_x 
&=& \overline{\underbrace{C \times_{\pi} \dots \times_{\pi} C}_s
\setminus \De} \\
&=&  \Spec \com [[ a,b_1, \dots, b_s
]] / (a = b_1^{m^1}= \dots = b_s^{m^s}).
\end{eqnarray*}
Hence the deformation space $\Def_{A_1}$ is isomorphic to 
$$
\Spec \com [[ a,b_1, \dots, b_s ]] / (a = b_1^{m^1}= \dots = b_s^{m^s}).
$$
The divisor $D=(y=x)$ is given by $(a=0)$.  

{\em Step 7${}'$:  Returning to the original problem.}
For $j>1$, $A_j$ is either a marked point $p^j_m$ with a ramification of
order $m$ over $z$ or a simple ramification.  In these cases we have $\Def
A_j = 0$ or $\Def A_j = \Spec \com[[c]]$ respectively.  Hence the deformation space
$\Def_{(C, \{ p^j_m \}, q, \pi)}$ is isomorphic to 
$$
\Spec \com [[a, b_1,\dots,b_s,c_1,\dots,c_{\dim X-1}]] / (a = b_1^{m^1} = \dots =
b_s^{m^s}).
$$
\epf

Proposition \ref{imultg0} and Corollary \ref{ilocalst} above are statements
about varieties, so long as $d \neq 2$.

In order to extend these results to components for which $p_a(C(0)) = 1$,
we  will need the following result.
\begin{pr}
Let $Y$ be a component of $D$, with $(C,\{p^j_m \}, q, \pi)$ the map
corresponding to a general point of $Y$, $\{ r^1, ..., r^s \} = C(0)
\cap \tilde{C}$, and 
$m^k$ the multiplicity of $\pi^*(z)$ along $\tilde{C}$ at $r^k$.  Then
$$
\oh_{C(0)}\left(\sum_{m,j} m p^j_m(0)\right) \cong
\oh_{C(0)}\left(\sum_{k=1}^s m^k r^k \right)
$$
where $\{ p^j_m(0) \}^{h_m(0)}_{j=1} \subset \{ p^j_m \}_{j=1}^{h_m}$ are
the marked points whose limits lie in $C(0)$.
\end{pr}
\bpf
For a map $(C, \{ p^j_m \}, q ,\pi)$ corresponding to a general point in $X$, we
have the following relation in the Picard group of $C$:
$$
\pi^*(\oh_{\proj^1}(1)) \cong \oh_C(\sum_{m,j} m p^j_m).
$$
Thus for the curve corresponding to a general point of our component of $D$
the invertible sheaf $\oh_C (\sum_{m,j} m p^j_m)$ must be a possible limit of
$\pi^*(\oh_{\proj^1}(1))$.  The statement of the lemma depends only on an analytic
neighborhood of $C(0)$, so we may assume (as in Step 2 of the proof of
Proposition \ref{imultg0}) that $\vh = \vh(0)$, $p_a(C) = 1$,
and $\tilde{C}$ consists of $k$ rational tails $C(1)$, \dots, $C(k)$ each totally ramified where
they intersect $C(0)$.  As the dual graph of $C$ is a tree, $C$ is of
compact type (which means that $\Pic C$ is compact).  One possible limit of
$\pi^*(\oh_{\proj^1}(1))$ is the line bundle that is trivial on $C(0)$ and degree
$m^k$ on $C(k)$.  If a curve $C'$ is
the central fiber of a one-dimensional family of curves, and $C' = C_1 \cup
C_2$, and a line bundle $\cll$ is the limit of a family of line bundles,
then the line bundle $\cll'$ whose restriction to $C_i$ is $\cll|_{C_i}(
(-1)^i C_1 \cap C_2)$ is another possible limit.  Thus the line bundle that
is trivial on $\tilde{C}$ and $\oh_C(\sum m^k r^k)$ on $C(0)$ is a possible
limit of $\pi^*(\oh_{\proj^1}(1))$.

If two line bundles on a nodal curve $C$ of compact type are possible
limits of the same family of line bundles, and they agree on all components
but one of $C$, then they must agree on the remaining component.  But
$\oh_C(\sum m p^j_m)$ is another limit of $\pi^*(\oh_{\proj^1}(1))$ that is
trivial on $\tilde{C}$, so the result follows.
\epf

We are now in a position to determine all components $Y$ of $D$ satisfying
$p_a(C(0)) = 1$.  For each choice of a partition of $\{ p^j_m \} $ into $\{
p^j_m(0) \} \cup \{ \tilde{p}^j_m \}$, a positive integer $s$, and $(m^1,
\dots, m^s)$ satisfying $\sum m^s + \sum m \tilde{h}^j_m = d$, there is a
variety (possibly reducible) that is the closure
in $\mbar_{g, \sum h_m + 1}(\proj^1,d)$ of points corresponding to maps
$$
(C(0) \cup \tilde{C}, \{ p^j_m \}, q, \pi)
$$
where 
\begin{enumerate}
\item[D1.] The curve $C(0)$ is a smooth elliptic curve with labeled points 
$$
\{ p^j_m(0) \} \cup \{ r^k \} \cup \{ q \},
$$
where $\oh_{C(0)}(\sum m p^j_m(0)) \cong \oh_{C(0)}( \sum m^k r^k )$,
and $\pi(C(0)) = z$.
\item[D2.] The curve $\tilde{C}$ is smooth of arithmetic genus $g-s$ with
marked points $\{ \tilde{p}^j_m \} \cup \{ r^k \}$.  The map $\pi$ is
degree $d$ on $\tilde{C}$, and
$$
( \pi|_{\tilde{C}})^*(z) = \sum m \tilde{p}^j_m + \sum m^k r^k.
$$
\item[D3.] $C(0) \cup \tilde{C}$ is nodal, and $C(0)$ and $\tilde{C}$ are glued
at the points $\{ r^k \}$.
\end{enumerate}

Let $U$ be the union of these varieties.  The divisorial
condition $\oh_{C(0)}(\sum m p^j_m(0)) \cong \oh_{C(0)}( \sum m^k r^k )$  defines a
substack $\cm'$ of pure codimension 1 in $\cm_{1,\sum h_m + 1 + s}$:  for
any 
$$
(C, \{ p^j_m \}, q, \{ r^k \}_{k>1})\in \cm_{1, \sum h_m + 1 + (s-1)}
$$
the subscheme of points $r^1 \in C$ satisfying 
$$
\oh(m^1 r^1) \cong \oh \left( 
\sum m p^j_m(0) -  \sum_{k>1} m^k r^k \right)
$$
is reduced of degree $(m^1)^2$.  The stack $\cm'$ is a degree $(m^1)^2$
\'{e}tale cover of $\cm_{1, \sum h_m + 1 + (s-1)}$.  By this observation
and (\ref{inaive}), $U$ has pure dimension $\dim X - 1$.

An irreducible component $Y$ of the divisor $D$ satisfying $p_a(C(0)) = 1$
has dimension $\dim X - 1$ and is a subvariety of $U$, which also has
dimension $\dim X - 1$.  Hence $Y$ must be a component of $U$ and the
stable map corresponding to a general point of $Y$ satisfies 
properties D1--D3 above.

The determination of multiplicity and local
structure is identical to the genus 0 case.
\begin{pr}
\label{imultg1}
Fix such a component $\cy$ with $p_a(C(0)) = 1$.  If $m^1$, \dots, $m^s$ are
the multiplicities of $\pi^* (z)$ along $\tilde{C}$ at the $s$ points
$C(0) \cap \tilde{C}$, then this divisor appears with multiplicity
$\prod_k m^k$.
\end{pr}
\bpf
The proof is identical to that of Proposition
\ref{imultg0}.  We summarize the steps here.

{\em Step 1.}  If $A^1$ is the special locus of $\pi$
containing $q$, then it suffices to analyze $\Def A_1$.

{\em Step 2.} We may assume that the map corresponding to a general point
in $\cy$ consists of $C(0)$ and $s$ rational tails ramifying completely
over $z$.

{\em Step 3.}
We require the $s$ rational tails to ramify completely over another point
$y$, and we label these ramifications $y^1$, \dots, $y^s$.

{\em Step 4.} Let $\cmbar'_{1, \sum h_m + 1 + s}$ be the
substack of $\cmbar_{1, \sum h_m + 1 + s}$ that is the closure
of the set of points representing smooth stable curves where
$\oh(\sum m p^j_m) \cong \oh(\sum m^k y^k)$.  If $\al$ is defined by 
$$
\al: \cx' \rightarrow \cmbar'_{1, \sum
h_m + 1 + s}
$$
then $\al |_{\cy}$ 
is dominant, so we may consider a fixed stable curve
$$
(C, \{ p^j_m \}, q, \{ y^k \})
\in 
\cmbar'_{1, \sum h_m + 1 + s}.
$$

{\em Step 5.}  The variety $X''$ is $\proj^1$.

{\em Step 6.}  The multiplicity calculation is identical.
\epf

Once again, we get the \'{e}tale or formal local structure of $\cx$ as a
corollary.  The statement and proof are identical to those of Corollary \ref{ilocalst}.
\begin{co}
Let $\cy$ be the same component as in Proposition \ref{imultg1}.  
In an \'{e}tale neighborhood of a general point of $\cy$, the
stack $\cx$ is isomorphic to
$$
\Spec \com [[a, b_1,\dots,b_s,c_1,\dots,c_{\dim X - 1}]] / (a = b_1^{m^1} = \dots =
b_s^{m^s})
$$
with $D$ given by $\{ a = 0 \}$, and $\cy$ given set-theoretically by the
same equation.
\end{co}

If $g=0$ or $1$, then $p_a(C(0)) = 0$ or $1$, so we have found all
components of $D=\{ \pi(q) = z \}$ and the multiplicity of $D$ along each
component.  We summarize this in two theorems which will be invoked in
Sections \ref{rational} and \ref{elliptic}.

\begin{tm}
\label{igenus0}
If $g=0$, the components of $D = \{ \pi(q) = z \}$ on $X^{d,g}(\vh)$
are of the following form.  Fix a positive integer $l$, $\{ d(k)
\}_{k=1}^l$ with $\sum_{k=0}^l d(k) = d$, and a partition of the
points $\{ p^j_m \}$ into $l$ subsets $\cup_{k=1}^l \{ p^j_m(k) \}$.
This induces a partition of $\vh$ into $\sum_{k=0}^l \vh(k)$.  Let
$m^k = d(k) - \sum m h_m(k)$.  Then the general member of the
component is a general map from $C(0) \cup
\dots \cup C(l)$ to $\proj^1$, where:
\begin{itemize}
\item The irreducible components $C(0)$, \dots, $C(l)$ are rational.
\item The curve $C(0)$ contains the marked points $\{ p^j_m(0) \}$, $q$,
$\{ r^k \}$, and $\pi(C(0)) = z$.
\item For $k>0$, $C(k)$ maps to $\proj^1$ with degree $d(k)$ and 
$$
\pi^*(z) |_{C(k)} = \sum_{m,j} m p^j_m(k) + m^k r^k.
$$
The curves $C(0)$ and $C(k)$ intersect at $r^k$.
\end{itemize}
This component appears with multiplicity $\prod_k m^k$.
\end{tm}

\begin{tm}
\label{igenus1}
If $g=1$, the components of $D= \{ \pi(q) = z \}$ on $X^{d,g}(\vh)$ are of the following
form.
Fix a positive integer $l$, $\{ d(k) \}_{k=1}^l$ with $\sum_{k=0}^l d(k) =
d$, and a partition of the points $\{ p^j_m \}$ into $l$ subsets
$\cup_{k=1}^l \{ p^j_m(k) \}$.  This induces a partition of $\vh$ into
$\sum_{k=0}^l \vh(k)$.  Let $m^k = d(k) - \sum m h_m(k)$. 
Then the general member of the component is a general map from $C(0) \cup
\dots \cup C(l)$ to $\proj^1$, where:
\begin{itemize}
\item The curve $C(0)$ contains the marked points $\{ p^j_m(0) \}$ and $q$,
and $\pi(C(0)) = z$.
\item For $k>0$, $C(k)$ maps to $\proj^1$ with degree $d(k)$.
\end{itemize}
Furthermore, one of the following cases holds:
\begin{enumerate}
\item[a)]The curve $C(1)$ is elliptic and the other components are
rational.  When 
$k>0$, 
the curves
$C(0)$ and $C(k)$ intersect at the point $r^k$, and 
$\pi^* (z) |_{C(k)} = \sum_{m,j} p^j_m(k) + m^k r^k$.  
\item[b)] All components are rational.  When $k>1$, 
the curves $C(0)$ and $C(k)$ intersect
at the point $r^k$, and
$\pi^* (z) |_{C(k)} =
\sum_{m,j} m p^j_m(k) + m^k r^k$.
The curves $C(1)$ and $C(0)$ intersect at two points $r^1_1$ and $r^1_2$, and the ramifications
$m^1_1$ and $m^1_2$ at these two points sum to $m^1$.  $\pi^*(z) |_{C(1)} =
\sum_{m,j} m p^j_m(k) + m^1_1 r^1_1 + m^1_2 r^1_2$.
\item[c)]The curve $C(0)$ is elliptic and the other components are
rational.  When $k>0$, the curves $C(0)$ and $C(k)$ intersect at $r^k$, and
$\pi^* (z) |_{C(k)} = \sum_{m,j} m p^j_m(k) + m^k r^k$.  Also,
$\oh_{C(0)}( \sum_{m,j} m p^j_m(0)) \cong \oh_{C(0)}( \sum m^k r^k)$.
\end{enumerate}
The components of type a) and c) appear with multiplicity $\prod_{k=1}^l  m^k$ and
those of type b) appear with multiplicity $m^1_1 m^1_2 \prod_{k=2}^l m^k$.
\end{tm}

In all three cases, the
multiplicity is the product of the ``new ramifications'' of
the components not mapped to $z$. 

For general $g$, the above argument identifies some of the components, but further
work is required to determine what happens when $p_a(C(0)) > 1$.

\subsection{Pathological behavior of $\mbgd$}

When $d>2$, $g>0$, $\mbgd$ has more than one component.  The most
interesting one consists (generically) of irreducible genus $g$ curves.
Call this one $\mbgd^o$.  A second consists (generically) of two
intersecting components, one of genus $g$ and mapping to a point, and the
other rational and mapping to $\proj^1$ with degree $d$.  The first has
dimension $2d+2g-2$, and the second has dimension $2d+3g-3$, so the second
is not in the closure of the first.

The local structure of $\mbgd^o$ may be complicated
where it intersects the other components.  One might hope that $\mbgd^o$ is
smooth (at least as a stack, or equivalently on the level of deformation
spaces). This is not the case; $\mbgd^o$ can be singular and even fail to
be unibranch.

Let $g=3$ (so a general degree 3 divisor has one section, but all degree 4
divisors have two) and $d=4$.  Then $\mbar_3(\proj^1,4)^o$ has dimension
12.  Consider the family $Y$ of stable maps whose general element
parametrizes a smooth genus 3 curve $C(0)$ with two rational tails $C(1)$
and $C(2)$.  The curve
$C(0)$ maps with degree 0 to $\proj^1$, and the rational tails each map with
degree 2 to $\proj^1$, both ramifying at their intersection
with $C$.

The subvariety $Y$ has dimension 11:  8 for the two-pointed genus 3 curve $C(0)$,
1 for the image of $C(0)$ in $\proj^1$, 
and 2 for the other ramification points of $C(1)$ and $C(2)$.  Thus if $Y$
is contained in $\mbar_3(\proj^1,4)^o$, it is a Weil divisor. 

\begin{pr}
$\mbar_3(\proj^1,4)^o$ has two smooth branches along $Y$, intersecting transversely.
\end{pr}
This gives an example of a map that could be smoothed in two
different ways.  

\bpf
For convenience, we use the language of stacks.  Let $\cx$ be the stack
corresponding to $X$, and $\cy$ the stack corresponding to $Y$.

Let $(C, \pi)$ be the map corresponding to a general point of $\cy$ (where
$C = C(0) \cup C(1) \cup C(2)$ where $C(0)$, $C(1)$, $C(2)$ are as
described above), and let $A_1$ be the special locus consisting of the
collapsed genus 3 curve.  By Proposition
\ref{ilocal}, it suffices to consider $\Def_{A_1}$.

As in Step 3 of Proposition \ref{imultg0}, we fix the other special loci.  Fix a
point $y \in \proj^1$.  We may restrict attention to maps ramifying at two
points above $y$, labeled $y^1$ and $y^2$.  Denote by $\cx'$ and $\cy'$ the
substacks of $\cmbar_{3,2}(\proj^1,4)$ (with ramification above $y$ at
$y^1$ and $y^2$) corresponding to $\cx$ and $\cy$.  It suffices to prove the
corresponding result for $\cx'$ and $\cy'$.

Next, fix another point $z \in \proj^1$ ($z \neq y$) and mark the four
pre-images of $z$ with the labels $p^1$, \dots, $p^4$.  Denote by $\cx''$
and $\cy''$ the substacks of $\cmbar_{3,6}(\proj^1,4)$ corresponding to
$\cx'$ and $\cy'$.  It suffices to prove the corresponding result for $\cx''$
and $\cy''$.

Let $\cm$ be the substack of $\cmbar_{3,6}$ that is the closure of the
points representing stable marked curves $(C, \{ p^i \}, y^1, y^2)$ where
$C$ is smooth and $\oh_C(p^1 + \dots + p^4) \cong \oh_C(2 y^1 + 2 y^2)$.
If $\al$ is the natural map $\al:  \cmbar_{3,6}(\proj^1,4) \rightarrow
\cmbar_{3,6}$ then  $\al(\cx'') \subset \cm$ as for a general map $(C, \{
p^i \}, y^1, y^2, \pi) \in \cx''$,
$$
\oh_C(p^1 + \dots + p^4) \cong \pi^* ( \oh_{\proj^1}(1)) \cong \oh_C(2 y^1
+ 2 y^2).
$$
Moreover $\al |_{\cy''}$ surjects onto $\cm$:  for any curve 
$(C, \{ p^i \}, r^1, r^2 ) \in \cm$ consider the map
$(C(0) \cup C(1) \cup C(2) , \{ p^i \}, y^1, y^2, \pi)$ where $C(0)$ and
$C(k)$ intersect at $r^k$ ($k=1,2$),
$$
(C(0), \{ p^i \}, r^1, r^2) \cong (C, \{ p^i \}, r^1, r^2 ) 
$$
as marked curves, $\pi(C(0)) = z$, and $C(k)$ ($k = 1,2$) maps to $\proj^1$
with degree 2 ramifying over $y$ (at $y^k$) and $z$ (at $r^k$).  The stable
model of the underlying curve of such a map is indeed isomorphic to 
$(C, \{ p^i \}, r^1, r^2 )$.

Thus we may restrict attention to a fixed (general) marked curve
$(C, \{ p^i \}, y^1, y^2)$ in $\cm$.  

We may now directly follow steps 5${}'$ and $6'$ in the proof of
Corollary \ref{ilocalst}.  The steps are summarized here.

{\em Step 5${}'$.}  Rigidify the moduli problem by eliminating automorphisms.
Fix a point $x \in \proj^1$ distinct from $y$ and $z$, and mark the points
of 
$\pi^{-1}(x)$ with the labels $\{ x^1, \dots, x^d \}$.  Call the resulting
stack $\cx''_x$.

{\em Step 6${}'$.}  Observe that 
$$
\cx''_x = \overline{\underbrace{C \times_{\pi} \dots \times_{\pi} C}_{d} \setminus \De}
$$
where $\De$ is the big diagonal, and the closure is in
$$
\underbrace{C \times_{\pi} \dots \times_{\pi} C}_{d}.
$$
Then show that 
$$
\Def_{A_1} = \Spec \com [[ a, b_1, b_2]] / (a = b_1^2 = b_2^2).
$$
\epf

By a similar argument, we can find a codimension 1 unibranch singularity of
$\mbar_4(\proj^1,5)^o$ and singularities of 
$\mbar_6(\proj^1,7)^o$) with 
several codimension 1 singular branches.

\subsection{Possible applications of these methods}

Ideas involving degenerations and the space of stable maps can be
used in other cases besides those dealt with in this article, and we will
briefly mention them here.  In [V1] the same methods are used
for genus $g$ curves in a divisor class $D$ on the rational ruled surface
${\mathbb F}_n$: the curves through $-K_{{\mathbb F}_n}
\cdot D + g - 1$ general points are enumerated; these also repreoduce the 
calculations in [CH3] of genus $g$ Gromov-Witten invariants of the
plane (i.e. the enumerative geometry of plane curves) in the language
of maps.  The genus $g$ Gromov-Witten invariants of the plane blown up
at up to five points are calculated in [V2].  In [V3], the analogous
question for plane curves with certain allowed singularities is
incompletely addressed.  It also seems possible that the classical
question of characteristic numbers of rational and elliptic curves in
$\proj^n$ would be susceptible to such an approach.  (The rational
case has already been settled by L. Ernstr\"{o}m and G. Kennedy by
different methods in [EK1].)  

\section{Rational Curves in Projective Space}
\label{rational}


In this section, we use the ideas and results of Section \ref{intro} to study the
geometry of varieties parametrizing degree $d$ rational curves in $\proj^n$
intersecting fixed general linear spaces and tangent to a fixed hyperplane
$H$ with fixed multiplicities along fixed general linear subspaces of $H$.

We employ two general ideas.  First, we specialize a linear space (that the
curve is required to intersect) to lie in the hyperplane $H$, and analyze
the limit curves.  It turns out that the limit curves are of the same form,
and are in some sense simpler.  Enumerative results have been proved using such
specialization ideas since the nineteenth century
(see [PZ], for example, especially pp. 268--275 and pp. 313--319).

The second general idea we use is Kontsevich's moduli spaces of stable
maps, particularly the spaces $\mbar_{0,m}(\proj^n,d)$ and
$\mbar_{0,m}(\proj^1,d)$ (and the stacks $\cmbar_{0,m}(\proj^n,d)$ and
$\cmbar_{0,m}(\proj^1,d)$).  The calculations in Section \ref{intro} on the
space $\mbar_{0,m}(\proj^1,d)$ will give the desired
results in $\proj^n$.  L. Caporaso and J. Harris' results on plane curves of
any degree and genus (cf. [CH3]) can also be reinterpreted in this light.  
The reader can verify that the
argument here for $n=2$ is in essence the same as that in [CH3] for genus 0
curves.

\subsubsection{Example:  2 lines through 4 general lines in $\proj^3$}
\begin{figure}
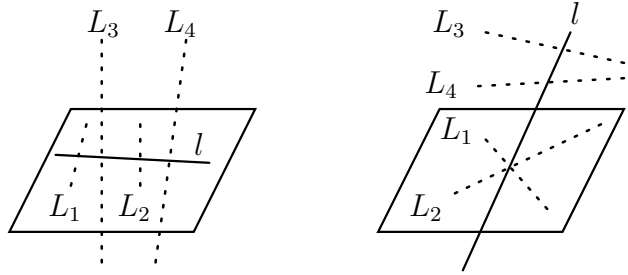

\begin{center}

	   \setlength{\unitlength}{.15\baseunit}
	    \input r2lines.tex 
\end{center}
\caption{Possible positions of $l$ after $L_1$ and $L_2$ have
degenerated to $H$}
\label{r2lines}
\end{figure}

We can follow through these ideas in a classical special case.  Fix four
general lines $L_1$, $L_2$, $L_3$, $L_4$ in $\proj^3$, and a hyperplane
$H$.  There are a finite number of lines in $\proj^3$ intersecting $L_1$,
$L_2$, $L_3$, $L_4$.  Call one of them $l$.  We will specialize the lines
$L_1$, $L_2$, $L_3$, and $L_4$ to lie in $H$ one at a time and see what
happens to $l$.  First, specialize the line $L_1$ to (a general line in)
$H$, and then do the same with $L_2$ (see Figure \ref{r2lines}).  If $l$
doesn't pass through the intersection of $L_1$ and $L_2$, it must still
intersect both $L_1$ and $L_2$, and thus lie in $H$.  Then $l$ is uniquely
determined: it is the line through $L_3 \cap H$ and $L_4
\cap H$.  Otherwise, if $l$ passes through the point $L_1 \cap L_2$, it is
once again uniquely determined (as only one line in $\proj^3$ can pass
through two general lines and one point --- this can also be seen through
further degeneration).

The above argument can be tightened to rigorously show that there are two
lines in $\proj^3$ intersecting four general lines.  The only information
one needs to know in advance is that there is one line through two distinct
points.  This is the same as the seed data for Kontsevich's recursive
formula in [KM], and it is all we will need in this section.

\subsubsection{Example:  92 conics through 8 general lines in $\proj^3$}  
\label{r92subsection}
\begin{figure}
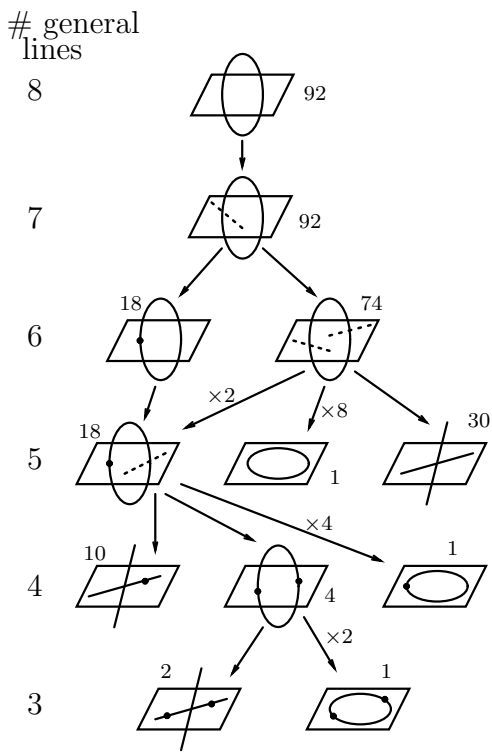

\begin{center}

	   \setlength{\unitlength}{.1\baseunit}
	    \input r92conics.tex 
\end{center}
\caption{Counting 92 conics in $\proj^3$ through 8 general lines}
\label{r92conics}
\end{figure}

The example of conics in $\proj^3$ is a simple extension of that of lines
in $\proj^3$, and gives a hint as to why stable maps are the correct way to
think about these degenerations.  Consider the question: How many conics
pass through 8 general lines $L_1$, \dots, $L_8$?  (For another discussion
of this classical problem, see [H1] p. 26.)  We introduce a pictorial
shorthand that will allow us to easily follow the degenerations (see
Figure \ref{r92conics}); the plane $H$ is represented by a parallelogram.

We start with the set of conics through 8 general lines (the top row of the
diagram --- the label 92 indicates the number of such conics, which we will
calculate last) and specialize one of the lines $L_1$ to $H$ to get row 7.
(The line $L_1$ in $H$ is indicated by the dotted line in the figure.)
When we specialize another line $L_2$, one of two things can happen: the
conic can intersect $H$ at the point $L_1 \cap L_2$ and one other (general)
point, or it can intersect $H$ once on $L_1$ and once on $L_2$ (at general
points).  (The requirement that the conics must pass through a fixed point
in the first case is indicated by the thick dot in the figure.)

In this second case (the picture on the right in row 6), if we specialize
another line $L_3$, one of three things can happen.
\begin{enumerate}
\item The conic can stay smooth, and not lie in $H$, in which case it must
intersect 
$H$ at $\{ L_1 \cap L_3, L_2 \}$ or $\{L_1, L_2 \cap L_3 \}$ (hence the
``$\times 2$'' in the figure).
\item The conic could lie in $H$.  In this case, there are eight conics through
five fixed points $L_4 \cap H$, \dots, $L_8 \cap H$ with marked points on
the lines $L_1$, $L_2$, and $L_3$.
\item The conic can degenerate into the union of two intersecting lines,
one ($l_0$) in $H$ and one ($l_1$) not.  These lines must intersect $L_4$,
\dots, $L_8$.  (The line $l_0$ already intersects $L_1$, $L_2$, $L_3$, so we don't
have to worry about these conditions.)  Either three or four of $\{ L_4,
\dots, L_8 \}$ intersect $l_1$.  In the first case, there are $\binom 5 3$
choices of the three lines, and two configurations $(l_0,l_1)$ once the
three lines are chosen.  In the second case there are a total of $\binom 5
4 \times 2$ configurations by similar reasoning.  Thus the total number of such configurations
is 30.
\end{enumerate}
We fill out the rest of the diagram in the same way.  Then, using the
enumerative geometry of lines in $\proj^3$ and conics in $\proj^2$ we can
work our way up the table, attaching numbers to each picture, finally
deducing that there are 92 conics through 8 general lines in $\proj^3$.
To make this argument rigorous, precise dimension counts and multiplicity
calculations are needed.

\begin{figure}
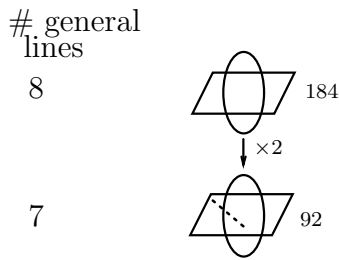

\begin{center}

	   \setlength{\unitlength}{.1\baseunit}
	    \input r92conics2.tex 
\end{center}
\caption{Counting 184 conics with two marked points on $H$ through 8
general lines}
\label{r92conics2}
\end{figure}

The algorithm described in this section is slightly different:  we will
parametrize rational curves with various conditions {\em and marked
intersections with $H$}.  In the case of conics through 8 lines, for
example, we would count 184 conics through 8 lines with 2 marked points on
$H$, and then divide by 2.  The argument will then be cleaner.  The
resulting pictorial table is almost identical to
Figure \ref{r92conics}; the only difference is in the first two lines (see
Figure \ref{r92conics2}).

\subsubsection{Example:  Cubics in $\proj^3$}
The situation in general is not much more complicated than our
calculations for conics in $\proj^3$.  Two additional twists come up,
which are illustrated in the case of the $80,160$ twisted cubics
through 12 general lines in $\proj^3$.  This calculation is indicated
pictorially in Figure \ref{rcubics} at the beginning of the
introduction.  The third figure in row 8 represents a nodal (rational)
cubic in $H$.  There are 12 nodal cubics through 8 general points in
$\proj^2$ (well-known, e.g. [DI] p. 85).  The algorithm of this
section will actually calculate $80,160 \times 3!$ cubics with marked
points on $H$ through 12 general lines.

On the left side of row 8 we see a new degeneration (from twisted cubics
through nine general lines intersecting $H$ along three fixed general lines
in $H$):  a conic tangent to $H$, intersecting a line in $H$. (The tangency
of the conic is indicated pictorially by drawing its lower horizontal
tangent inside the parallelogram representing $H$.)  We also have an
unexpected multiplicity of 2 here.  

The appearance of these  new degenerations indicate why, in order to
enumerate rational curves through general linear spaces by these
degeneration methods, we must expand the set of curves under consideration to
include those required to intersect $H$ with given multiplicity, along linear
subspaces.

\subsubsection{The algorithm (informally)}

The algorithm in general is (informally) as follows.  
(Theorem \ref{rrecursiveX2} describes the algorithm rigorously.)

Fix positive integers $n$ and $d$, and fix a hyperplane $H$ in $\proj^n$.
Let $\vh = ( h_{m,e} )_{m \geq 1, e \geq 0}$ and $\vi = ( i_e )_{e \geq
0}$ be sets of non-negative integers.  Fix a set of general linear spaces
$\Ga = \{ \Ga^j_{m,e} \}_{m,e,1 \leq j \leq h_{m,e}}$ in $H$ where $\dim
\Ga^j_{m,e} = e$.  Fix a set of general linear spaces
$\De = \{ \De^j_e \}_{e,1 \leq j \leq i_e}$ in $\proj^n$ where $\dim
\De^j_{e} = e$.  Let $X_n(d, \Ga, \De)$ be the closure in
$\mbar_{0, \sum h_{m,e} + \sum i_e}(\proj^n,d)$ of points representing
maps
$$
(\proj^1, \{ p^j_{m,e} \}_{m,e,1 \leq j \leq h_{m,e}}, \{ q^j_e \}_{e, 1
\leq j \leq i_e }, \pi)
$$
where
\begin{itemize}
\item $\pi(p^j_{m,e}) \in \Ga^j_{m,e}$, $\pi(q^j_e) \in \De^j_e$.
\item $\pi^{-1}H$ is a set of points, and as divisors $\pi^* H =
\sum_{m,e,j} m p^j_{m,e}$.
\end{itemize}

Assume $X(d,\Ga,\De)$ is a finite set.  We will count the points of
$X(d,\Ga, \De)$.  Specialize one of the linear spaces
$\De^{j_1}_{e_1}$ to lie in $H$, and consider the limits of the $\#
X(d,\Ga, \De)$ stable maps.  One of the two following types of limits will appear.

\begin{enumerate}
\item
The limit map is of the form 
$$
(C(0) \cup C(1),\{ p^j_{m,e} \}_{m,e,j}, \{ q^j_e \}_{e, j}, \pi)
$$
where the curves $C(0)$ and $C(1)$ are both smooth and
rational, $\pi(C(0))$ is a point, for some $(m_0, e_0, j_0)$ the curve
$C(0)$ contains the marked points $q^{j_1}_{e_1}$ and $p^{j_0}_{m_0, e_0}$
(and $C(1)$ contains the other marked points), and
$$
(\pi \mid_{C(1)})^* H = \sum_{ \substack{  {m,e,j} \\ {(m,e,j) \neq (m_0, e_0,
j_0)}}} m p^j_{m,e} + m_0 (C(1) \cap C(0)).
$$
Also, $\pi ( p^j_{m,e}) \in \Ga^j_{m,e}$, $\pi(q^j_e) \in \De^j_e$, and (as
a consequence) $\pi(C(0))$ is contained in $\De^{j_1}_{e_1} \cap \Ga^{j_0}_{m_0,e_0}$.
We can ignore the rational tail, replacing it
with another marked point, and continue the process.

There are $m_0$ curves of $X(d,\Ga,\De)$ tending to this limit.

\item The limit map is of the form 
$$
(C = C(0) \cup C(1) \cup \dots \cup C(l),\{ p^j_{m,e} \}_{m,e,j}, \{ q^j_e \}_{e, j }, \pi)
$$
where $C(k)$ ($0 \leq k \leq l$) is smooth and rational.  The points $
\{ p^j_{m,e} \}$, $\{ q^j_e \}$ are  partitioned into sets $\{ p^j_{m,e}(k)
\}$, $\{ q^j_e(k) \}$, where the $k^{\text{th}}$ subset lies in $C(k)$;
this induces partitions $\vh = \sum_{k=0}^l \vh(k)$ and $\vi = \sum_{k=0}^l
\vi(k)$.   The marked point $q^{j_1}_{e_1}$ lies on $C(0)$; that is, 
$q^{j_1}_{e_1} \in \{ q^j_e(0) \}$.  The component $C(0)$ intersects all
other components $C(1)$, \dots, $C(l)$.  The map $\pi$ sends $C(0)$ to $H$ with
positive degree, and sends no other component of $C$ to $H$.  If $m^k =
\deg (\pi \mid_{C(k)}) - \sum_{m,e} m h_{m,e}(k)$, then
$$
( \pi \mid_{C(k)})^* H = \sum_{m,e,j} m p^j_{m,e}(k) + m^k (C(0) \cap C(k))
$$
as divisors for $1 \leq k \leq l$.  Finally, $\pi(p^j_{m,e}) \in
\Ga^j_{m,e}$ for all $m$, $e$, $j$, and $\pi(q^j_e) \in \De^j_e$ for all
$e$, $j$.  There are $\prod_{k=1}^l m^k$ curves in $X(d,\Ga, \De)$ tending
to this limit.
\end{enumerate}

Examples of both types of limits can be seen in Figure \ref{rcubics}.
Given the results of Subsection \ref{ikey}, the
algorithm and multiplicities are not completely unexpected.  

\subsection{Notation and summary}

For convenience, let $\vep_e$, $\vep_{m,e}$ be the natural basis vectors:
$(\vep_e)_{e'} = 1$ if $e = e'$ and 0 otherwise; and $(\vep_{m,e})_{m',e'}
= 1$ if $(m,e)=(m',e')$, and 0 otherwise.  Fix a hyperplane $H$ in
$\proj^n$, and a hyperplane $A$ of $H$. 

\subsubsection{The schemes $X(\ce)$}
\label{r21}
We now define the primary objects of interest to us.

Let $n$ and $d$ be positive integers, and let $H$ be a hyperplane in $\proj^n$.
Let $\vh = ( h_{m,e} )_{m \geq 1, e \geq 0}$ and $\vi = ( i_e )_{e \geq
0}$ be sets of non-negative integers.  Let 
$\Ga = \{ \Ga^j_{m,e} \}_{m,e,1 \leq j \leq h_{m,e}}$ be a set of linear
spaces in $H$ where $\dim 
\Ga^j_{m,e} = e$.   Let 
$\De = \{ \De^j_e \}_{e,1 \leq j \leq i_e}$ be a set of linear spaces in
$\proj^n$ where $\dim
\De^j_{e} = e$.

\begin{defn}
\label{defnX}
The scheme $X_n(d,\Ga,\De)$ is the (scheme-theoretic) closure of
the locally closed subset of $\mbar_{0,\sum h_{m,e} + \sum i_e}(\proj^n,d)$ (where the
points 
are labeled $\{ p_{m,e}^j \}_{1 \leq j \leq h_{m,e}}$ and $\{ q_e^j \}_{1 \leq
j \leq i_e}$) representing stable maps $(C, \{ p^j_{m,e} \}, \{
q^j_e \}, \pi)$ satisfying
$\pi(p_{m,e}^j) \in \Ga_{m,e}^j$, $\pi(q_e^j) \in \De_e^j$,  and $\pi^* H =
\sum_{m,e,j} m p^j_{m,e}$. 
\end{defn}
In particular, $\sum_{m,e} m
h_{m,e} = d$, and no component of $C$ is contained in $\pi^{-1}H$.
The incidence conditions define closed subschemes of $\mbar_{0,\sum h_{m,e} + \sum
i_e}(\proj^n,d)$, so the union of these conditions indeed defines a closed
subscheme of $\mbar_{0,\sum h_{m,e} + \sum i_e}(\proj^n,d)$.  

Define $\cx_n(d,\Ga,\De)$ in the same way as a substack of
$$
\cmbar_{0,\sum h_{m,e} + \sum i_e}(\proj^n,d).
$$ 
When we speak of properties that are constant for general $\Ga$ and
$\De$ (such as the dimension of $X(d,\Ga,\De)$), we will write
$X_n(d,\vh, \vi)$.  For convenience, write $\ce$ (for $\ce$verything)
for the data $d, \vh, \vi$, so $X_n(\ce) = X_n(d,\vh,\vi)$.  Also, the
$n$ will often be suppressed for convenience.

The variety $X(d, \Ga, \De)$ can be loosely thought of as parametrizing 
degree $d$ rational curves in projective space intersecting
certain linear subspaces of $\proj^n$, and intersecting $H$ with different
multiplicities along certain linear subspaces of $H$.
For example, if $n=3$, $d=3$, $h_{2,0}=1$, $h_{1,2}=1$, $X$ parametrizes
twisted cubics in $\proj^3$ tangent to $H$ at a fixed point.

In the special case where $h_{m,e}=0$ when $e<n-1$ and $\vi=\vep_n$, define
$\cehat$ by $\hat{d} = d$, $\hat{i}_1 = 1$, $\hat{h}_{m,0} = h_{m,n-1}$.
We will relate the geometry of $\cx_n(\ce)$ to that of $\cx_1(\cehat)$,
which was studied in Subsection
\ref{ikey}.  (The general point of $\cx_n(\ce)$ corresponds to a general
degree $d$ map from $\proj^1$ to $\proj^n$ with $\pi^*H$ consisting of
points with multiplicity given by the partition $(h_{1,n-1}, h_{2,n-1},
\dots)$.  The general point of $\cx_1(\cehat)$ corresponds to a general
degree $d$ map from $\proj^1$ to $\proj^1$ with $\pi^*z$ consisting of
points with multiplicity given by the partition $(h_{1,n-1}, h_{2,n-1},
\dots)$.)  The geometry of $\cx_n(\ce)$ for general $\ce$ can be understood
from this special case.  For example, consider $\cx = \cx_3(d = 2, h_{2,0}
= 1, i_1 = 2)$, the stack parametrizing conics in $\proj^3$ through two
general lines, tangent to $H$ at a fixed point of $H$.  To analyze $\cx$,
we study the stack $\cx_1 \subset \cmbar_{0,2}(\proj^3,2)$ parametrizing
conics tangent to $H$ (where the tangency is labeled $p^1_{2,0}$) with a
marked point $q^1_1$ (with no other incidence conditions).  We take the
universal curve over this stack $\cx_2$ (which can be seen as a substack of
$\cmbar_{0,3}(\proj^3,2)$), and label the point of the universal curve
$q^2_1$.  Then we require $\pi(p^1_{2,0})$ to lie on two general
hyperplanes $H_1$ and $H_2$ (thus requiring $\pi(p^1_{2,0})$ to be a fixed
general point $H \cap H_1 \cap H_2$ of $H$), $\pi(q^1_1)$ to lie on two
general hyperplanes $H_3$ and $H_4$ (thus requiring $\pi(q^1_1)$ to lie on
a fixed general line $H_3 \cap H_4$ of $\proj^3$), and $\pi(q^2_1)$ to lie
on two general hyperplanes $H_5$ and $H_6$ (thus requiring $\pi(q^2_1)$ to
lie on a fixed general line $H_5 \cap H_6$ of $\proj^3$). We shall prove
(in the next section) that the resulting stack is indeed $\cx$.

By these means we show that if the linear spaces $\Ga$, $\De$ are general,
these varieties have the dimension one would naively expect.  The family of
degree $d$ rational curves in $\proj^n$ has dimension $(n+1)d+(n-3)$.
Requiring the curve to pass through a fixed $e$-plane should be a codimension
$(n-1-e)$ condition.  Requiring the curve to be $m$-fold tangent to $H$
along a fixed $e$-plane of $H$ should be a codimension $(m-1)+(n-1-e)$ condition.
Thus we will show (Theorem \ref{rdimX}) that when the linear spaces in
$\Ga$, $\De$ are general, each component of $X(\ce, \Ga, \De)$ has dimension
$$
(n+1) d + (n-3) - \sum_{m,e} (n+m-e-2)h_{m,e} - \sum_e  (n-1-e)i_e.
$$
Moreover, $X(\ce)$ is reduced.  When the dimension is 0, $X(\ce)$ consists
of a finite number of reduced points.  We call their number $\# X(\ce)$ --- these are the
numbers we want to calculate.  Define $\# X(\ce)$ to be zero if $\dim
X(\ce) > 0$.

For example, when $n=3$, $d=3$, $h_{1,2} = 3$, $i_1 = 12$, $\# X(\ce)$
is 3! times the number of twisted cubics through 12 general lines.
(The ``3!'' arises from the markings of the three points of intersection of the
cubic with $H$.)  

\subsubsection{The schemes $Y ( \ce(0); \dots;  \ce(l) )$} 
We will be naturally led to consider subvarieties of $X(\ce,\Ga,\De)$ which
are similar in form.  
Fix $n$, $d$, $\vh$, $\vi$, $\Ga$, $\De$, and a non-negative integer $l$.  Let
$\sum_{k=0}^l d(k)$ be a partition of $d$.  Let the points $\{ p^j_{m,e}
\}_{m,e,j}$ be partitioned into $l+1$ subsets $\{ p^j_{m,e}(k) \}_{m,e,j}$
for $k= 0$, \dots, $l$.  This induces a partition of $\vh$ into
$\sum_{k=0}^l \vh(k)$ and a partition of the set $\Ga$ into
$\coprod_{k=0}^l \Ga(k)$.  
Let the points $\{ q^j_e
\}_{e,j}$ be partitioned into $l+1$ subsets $\{ q^j_{e}(k) \}_{e,j}$
for $k= 0$, \dots, $l$.  This induces a partition of $\vi$ into
$\sum_{k=0}^l \vi(k)$ and a partition of the set $\De$ into
$\coprod_{k=0}^l \De(k)$.    Define $m^k$ by $m^k = d(k) - \sum_m m
h_m(k)$, and assume $m^k>0$ for all $k = 1$, \dots, $l$.

\begin{figure}
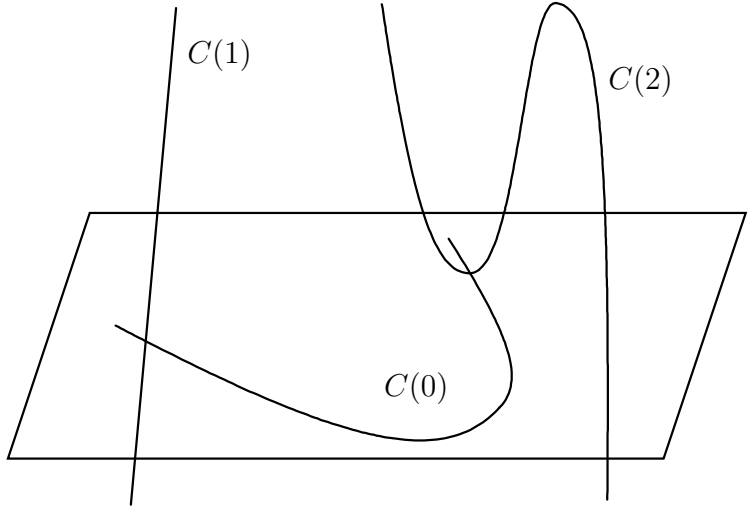

\begin{center}

	   \setlength{\unitlength}{.1\baseunit}
	    \input rtype2eg.tex 
\end{center}
\caption{An example of a map corresponding to a general point of some $Y(\ce(0); \ce(1); \ce(2))$}
\label{rtype2eg}
\end{figure}

\begin{defn}
\label{rdefY}
The scheme 
$$
Y_n(d(0),\Ga(0),\De(0);
\dots; d(l),\Ga(l),\De(l))
$$
is the (scheme-theoretic) closure of
the locally closed subset of $\mbar_{0,\sum h_{m,e} + \sum i_e}(\proj^n,d)$
(where the 
points are labeled $\{ p_{m,e}^j \}_{1 \leq j \leq h_{m,e}}$ and $\{ q_e^j
\}_{1 \leq j \leq i_e}$) representing stable maps $(C, \{ p^j_{m,e}
\}, \{ q^j_e \}, \pi)$ satisfying the following conditions
\begin{enumerate}
\item[Y1.] The curve $C$ consists of $l+1$ irreducible components $C(0)$,
\dots, $C(l)$ with all components intersecting $C(0)$.  The map $\pi$ has
degree $d(k)$ on curve $C(k)$ ($0 \leq k \leq l$).
\item[Y2.] The points $\{ p^j_{m,e}(k) \}_{m,e,j}$ and $\{
q^j_e(k)\}_{e,j}$ lie 
on $C(k)$, and $\pi(p^j_{m,e}(k)) \in \Ga^j_{m,e}(k)$, $\pi(q^j_e(k)) \in
\De^j_e(k)$.  
\item[Y3.] As sets, $\pi^{-1}H = C(0) \cup \{ p^j_{m,e} \}_{m,e,j}$, and for $k>0$,
$$
( \pi \mid_{C(k)} )^* H = \sum_{m,e,j} m p^j_{m,e}(k) + m^k (C(0) \cap C(k)).
$$
\end{enumerate}
\end{defn}
A pictorial representation of such a map is given in Figure \ref{rtype2eg}.
Note that $d(k)>0$ for all positive $k$ by the last condition.

When discussing properties that hold for general $\{ \Ga^j_{m,e}
\}_{m,e,j}$, $\{ \De^j_e \}_{e,j}$, we will write 
$$
Y(\ce(0); \dots;
\ce(l))
= Y(d(0),\vh(0),\vi(0); \dots; d(l), \vh(l),\vi(l)).
$$
If $\vh(k) + \vi(k) \neq \vec{0}$ for all $k>0$,
$Y(d(0),\Ga(0),\De(0); 
\dots; d(l),\Ga(l),\De(l))$ is isomorphic to a closed subscheme of
$$
\mbar_{0,\sum h_{m,e}(0) + \sum i_e(0)+l}(H,d(0)) \times \prod_{k=1}^l 
X(d(k),Ga'(k),\De(k)).
$$
where $\vec{h'}(k) = \vh(k) + \vep_{m^k,n-1}$ and $\Ga'(k)$ is the same as $\Ga(k)$ except $\Ga^{h_{m^k,n-1}}_{m^k,n-1}=H$
Define $\cy(d(0), \De(0), \Ga(0); \dots; d(l), \De(l), \Ga(l))$
as the analogous stack. 

\subsubsection{Enumeratively meaningful subvarieties of
$\mbar_{g,m}(\proj^n,d)$}
\label{rgeomean}
For any irreducible substack $\cv$ of $\cmbar_{g,m}(\proj^n,d)$, there is an
open subset $\cu$ such that for the stable maps $(C, \{ p_i \}, \pi)$
corresponding to points of $\cu$, the reduced image curve
$\pi(C)^{\text{red}}$ has constant Hilbert polynomial $f(t)$.  This gives a
morphism $\xi$ from $\cu$ to the Hilbert scheme $H_{f(t)}$.
\begin{defn}
The substack $\cv$ is {\em enumeratively meaningful} if the dimension of
$\xi(\cu)$ is the same as that of $\cu$ (i.e. $\xi$ is generically finite onto
its image).
\end{defn}
Define {\em enumeratively meaningful subvarieties} of
$\mbar_{g,m}(\proj^n,d)$ in the same way.  When making enumerative
calculations, we are counting reduced points of $\mbar_{g,m}(\proj^n,d)$,
which are obviously enumeratively meaningful.  This definition will be
important in Section 3.

\subsection{Preliminary results}

The following proposition is an analog of Bertini's theorem.

\begin{pr}
\label{rgeneral}
Let $\ca$ be a reduced irreducible substack of $\cmbar_{g,m}(\proj^n,d)$,
and let $p$ be one of the labeled points.  Then there is a Zariski-open
subset $U$ of the dual projective space $(\proj^n)^*$ such that for all
$[H'] \in U$ the intersection $\ca \cap \{ \pi(p) \in H' \}$, if nonempty,
is reduced of dimension $\dim \ca - 1$.  
\end{pr}
Loosely, this result states that the requirement that a marked
point lie on a general hyperplane imposes one transverse condition on  an
irreducible substack of $\cmbar_{g,m}(\proj^n,d)$.  To prove this
proposition, we will invoke Theorem 2 of [Kl].

\begin{tm}[Kleiman]
Let $G$ be an integral algebraic group scheme, $X$ an integral algebraic
scheme with a transitive $G$-action.  Let $f:  Y \rightarrow X$ and $g: Z
\rightarrow X$ be two maps of integral algebraic schemes.  For each
rational point $s$ of $G$, let $sY$ denote $Y$ considered as an $X$-scheme
via the map $y \mapsto sf(y)$.
\begin{enumerate}
\item[(i)]  Then, there exists a dense open subset $U$ of $G$ such that,
for each rational point $s$ in $U$, either the fibered product,
$(sY)\times_X Z$ is empty or it is equidimensional and its dimension is
given by the formula,
$$
\dim(( sY) \times_X Z) = \dim(Y) + \dim(Z) - \dim(X).
$$
\item[(ii)]  Assume the characteristic is zero, and $Y$ and $Z$ are
regular.  Then, there exists a dense open subset $U$ of $G$ such that, for
each rational point $s$ in $U$, the fibered product $(sY) \times_X Z$, is
regular.
\end{enumerate}
\end{tm}

The proof of Kleiman's theorem carries through without change if $Z$ is an
algebraic stack.

\noindent {\em Proof of Proposition \ref{rgeneral}.}  
Let $G = PGL(n)$, $X= \proj^3$.  Let $Y$ be a hyperplane of $X$ with $f: Y
\rightarrow X$ the immersion.  Let $Z$ be the smooth points of $\ca$, with
$g:  Z \rightarrow X$ given by evaluation at $p$.  Then the result follows
immediately from Kleiman's theorem.
\epf

The next proposition is a variation of Proposition \ref{rgeneral}.  
\begin{pr}
\label{rgeneralproper}
With the hypotheses of Proposition \ref{rgeneral}, let $\cb$ be a proper
closed substack of $\ca$.  Then there is a Zariski-open subset 
$U$ of the dual projective space $(\proj^n)^*$ such that for all
$[H'] \in U$, each component of $\cb \cap \{ \pi(p) \in H' \}$ is a proper
closed substack of a component of $\ca   \cap \{ \pi(p) \in H' \}$.
\end{pr}
\bpf
The components of $\ca   \cap \{ \pi(p) \in H' \}$ are each of dimension
$\dim \ca - 1$ (by Proposition \ref{rgeneral}), and the components of 
$\cb \cap \{ \pi(p) \in H' \}$ are each of dimension less than $\dim \ca -
1$ (by Proposition \ref{rgeneral} applied to the irreducible components of
$\cb$).  
\epf

We now summarize the results of the remainder of this section.  There are
many relations among the various spaces $X(\ce)$ as $\ce$ varies.  Some are
universal curves over others (Proposition \ref{runiversal}).  Some are the
intersections of others with a divisor (Proposition \ref{rgeneral2}).  The
variety $X(\ce)$ can be identified with $Y(\ce(0);\ce(1))$ for
appropriately chosen $\ce(0)$ and $\ce(1)$ (Proposition \ref{rXY}).  A
smoothness result (Proposition \ref{rbig}) allows us to use results about
stable maps to $\proj^1$ proved in the previous section.  A first
application of Propositions \ref{runiversal}, \ref{rgeneral2}, and
\ref{rbig} is a calculation of the dimension of $X(\ce)$ and $Y(\ce(0);
\dots; \ce(l))$ (Proposition \ref{rdimX}); more will follow in subsequent
sections.  Finally, Proposition \ref{rXnice} ensures that the image of the
stable map corresponding to a general point of $X(\ce)$ is smooth.

Let $A$ be a general $(n-2)$-plane in $H$.  The projection $p_A$ from $A$
induces a rational map $\rho_A:
\cmbar_{0,m}(\proj^n,d) \dashrightarrow \cmbar_{0,m}(\proj^1,d)$, that is a
morphism (of stacks) at points representing maps $(C, \{ p_i \},
\pi)$  whose image $\pi(C)$ does not intersect $A$.  Via
$\cmbar_{0,m}(\Bl_A \proj^n,d)$, the morphism can be extended over the set
of maps $(C, \{ p_i \}, \pi)$ where $\pi^{-1} A$ is a union of reduced
points distinct from the $m$ marked points $\{ p_i \}$.  The image of such curves in
$\cmbar_{0,m}(\proj^1,d)$ is a stable map
$$
(C \cup C_1 \cup \dots \cup C_{\# \pi^{-1} A }, \{  p_i \}, \pi')
$$
where $C_1$, \dots, $C_{\# \pi^{-1} A }$ are rational tails attached to $C$
at the points of $\pi^{-1} A$, 
$$
\pi' \mid_{ \{ C \setminus \pi^{-1} A \} } = ( p_A
\circ \pi ) \mid_{ \{ C \setminus \pi^{-1} A \} }
$$
(which extends to a morphism from all of $C$) and $\pi' \mid_{C_k}$ is a
degree 1 map to $\proj^1$ ($1 \leq k \leq \# \pi^{-1} A
$).

\begin{pr}
If $(C, \{ p_i \}, \pi) \in \cmbar_{0,m}(\proj^n,d)$ and $\pi^{-1}
A$ is a union of reduced points disjoint from the marked points, then at
the point 
$(C,\{ p_i \}, \pi)$, $\rho_A$ is a smooth morphism of stacks of
relative dimension $(n-1)(d+1)$.
\label{rbig}
\end{pr}
\bpf
To show that a morphism of stacks $\ca \rightarrow \cb$ is smooth at a
point $a \in \ca$, where $\cb$ is smooth and $\ca$ is equidimensional, it
suffices to show that the fiber is smooth at $a$, or equivalently that the
Zariski tangent space to the fiber at $a$ is of dimension $\dim \ca - \dim
\cb$.

Recall that $\cmbar_{0,m}(\proj^n,d)$ is a smooth stack of dimension
$(n+1)d+m-1$ (see Subsubsection \ref{itmsosm}).  The first
order deformations of $(C, \{ p_i \},
\pi)$ in the fiber of $\rho_A$ can be identified with sections of the
vector bundle $\pi^*(\oh(1)^{n-1})$: these are deformations of a map
$$
(C, \{ p_i \},\pi)
 {\stackrel {(s_0, s_1, \dots, s_n)} \longrightarrow} \proj^n 
$$
keeping the marked curve $(C, \{ p_i \}, \pi)$ and the 
sections $(s_0,s_1)$ constant.
But
\begin{eqnarray*}
h^0(C, \pi^* \oh_{\proj^n}(1)^{n-1}) &=& (n-1) h^0(C, \pi^*
\oh_{\proj^n}(1)) \\
&=& (n-1) (d+1) \\
&=& \dim \cmbar_{0,m}(\proj^n,d) - \dim \cmbar_{0,m} (\proj^1,d)
\end{eqnarray*}
as desired.
\epf

The following two propositions give relationships among the spaces $X(\ce)$
as $\ce$ varies.

\begin{pr}
\label{runiversal}
Given $d$, $\vh$, $\vi$, $\Ga$, $\De$, let $\vec{i'} = \vi + \vep_n$, and
define $\De'$ to be the same as $\De$ except $\De^{i'_n}_n = \proj^n$.
Then $\cx(d, \Ga,
\De')$ is the universal curve over $\cx(d, \Ga, \De)$.
\end{pr}

\bpf  
The moduli stack $\cmbar_{0, \sum h_{m,e} + \sum i_e + 1} (\proj^n,d)$ is the
universal curve over $\cmbar_{0, \sum h_{m,e} + \sum i_e } (\proj^n,d)$
(see Subsubsection \ref{itmsosm}).  The
proposition is a consequence of the commutativity of the following diagram:
$$\begin{CD}
\cx(d, \Ga, \De')     @>>> \cmbar_{0, \sum h_{m,e} + \sum i_e + 1}(\proj^n,d) \\
@V{p}VV    @VV{p}V \\
\cx(d, \Ga, \De)   @>>> \cmbar_{0,\sum h_{m,e} + \sum i_e}(\proj^n,d)
\end{CD}$$
\epf

\begin{pr}
\label{rgeneral2}
Let $H'$ be a general hyperplane of $\proj^n$.  
\begin{enumerate}
\item[a)] The divisor 
$$
\{ \pi(p^{j_0}_{m_0,e_0}) \in H' \} \subset 
\cx(d, \Ga, \De)
$$
is $\cx(d, \Ga', \De)$ where
\begin{itemize}
\item $\vec{h'}=\vh - \vep_{m_0,e_0} +\vep_{m_0,e_0-1}$
\item For $(m,e) \neq (m_0,e_0), (m_0,e_0-1)$, ${\Ga'}^j_{m,e} = \Ga^j_{m,e}$.
\item $\{ {\Ga'}^j_{m_0,e_0} \}_j = \{ \Ga^j_{m_0,e_0} \}_j  \setminus \{
\Ga^{j_0}_{m_0,e_0} \}$, $\{ {\Ga'}^j_{m_0,e_0-1} \} = \{ \Ga^j_{m_0,e_0-1}
\}_j \cup \{ \Ga^{j_0}_{m_0,e_0} \cap H' \} $
\end{itemize}
\item[b)] The divisor 
$$
\{ \pi(q^{j_0}_{e_0}) \in H' \} \subset 
\cx(d, \Ga, \De)
$$ 
is $\cx(d, \Ga, \De')$ where
\begin{itemize}
\item $\vec{i'} = \vi - \vep_{e_0} + \vep_{e_0-1}$
\item For $e \neq e_0, e_0-1$, ${\De'}^j_{e} = \De^j_{e}$.
\item $\{ {\De'}^j_{e_0} \}_j = \{ \De^j_{e_0} \}_j  \setminus \{
\De^{j_0}_{e_0} \}$, $\{ {\De'}^j_{e_0-1} \}_j = \{ \De^j_{e_0-1} \}_j
\cup \{ \De^{j_0}_{e_0} \cap H' \} $ 
\end{itemize}
\end{enumerate}
\end{pr}
\bpf
We prove a) first.  Every point of 
$\{ \pi(p^{j_0}_{m_0,e_0}) \in H' \}$ represents a map where
$\pi(p^j_{m,e}) \subset \Ga^j_{m,e}$, $\pi(q^j_e) \subset \De^j_e$,
$\pi(p^{j_0}_{m_0,e_0}) \in H'$.  Clearly 
$$
\cx(d, \Ga', \De) \subset \{ \pi(p^{j_0}_{m_0,e_0}) \in H' \};
$$
each component of $\cx(d, \Ga', \De)$ appears with
multiplicity one by Proposition \ref{rgeneral}.  The only other possible
components of $\{ \pi(p^{j_0}_{m_0,e_0}) \in H' \}$ are those whose general
point represents a map where $\pi^{-1} H$ is not a union of points
(i.e. contains a component of $C$).  But such maps form a union of proper
subvarieties of components of $\cx(d,  \Ga, \De)$, and by
Proposition \ref{rgeneralproper} such maps cannot form a component of 
$$
\{ \pi(p^{j_0}_{m_0,e_0}) \in H' \} \cap \cx(d,  \Ga, \De).
$$

Replacing $p^{j_0}_{m_0,e_0}$ with $q^{j_0}_{e_0}$ in the previous paragraph gives
a proof of b).
\epf

The next observation is analogous to a well-known fact about the moduli
space of stable marked curves.  Consider the stack $\cmbar_{g,m}$ where the $m$
marked points are labeled $a_1$, $a_2$, $b_1$, \dots, $b_{m-2}$.  Then the
closed substack $\cv$ of $\cmbar_{g,m}$ parametrizing marked curves $A \cup
B$ with $a_i \in A$, $b_i \in B$, $p_a(A) = 0$, $p_a(B) = g$ is isomorphic
to $\cmbar_{g,m-1}$ where the $m-1$ marked points are labeled $c$, $b_1$,
\dots, $b_{m-2}$.  The isomorphism $\cmbar_{g,m-1} \rightarrow \cv$
involves gluing a rational tail (with marked points $a_1$, $a_2$) at $c$.

Fix $\ce$, integers $m_0$, $e_0$, $e_1$, and general $\Ga$, $\De$.  Let
$j_0 = h_{m_0,e_0}$, $j_1 = i_{e_1}$, $e' = e_0 + e_1 - n$, and $j' =
h_{m_0,e'} + 1$.  There is a subvariety $Y$ of $X(d,\Ga,\De)$ where $\pi(
p^{j_0}_{m_0,e_0}) = \pi(q^{j_1}_{e_1})$.  The general point of $Y$
represents a map $\pi: C(0) \cup C(1)
\rightarrow \proj^n$ where $C(0)$ and $C(1)$ are both isomorphic to
$\proj^1$, $\pi$ collapses $C(0)$ to a point, $p^{j_0}_{m_0,e_0}$ and
$q^{j_1}_{e_1}$ are on $C(0)$, and the rest of the marked points are on
$C(1)$.  Necessarily $\pi(C(0)) \subset \Ga^{j_0}_{m_0,e_0} \cap
\De^{j_1}_{e_1}$.  Such maps form a dense open subset of
$Y(d(0),\Ga(0),\De(0); d(1),\Ga(1),\De(1))$ where
\begin{itemize}
\item $(d(0),\vh(0),\vi(0)) = (0,\vep_{m_0,e_0},\vep_{e_1})$, $\Ga(0)
= \{ \Ga^{j_0}_{m_0,e_0} \}$, $\De(0) = \{ \De^{j_1}_{e_1} \}$.
\item $\ce(1) = \ce - \ce(0)$, $\Ga(1) = \Ga \setminus \Ga(0)$, $\De(1) =
\De \setminus \De(0)$.  
\end{itemize}

Now let $d' = d$, $\vec{h'} = \vh(1)
+ \vep_{m_1,e'}$, $\vec{i'} = \vi(1)$, $\Ga' = \Ga(1) \cup \{
\Ga^{j_0}_{m_0,e_0} \cap \De^{j_1}_{e_1} \}$, and $\De' =
\De(1)$.  The stable map corresponding to a general point of $Y$ can also be
identified with a stable map $(C(1), \{ p^j_{m,e} \}, \{ q^j_e \}, \pi)$ in
$X(d',\Ga',\De')$ where $C(1)$ is smooth and not contained in
$\pi^{-1}H$, by attaching a rational tail $C(0)$ (with two marked points
$p^{j_0}_{m_0,e_0}$ and $q^{j_1}_{e_1}$) at the point $p^{j'}_{m_0,e'}$ of
$C(1)$.  

In this way we get an isomorphism of $X(d',\Ga',\De')$ with 
$$
Y(d(0),\Ga(0),\De(0);d(1),\Ga(1),\De(1)):
$$

\begin{pr}
There is a natural isomorphism
$$
\phi:  X(d', \Ga', \De') \rightarrow Y(d(0), \Ga(0), \De(0);
d(1),\Ga(1),\De(1)). 
$$
\label{rXY}
\end{pr}
\bpf
The points in a dense open set of $X(d', \Ga', \De')$ represent
degree $d$ stable maps $\pi$ from a smooth curve $C$ to $\proj^n$ with
incidences $\pi({p'}^j_{m,e}) \in \Ga^j_e$, $\pi( {q'}^j_e) \in \De^j_e$
and an equality of divisors $\pi^* H = \sum m {p'}^j_{m,e}$ on $C$.  To
each such map, consider the map to $\proj^n$ where the marked point
${p'}^j_{m,e}$ is replaced by $p^j_{m,e}$ for $(m,e,j) \neq (m_0, e',
h'_{m_0,e'})$, ${q'}^j_e$ is replaced by $q^j_e$, and
$p^{h'_{m_0,e'}}_{m_0,e'}$ is replaced by a rational tail with additional
marked points $p^{j_0}_{m_0,e_0}$ and $q^{j_1}_{e_1}$.  The resulting
stable maps corresponds to points in a dense open set of $Y(d(0), \Ga(0),
\De(0);
d(1),\Ga(1),\De(1))$.
\epf

\begin{pr}
\label{rdimX}
Every component of $X(\ce)$ is reduced of dimension
$$
(n+1) d + (n-3) - \sum_{m,e} (n+m-e-2)h_{m,e} - \sum_e
(n-1-e)i_e.
$$
The general element of each component is (a map from) an
irreducible curve.

If $\sum_{k=0}^l \ce(k) = \ce$, then every component of $Y(\ce(0); \dots;
\ce(l)) $ is reduced of dimension $\dim X(\ce) - 1$.
\end{pr}
\bpf
We will prove the result about $\dim X(\ce)$ in the special case $\vi =
\vec{0}$ and $h_{m,e}=0$ when $e<n-1$.  Then the result holds when $\vi =
i_n \vep_n$ by Proposition \ref{runiversal} (applied $i_n$ times), and we can invoke Proposition
\ref{rgeneral2} repeatedly to obtain the result in full generality.  (This
type of reduction will be used often.)  In this special case, we must
prove that each component of $X(\ce)$ is reduced of dimension
$$
(n+1) d + (n-3) - \sum_m (m-1) h_{m,n-1}.
$$
The natural map $\cx(\ce) \dashrightarrow \cx_1(\cehat)$ induced by 
$$
\rho_A: \cmbar_{0,\sum h_{m,n-1}}(\proj^n,d) \dashrightarrow \cmbar_{0,\sum
h_{m,n-1}}(\proj^1,d)
$$ 
is smooth of relative dimension $(n-1)(d+1)$ at a
general point of any component of $\cx(\ce)$ by Proposition
\ref{rbig}.  The stack $\cx_1(\cehat)$ is reduced of dimension $2d-1 - \sum(m-1)
h_{m,n-1}$ by Subsection \ref{ikey}, so $\cx(\ce)$ is
reduced of dimension 
$$
(n-1)(d+1) + \dim  \cx_1(\cehat) = (n+1)d + (n-3) - \sum_m (m-1) h_{m,n-1}
$$
as desired.
As the general element of $\cx_1(\cehat)$ is (a map from) an irreducible
curve, the same is true of $\cx(\ce)$, and thus $X(\ce)$.

The same argument works for $Y$, as in Subsection
\ref{ikey} it was shown that $\dim Y_1(\cehat) = \dim X_1 (\cehat) - 1$.
\epf 

The following proposition is completely irrelevant to the rest of the
argument.  It is included to ensure that we are actually counting what we
want.

\begin{pr}
\label{rXnice}
If $n \geq 3$, and $(C, \{ p^j_{m,e} \}, \{ q^j_e \}, \pi)$ is the stable
map corresponding to a general point of a component of $X(\ce)$, then $C
\cong
\proj^1$ (with distinct marked points), and $\pi$ is a closed immersion.
\end{pr}
\bpf
By Propositions \ref{runiversal} and \ref{rgeneral2} again, we can assume
$\vi = \vec{0}$ and 
$h_{m,e}=0$ when $e<n-1$.  By the previous proposition, 
the curve $C$ is irreducible.  We need only check that $\pi$ is a closed
immersion.   The line bundle $\oh_C(d)$ is very ample, so a given non-zero
section $s_0$ and three general sections $t_1$, $t_2$, $t_3$ will separate
points and tangent vectors.  If $\pi = (s_0, s_1, s_2,s_3, \dots)$ then the
infinitesimal deformation $(s_0, s_1 + \varepsilon t_1, s_2 + \varepsilon t_2,
s_3 + \varepsilon t_3, s_4, \dots)$ will separate points and tangent vectors and still
lie in $X(\ce)$.  As
$(C,\{ p^j_{m,e} \}, \{ q^j_e \},\pi)$ corresponds to a general point in
$X(\ce)$, the map $\pi$ must be an immersion at this point.
\epf

\subsection{Degenerations set-theoretically}
\label{rdst}
Fix $\ce=(d,\vh,\vi)$ and a non-negative integer $E$, and let $\Ga$ and $\De$ be sets
of general linear spaces of $\proj^n$ (as in the definition of
$X(d,\Ga,\De)$).  Let $q$ be the marked point corresponding to one of the
(general) $E$-planes $Q$ in $\De$.

Let $D_H = \{ \pi(q) \in H \}$ be the divisor on $X(\ce)$ that corresponds
to requiring $q$ to lie on $H$.  In this
section, we will determine the components of $D_H$.  That is, we will give a list of subvarieties, and show
that the components of $D_H$ are a subset of this
list.  In the next section, we will determine the multiplicity with which
each component appears.  In particular, we will
see that the multiplicity of each component is at least one, so each
element of the list is indeed a component of $D_H$.  

\begin{figure}
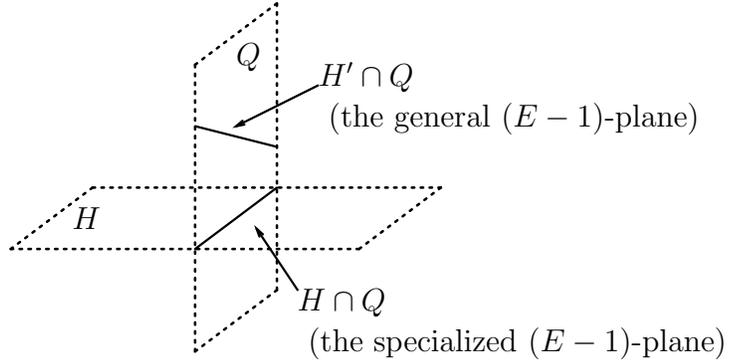

\begin{center}

	   \setlength{\unitlength}{.1\baseunit}
	    \input rspecialize.tex 
\end{center}
\caption{Specializing $H' \cap Q$ to lie in $H$}
\label{rspecialize}
\end{figure}

But first, let us relate this result to the enumerative problem we wish to
solve.  If $X(\ce^-)$ is a finite set of reduced points, we can determine
$\# X(\ce^-)$ by specializing one of the linear spaces of $\De$, of
dimension $E-1$, to lie in the hyperplane $H$ (see Figure \ref{rspecialize}).
Define $\ce$ by 
$$
(d, \vh, \vi) = (d^-, \vec{h^-}, \vec{i^-} + \vep_E - \vep_{E-1}).
$$
Then $X(\ce)$ is a dimension 1 variety by Proposition \ref{rdimX}. Let $q$
be the marked point on one of the $E$-planes $Q$ in $\De^-$.  The
conditions of $\ce$ are weaker than those of $\ce^-$: in $\ce$ we allow $q$
to lie on a linear space of dimension one more than in $\ce^-$.  Let
$D_{H'}$ be the divisor on $X(\ce)$ that corresponds to requiring $q$ to
lie on a fixed general hyperplane $H'$ in $\proj^n$.  Then the divisor
$D_{H'}$ on $X(\ce)$ is $X(\ce^-)$ by Proposition \ref{rgeneral2}.  If we
specialize $H'$ to $H$, $H' \cap Q$ will specialize to a general
$(E-1)$-plane $H \cap Q$ in $H$.  As $(D_H)
\sim (D_{H'})$ as divisor classes on the complete curve $X(\ce)$, $\deg D_H
=
\deg D_{H'}$.  So to calculate $\# X(\ce^-)$, we can simply enumerate the
points $D_H$ on $X(\ce)$, with the appropriate multiplicity.  Only
enumeratively meaningful divisors on $X(\ce)$ are relevant to such
enumerative calculations:  we are counting points on $X(\ce)$, which are
obviously enumeratively meaningful.

The components of $D_H$ on $X(\ce)$ are given by the following result.

\begin{tm}
\label{rlist1}
If $\Ga$ and $\De$ are general, 
each component of $D_H$ (as a divisor on $X(d,\Ga,\De)$) is a component of
$$
Y(d(0), \Ga(0), \De(0); \dots; d(l), \Ga(l), \De(l))
$$ 
for some 
$l$, $\ce(0)$, \dots, $\De(l)$, with $\ce = \sum_{k=0}^l \ce(k)$, 
$\Ga = \cup_{k=0}^l \Ga(k)$, $\De =
\cup_{k=0}^l \De(k)$, $Q \in \Ga(0)$. 
\end{tm}

\bpf
We may assume that $h_{m,e} = 0$ unless $e=n-1$, and that $\vi = \vep_n$
(and that $E=n$ and $q=q^1_n$).  The general case follows by adding more marked
points (Proposition \ref{runiversal}) and requiring each marked point to lie on
a certain number of general hyperplanes (Proposition \ref{rgeneral2}).

With these assumptions, the result becomes much simpler.  The stack
$\cx(d,\vh,\vi)$ is the universal curve over $\cx(d,\vh,\vec{0})$, and we
are asking which points of the universal curve lie in $\pi^{-1}H$.

Let $(C, \{ p^j_{m,n-1} \}, q, \pi)$ be the stable map corresponding to a
general point of a component of $D_H$.  Choose a general $(n-2)$-plane $A$
in $H$.  The set $\pi^{-1} A$ is a union of reduced points on $C$, so
by Proposition \ref{rbig} $\rho_A$ is smooth (as a morphism of stacks) at the point representing
$(C, \{ p^j_{m,n-1} \}, q, \pi)$ (by Proposition \ref{rbig}).  As a set,
$D_H$ contains the entire fiber of $\rho_A$ above $\rho_A(C, \{ p^j_{m,n-1}
\}, q, \pi)$, so $\rho_A(D_H)$ is a Weil divisor on $X_1(\cehat)$ that is a
component of $\{\pi(q) = z\}$ where $z = p_A(H)$.  By 
Theorem \ref{igenus0}, the curve $C$ is a union of irreducible components
$C(0) \cup \dots \cup C(l')$ with $\rho_A \circ \pi(C(0)) = z$
(i.e. $\pi(C(0)) \subset H$), $C(0) \cap C(k) \neq \phi$, and the marked
points split up among the components: $\vh = \sum_{k=0}^{l'} \vh(k)$.  If
$d(0) = \deg \pi |_{C(0)}$, then $d(0)$ of the curves $C(1)$, \dots,
$C(l')$ are rational tails that are collapsed to the $d(0)$ points of $C(0)
\cap A$; they contain no marked points.  Let $l = l' - d(0)$.  Also,
$\vi(k) = \vec{0}$ for $k>0$, as the only incidence condition in $\vi$ was
$q \in Q$, and $q \in C(0)$. 

Therefore this component of $D_H$ is
contained in 
$$
Y = Y(d(0),
\Ga(0), \De(0); \dots; d(l), \Ga(l), \De(l) ).
$$
  But $\dim Y = \dim X(\ce) -
1$ (by Proposition \ref{rdimX}), so the result follows.
\epf

For enumerative calculations, we need only consider enumeratively
meaningful components.  With this in mind, we restate Theorem \ref{rlist1}
in language reminiscent of [CH3].  Let $\phi$ be the isomorphism of
Proposition \ref{rXY}.    The following theorem will be more convenient for
computation. 

\begin{tm}
\label{rlist}
If $\Ga$ and $\De$ are general,
each enumeratively meaningful component of $D_H$ (as a divisor on
$X(d,\Ga,\De)$) is one of the following.
\begin{enumerate}
\item[(I)] A component of $\phi(X(d',\Ga', \De'))$, where, for some $m_0, e_0$,
$1 \leq j_0 \leq h_{m_0, e_0}$, $e' := e_0 + E - n \geq 0$:
\begin{itemize}
\item $d' = d$, $\vec{h'} = \vh - \vep_{m_0,e_0} + \vep_{m_0, e'}$, 
$\vec{i'} = \vi - \vep_E$
\item ${\Ga'}^j_{m,e} = \Ga^j_{m,e} \quad \text{if $(m,e) \neq (m_0,
e_0)$}$
\item $ \{ {\Ga'}^j_{m_0, e_0} \}_j = \{ \Ga_{m_0,e_0}^j \}_j \setminus \{
\Ga^{j_0}_{m_0,e_0} \}$
\item ${\Ga'}_{m_0, e'}^{h'_{m_0, e'}} =
\Ga^{j_0}_{m_0,e_0} \cap Q$ 
\item ${\De'}^j_e = \De^j_e$ if $e \neq E$, and $\{ {\De'}_E^j \}_j = \{
\De^j_E \}_j \setminus \{ Q \}$. 
\end{itemize}
\item[(II)]  A component of
$Y(d(0),\Ga(0),\De(0);\dots;d(l),\Ga(l),\De(l))$ for some
$l$, $\ce(0)$, \dots, $\De(l)$, with $\ce = \sum_{k=0}^l \ce(k)$, 
$\Ga = \cup_{k=0}^l \Ga(k)$, $\De =
\cup_{k=0}^l \De(k)$, $Q \in \Ga(0)$, and $d(0)>0$. 
\end{enumerate}
\end{tm}

Call these components {\em Type I components} and {\em Type II components}
respectively. 

\bpf
Consider a component $Y$ of $D_H$ that is not a Type II component (so $d(0)
= 0$).  Let $\{ C(0) \cup \dots \cup C(l), \{ p^j_{m,e} \}, q, \pi \}$ be
the stable map corresponding to a general point of $Y$.  The curve $C(0)$
has at least 3 special points: $q$, one of $\{ p^j_{m,e}\}$ (call it
$p^{j_0}_{m_0,e_0}$), and $C(0) \cap C(1)$.  If $C(0)$ had more than 3
special points, then the component would not be enumeratively meaningful,
due to the moduli of the special points of $C(0)$.  Thus $l=1$, and $Y$ is
a Type I component.
\epf

\subsection{Multiplicity calculations}
\label{rmultgen}
In Subsection \ref{rdst}, we saw that, in a neighborhood of a general point of
a component of the divisor $D_H$, there was a smooth morphism $\rho_A$ to
$X_1(\cehat)$ whose behavior at the corresponding divisor we understood
well.  For this reason, the multiplicity will be easy to calculate.  We
will see that the multiplicity with which the component $Y(\ce(0);
\dots; \ce(l))$ appears is $\prod_{k=1}^l m^k$, where $m^k = d(k) -
\sum_{m,e} m h_{m,e}(k)$ as defined earlier.  As usual, Propositions
\ref{runiversal} and \ref{rgeneral2} allow us to assume that $h_{m,e} = 0$
unless $e=n-1$, and $\vi = \vep_n$.


Recall that $X_1(\cehat)$ is the closure of
the (locally closed) subvariety of $\mb(\proj^1, d)$ parametrizing stable
maps from pointed rational curves with points $\{ p^j_{m,n-1} \}_{1 \leq j
\leq h_m}$, $q$ such that $\pi^* z = \sum_m m p_{m,n-1}^j$.  By (\ref{idimX}) (or a quick count of ramification points
away from $z$), $\dim X_1(\cehat) = d-1 + \sum_m h_{m,n-1}$.

Consider $Y_1 = Y_1(\cehat(0),\dots,\cehat(l),\clehat(1),\dots,\clehat(d(0)))$
where $\clehat(i) = (1,\vec{0},\vec{0})$.  The general point of $Y$ is a map
from a tree of rational curves $A^0$, $A^1$, \dots, $A^l$, $B^1$, \dots,
$B^{d(0)}$ with $A^0$ intersecting the other components and mapping to
$z$, $B^1, \dots, B^{d(0)}$ mapping to $\proj^1$ with degree 1, and $A^k$
mapping to $\proj^1$ with degree $d(k)$ (for $k>0$), with
$$
(\pi\mid_{A^k}) ^* z  = \sum_m \left( \sum_j m p^j_{m,n-1}(k) \right)
 + m^k (A^k \cap A_1).
$$
The rest of the marked
points (including $q$) are on $A^0$.  By Theorem \ref{igenus0}, $Y_1$ is a Weil divisor on $X_1(\cehat)$.

Let $D$ be the Cartier divisor on $X_1(\cehat)$ defined by $\{ \pi(q) = z
\}$.


Choose a general point $(C, \{ p^j_{m,n-1} \}, q,
\pi)$ of our component.  Note that that the image $\rho_A(C, \{
p^j_{m,n-1} \}, q,  \pi)$ is a general point of $Y_1$, with $A^k = C(k)$
(for $0 \leq k \leq l$).  The additional components $B^1$, \dots,
$B^{d(0)}$ come from the $d(0)$ intersections of 
$C(0)$ with the general $(n-2)$-plane $A$ of $H$.

\begin{lm}
In a neighborhood of $(C, \{ p^j_{m,n-1} \}, q,\pi)$:
\begin{enumerate}
\item[(a)] $\cmb(\proj^n,d) \rightarrow \cmb(\proj^1,d)$ is a smooth morphism of smooth stacks
or algebraic spaces.  Equivalently, the morphism is smooth on the level of
deformation spaces.
\item[(b)] The diagram  
$$\begin{CD}
\cx(\ce)     @>>> \cmb(\proj^n,d) \\
@V{p}VV    @VV{p}V \\
\cx_1(\cehat)   @>>> \cmb(\proj^1,d)
\end{CD}$$
is a fiber square.
\item[(c)] The component $\cy(\ce(0); \dots; \ce(l))$ is $p^*
\cy_1$ on $\cx(\ce)$. 
\item[(d)] As Cartier divisors, $p^* D = D_H$.
\end{enumerate}
\end{lm}

\bpf
Part (a) is Proposition \ref{rbig}.  Both (b) and (c) are clearly true
set-theoretically, and the fact that they are true stack-theoretically
follows from (a).  Part (d) is clear: the divisor $D_H$ is $\{ \pi(q) \in H
\}$, 
and the divisor $D$ is $\{ p_A(\pi(q)) = z \}$, where $p_A$ is the
projection from $A$ in $\proj^n$.
\epf

Combining these four statements with Theorem \ref{igenus0} and
Corollary \ref{ilocalst}, we have:
\begin{tm}
\label{rmult2}
The multiplicity of $D_H$ along the component $Y' = Y(\ce(0);
\dots; \ce(l))$ is the multiplicity of $D$ along $Y$, which is
$\prod_{k=1}^l m^k$.  In an 
\'{e}tale or formal neighborhood of a general point of $Y'$, $X(\ce)$ is
isomorphic to
$$
\Spec \com[[a, b_1,\dots,b_l,c_1,\dots,c_{\dim X(\ce) -1 } ]] / ( a =  b_1^{m^1} = \dots =
b_l^{m^l})
$$
where $D_H$ is given by $a=0$.
\end{tm}
Thus if $\lambda = \lcm(m^1,\dots,m^l)$, then $X(\ce)$ has $\prod m^k /
\lambda$ distinct reduced branches in an \'{e}tale neighborhood of a
general point of $Y$, all smooth 
if and only if $\lambda = m^k$ for some $k$.

\subsubsection{Multiplicity of $D_H$ along Type I components}
\label{rmultI}
Recall that a Type I component parametrizes those stable maps in $X(\ce)$
where one of the marked points $p^{j_0}_{m_0,e_0}$ is mapped to the
linear space $Q$; call this component $Z =
Z(m_0,e_0,j_0)$.  By the above argument, 
$Z$ appears with multiplicity $m_0$.  But the following
argument is more direct.

By Propositions \ref{runiversal} and \ref{rgeneral2}, we may assume
$h_{m,e} = 0$ 
unless $e=n-1$, and $\vi = \vep_n$.  The stack $\cx(\ce)$ is the universal
curve over $\cx(d,\vh,\vec{0})$, and the Type I component $Z(m_0,e_0,j_0)$
corresponds to the section $p^{j_0}_{m_0,e_0}$ of the universal curve.  On
the general fiber $C$ of the family $\cx(\ce) \rightarrow
\cx(d,\vh,\vec{0})$, $D_H = \sum m p^j_{m,e}$.  Hence $D_H$
contains $Z(m_0,e_0,j_0)$ with multiplicity $m_0$.

\subsection{Recursive formulas}
\label{rrecursive}

\subsubsection{The enumerative geometry of $Y$ from that of $X$}

Now that we inductively understand the enumerative geometry of varieties of
the form $X(\ce)$, we can compute $\#Y(\ce(0); \dots; \ce(l))$.

\begin{figure}
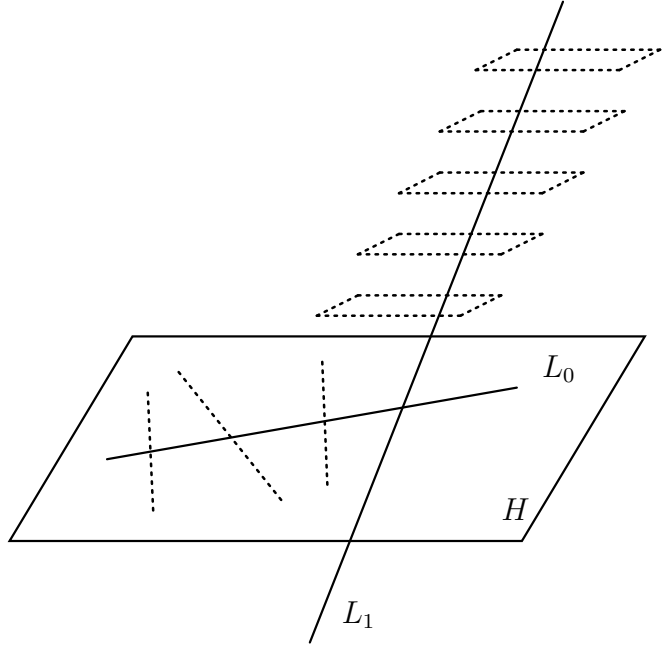

\begin{center}

	   \setlength{\unitlength}{.1\baseunit}
	    \input rYeg.tex 
\end{center}
\caption{How many $(L_0,L_1)$ satisfy the desired conditions?}
\label{rYeg}
\end{figure}

The method can be seen through a simple example.  Fix a hyperplane $H
\subset \proj^4$.  In $\proj^4$ the number
of ordered pairs of lines $(L_0, L_1)$ consisting of lines $L_0 \subset H$
and $L_1 \subset \proj^4$, with $L_0$ intersecting 3 fixed general lines
$a_1$, $a_2$, $a_3$ in $H$, $L_1$ intersecting 5 fixed general 2-planes
$b_1$, \dots, $b_5$ in $\proj^4$, and $L_0$ intersecting $L_1$ (see
Figure \ref{rYeg}), can be determined as follows.

There is a one-parameter family of lines $L_0$ in $H$ intersecting the
general lines $a_1$, $a_2$, $a_3$, 
and this family sweeps out a surface $S \subset H$ of
some degree $d_0$.  The degree $d_0$ is the number of lines $l_0$
intersecting the lines $a_1$, $a_2$, and $a_3$
{\em and another general line in $H$}, so this is $\#
X_3(\ce'(0))$ for $d'(0) = 1$, $\vec{h'}(0) = \vep_{1,2}$, $\vec{i'}(0) = 4
\vep_1$.  There is also a
one-parameter family of lines $L_1$ intersecting the general 2-planes $b_1$,
\dots, $b_5$, and
the intersection point of such $L_1$ with $H$ sweeps out a curve $C \subset
H$ of some degree $d_1$.  The degree $d_1$ is the number of lines
intersecting the 2-planes $b_1$, \dots, $b_5$ in $\proj^4$ {\em and another
general 2-plane in $H$}.  Thus $d_1 = \# X_4(\ce'(1))$ for $d'(1) = 1$,
$\vec{h'}(1) = \vep_{1,2}$, $\vec{i'}(1) = 5 \vep_1$.  The answer we seek
is $\# (C \cap S) = d_0 d_1$. 

The same argument in general yields:

\begin{pr}
\label{rrecursiveY}
$$
\# Y_n ( \ce(0); \dots; \ce(l) ) = \frac { \# X_{n-1} ( \ce'(0)) }{ d(0)!} \prod_{k=1}^l \# X_n(\ce'(k))
$$
where
\begin{itemize}
\item $h'(0) = d(0) \vep_{1,n-2}$
\item $i'_e(0) = i_{e+1}(0) + \# \{ \dim X_n(\ce(k)) = e \}_{1 \leq k \leq l} +
\sum_m h_{m,e}(0)$
\item for $1 \leq k \leq l$, $\vec{i'}(k) = \vi(k)$ and 
$\vec{h'}(k) = \vh(k) + \vep_{m^k, n-1-\dim X_n(\ce(k))}$.
\end{itemize}
\end{pr}
The $d(0)!$ is included to account for the $d(0)!$ possible labelings of
the intersection points of a degree $d(0)$ curve in $H$ with a
fixed general hyperplane $H'$ of $H$.

The following result is trivial but useful.

\begin{pr}
\label{rtrivial}
If the data $\ce'$ is the same as $\ce$ except $\vec{i'} = \vi+
\vep_{n-1}$, then $\# X(\ce') = d \cdot \# X(\ce)$.
\end{pr}
\bpf
The stable maps in $X(\ce')$ are just the stable maps in $X(\ce)$
along with a marked point mapped to a fixed general hyperplane.  There
are $d$ choices of this marked point.  (This is analogous to the
divisorial axiom for Gromov Witten invariants, cf. [FP] p. 35 (III).)
\epf

We now summarize the results of Subsections \ref{rdst} and 
\ref{rmultgen}.  Along with Propositions \ref{rrecursiveY} and
\ref{rtrivial}, this will give an algorithm to compute $\# X(\ce)$ for any
$\ce$.  (Proposition
\ref{rtrivial} isn't strictly necessary, but will make the algorithm
faster.)  The only initial data
needed is the ``enumerative geometry of $\proj^1$'': the number of stable
maps to $\proj^1$ of degree 1 is 1.  

Given $\ce$, fix an $E$ such that $i_E>0$.  Partitions of $\ce$
are simultaneous partitions of $d$, $\vh$, and $\vi$.  Define
multinomial coefficients with vector arguments as the product of the
multinomial coefficients of the components of the vectors:
$$
 \binom {\vh }{ { \vh(0), \dots, \vh(l)} } = \prod_{m,e} \binom {
h_{m,e} }{ {h_{m,e}(0), \dots, h_{m,e}(l)}}, 
$$
$$
\binom { \vi }{ { \vi(0), \dots, \vi(l)} } = \prod_e \binom {
i_e }{{i_e(0), \dots, i_e(l)}}.
$$
Define $\ce^-$ by $(d^-, \vec{h^-}, \vec{i^-}) = ( d, \vh, \vi - \vep_E +
\vep_{E-1})$, and let $\Ga^- = \Ga$ and 
$\De^- = \De \cup \{ \De^{i_E}_E \cap H' \}  \setminus \{ \De^{i_E}_E \}$
where $H'$ is a general hyperplane.  (This notation was used earlier, in
Subsection \ref{rdst}.)

\begin{tm}
\label{rrecursiveX1}
In $A^1(X_n(d, \Ga, \De))$, the cycle 
$X_n(d^-, \Ga^-, \De^-)$ is rationally equivalent to 
$$
\sum \left( \prod_{k=1}^l m^k \right) Y_n(d(0),\Ga(0),\De(0); \dots;
d(l),\Ga(l),\De(l)) 
$$
where the sum is over all $l$, 
$\ce(0)$, \dots, $\De(l)$, with $\ce = \sum_{k=0}^l \ce(k)$, 
$\Ga = \coprod_{k=0}^l \Ga(k)$, $\De =
\coprod_{k=0}^l \De(k)$, $\Ga_{E}^{i_E} \in \Ga(0)$. 
\end{tm}
\bpf
The left side is rationally equivalent (in $\cx(\ce)$) to $D_H = \{
\pi(q^{i_E}_E) \in H \}$.  The right side is set-theoretically $D_H$ by
Theorem \ref{rlist1}, and the multiplicity $\prod m^k$ was determined in
Subsection \ref{rmultgen}.
\epf

If $\# X_n(\ce^-)$ is finite, the following statement is more suitable for computation.

\begin{tm}
\label{rrecursiveX2}
$$
\# X_n(\ce^-) = \sum_{m,e} m h_{m,e} \cdot \# X_n(\ce'(m,e)) 
$$
$$
+ \sum \left(
\prod_{k=1}^l m^k \right) \binom { {\vh} }{ { \vh(0), \dots, \vh(l)} }
\binom {
{\vi - \vep_E} }{ { \vi(0)-\vep_E,\vi(1), \dots,
\vi(l)} } 
$$
$$
\quad \quad \quad \cdot \frac { \# Y_n(\ce(0); \dots; 
\ce(l)) }{ \Aut( \ce(1),\dots,\ce(l))}
$$
where, in the first sum, $\ce'(m,e) = (d,\vh -\vep_{m,e}+\vep_{m,e+E-n}, \vi -
\vep_E)$; the second sum is over all $l$ and all partitions 
$\ce(0)$,\dots, $\ce(l)$ of $\ce$ with $d(0)>0$.
\end{tm}
This follows from Theorem \ref{rlist} and the multiplicity calculations of
Subsection \ref{rmultgen}.  The only new points
requiring 
explanation are the combinatorial aspects: the $h_{m,e}$ in the first sum,
and the ``$\Aut$'' and various multinomial coefficients in the second.
In Theorem \ref{rlist}, the Type I components were indexed by
$(m_0, e_0, j_0)$.  But for fixed $(m_0,e_0)$, $\# X(\ce'(m_0,e_0))$ is
independent of $j_0$, so the above formula eliminates this redundancy.
Similarly, in Theorem \ref{rlist}, the Type II components were indexed by
partitions of the points $\{ p^j_{m,e} \}_{m,e,j}$ and $\{ q^j_e \}_{e,j}
\setminus \{ q \}$, but the value of $\# Y_n(\ce(0); \dots; \ce(l))$
depends only on $\{ \vh(k), \vi(k) \}_{k=0}^l$ and not on the actual
partitions.  The multinomial coefficients in the second line eliminate
this redundancy.  Finally, we divide the last term by $\Aut(\ce(1);
\dots; \ce(l))$ to ensure that we are counting each Type II component
once.  

\subsubsection{Transposing these results to subvarieties of the Hilbert scheme}

Our original result was an equality of divisors on $X(\ce)$.  We will
briefly sketch the analogous equality in the Chow ring of the Hilbert
scheme.  

Assume for convenience that $n \geq 3$, and that $h_{m,e} = i_e = 0$ when
$e > n-2$.  By Proposition \ref{rXnice}, there is a dense open subset $U$ of
$X(\ce)$ such that the image of the map corresponding to a point on $U$ is
smooth.  We can take a smaller $U$ such that the images of the
corresponding maps intersect each $\Ga^j_{m,e}$ and $\De^j_e$ in one point.

Define the closed subscheme $X^{\Hilb}(\ce)$ of the Hilbert scheme to be the
closure of the points $U^{\Hilb}$ representing the images of the maps
corresponding to points of $U$.  The subvarieties $Y^{\Hilb}(\ce(0); \dots;
\ce(l))$ can be defined analogously.  There is a rational map
$$
\psi:  X^{\Hilb}(\ce) \dashrightarrow X(\ce)
$$ that restricts to an isomorphism from $U^{\Hilb}$ to $U$.  (The map
$\psi$ is the inverse of the rational map $\xi$ defined in Subsubsection
\ref{rgeomean}.)  Let $\Phi_1$,
$\Phi_2$ be the projection of the graph of $\psi$ to $X^{\Hilb}(\ce)$ and
$X(\ce)$ respectively.    The exceptional divisors of $\psi$ are defined to
be the image under $\Phi_1$ of the divisors on the graph collapsed by
$\psi_2$.  (It is not clear to the author if such divisors exist.)
Then Theorem \ref{rrecursiveX1} can be reinterpreted as follows.

\begin{tm}  
\label{rrecHilb}
In $A^1( X^{\Hilb}(d,\Ga,\De))$, modulo the exceptional divisors of
$\psi$,
$$
X^{\Hilb}(d^-,\Ga^-,\De^-) = \sum \left( \prod_{k=1}^l m^k \right)
Y^{\Hilb}(d(0),\Ga(0),\De(0); \dots; d(l),\Ga(l),\De(l))
$$
where the sum is over all
$l$, $\ce(0)$, \dots, $\De(l)$, with $\ce = \sum_{k=0}^l \ce(k)$, 
$\Ga = \coprod_{k=0}^l \Ga(k)$, $\De =
\coprod_{k=0}^l \De(k)$, $\Ga_{E}^{i_E} \in \Ga(0)$. 
\end{tm}

This result follows from Theorem \ref{rrecursiveX1} and the multiplicity
calculations of Subsection \ref{rmultgen}.

\section{Elliptic Curves in Projective Space}
\label{elliptic}


In this section, we extend our methods to study the geometry of varieties
$W(\ce)$ parametrizing degree $d$ elliptic curves in $\proj^n$ intersecting
fixed general linear spaces and tangent to a fixed hyperplane $H$ with
fixed multiplicities along fixed general linear subspaces of $H$.  We use
the same general ideas as in the preceding section: we work with the
variety $\mbar_{1,m}(\proj^n,d)$ (and the stack $\cmbar_{1,m}(\proj^n,d)$)
and specialize linear spaces (which the curve is required to intersect) to
lie in $H$ one at a time.  Many arguments will carry over wholesale.  The
main additions deal with new types of degenerations.

\subsubsection{Example:  Cubic elliptic space curves}
\label{ecubics}
\begin{figure}
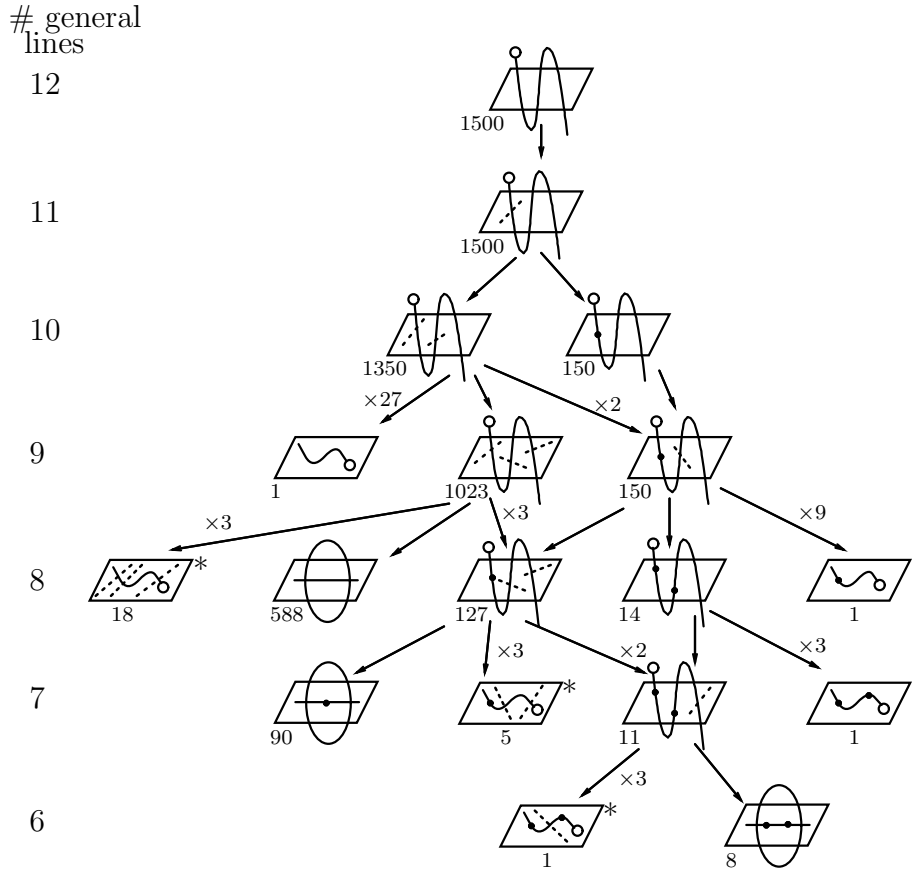

\begin{center}

	   \setlength{\unitlength}{.1\baseunit}
	    \input e1500cubics.tex 
\end{center}
\caption{Counting 1500 elliptic cubics through 12 general lines in $\proj^3$}
\label{e1500cubics}
\end{figure}

The example of smooth elliptic cubics in $\proj^3$ illustrates some of the
degenerations we will see, and shows a new complication.  There are 1500
smooth elliptic cubics in $\proj^3$ through 12 general lines, and we can
use the same degeneration ideas to calculate this number.  Figure
\ref{e1500cubics} is a pictorial table of the degenerations; smooth
elliptic curves are indicated by an open circle.

The degenerations marked with an asterisk have a new twist.  For example,
consider the cubics through 9 general lines $L_1$, \dots, $L_9$ and 3 lines
$L_{10}$, $L_{11}$, $L_{12}$ in $H$ (in row 9) and specialize $L_9$ to lie in
$H$.  The limit cubic could be a smooth plane curve in $H$ (the left-most
picture of row 8 in the figure). In this case, it must pass through the
eight points $L_1 \cap H$, \dots, $L_8 \cap H$.  But there is an additional
restriction.  The cubics (before specialization) intersected $L_{10}$, $L_{11}$,
$L_{12}$ in three points $p_{10}$, $p_{11}$, $p_{12}$ ($p_i \in L_i$), and as
elliptic cubics are planar, these three points must have been collinear.
Thus the possible limits are those curves in $H$ through $L_1 \cap H$,
\dots, $L_8 \cap H$ and passing through collinear points $p_{10}$, $p_{11}$,
$p_{12}$ (with $p_i \in L_i$).  (There is also a choice of a marked point of
the curve on $L_9$, which will give a multiplicity of 3.)  This
collinearity condition can be written as $\pi^*(\oh(1)) = p_{10} + p_{11} +
p_{12}$ in the Picard group of the curve.  The existence of such
degenerations is analogous to the divisorial condition of Theorem
\ref{igenus1}.  

We will have to count elliptic curves with such a divisorial condition
involving the marked points; this locus forms a divisor on a variety of the
form $W(\ce)$. Fortunately, we can express this divisor in terms of
divisors we understand well (Subsubsection \ref{eevalZ}).  Thus as a side
benefit, we get enumerative data about elliptic curves in $\proj^n$ with
certain incidence and tangency conditions, and a divisorial condition as
well.

\subsection{Notation and summary}

For convenience, let $\vep_e$, $\vep_{m,e}$ be the natural basis vectors:
$(\vep_e)_{e'} = 1$ if $e = e'$ and 0 otherwise; and $(\vep_{m,e})_{m',e'}
= 1$ if $(m,e)=(m',e')$, and 0 otherwise.  Fix a hyperplane $H$ in
$\proj^n$, and a hyperplane $A$ of $H$. 
From the previous section, recall the definitions of ``enumeratively
meaningful'', $X(\ce)$, $\cx(\ce)$, $Y(\ce(0);\dots;\ce(l))$, and
$\cy(\ce(0);\dots;\ce(l))$. 

Motivated by the analysis in Subsection
\ref{ikey} of divisors on subvarieties of
$\mbar_{1,m}(\proj^1,d)$, we define five new classes of varieties, labeled
$W$, $Y^a$, $Y^b$, $Y^c$, and $Z$ and corresponding stacks, labeled $\cw$,
$\cy^a$, $\cy^b$, $\cy^c$, and $\cz$.

\subsubsection{The schemes $W(\ce)$}

The objects of primary interest to us are smooth degree $d$ elliptic curves
in $\proj^n$ ($n \geq 2$) intersecting a fixed hyperplane $H$ with various
multiplicities along various linear subspaces of $H$, and intersecting
various general linear spaces in $\proj^n$.  We will examine these objects
as stable maps from marked curves to $\proj^n$ (where the markings will be
the various intersections with $H$ and incidences). 

Let $n$ and $d$ be positive integers, and let $H$ be a hyperplane in $\proj^n$.
Let $\vh = ( h_{m,e} )_{m \geq 1, e \geq 0}$ and $\vi = ( i_e )_{e \geq
0}$ be sets of non-negative integers.  Let 
$\Ga = \{ \Ga^j_{m,e} \}_{m,e,1 \leq j \leq h_{m,e}}$ be a set of linear
spaces in $H$ where $\dim 
\Ga^j_{m,e} = e$.   Let 
$\De = \{ \De^j_e \}_{e,1 \leq j \leq i_e}$ be a set of linear spaces in
$\proj^n$ where $\dim
\De^j_{e} = e$.

\begin{defn}
The scheme $W_n(d,\Ga,\De)$ is the (scheme-theoretic) closure of
the locally closed subset of $\mbar_{1,\sum h_{m,e} + \sum i_e}(\proj^n,d)$
(where the 
marked points are labeled $\{p_{m,e}^j\}_{1 \leq j \leq h_{m,e}}$ and $\{
q_e^j \}_{1 \leq j
\leq i_e}$) representing stable maps $(C, \{ p^j_{m,e} \}, \{
q^j_e \}, \pi)$ satisfying $\pi(p_{m,e}^j) \in \Ga_{m,e}^j$, $\pi(q_e^j) \in
\De_e^j$, $\pi^* H = \sum_{m,e,j} m p^j_{m,e}$, and where no components
of $C$ are collapsed by $\pi$.
\end{defn}

In particular, $\sum_{m,e} m h_{m,e} = d$, and no component of $C$ is
contained in $\pi^{-1}H$.  The incidence conditions define closed
subschemes of $\mbar_{1,\sum h_{m,e} + \sum i_e}(\proj^n,d)$, so the union
of these conditions indeed defines a closed subscheme of $\mbar_{1,\sum
h_{m,e} + \sum i_e}(\proj^n,d)$.

Define $\cw_n(d,\Ga,\De)$ in the same way as a substack of 
$\cmbar_{1,\sum h_{m,e} + \sum i_e}(\proj^n,d)$.  When we speak of propertis that are constant for general $\Ga$ and $\De$ (such as the dimension), we will write $W_n(d,\vh,\vi)$.  For convenience, 
write $\ce$ (for $\ce$verything) for the data $d, \vh, \vi$, so
$W_n(\ce) = W_n(d,\vh,\vi)$.  Also, the $n$ will often be suppressed
for convenience.

The variety $W(d, \Ga, \De)$ (analogous to $X(d,\Ga,\De)$ defined in
Section \ref{rational}) can be loosely thought of as parametrizing
degree $d$ elliptic curves in projective space intersecting certain linear
subspaces of $\proj^n$, and intersecting $H$ with different multiplicities
along certain linear subspaces of $H$.  For example, if $n=3$, $d=3$,
$h_{2,0}=1$, $h_{1,2}=1$, $W$ parametrizes elliptic cubics in $\proj^3$
tangent to $H$ at a fixed point.

In the special case where $h_{m,e}=0$ when
$e<n-1$ and $\vi=\vep_n$, define $\cehat$ by $\hat{d} = d$, $\hat{i}_1 =
1$, $\hat{h}_{m,0} = h_{m,n-1}$.  We will relate the geometry of
$\cw_n(\ce)$ to that of $\cw_1(\cehat)$, which was studied in Subsection
\ref{ikey}.  The geometry of $\cw_n(\ce)$ for general $\ce$ can be
understood from this special case.

If the linear spaces $\Ga$, $\De$ are general,
these varieties have the dimension one would naively expect.  The family of
degree $d$ elliptic curves in $\proj^n$ has dimension $(n+1)d$.
Requiring the curve to pass through a fixed $e$-plane should be a codimension
$(n-1-e)$ condition.  Requiring the curve to be $m$-fold tangent to $H$
along a fixed $e$-plane of $H$ should be a codimension $(m-1)+(n-1-e)$ condition.
Thus we will show (Theorem \ref{edimW}) that when the linear spaces in
$\Ga$, $\De$ are general, each component of $W(d, \Ga, \De)$ has dimension
$$
(n+1) d - \sum_{m,e} (n+m-e-2)h_{m,e} - \sum_e  (n-1-e)i_e.
$$
Moreover, $W(\ce)$ is reduced.  When the dimension is 0, $W(\ce)$ consists
of reduced points.  We call this number $\# W(\ce)$ --- these are the
numbers we want to calculate.  Define $\# W(\ce)$ to be zero if $\dim
W(\ce) > 0$.

For example, when $n=3$, $d=3$, $h_{1,2} = 3$, $i_1 = 12$, $W(\ce)$
consists of a certain number of reduced points: 3! times the number of
elliptic cubics through 12 general lines.  (The 3! arises from the markings
of the three intersections of the cubic with $H$.)

\subsubsection{The schemes $Y^a$,
$Y^b$, and $Y^c$} We
will be naturally led to consider subvarieties of $W(d,\Ga,\De)$ which
are similar in form to the varieties 
$$
Y(d(0), \Ga(0), \De(0); \dots; d(l), \Ga(l), \De(l) )
$$ 
of the previous section.
Fix $n$, $\ce$, $\Ga$, $\De$, and a non-negative integer $l$.  Let
$\sum_{k=0}^l d(k)$ be a partition of $d$.  Let the points $\{ p^j_{m,e}
\}_{m,e,j}$ be partitioned into $l+1$ subsets $\{ p^j_{m,e}(k) \}_{m,e,j}$
for $k= 0$, \dots, $l$.  This induces a partition of $\vh$ into
$\sum_{k=0}^l \vh(k)$ and a partition of the set $\Ga$ into
$\coprod_{k=0}^l \Ga(k)$.  
Let the points $\{ q^j_e
\}_{e,j}$ be partitioned into $l+1$ subsets $\{ q^j_{e}(k) \}_{e,j}$
for $k= 0$, \dots, $l$.  This induces a partition of $\vi$ into
$\sum_{k=0}^l \vi(k)$ and a partition of the set $\De$ into
$\coprod_{k=0}^l \De(k)$.    Define $m^k$ by $m^k = d(k) - \sum_m m
h_m(k)$, and assume $m^k>0$ for all $k = 1$, \dots, $l$.

\begin{defn}
The scheme 
$$
Y^a_n(d(0),\Ga(0),\De(0);
\dots; d(l),\Ga(l),\De(l))
$$
is the (scheme-theoretic) closure of
the locally closed subset of $\mbar_{1,\sum h_{m,e} + \sum i_e}(\proj^n,d)$ (where the
points are labeled $\{ p_{m,e}^j \}_{1 \leq j \leq h_{m,e}}$ and $\{ q_e^j
\}_{1 \leq j \leq i_e}$) representing stable maps $(C, \{ p^j_{m,e}
\}, \{ q^j_e \}, \pi)$ satisfying the following conditions
\begin{enumerate}
\item[Y1.] The curve $C$ consists of $l+1$ irreducible components $C(0)$, \dots,
$C(l)$ with all components intersecting $C(0)$.  The map $\pi$ has
degree $d(k)$ on curve $C(k)$ ($0 \leq k \leq l$).
\item[Y2.] The points $\{ p^j_{m,e}(k) \}_{m,e,j}$ and $\{ q^j_e(k)\}_{e,j}$ lie
on $C(k)$, and $\pi(p^j_{m,e}(k)) \in \Ga^j_{m,e}(k)$, $\pi(q^j_e(k)) \in
\De^j_e(k)$.  
\item[Y3.] As sets, $\pi^{-1}H = C(0) \cup \{ p^j_{m,e} \}_{m,e,j}$, and for $k>0$,
$$
( \pi \mid_{C(k)} )^* H = \sum_{m,e,j} m p^j_{m,e}(k) + m^k (C(0) \cap C(k)).
$$
\item[Y${}^{\text{a}}$4.] The curve $C(1)$ is elliptic and the other
components are rational.
\end{enumerate}
\end{defn}
Conditions Y1--Y3 appeared in the definition of $Y$ (Definition \ref{rdefY}). 
Note that $d(k)>0$ for all positive $k$ by condition Y3.

When discussing properties that hold for general $\{ \Ga^j_{m,e}
\}_{m,e,j}$, $\{ \De^j_e \}_{e,j}$, we will write 
$$
Y^a_n(\ce(0); \dots;
\ce(l)) =
Y^a_n(d(0),\vh(0),\vi(0); \dots;
d(l),\vh(l),\vi(l)) .
$$
  The $n$ will often be suppressed for convenience.  
If $\vh(k) + \vi(k) \neq \vec{0}$ for all $k>0$, 
$$
Y^a(\ce(0),\Ga(0),\De(0); \dots; \ce(l),\Ga(0),\De(0))
$$ 
is isomorphic to a closed subscheme of
$$
\mbar_{0,\sum h(0) + \sum i(0)+l}(H,d(0)) \times
W(d(1),\Ga(1),\De(1))
$$
$$
\times \prod_{k=2}^l 
X(d(k),\Ga'(k),\De(k)),
$$
where for $k= 1, \dots, l$, $\vec{h'}(k) = \vh(k) + \vep_{m^k,n-1}$ and
$\Ga'(k)$ is the same as $\Ga(k)$ except
$\Ga^{h_{m^k,n-1}+1}_{m^k,h_{m^k,n-1}+1} = H$.

Define $\cy^a(\ce(0), \De(0), \Ga(0); \dots; \ce(l), \De(l), \Ga(l))$
as the analogous stack. 

\begin{defn}
The scheme 
$$
Y^b_n(d(0),\Ga(0),\De(0);
\dots; d(l),\Ga(l),\De(l))
$$ 
is the (scheme-theoretic) closure of
the locally closed subset of $\mbar_{1,\sum h_{m,e} + \sum i_e}(\proj^n,d)$
(where the points are labeled $\{ p_{m,e}^j \}_{1 \leq j \leq h_{m,e}}$ and
$\{ q_e^j \}_{1 \leq j \leq i_e}$) representing stable maps $(C, \{
p^j_{m,e}  \}, \{ q^j_e \}, \pi)$ satisfying the conditions Y1--Y2 above,
and 
\begin{enumerate}
\item[Y${}^{\text{b}}$3.] As sets, $\pi^{-1}H = C(0) \cup \{ p^j_{m,e} \}_{m,e,j}$, and for $k>1$,
$$
( \pi \mid_{C(k)} )^* H = \sum_{m,e,j} m p^j_{m,e}(k) + m^k (C(0) \cap C(k)).
$$
\item[Y${}^{\text{b}}$4.] All components of $C$ are rational.  The curves
$C(0)$ and $C(1)$ intersect at two distinct points $\{ a_1, a_2 \}$.
(These points are not marked; monodromy may exchange them.)  Also,
$$
(\pi \mid_{C(1)})   ^* H  = \sum_{m,e}
\sum_{j=1}^{h^k_{m,e}} m p^j_{m,e} + m^1_1 a_1 + m^1_2 a_2 
$$
where $m^1_1 + m^1_2 = m^1$.  
\end{enumerate}
\end{defn}
Thus $Y^b(\ce(0); \dots; \ce(l) )$ is naturally the union of $[m^1/2]$
(possibly reducible) schemes (where $[\cdot]$ is the greatest-integer
function), indexed by $m^1_1$.  For convenience, label these varieties 
$$
\{
Y^b(\ce(0); \dots; \ce(l))_{m_1^1}\}_{1 \leq m^1_1 < m^1},
$$
so $\{ Y^b(\ce(0); \dots; \ce(l))_{m_1^1}\}_{m^1_1} =
\{ Y^b(\ce(0); \dots; \ce(l))_{m_1^1}\}_{m^1 - m^1_1}$.

For enumerative reasons, we define a slightly different variety.
\begin{defn}
The scheme 
$$
\tilde{Y}^b_n(d(0),\Ga(0),\De(0);
\dots; d(l),\Ga(l),\De(l))_{m^1_1} 
$$
is the
(scheme-theoretic) closure of 
the locally closed subset of the universal curve over $\mbar_{1,\sum h_{m,e} + \sum
i_e}(\proj^n,d)$ (where the 
points are labeled $\{ p_{m,e}^j \}_{1 \leq j \leq h_{m,e}}$ and $\{ q_e^j
\}_{1 \leq j \leq i_e}$, and the point on the universal curve is labeled
$a_1$) representing stable maps $(C, \{ 
p^j_{m,e} 
\}, \{ q^j_e \}, \pi)$ (with additional point $a_1$) satisfying the
conditions Y1, Y2, Y${}^{\text{b}}$3, and Y${}^{\text{b}}$4 above (for some
other point $a_2$).
\end{defn}
There is a morphism $\tilde{Y}^b_n(d(0), \dots, \De(l))_{m^1_1}
\rightarrow Y^b_n(d(0), \dots, \De(l))_{m^1_1}$
corresponding to forgetting the point $a_1$.  This morphism is an
isomorphism if $m^1_1 \neq m^1_2$ and it is generically two-to-one when
$m^1_1 = m^1_2$.

When discussing properties that hold for general $\{ \Ga^j_{m,e}
\}_{m,e,j}$, $\{ \De^j_e \}_{e,j}$, we will write 
$$
Y^b_n(\ce(0); \dots; \ce(l)) =Y^b_n(d(0),\vh(0),\vi(0); \dots; d(l),\vh(l),\vi(l))
$$ 
and 
$$
\tilde{Y}^b_n(\ce(0); \dots; \ce(l)) =\tilde{Y}^b_n(d(0),\vh(0),\vi(0); \dots; d(l),\vh(l),\vi(l))
$$
  The $n$ will often be suppressed for convenience.  If $\vh(k) +
\vi(k) \neq \vec{0}$ for all $k>0$, $\tilde{Y}^b(d(0),\Ga(0),\De(0);
\dots; d(l),\Ga(l), \De(l))$ is isomorphic to a closed subscheme of
$$
\mbar_{0,\sum h(0) + \sum i(0)+l+1}(H,d(0)) \times
\prod_{k=1}^l 
X(d(k),  \Ga'(k),\De(k))
$$
for appropriately chosen $\Ga'(k)$, $k = 1$, \dots, $l$.

Define $\cy^b(d(0), \De(0), \Ga(0); \dots; d(l), \De(l), \Ga(l))$
as the analogous stack. 

\begin{defn}
\label{eYcdef}
The scheme 
$$
Y^c_n(
d(0),\Ga(0),\De(0); \dots; d(l),\Ga(l),\De(l))
$$ 
is the (scheme-theoretic) closure of
the locally closed subset of $\mbar_{1,\sum h_{m,e} + \sum i_e}(\proj^n,d)$ (where the
points are labeled $\{ p_{m,e}^j \}_{1 \leq j \leq h_{m,e}}$ and $\{ q_e^j
\}_{ 1 \leq j \leq i_e}$) representing stable maps $(C, \{ p^j_{m,e}
\}, \{ q^j_e \}, \pi)$ satisfying conditions Y1--Y3, and
\begin{enumerate}
\item[Y${}^{\text{c}}$4.] The curve $C(0)$ is elliptic and the other
components are rational.  The morphism $\pi$ has positive degree on every
component.  
\item[Y${}^{\text{c}}$5.] 
In $\Pic(C(0))$,
\begin{eqnarray*}
\pi^*(\oh_{\proj^n}(1)) & \otimes & \oh_{C(0)} \left( \sum_{k=1}^l m^k (C(0)
\cap C(k)) \right) \\
& \cong &
\oh_{C(0)} \left( \sum_{m,e} \sum_{j=1}^{h_{m,e}(0)} m p_{m,e}^j(0) \right).
\end{eqnarray*}
\end{enumerate}
\end{defn}
The divisorial condition Y${}^{\text{c}}$5 is motivated by the ideas of
Subsection \ref{ikey}.  If $\ce = \sum_{k=0}^l \ce(k)$,
and $\Ga$ and $\De$ are general with $\Ga = \coprod_{k=0}^l \Ga(k)$, $\De =
\coprod_{k=0}^l \De(k)$, then the variety 
$$
Y^c_n(d(0), \Ga(0),\De(0); \dots; d(l),\Ga(l),\De(l))
$$
will turn out to be a Weil divisor on $W(d, \Ga, \De)$.  The stable
map $(C, \{ p^j_{m,e} \}, \{ q^j_e \}, \pi)$ corresponding to a general
point of $W(d,\Ga,\De)$ satisfies $\pi^* (\oh_{\proj^n}(1)) \cong \oh_C (
\sum_{m,e,j} m 
p^j_{m,e})$, and this condition must in some sense be inherited by 
the map corresponding to a general
point on the Weil divisor.

This condition was actually present in $Y^a$ and $Y^b$ (and the Type II
component $Y$ of the previous section), but as $C(0)$ was rational in each
of these cases, the
requirement reduced to
$$
d(0) +  \sum_{k=1}^l m^k =  \sum_{m,e} m h_{m,e}(0)
$$
which was always true.

When discussing properties that hold for general $\{ \Ga^j_{m,e}
\}_{m,e,j}$, $\{ \De^j_e \}_{e,j}$, we will write $Y^c_n(\ce(0); \dots;
\ce(l))$.  The $n$ will often be suppressed for convenience.  If $\vh(k) +
\vi(k) \neq \vec{0}$ for all $k>0$,  
$$
Y^c(d(0), \Ga(0), \De(0); 
\dots; d(l), \Ga(l), \De(l))
$$ 
is isomorphic to a closed subscheme of
$$
\mbar_{1,\sum h_{m,e}(0) + \sum i_e(0)+l}(H,d(0)) \times \prod_{k=1}^l 
X(d(k), \Ga'(k), \De(k))
$$
for appropriately chosen $\Ga'(k)$.

Define $\cy^c(d(0), \De(0), \Ga(0); \dots; d(l), \De(l), \Ga(l))$
as the analogous stack. 

\begin{figure}
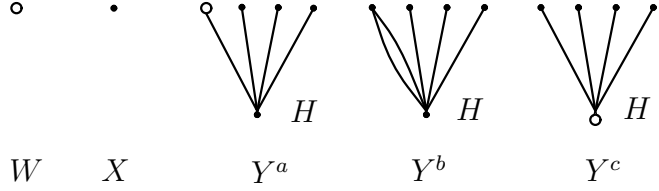

\begin{center}

	   \setlength{\unitlength}{.1\baseunit}
	    \input ewxy.tex 
\end{center}
\caption{Five classes of varieties}
\label{ewxy}
\end{figure}

The five classes of varieties $W$, $X$, $Y^a$, $Y^b$, $Y^c$ are illustrated
in Figure \ref{ewxy}.  In the figure, the dual graph of the curve
corresponding to a general point of the variety is given.  Vertices
corresponding to components mapped to $H$ are labeled with an $H$, and
vertices corresponding to elliptic components are open circles.

\subsubsection{The scheme $Z(d,\vi)_{\cd}$}
\label{esubZ}
Because of the divisorial condition Y${}^{\text{c}}$5 in the definition of
$Y^c(\ce)$, we
will also be interested in the variety parametrizing smooth degree $d$
elliptic curves in $\proj^n$ ($n \geq 2$) with a condition in the Picard
group of the curve involving the marked points and
$\pi^*(\oh_{\proj^n}(1))$.  Let $d$ and $n$ be positive integers and $\vi =
(i_e)_{e \geq 0}$ a set of non-negative integers.  Let $\De = \{ \De^j_e
\}_{e, 1 \leq j \leq i_e}$ be a set of linear spaces in $\proj^n$ where
$\dim \De^j_e = e$.  Let $\cd$ be a linear equation in formal variables $\{
q^j_e \}_{e,j}$ with integral coefficients summing to $d$.

\begin{defn}
The scheme $Z_n(d,\De)_{\cd}$ is the (scheme-theoretic) closure of the
locally closed subset of $\mbar_{1, \sum i_e}(\proj^n,d)$ (where the points are
labeled $\{q^j_e \}_{1 \leq j \leq i_e}$) representing stable maps
$(C, \{ q^j_e \}, \pi)$ satisfying the following conditions:
\begin{enumerate}
\item[(i)] The curve $C$ is smooth,
\item[(ii)] $\pi(q^j_e) \in \De^j_e$, and
\item[(iii)] in $\Pic(C)$, $\pi^* (\oh_{\proj^n}(1)) \cong \oh_C(\cd)$.
\end{enumerate}
\end{defn}
When discussing properties that hold for general $\{ \De^j_e \}_{e,j}$, we
will write $Z_n(d,\vi)_{\cd}$.  The $n$ will often be suppressed for
convenience.  Define $\cz(d,\vi)_{\cd}$ as the analogous 
stack.

For example,
$$
Z_2(d=4,i_0 = 11)_{q^1_0 + q^2_0 + q^3_0 + q^4_0}.
$$
parametrizes the finite number of smooth two-nodal quartic plane curves
through 11 fixed general points $\{ q^j_0 \}_{1 \leq j \leq 11}$ satisfying
$$
\pi^*(\oh_{\proj^2}(1)) \cong \oh(q^1_0 + q^2_0 + q^3_0 + q^4_0)
$$
in the Picard group of the normalization of the curve.

When $\Ga$ and $\De$ are general, all of the varieties $W$, $Y^a$, $Y^b$,
$Y^c$, $Z$ defined above will be seen to be 
reduced (Propositions \ref{edimW} and \ref{edimZ}).  When the dimension is 0 (and, as before, $\Ga$ and $\De$ are
general), they consist of reduced points, and the number of points is
independent of $\Ga$ and $\De$.  We call this number $\# W(\ce)$, $\#
Y^a(\ce)$, etc.  We will calculate all of these values for all $n$ and
$\ce$.

\subsection{Preliminary results}

In this section, we prove preliminary results we will need.

Recall Propositions \ref{rgeneral} and
\ref{rgeneralproper}, which are collected in the following proposition: 

\begin{pr}
\label{egeneral}
Let $\ca$ be a reduced irreducible substack of $\cmbar_{g,m}(\proj^n,d)$,
and let $p$ be one of the labeled points.  Then there is a Zariski-open
subset $U$ of the dual projective space $(\proj^n)^*$ such that for all
$[H'] \in U$ the intersection $\ca \cap \{ \pi(p) \in H' \}$, if nonempty,
is reduced of dimension $\dim \ca - 1$.

Let $\cb$ be a proper closed substack of $\ca$.  Then there is a
Zariski-open subset $U'$ of the dual projective space $(\proj^n)^*$ such
that for all $[H'] \in U'$, each component of $\cb \cap \{ \pi(p) \in H'
\}$ is a proper closed substack of a component of $\ca \cap \{ \pi(p) \in
H' \}$.
\end{pr}

The following two propositions are variants of Propositions
\ref{runiversal} and \ref{rgeneral2} of the previous section.  The proofs
are identical once $\cx$ is replaced with $\cw$.

\begin{pr}
\label{euniversal}
Given $d$, $\vh$, $\vi$, $\Ga$, $\De$, define $\vec{i'} = \vi + \vep_n$ and $\De'$ the same as
$\De$ except $\De^{i'_n}_n = \proj^n$.  Then $\cw(d, \Ga,
\De')$ is the universal curve over $\cw(d,\Ga, \De)$.
\end{pr}

\begin{pr}
\label{egeneral2}
Let $H'$ be a general hyperplane of $\proj^n$.  
\begin{enumerate}
\item[a)] The divisor 
$$
\{ \pi(p^{j_0}_{m_0,e_0}) \in H' \} \subset
\cw(d, \Ga, \De)
$$
is $\cw(d, \Ga', \De)$ where
\begin{itemize}
\item $\vec{h'}= \vh - \vep_{m_0,e_0} + \vep_{m_0,e_0-1}$ 
\item For $(m,e) \neq (m_0,e_0), (m_0,e_0-1)$, ${\Ga'}^j_{m,e} = \Ga^j_{m,e}$.
\item $\{ {\Ga'}^j_{m_0,e_0} \}_j = \{ \Ga^j_{m_0,e_0} \}_j  \setminus \{
\Ga^{j_0}_{m_0,e_0} \}$
\item $\{ {\Ga'}^j_{m_0,e_0-1} \} = \{ \Ga^j_{m_0,e_0-1}
\}_j \cup \{ \Ga^{j_0}_{m_0,e_0} \cap H' \} $
\end{itemize}
\item[b)] The divisor 
$$
\{ \pi(q^{j_0}_{e_0}) \in H' \} \subset
\cw(d, \Ga, \De)
$$ 
is $\cw(d, \Ga, \De')$ where
\begin{itemize}
\item $\vec{i'} = \vi - \vep_{e_0} + \vep_{e_0-1}$
\item For $e \neq e_0, e_0-1$, ${\De'}^j_{e} = \De^j_{e}$.
\item $\{ {\De'}^j_{e_0} \}_j = \{ \De^j_{e_0} \}_j  \setminus \{
\De^{j_0}_{e_0} \}$, $\{ {\De'}^j_{e_0-1} \}_j = \{ \De^j_{e_0-1} \}_j
\cup \{ \De^{j_0}_{e_0} \cap H' \} $ 
\end{itemize}
\end{enumerate}
\end{pr}

An analogous proposition holds for $\cz_n(d,\De)$.  The proof is
essentially the same, and is omitted.

\begin{pr}
\label{ezgeneral2}
Let $H'$ be a general hyperplane of $\proj^n$.  
The divisor $\{ \pi(q^{j_0}_{e_0}) \in H' \}$ on
$\cz(d,\De)_{\cd}$ is $\cz(d, \De')_{\cd}$ where
\begin{itemize}
\item $\vec{i'} = \vi - \vep_{e_0} + \vep_{e_0-1}$
\item For $e \neq e_0, e_0-1$, ${\De'}^j_{e} = \De^j_{e}$.
\item $\{ {\De'}^j_{e_0} \}_j = \{ \De^j_{e_0} \}_j  \setminus \{
\De^{j_0}_{e_0} \}$, $\{ {\De'}^j_{e_0-1} \}_j = \{ \De^j_{e_0-1} \}_j
\cup \{ \De^{j_0}_{e_0} \cap H' \} $ 
\end{itemize}
\end{pr}

For a given $\ce'$, $\Ga'$, $\De'$, we also have an isomorphism $\phi_1$
between $W(d', \Ga', \De')$ and
$Y^a(d(0),\Ga(0),\De(0);d(1),\Ga(1),\De(1))$ (for appropriately chosen
$\ce(0)$, \dots, $\De(1)$) that is similar to the isomorphism $\phi$ of
Proposition \ref{rXY}.  The notation used in this
proposition is the same, and the proof is also identical: the
morphism involves attaching a rational tail with two marked
points.

\begin{pr}
Fix $\ce$, integers $m_0$, $e_0$, $e_1$, and general $\Ga$, $\De$.
Let $j_0 = h_{m_0,e_0}$, $j_1 = i_{e_1}$, $e' = e_0+e_1-n$, and $j' =
h_{m_0,e'}+1$. 
Let $(d(0), \vh(0), \vi(0)) = ( 0, \vep_{m_0,e_0},  \vep_{e_1})$, $\Ga(0)
= \{ \Ga^{j_0}_{m_0,e_0} \}$, $\De(0) = \{ \De^{j_1}_{e_1} \}$,
$\ce(1) = \ce - \ce(0)$, $\Ga(1) = \Ga \setminus \Ga(0)$, $\De(1) =
\De \setminus \De(0)$.  
Let $(d', \vec{h'}, \vec{i'}) = ( d, \vh(1)+ \vep_{m_1,e'}, \vi(1))$, $\Ga'
= \Ga(1) \cup \{ 
\Ga^{j_0}_{m_0,e_0} \cap \De^{j_1}_{e_1} \}$, and $\De' =
\De(1)$.  
Then there is a natural isomorphism
$$
\phi_1:  W(d', \Ga', \De') \rightarrow Y^a(d(0), \Ga(0), \De(0);
d(1),\Ga(1),\De(1)). 
$$
\label{eWY}
\end{pr}

Proposition \ref{ebig} is a variation of the smoothness result (Proposition
\ref{rbig}) that was so useful in Section \ref{rational}.  To prove it, we
will need some preliminary results about stable maps from elliptic curves
to $\proj^1$.

\begin{lm}
\label{ecomb}
Let $C$ be a complete reduced nodal curve of arithmetic genus 1.  Let
$\pi$ be a morphism $\pi: C \rightarrow \proj^1$ contracting no component
of $C$ of arithmetic genus 1.
Then 
$$
H^1(C,\pi^*(\oh_{\proj^1}(1))) = H^1(C,\pi^*(\oh_{\proj^1}(2))) = 0.
$$
\end{lm}
By ``contracting no component of $C$ of arithmetic genus 1'' we mean that
all connected unions of contracted irreducible components of $C$ have arithmetic genus
0.

\bpf
By Serre duality, it suffices to show that 
$$
H^0(C,K_C \otimes \pi^*(\oh(-1))) = 0.
$$
Assume otherwise that such $(C,\pi)$ exists, and choose one with the fewest
components, and choose a nonzero global section $s$ of $K_C
\otimes \pi^*(\oh(-1))$.  If $C = C' \cup R$ where $R$ is a rational
tail (intersecting $C'$ at one point), then $s=0$ on $R$ as
$$
\deg_R ( K \otimes \pi^* ( \oh(-1))) = -1 - \deg_{\pi} R < 0.
$$
Then $s|_{C'}$ is a section of $(K_C \otimes \pi^*( \oh(-1)))|_{C'}$ that
vanishes on $C' \cap R$.  But $K_{C'} = K_C(-C' \cap R) |_{C'}$, so this
induces a non-zero section of $K_{C'} \otimes (\pi|_{C'})^*(\oh(-1))$,
contradicting the minimality of the number of components.  Thus $C$ has no
rational tails, and $C$ is either an irreducible elliptic curve or a cycle
of rational curves.  If $C$ is an irreducible elliptic curve, then $C$
isn't contracted by hypothesis, so $K_C \otimes \pi^*(\oh(-1))$ is negative
on $C$ as desired.  If $C$ is a cycle $C_1 \cup
\dots \cup C_s$ of $\proj^1$'s, then 
$$
\deg_{C_i} ( K_C \otimes \pi^* (\oh(-1) ))  = - \deg C_i \leq 0.
$$
As one of the curves has positive degree, there are no global sections of $K_C
\otimes \pi^* (\oh(-1))$.
\epf

\begin{lm}
\label{esm1}
Let $(C, \{ p_i \}_{i=1}^m, \pi)$ be a stable map in
$\cmbar_{1,m}(\proj^1,d)$ having no contracted component of arithmetic genus
1.  Then $\cmbar_1(\proj^1,d)$ is smooth of dimension $2d+m$ at
$(C,\{ p_i \},\pi)$.
\end{lm}
\bpf
As $H^1(C,\pi^* T_{\proj^1}) = 0$ by the previous lemma,
$\cmbar_{1,m}(\proj^1,d)$ is smooth of dimension $\deg \pi^*T_{\proj^1} + m =
2d+m$.  The argument is well-known, but for completeness we give it here.

From the exact sequence for infinitesimal deformations of stable maps (see
Subsubsection \ref{itmsosm}), we
have
\begin{eqnarray}
\label{edefsm}
0 &\longrightarrow& \Aut (C, \{ p_i \}) \longrightarrow H^0(C,\pi^*
T_{\proj^1}) \\
\nonumber
 \longrightarrow  \Def (C, \{ p_i \}, \pi ) & \longrightarrow &
\Def (C, \{ p_i \}) \longrightarrow H^1(C,\pi^*
T_{\proj^1} ) \\
 \longrightarrow  \Ob (C, \{ p_i \}, \pi ) &\longrightarrow & 0
\nonumber
\end{eqnarray}
where $\Aut (C,\{ p_i \}) = \Hom (\Omega_C(p_1 + \dots + p_m),\oh_C)$
(resp. $\Def (C, \{ p_i \}) = \Ext^1(\Omega_C(p_1 + \dots + p_m),\oh_C )$) are the
infinitesimal 
automorphisms (resp. infinitesimal deformations) of the marked curve, and
$\Def (C, \{ p_i \}, \pi)$ (resp. $\Ob  (C, \{ p_i \}, \pi)$)
are the infinitesimal deformations (resp. obstructions) of the stable map.
As $H^1(C,\pi^* T_{\proj^1}) = 0$, $\Ob(C, \{ p_i \}, \pi) = 0$ from 
(\ref{edefsm}).  Thus the deformations of $(C,\{ p_i \} , \pi)$
are unobstructed, and the dimension follows from:
\begin{eqnarray*}
\dim \Def (C, \{ p_i \}, \pi) &-& \dim \Ob (C, \{ p_i \}, \pi) \\
&=& ( \dim \Def (C, \{ p_i \}) - \dim \Aut (C, \{ p_i \})) \\
& & + (h^0(C,\pi^* T_{\proj^1})  - h^1(C,\pi^* T_{\proj^1})  ) \\
&=& m + 2d.
\end{eqnarray*}
\epf

The next lemma will be useful for studying the behavior of
the space $\mbar_{1,m}(\proj^n,d)$ at points representing maps with contracted
elliptic components.
\begin{lm}  
\label{etancondition}
Let $C$ be a complete reduced nodal curve of arithmetic genus 1, and let
$\pi: C \rightarrow \proj^n$.  Assume
$(C,\pi)$ can be smoothed.  If $B$ is a connected union of contracted
components of $C$ of arithmetic genus 1, intersecting $\overline{C
\setminus B}$ in $k$ points, and $T_1$, \dots, $T_k$ are the tangent vectors to
$\overline{C \setminus B}$ at those points, then $\{ \pi( T_i ) \}_{i=1}^k$
are linearly dependent in $T_{\pi(B)} \proj^n$.
\end{lm}
More generally, this result will hold whenever $\pi$ is a map to an
$n$-dimensional variety $X$, and $B$ is contracted to a smooth point of $X$.
It is a variation of [V2] Theorem 1 in
higher dimensions.

\bpf
Let $\De$ be a smooth curve parametrizing maps $(\cc_t, \pi)$ (with
total family $(\cc,\pi)$) to $\proj^n$, with $(\cc_0,\pi) = (C,\pi)$
and general member a map from a smooth curve.  Blowing up points of the
central fiber changes $C$, but 
does not change the hypotheses of the proposition, so we may assume
without loss of generality that the total family $\cc$ is a smooth
surface.  The following diagram is commutative.
$$\begin{array}{rcccl}
\cc & \; & {\stackrel \pi \longrightarrow} & \; & \proj^n \times \De \\
\; & \searrow & \; & \swarrow & \; \\
\; & \; & \De & \; & \;
\end{array}$$
There is an open neighborhood $U$ of $B \subset \cc$ such that $\pi \mid_{U
\setminus B}$ is an immersion.  Thus $\pi$
factors through a family $\cc'$ that is the same as $\cc$ except $B$ is
contracted.  Let $\pi'$ be the contraction $\pi': \cc \rightarrow \cc'$.
The family $\cc'$ is also flat, and its general fiber has genus 1.  The
central fiber is a union of rational curves, at most nodal away from the
image of $B$.  If the images of $T_1$, \dots, $T_k$ in $\cc'_0$ are
independent, the reduced fiber above 0 would have arithmetic genus 0, so
the central fiber (reduced away from the image of $B$) would have
arithmetic genus at most zero, contradicting the constancy of arithmetic
genus in flat families.  Thus the images of $T_1$, \dots, $T_k$ in
$T_{\pi'(B)} \cc'_0$
must be dependent, and hence their
images in $T_{\pi(B)}\proj^n$ must be dependent as well.
\epf

In Lemma \ref{esm2}, we will prove that the moduli stack
$\cmbar_{1,m}(\proj^1,d)$ is smooth even at some points with contracted
components of arithmetic genus 1.  Let $(C,\{ p_i \}_{i=1}^m, \pi )$ be a
stable map in $\cmbar_{1,m}(\proj^1,d)$ with $\pi^{-1}(z)$ containing (as a
connected component) a curve $E$ of arithmetic genus 1, where $E$
intersects the rest of the components $R$ at two points $p$ and $q$ (and
possibly others) with the $\pi|_R$ \'{e}tale at $p$.  (This result should
be true even without the
\'{e}tale condition.)
\begin{lm}
\label{esm2}
The moduli stack $\cmbar_{1,m}(\proj^1,d)$ is smooth at $(C, \{ p_i \},
\pi)$ of dimension $2d+m$.
\end{lm}
\bpf
For convenience (and without loss of generality) assume $m=0$.  The
calculations of Lemma \ref{ecomb} show that $h^1(C,\pi^* T_{\proj^1})
= 1$, so our proof of Lemma \ref{esm1} will not carry through.
However, $\Def (C,\pi)$ does not surject onto $\Def(C)$ in long exact
sequence (\ref{edefsm}), as it is not possible to smooth the nodes
independently: one cannot smooth the node at $p$ while preserving the
other nodes even to first order.  (This is well-known; one argument,
due to M. Thaddeus, is to consider a stable map $(C, \pi)$ in
$\cmbar_1(\proj^1,1)$ and express the obstruction space $\hyperExt^2 (
\underline{\Omega}_{\pi}, \oh_C)$ as the dual of $H^0(C, \cf)$ for a
certain sheaf $\cf$.)  Thus the map $\Def(C) \rightarrow H^1(C,\pi^*
T_{\proj^1})$ is not the zero map, so $\Def(C)$ surjects onto
$H^1(C,\pi^* T_{\proj^1})$.  Therefore $\Ob(C,
\pi) = 0$, so the deformations are unobstructed.

The rest of the proof is identical
to that of Lemma \ref{esm1}.
\epf

With these lemmas in hand we are now ready to prove an important smoothness
result.  Let $A$ be a general $(n-2)$-plane in $H$.  Projection
from $A$ induces a rational map $\rho_A: \cmbar_{1,m}(\proj^n,d)
\dashrightarrow \cmbar_{1,m}(\proj^1,d)$, that is a morphism (of stacks) at
points representing maps $(C, \{ p_i \}, \pi)$ whose 
image $\pi(C)$ doesn't intersect $A$.  Via $\cmbar_{1,m}(\Bl_A \proj^n,d)$, the
morphism can be extended over the set of maps $(C,\{ p_i \}, \pi)$ where
$\pi^{-1} A$ is a union
of reduced points distinct from the $m$ marked points $\{ p_i \}$.
The image of such curves in
$\cmbar_{1,m}(\proj^1,d)$ is a stable map
$$
(C \cup C_1 \cup \dots \cup C_{\# \pi^{-1} A }, \{  p_i \}, \pi')
$$
where $C_1$, \dots, $C_{\# \pi^{-1} A }$ are rational tails attached to $C$
at the points of $\pi^{-1} A$, 
$$
\pi' \mid_{ \{ C \setminus \pi^{-1} A \} } = ( p_A
\circ \pi ) \mid_{ \{ C \setminus \pi^{-1} A \} }
$$
(which extends to a morphism from all of $C$) and $\pi' \mid_{C_k}$ is a
degree 1 map to $\proj^1$ ($1 \leq k \leq \# \pi^{-1} A$).

\begin{pr}
If $(C,\{ p_i \}, \pi) \subset \cmbar_{1,m}(\proj^n,d)$, the scheme
$\pi^{-1} A$ is a
union of reduced points disjoint from the marked points, and $\pi$
collapses no components of arithmetic genus 1, then at $(C,\{ p_i \},
\pi)$, $\rho_A$ is a smooth morphism of stacks of relative dimension
$(n-1)d$.
\label{ebig}
\end{pr}
\bpf
If no components of $C$ of arithmetic genus 1 are mapped to $H$, then
$\rho_A(C,\pi)$ is a smooth point of $\cmbar_{1,m}(\proj^1,d)$ by Lemma
\ref{esm1}.  If a component of $C$ of arithmetic genus 1 is mapped to $H$,
it must intersect $A$ in at least two points.  In this case $\rho_A(C,\pi)$
consists of a curve with a contracted elliptic component, and this elliptic
component has at least two rational tails that map to $\proj^1$ with degree
1.  Thus by Lemma \ref{esm2}, $\rho_A(C,\pi)$ is a smooth point of
$\cmbar_{1,m}(\proj^1,d)$ as well.

By Lemma \ref{ecomb}, $H^1(C,\pi^*(\oh(1))) = 0$, so
$h^0(C,\pi^*(\oh(1))) = d$ by Riemann-Roch.  The proof is then identical to
that of Proposition \ref{rbig}.
\epf

We now calculate the dimension of the varieties $W$, $Y^a$,
$Y^b$, $Y^c$, and $Z$.

\begin{pr}
\label{edimW}
Every component of $W(\ce)$ is reduced of dimension
$$
(n+1) d - \sum_{m,e} (n+m-e-2)h_{m,e} - \sum_e (n-1-e)i_e.
$$
The general element of each component is (a map from) a smooth curve.

If $\sum_{k=1}^l \ce(k) = \ce$, then every component of $Y^a =
Y^a(\ce(0);\dots;\ce(l))$ (respectively $Y^b$, $Y^c$) is reduced of
dimension $\dim W(\ce)-1$.
\end{pr}
\bpf
We will prove the result about $\dim W(\ce)$ in the special case $\vi =
\vec{0}$ and $h_{m,e}=0$ when $e<n-1$.  
Then the result holds when $\vi =
i_n \vep_n$ by Proposition \ref{euniversal} (applied $i_n$ times), and we can invoke Proposition
\ref{egeneral2} repeatedly to obtain the result in full generality.  (As in
the previous section, this type of reduction will be used often.)  In
this special case, we must prove that each component of $W(\ce)$ is reduced
of 
dimension
$$
(n+1) d - \sum_m (m-1) h_{m,n-1}.
$$
Consider any point $(C, \{ p^j_{m,e} \}, \{ q^j_e \}, \pi)$ on $W(\ce)$
where no component maps to $H$ and $\pi$ collapses no component of
arithmetic genus 1.  The natural map $\cw(\ce) \dashrightarrow
\cw_1(\cehat)$ induced by $\rho_A:
\cmbar_{1,\sum h_{m,n-1}}(\proj^n,d) \dashrightarrow \cmbar_{1,\sum
h_{m,n-1}}(\proj^1,d)$ is smooth of relative dimension $(n-1)d$ at the
point $(C, \{
p^j_{m,e} \}, \{ q^j_e \}, \pi)$ by Proposition
\ref{ebig}.  The stack $\cw_1(\cehat)$ is reduced of dimension $2d+1 - \sum(m-1)
h_{m,n-1}$ by Subsection \ref{ikey}, so $\cw(\ce)$ is
reduced of dimension 
$$
(n-1)d + \dim ( \cw_1(\cehat)) = (n+1)d  - \sum_m (m-1) h_{m,n-1}
$$
as desired.
As the general element of $\cw_1(\cehat)$ is (a map from) an irreducible
curve, the same is true of $\cw(\ce)$, and thus $W(\ce)$.

The same argument works for $Y^a$, $Y^b$, and $Y^c$, as in 
Subsection \ref{ikey}, it was shown that 
$Y^a_1(\cehat)$, $Y^b_1(\cehat)$, and $\dim Y^c_1(\cehat)$ are reduced
divisors of $W_1(\cehat)$. 
\epf

\begin{pr}
\label{edimZ}
Every component of $Z_n(d,\vi)_\cd$ is reduced of dimension
$$
(n+1)d - \sum_e (n-1-e) i_e - 1.
$$
\end{pr}
\bpf
It suffices to prove the result for the generically degree $d!$ cover
$Z'_n(d,\vi)_\cd$ obtained by marking the points of intersection with a
fixed general hyperplane $H$.  This is a subvariety of
$W(d,d\vep_{1,n-1},\vi)$, and as
$$
\dim W_n(d,d\vep_{1,n-1},\vi) = (n+1)d - \sum_e (n-1-e) i_e,
$$
we wish to show that $Z'_n(d,\vi)_\cd$ is a reduced Weil divisor of
the variety $W_n(d,d\vep_{1,n-1},\vi)$.  By Proposition
\ref{ezgeneral2}, we may assume that $i_e = 0$ unless $e=n$.

By relabeling if necessary, assume $q^{i_n}_n$ appears in $\cd$ with
non-zero coefficient $\al$ (so $\cd - \al q^{i_n}_n$ is a sum of integer
multiples of $q^1_n$, \dots, $q^{i_n-1}_n$).  Let $\cw(d_1, d
\vep_{1,n-1}, \vi - \vep_n)^o$ be the open subset of $\cw(d_1, d
\vep_{1,n-1}, \vi - \vep_n)$ representing maps from smooth elliptic curves.
On the 
universal curve over $\cw(d_1, d \vep_{1,n-1}, \vi - \vep_n)^o$ there is a
reduced divisor $\cz$ corresponding to points $q$ such that
$$
\al q = (\cd - \al q^{i_n}_n) -
\pi^*(\oh(1))
$$ 
in the Picard group of the fiber.  The universal curve over the stack
$\cw(d_1, d \vep_{1,n-1}, \vi - \vep_n)$ is $\cw(d_1, d \vep_{1,n-1}, \vi)$
by Proposition \ref{euniversal}, so by definition the closure of $\cz$ in
$\cw(d_1, d \vep_{1,n-1}, \vi)$ is $\cz(d,\vi)_{\cd}$.
\epf

We will need to avoid the locus on $W(\ce)$ where an elliptic component is
contracted.  Lemma \ref{etancondition} identifies which such stable maps
could lie in $W(\ce)$.
It is likely that every stable map of the form described in the lemma can
be smoothed, which would suggest (via a dimension estimate) that when $k
\leq n+1$ those maps with a collapsed elliptic component intersecting $k$
noncontracted components (with linearly dependent images of tangent
vectors) form a Weil divisor of $W(\ce)$.  Because of the moduli of
$\mbar_{1,k}$, none of these divisors would be enumeratively meaningful.
Thus the following result is not surprising.
\begin{pr}
\label{ecodim2}
If $W'$ is an irreducible subvariety of $W(\ce)$ whose general map has a
contracted elliptic component (or more generally a contracted connected
union of components of arithmetic genus 1) and $W'$ is of codimension 1,
then $W'$ is not enumeratively meaningful.
\end{pr}
\bpf
By Proposition \ref{egeneral2}, we
may assume $i_e = 0$ unless $e=n$,  and $h_{m,e} = 0$ unless
$e=n-1$.  We could proceed naively by using the previous lemma and simply
counting dimensions, but the following argument is slightly cleaner.

Let $(C, \{ p^j_{m,n-1} \}, \{ q^j_n \}, \pi)$ be a general point of $W'$,
and let $E$ be the contracted component of $C$.  Say $E$ has $s$ special
points (markings or intersections with noncontracted components) including
$k$ intersections with noncontracted components.  Replace $E$ by a rational
$R=\proj^1$, with the $s$ special points distinct, to obtain a new stable
map $(C', \{ p^j_{m,n-1} \}, \{ q^j_n \}, \pi') \in X(\ce)$.  The family of such $(C',\pi')$ forms a
subvariety $X'$ of $X(\ce)$, and $X'$ is contained in $X''$ where in the
latter we don't impose the dependence of tangent vectors required by the
previous lemma.  Let $\xi$ be the natural rational map to the
Hilbert scheme of Subsubsection \ref{rgeomean}.

If $s \geq 2$, $X''$ is codimension at least 1 in $X(\ce)$.  Due to the
moduli of $s$ points on $R$, $\xi(X'')$ is codimension at least $1
+ (s-3) = s-2$ in $\xi(X(\ce))$.  The previous lemma imposes an
additional $\max(n+1-k,0)$ conditions, which are independent as the
rational curves intersecting $R$ can move freely under automorphisms of
$\proj^n$ preserving $H$.  Thus the codimension of $\xi(X')$ in
$\xi(X(\ce))$ is at least $n-1+(s-k) \geq n-1$.  But $\dim X(\ce) -
\dim W(\ce) = n-3$, so $\dim \xi(W') < \dim W(\ce)-2 = \dim W'$, as desired

Otherwise, $k=s=1$.  By Proposition \ref{euniversal}, we may assume that
$\vi = \vec{0}$ as there are no marked points on the contracted component
$E$.  Then $X''$ can be identified with the subvariety of $X(d, \vh,
\vep_n)$ where the corresponding map $\pi:  (C, \{ p^j_{m,n-1} \}, q^1_n)
\rightarrow \proj^n$ is singular at $q^1_n$.  As the singularity condition
imposes $n$ conditions,
\begin{eqnarray*}
\dim \xi(W') & \leq & \dim X(d, \vh, \vep_n) - n \\ 
&=& \dim X(\ce)+1-n \\
&=& (\dim W(\ce) + n - 3 ) + 1-n \\
&=& \dim W(\ce) - 2 \\
&=& \dim W' - 1
\end{eqnarray*}
as desired.
\epf

\subsection{Degenerations set-theoretically}

The theorem listing the possible degenerations follows the same pattern as
the corresponding results (Theorems \ref{rlist1} and \ref{rlist}) of the
previous section.  Fix $\ce$ and a non-negative integer $E$, and let $\Ga$
and $\De$ be sets of general linear spaces of $\proj^n$ (as in the
definition of $W(d,\Ga,\De)$).  Let $q$ be the marked point corresponding
to one of the (general) $E$-planes $Q$ in $\De$.

Let $D_H = \{ \pi(q) \in H \}$ be the divisor on $W(d, \Ga,\De)$ that corresponds
to requiring $q$ to lie on $H$.  In this section, we will determine the
enumeratively meaningful components of $D_H$.  That is, we will give a list
of subvarieties, and show that the enumeratively meaningful components of
$D_H$ are a subset of this list.  In the subsequent section, we will
determine the multiplicity with which each enumeratively meaningful
component appears.  In particular, we will see that the multiplicity of
each component on the list is at least one, so each element of the list is
indeed a component of $D_H$.

As before, we can relate this result to the enumerative problem we wish to
solve.  If $W(\ce^-)$ is a union of points, we can determine $\# W(\ce^-)$ by
specializing one of the linear spaces of $\De$, of dimension $E-1$, to the
hyperplane $H$.  Define $\ce$ by
$(d,\vh,\vi) = (d^-,\vec{h^-},\vec{i^-} + \vep_E -
\vep_{E-1})$.  To calculate $\# W(\ce^-)$, we simply enumerate the
points $D_H$ on $W(\ce)$, with the appropriate multiplicity.  Let $\phi_1$ be
the isomorphism of Proposition \ref{eWY}.

\begin{tm}
\label{elist}
If $\Ga$ and $\De$ are general, each enumeratively meaningful component of $D_H$ (as a divisor on
$W(d,\Ga,\De)$) is one of the following.
\begin{enumerate}
\item[(I)] A component of $\phi_1(W(d',\Ga', \De'))$, where, for some $m_0, e_0$,
$1 \leq j_0 \leq h_{m_0, e_0}$, $e' := e_0 + E - n \geq 0$:
\begin{itemize}
\item $(d',\vec{h'}, \vec{i'}) = (d,\vh - \vep_{m_0, e_0} + \vep_{m_0,
e'},  \vi - \vep_E    )$
\item ${\Ga'}^j_{m,e} = \Ga^j_{m,e} \quad \text{if $(m,e) \neq (m_0,
e_0)$}$
\item $ \{ {\Ga'}^j_{m_0, e_0} \}_j = \{ \Ga_{m_0,e_0}^j \}_j \setminus \{
\Ga^{j_0}_{m_0,e_0} \}$
\item ${\Ga'}_{m_0, e'}^{h'_{m_0, e'}} = \Ga^{j_0}_{m_0,e_0} \cap Q$
\item ${\De'}^j_e = \De^j_e \quad \text{if $e \neq E$}$,
$\{ {\De'}_E^j \}_j = \{ \De^j_E \}_j \setminus \{ Q \}$.
\end{itemize}
\item[(II)]  A component of
$Y^a(d(0),\dots, \De(l))$,
$Y^b(d(0),\dots,\De(l))$, or
$Y^c(d(0),\dots,\De(l))_{m^1_1}$ for some $l$,
$\ce(0)$, \dots, $\De(l)$, with $\ce = \sum_{k=0}^l \ce(k)$,
$\Ga = \cup_{k=0}^l \Ga(k)$, $\De =
\cup_{k=0}^l \De(k)$, $Q \in \Ga(0)$, and $d(0)>0$. 
\end{enumerate}
\end{tm}
Call the components of (I) Type I components, and call the three types of
components of (II) Type IIa, IIb, and IIc components
respectively.

\bpf  We follow the proofs of Theorems \ref{rlist1}
and \ref{rlist}.  By Propositions \ref{euniversal} and \ref{egeneral2}, we
may assume that $\vi = \vep_n$, $E=n$, and $h_{m,e} = 0$ unless $e=n-1$.
With these assumptions, the result becomes much simpler.
The stack $\cw(d,\vh,\vi)$ is the universal curve over
$\cw(d,\vh,\vec{0})$, and we are asking
which points of the universal curve lie in $\pi^{-1}H$.

Let $(C, \{ p^j_{m,n-1} \}, q, \pi)$ be the map corresponding to a general
point of a enumeratively 
meaningful component $Y$ of $D_H$.  By Proposition \ref{ecodim2}, the
morphism $\pi$
doesn't contract any component of $C$ of arithmetic genus 1.  Choose a
general $(n-2)$-plane $A$ in $H$.  The set $\pi^{-1} A$ is a union of
reduced points on $C$, so $\rho_A$ is smooth (as a morphism of stacks) at
$(C, \{ p^j_{m,n-1} \}, q, \pi)$ by Proposition \ref{ebig}.  As a set,
$D_H$ contains the entire fiber of $\rho_A$ above $\rho_A(C,\pi)$, so
$\rho_A(D_H)$ is a Weil divisor on $W_1(\cehat)$ that is a component of
$\{\pi(q) = z\}$.  By Theorem \ref{igenus1}, the curve
$C$ is a union of irreducible components $C(0) \cup \dots \cup C(l')$ with
$\rho_A \circ \pi(C(0)) = z$ (i.e. $C(0) \subset \pi^{-1}H$), $C(0)\cap C(k)
\neq \phi$, and the marked points split up among the components:
$\vh = \sum_{k=0}^{l'} \vh(k)$.  If $d(0) = \deg \pi |_{C(0)}$, then $d(0)$
of the curves $C(1)$, \dots, $C(l')$ are rational tails that are collapsed
to the $d(0)$ points of $C(0) \cap A$; they contain no marked points.  Let
$l = l' - d(0)$.  Also, $\vi(0) = \vep_n$, and $\vi(k) =
\vec{0}$ for $k>0$, as the only incidence condition was $q \in Q$, and
$q(0) \in C(0)$.

{\em Case $d(0)>0$.} By Theorem \ref{igenus1}, the
component $Y$ is contained in 
$$
Y^a(d,
 \dots, \De(l)), \quad Y^b (d(0), \dots, \De(l)), \quad \text{or} \quad Y^c 
(d(0), \dots,
\De(l)).
$$
As the dimensions of each of these three is $\dim W(d,\Ga,\De)
- 1 = \dim Y$, $Y$ must be a Type II component as described in the
statement of the theorem. 

{\em Case $d(0)=0$.}  As the morphism $\pi$ contracts no elliptic
components, the curve $C(0)$ is rational.  Also, $C(0)$ has at least 3
special points: $q$, one of $\{ p^j_{m,e} \}$ (call it
$p^{j_0}_{m_0,e_0}$), and $C(0) \cap C(1)$.  If $C(0)$ had more than 3
special points, then the component would not be enumeratively meaningful,
due to the moduli of the special points of $C(0)$.  Thus $l=1$, and $Y$ is
a Type I component.
\epf

In fact, this argument determines {\em all} the components of $D_H$ except
those representing maps with collapsed elliptic tails.  (It is not
clear whether such components exist.)

When $n=2$, the only enumeratively meaningful Type II divisors are Type IIa
and Type IIb with $d(0)=1$.  This agrees with the genus 1 case of [CH3].

\subsection{Multiplicity calculations}
\label{emultgen}
The proof of multiplicities of $D_H$ along the enumeratively meaningful components
are the 
same as in the genus 0 case (Subsection \ref{rmultgen}).  The Type I component
$\phi_1(W(d',\Ga',\De'))$ appears with multiplicity $m_0$, where $m_0$ was
defined in Theorem \ref{elist}.  (The argument of 
Subsubsection \ref{rmultI} also proves this.)
The Type IIa component $Y^a(\ce(0);\dots;\ce(l))$ appears with multiplicity
$\prod_{k=1}^l m^k$, where $m^k = d(k) - \sum_{m,e,j} m h^j_{m,e}(k)$ as defined
earlier.  The Type IIc component $Y^c(\ce(0);\dots;\ce(l))$ appears with
multiplicity $\prod_{k=1}^l m^k$ as well.  The Type IIb component
$Y^b(\ce(0);\dots;\ce(l))_{m^1_1}$ appears with multiplicity
$$
m^1_1 m^1_2 \prod_{k=2}^l m^k = 
m^1_1 (m^1- m^1_1) \prod_{k=2}^l m^k .
$$
By Corollary \ref{ilocalst} of Section \ref{intro}, we also get the same
results about the structure of $W(\ce)$ in a formal,
\'{e}tale, or analytic neighborhood of these components. 

\subsection{Recursive formulas}

We now enumerate the points of our varieties when the number is
finite.
The only initial data
needed is the ``enumerative geometry of $\proj^1$'': the number of stable
maps to $\proj^1$ of degree 1 is 1.  

\subsubsection{A recursive formula for $\# W$}

Given $\ce$, fix an $E$ such that $i_E>0$.  Partitions of $\ce$
are simultaneous partitions of $d$, $\vh$, and $\vi$.  Define
multinomial coefficients with vector arguments as the product of the
multinomial coefficients of the components of the vectors.

Define $\ce^-$ by $(d^-, \vec{h^-}, \vec{i^-}) = ( d, \vh, \vi - \vep_E +
\vep_{E-1})$, and let $\Ga^- = \Ga$ and 
$\De^- = \De \cup \{ \De^{i_E}_E \cap H' \}  \setminus \{ \De^{i_E}_E \}$
where $H'$ is a general hyperplane.

The following theorem is an analog of Theorem
\ref{rrecursiveX1}.

\begin{tm}
\label{erecursiveW1}
In $A^1 (W(d,\Ga, \De))$, modulo enumeratively meaningful divisors,
$$
W_n(d^-, \Ga^-, \De^-) = \sum   m_0 \cdot  W_n(d'(m_0,e_0,j_0), \Ga',
\De')
$$
$$
+ \sum \left(
\prod_{k=1}^l m^k \right)  Y^a_n(d(0), \dots, \De(l) )
$$
$$
+ \sum \left( \prod_{k=2}^l m^k \right) 
\left( \sum_{m^1_1 = 1}^{ [ m^1/2]}
m^1_1 (m^1 - m^1_1) Y^b_n(d(0), \dots, \De(l))_{m^1_1} \right)
$$
$$
+ \sum \left( \prod_{k=1}^l m^k \right) Y^c_n(d(0), \dots, \De(l) )
$$
where the first sum is over all $(m_0,e_0,j_0)$, and $(d', \Ga', \De')$
is as defined in Theorem \ref{elist}; and the last three sums are over
all $l$, $\ce(0)$, \dots, $\De(l)$ with $\ce = \sum_{k=0}^l \ce(k)$, $\Ga =
\coprod_{k=0}^l \Ga(k)$, $\De = \coprod_{k=0}^l \De(k)$, $\Ga_E^{i_E} \in
\Ga(0)$, and $d(0)>0$.
\end{tm}

\bpf
The left side is linearly equivalent (in $\cw(d,\Ga,\De)$) to $D_H = \{
\pi(q^{i_E}_E) \in H \}$.  The right side is set-theoretically $D_H$ by
Theorem \ref{elist}, and the multiplicities were determined in
Subsection \ref{emultgen}.
\epf

If $\# W_n(\ce^-)$ is finite, the following statement is more suitable for
computation.  It is an analog of Theorem \ref{rrecursiveX2}.

\begin{tm}
\label{erecursiveW}
{ \scriptsize
$$
\# W_n(\ce^-) = \sum_{m,e} h_{m,e} \cdot m \cdot \# W_n(\ce'(m,e)) 
$$
$$
+ \sum \left(
\prod_{k=1}^l m^k \right) { \binom {\vh} { \vh(0), \dots, \vh(l)} } {
\binom {\vi - \vep_E}  { \vi(0)-\vep_E,\vi(1), \dots,
\vi(l)} } \frac { \# Y^a_n(\ce(0); \dots; 
\ce(l)) }{ \Aut( \ce(2),\dots,\ce(l))}
$$
$$
+ \frac 1 2 \sum \left( \prod_{k=2}^l m^k \right){ \binom {\vh} { \vh(0),
\dots, \vh(l)} } { \binom {\vi - \vep_E} {\vi(0)-\vep_E,\vi(1), \dots,
\vi(l)} }  
$$
$$
\quad \quad \quad 
\cdot \left( \sum_{m^1_1 = 1}^{  m^1-1}
m^1_1 (m^1 - m^1_1) \frac { \# \tilde{Y}^b_n(\ce(0); \dots; \ce(l))_{m^1_1} }{ 
{\Aut( \ce(2),\dots,\ce(l))}} \right)
$$
$$
+ \sum \left( \prod_{k=1}^l m^k \right){ \binom {\vh}  { \vh(0), \dots,
\vh(l)} } { \binom {\vi - \vep_E}  {
\vi(0)-\vep_E,\vi(1), \dots, \vi(l)} } \frac { \# Y^c_n(\ce(0);
\dots; \ce(l)) }{ \Aut( \ce(1),\dots,\ce(l))}
$$
}
where, in the first sum, $\ce'(m,e) = (d,
\vh -\vep_{m,e}+\vep_{m,e+E-n}, \vi - \vep_E)$; the last
three sums are over all $l$ and all partitions $\ce(0)$, \dots, $\ce(l)$ of
$\ce$ with $d(0)>0$.
\end{tm}
\bpf
When $\#W_n(\ce^-)$ is finite, all components of $D_H$ are enumeratively
meaningful.  Take degrees of both sides of the equation in Theorem
\ref{erecursiveW1}.  As
$$
\# \tilde{Y}^b_n(\ce(0);
\dots; \ce(l))_{m^1_1} =  
\begin{cases}
2 (\# Y^b_n(\ce(0); \dots; \ce(l))_{m^1_1})  & \text{if $2 m^1_1 = m^1$,} \\
\# Y^b_n(\ce(0); \dots; \ce(l))_{m^1_1}  & \text{otherwise,}
\end{cases}
$$
it follows that
\begin{eqnarray*}
& & \sum_{m_1^1 = 1}^{[m^1/2]} m_1^1 ( m^1 - m^1_1) \# Y^b_n(\ce(0); \dots;
\ce(l))_{m^1_1} \\
& & \quad \quad \quad \quad \quad \quad = \frac 1 2 
\sum_{m_1^1 = 1}^{m^1-1} m_1^1 ( m^1 - m^1_1) \# \tilde{Y}^b_n(\ce(0);
\dots; \ce(l))_{m^1_1}. 
\end{eqnarray*}
The only additional points requiring explanation are the combinatorial aspects:
the $h_{m,e}$ in the first sum, and the ``$\Aut$'' and various multinomial
coefficients in the last three.  In Theorem \ref{elist}, the Type I
components were indexed by $(m_0, e_0, j_0)$.  But for fixed $(m_0,e_0)$,
$\# W(\ce'(m_0,e_0))$ is independent of $j_0$, so the above formula
eliminates this redundancy.  Similarly, in Theorem \ref{rlist}, the Type
IIa, IIb, and IIc components were indexed by partitions of the points $\{
p^j_{m,e} \}_{m,e,j}$ and $\{ q^j_e \}_{e,j}
\setminus \{ q \}$, but the values of $\# Y^a_n(\ce(0); \dots; \ce(l))$, 
$\# Y^b_n(\ce(0); \dots; \ce(l))$, and $\# Y^c_n(\ce(0); \dots;
\ce(l))$
depend only on $\{ \vh(k), \vi(k) \}_{k=0}^l$ and not on the actual
partitions of marked points.  The multinomial coefficients eliminate
this redundancy.  We divide by $\Aut(\ce(1); \dots; \ce(l))$
(resp.  $\Aut(\ce(2); \dots; \ce(l))$,  $\Aut(\ce(1); \dots; \ce(l))$)
to ensure that we are counting each Type IIa (resp. Type IIb, Type IIc)
component once.
\epf

Theorem \ref{erecursiveW1} can be strengthened to be true modulo those
divisors whose general map has a collapsed elliptic component.  There is a
similar statement in the Chow ring of the Hilbert scheme modulo exceptional
divisors of $\phi_1$ (cf. Theorem \ref{rrecHilb}).

\subsubsection{A recursive formula for $Y^a$}
\label{eYarec}
In the previous section, we found a formula for $\# Y$
in terms of $\# X$ (Proposition \ref{rrecursiveY}).  By the same argument, we
have 
$$
\# Y^a_n( \ce(0);\dots;\ce(l)) = \frac { \# X_{n-1}(\ce'(0)) }{ d(0)!} \cdot \#
W_n(\ce'(1)) \cdot \prod_{k=2}^l \# X_n(\ce'(k))
$$
where
\begin{itemize}
\item $d'(0) = d(0)$, $\vec{h'}(0) = d(0) \vep_{1,n-2}$
\item $\vec{i'}(0) = \sum_e i_{e+1}(0) \vep_e + \vep_{\dim W_n(\ce(1))} + \sum_{k=2}^l
\vep_{\dim X_n(\ce(k))} + \sum_{m,e} h_{m,e}(0) \vep_e$
\item $\vec{h'}(1) = \vh(1) + \vep_{m^1, n-1-\dim W_n(\ce(1))}$ and
$\vec{i'}(1) = \vi(1)$. 
\item For $2 \leq k \leq l$, 
$\vec{h'}(k) = \vh(k) + \vep_{m^k, n-1-\dim X_n(\ce(k))}$ and
$\vec{i'}(k) = \vi(k)$.
\end{itemize}

\subsubsection{Computing $\# Y^b$ and $\# \tilde{Y}^b$}  
The difficulty in computing $\# Y^b$ and $\# \tilde{Y}^b$ comes from
requiring the curves $C(0)$ 
and $C(1)$ to intersect twice (at marked points of each curve, on $H$), so
a natural object of study is the blow-up of
$H \times H$ along the diagonal $\De$, $\Bl_\De H \times H$.

But when $n=2$, the situation is simpler.  The curve $C(0)$ is $H$, and $C(1)$ will
always intersect it.  In this case, for $C(0) \cup \dots \cup C(l)$ to be
determined by the incidence conditions (up to a finite number of
possibilities), each of $C(1)$, \dots, $C(l)$ must also be determined (up
to a finite number).  The analogous formula to that for $\# Y^a$ (and $\#
Y$ in the previous section) is 
$$
\# \tilde{Y}^b_2(\ce(0); \dots; \ce(l))_{m^1_1} = \prod_{k=2}^l \# X_2(d(k),
\vec{h'}(k), \vi(k))
$$
where
$\vec{h'}(1) = \vh(1) + \vep_{m^1_1,1} + \vep_{m^1_2,1}$, and 
for $2 \leq k \leq l$, $\vec{h'}(k) = \vh(k) + \vep_{m^k,1}$.

This is in agreement with Theorem 1.3 of [CH3].

We now calculate $\# Y^b$ and $\# \tilde{Y}^b$
when $n=3$; the same method will clearly work for $n>3$.  As an
illustration of the method, consider the following enumerative problem.

{\em Fix seven general lines $L_1$, \dots, $L_7$ in $\proj^3$ and a point $p$
on a hyperplane $H$.  How many pairs of curves $(C(0),C(1))$ are there with
$C(0)$ a line in $H$ through $p$ and $C(1)$ a conic intersecting $L_1$,
\dots, $L_7$ and intersecting $C(0)$ at two distinct points, where the
intersections are labeled $a_1$ and $a_2$?}

The answer to this enumerative problem is
$$
\# Y^b = \# Y^b_3(1, \vep_{1,0},\vec{0}; 2, \vec{0}, 7 \vep_1);
$$
we will calculate instead $\# \tilde{Y}^b = 2 ( \# Y^b )$.

The space of lines in $H$ passing through $p$ is one-dimensional, and this
defines a three-dimensional locus in $\Bl_\De H \times H$ (which is the
Fulton-Macpherson configuration space of 2 points in $\proj^3$;
alternately, it is a degree 2 \'{e}tale cover of $\Sym^2 H$).
The space of conics in $\proj^3$ passing through 7 general lines is
one-dimensional, and thus defines a one-dimensional locus in $\Bl_\De H 
\times H$
parametrizing the points of intersection of the conic with $H$.  Then $\#
\tilde{Y}^b$ is the intersection of these two classes, and $\# Y^b$ (and
the answer to the enumerative problem) is half this.

Let $h_i$ be the class (in the Chow group) of the hyperplane on the
$i^{th}$ factor of $\Bl_{\De} H \times H$ ($i = 1,2$), and let
$e$ be the class of the exceptional divisor.  Then the Chow ring of
$\Bl_\De H \times H$ is generated (as a ${\mathbb Z}$-module) by the classes
listed below 
with the relations
$$
h_1^3 = 0, \quad h_2^3 = 0, \quad e^2 = 3 h_1 e - h_1^2 - h_1 h_2 - h_2^2.
$$   
\begin{center}
\begin{tabular}{c|c}
Codimension & Classes \\
\hline
0 & 1 \\
1 & $h_1$, $h_2$, $e$ \\
2 & $h_1^2$, $h_1 h_2$, $h_2^2$, $h_1e = h_2e$ \\
3 & $h_1^2 h_2$, $h_1 h_2^2$, $h_1^2 e = h_1 h_2 e = h_2^2 e$ \\
4 & $h_1^2 h_2^2$\\
\end{tabular}
\end{center}

Let the image of possible pairs of points on $C(0)$ be the class $\cc(0) 
= \al (h_1 + h_2) + \be e$.  Then $\cc(0) \cdot h_1^2 h_2 = \al$ and
$\cc(0) \cdot e h_1^2  =  - \be$.  But $\cc(0) \cdot h_1^2 h_2$ is the
number of lines in $H$ passing through $p$ (class $\cc(0)$) and another
fixed point (class $h_1^2$) with a marked point on a fixed general line
(class $h_2$), so $\al =1$.  Also, $\cc(0) \cdot e h_1^2$ is the number of
lines in the plane through $p$ (class $\cc(0)$) and another fixed point $q$
(class $h_1^2$) with a marked point mapping to $q$, so $\be = -1$.  Thus
$\cc(0) = h_1 + h_2 - e$.

Let the image of possible pairs of points on $C(1) \cap H$ be the class
$\cc(1) = \al ( h_1^2 h_2 + h_1 h_2^2) + \be h_1 h_2 e$.  Thus $\cc(1)
\cdot h_1 = \al$ and $\cc(1) \cdot e = - \be$.  Then $\cc(1) \cdot h_1$ is
the number of conics in $\proj^3$ through 7 general lines in $\proj^3$ and
a general line in $H$, which is 92 from Subsection \ref{r92subsection} of the
previous section.  Also, $\cc(1) 
\cdot e$ counts the number of conics in $\proj^3$ through 7 general lines
in $\proj^3$ and tangent to $H$, which is 116 using the methods of the
previous section.  Thus $\cc(1) = 92 (h_1^2 h_2 + h_1 h_2^2) - 116 e h_1
h_2$.

Finally, the number of pairs of curves is
$$
(h_1 + h_2 - e) \left( 92 (h_1^2 h_2 + h_1 h_2^2 \right) - 116 e h_1 h_2)
=  92 + 92 - 116 = 68.
$$
Thus $\# \tilde{Y}^b$ is 68, and $\# Y^b$ (and the answer to the enumerative
problem) is 34.

When $n=3$ in general, there are three cases to consider.  Let $C=C(0)
\cup C(1) \cup \dots \cup C(l)$ as usual.

{\em Case i).} If the incidence conditions $\ce(1)$ specify $C(1)$ up to a
finite number of possibilities, then:
$$
\# \tilde{Y}^b (\ce(0); \ce(1); \dots; \ce(l)) \quad \quad \quad
$$
$$
\quad \quad \quad
= \# X (d(1),\vec{h'}(1),\vi(1)) \cdot \# Y (d(0),\vec{h'}(0),\vi(0); \ce(2);\dots;\ce(l))
$$
where 
\begin{itemize}
\item $\vec{h'}(1) = \vh(1) +
\vep_{m^1_1,2} + \vep_{m^1_2,2}$ (the curve $C(1)$ intersects $H$ at two
points $a_1$ and $a_2$ with multiplicity $m^1_1$ and $m^1_2$ respectively;
these will be the intersections with $C(0)$)
\item $\vec{h'}(0) = \vh(0) + 2
\vep_{1,0}$ (the curve $C(0)$ must pass through two points of intersection
$a_1$ and $a_2$ of $C(1)$
with $H$).
\end{itemize}

{\em Case ii).} If the incidence conditions $\ce(1)$ specify $C(1)$ up to a
one-parameter family, we are in the same situation as in the enumerative
problem above. 
Then
$$
\# \tilde{Y}^b (\ce(0); \ce(1); \dots; \ce(l)) = \quad \quad \quad
$$
$$
 \Bigl( d(0) \left( \# X (d(1),\vec{h'},\vi(1)) + d(0) \# X
(d(1),\vec{h''},\vi(1)) \right)
$$
$$
\quad \quad \quad \quad  - 
\# X (d(1),\vec{h'''},\vi(1)) \Bigr) 
\cdot \# Y (d(0),\vec{h'}(0),\vi(0); \ce(2);\dots;\ce(l))
$$
where 
\begin{itemize}
\item $\vec{h'} = \vh(1) + \vep_{m^1_1,1} + \vep_{m^1_2,2}$ (the
curve $C(1)$ intersects $H$ with multiplicity $m^1_1$ at $a_1$ along a
fixed general line and with multiplicity $m^1_2$ at $a_2$ at another point
of $H$),
\item $\vec{h''} = \vh(1) + \vep_{m^1_2,1} + \vep{m^1_1,2}$ (the
curve $C(1)$ intersects $H$ with multiplicity $m^1_2$ at $a_2$ along a
fixed general line and with multiplicity $m^1_1$ at $a_1$ at another point
of $H$),
\item $\vec{h'''} = \vh(1) + \vep_{m^1,2}$ (the points $a_1$ and $a_2$
on the curve $C(1)$ coincide, and $C(1)$ is 
required to intersect $H$ at this point with multiplicity $m^1 = m^1_1 +
m^1_2$),
\item $\vec{h'}(0) = \vh(0) + \vep_{1,0}$ (the curve $C(0)$ is additionally required to
pass through a fixed point in $H$).
\end{itemize}

{\em Case iii).} If the incidence conditions on $C(0) \cup C(2) \cup \dots \cup
C(l)$ specify the union of these curves up to a finite number of
possibilities (and the incidence conditions on $C(1)$ specify $C(1)$ up to
a two-parameter family), a similar argument gives
$$
\# \tilde{Y}^b (\ce(0); \ce(1); \dots; \ce(l)) =  \quad \quad \quad
$$
$$
 \left( d(0) \# X (d(1), \vec{h'}(1),\vi(1)) - 
\# X (d(1), \vec{h''}(1), \vi(1)) \right)
$$
$$
\quad \quad \quad  \cdot \# Y (\ce(0); \ce(2);\dots;\ce(l))
$$
where 
\begin{itemize}
\item $\vec{h'}(1) = \vh(1) + \vep_{m^1_1,1} + \vep_{m^1_2,1}$ (the curve
$C(1)$ must  intersect $H$ along two fixed general lines at the points
$a_1$ and $a_2$ with multiplicity $m^1_1$ and $m^1_2$ respectively)
\item $\vec{h''}(1) = \vh(1) + \vep_{m^1,1}$
(the points $a_1$ and $a_2$ coincide on the curve $C(1)$, and $C(1)$ is
required to intersect $H$ at that point with multiplicity $m^1
= m^1_1 + m^1_2$ along a fixed
general line of $H$).
\end{itemize}

These three cases are illustrated pictorially in Figure \ref{eYb} for the
special case of conics in $\proj^3$ intersecting a line in $H$ at two
points, with the entire configuration required to intersect 8 general lines
in $\proj^3$.  One of the intersection points of the conic with $H$ is
marked with an ``$\times$'' to remind the reader of the marking $a_1$.  The
distribution of the line conditions (e.g. the number of line conditions on
the conic) is indicated by a small number.  The bigger number beside each
picture is the actual solution to the enumerative problem corresponding to
the picture.  For example, there are 116 conics in $\proj^3$ tangent to a
general hyperplane $H$ intersecting 7 general lines.

\begin{figure}
\begin{center}

	   \setlength{\unitlength}{.07\baseunit}
	    \input eYb.tex 
\end{center}
\caption{Calculating $\# \tilde{Y}^b$:  A pictorial example}
\label{eYb}
\end{figure}

\subsubsection{A recursive formula for $Y^c$}
By the same method as in Subsubsection \ref{eYarec} for $Y^a$ (and $Y$ in the
previous section), we have
$$
\# Y^c(\ce(0); \dots; \ce(l)) = \quad \quad \quad \quad
$$
$$
\# Z_{n-1} (d(0), \vec{i'}(0) )_{\sum_k
m^k {q''}^k - \sum_{m,e,j} m {q'}^j_{m,e}} \cdot \prod_{k=1}^l \# X_n(d(k),
\vec{h'}(k), \vi(k))
$$
where
\begin{itemize}
\item For $1 \leq k \leq l$, 
$\vec{h'}(k) = \vh(k) + \vep_{m^k, n-1-\dim X(\ce(k))}$.
\item $i'_e(0) = i_{e+1}(0) + \# \{ \dim X(\ce(k)) = e \}_{1 \leq k \leq l} +
\sum_m h_{m,e}(0)$
\item the marked points on the curves of $Z_{n-1}(d(0),\vec{i'}(0))$ have been relabeled
$$
\{ {q'}^1_{m,e}, \dots, {q'}^{h_{m,e}(0)}_{m,e} \}_{m,e} \cup
\{ {q''}^1, \dots, {q''}^l \} \cup
\{ {q'''}^1_e, \dots, {q'''}^{i_{e+1}(0)}_e \}_e
$$
where ${q'}^j_{m,e} = p^j_{m,e}(0)$, ${q''}^k = C(0) \cap C(k)$,
${q'''}^j_e = q^j_{e+1}(0)$.
\end{itemize}
This formula is merely a restatement of the divisorial condition
Y${}^{\text{c}}$5 in Definition \ref{eYcdef}.

\subsubsection{Evaluating $\#Z$}
\label{eevalZ}
We will use intersection theory on elliptic fibrations over a curve --- the
Chow ring modulo algebraic equivalence or numerical equivalence will
suffice.  Let $\cf$ be an elliptic fibration over a smooth curve whose
fibers are smooth elliptic curves, except for a finite number of fibers
that are irreducible nodal elliptic curves.  Let $F$ be the class of a
fiber.  Then the self-intersection of a section is independent of the
choice of section.  ({\em Proof:} $K_\cf$ restricted to the generic fiber
is trivial, so $K_\cf$ is a sum of fibers.  Let $S_1$, $S_2$ be two
sections.  Using adjunction, $S_1^2 + K_\cf \cdot S_1 = (K_\cf + S_1) \cdot
S_1 = 0$, so $S_1^2 = - K_\cf \cdot S_1 = - K_\cf
\cdot S_2 = S_2^2$.)

For convenience, call the self-intersection of a section $S^2$.  The
parenthetical proof above shows that $K_\cf = - S^2 F$. 
\begin{pr}
Let $S$ be a section, and $C$ a class on $\cf$ such that $S=C$ on the
general fiber.  Then $S = C + (\frac { S^2 - C^2 }{ 2})F$.
\end{pr}
\bpf
As all fibers are irreducible, $S = C + kF$ for some $k$.  By adjunction,
\begin{eqnarray*}
0 &=& S \cdot (K_\cf + S) \\
&=& (C + kF) (C + (k - S^2)F) \\
&=& C^2 + 2k-S^2.
\end{eqnarray*}
Hence $k = (S^2 - C^2)/2$.
\epf

If the dimension of $Z_n(d,\vi)_{\sum m^j_e q^j_e}$
is 0, then consider the universal family over the curve parametrizing maps
to $\proj^n$ with the incidence conditions of $\vi$.  This is $W_n(d, d
\vep_{1,n-1},\vi)$ modulo the symmetric group $S_d$.  The general point of
$W_n(d, d \vep_{1,n-1},\vi)$ represents a smooth elliptic curve.  The
remaining points of $W_n(d, d \vep_{1,n-1},\vi)$ represent curves that are
either irreducible and rational or elliptic with rational tails.  (This can
be proved by simple dimension counts on $W_1(d, d \vep_{1,0}, \vec{0})$.) 
Normalize the base
(which will normalize the family), and blow down (-1)-curves in fibers.
The curves blown down come from maps from nodal curves $C(0) \cup C(1)$,
where $C(0)$ is rational and $C(1)$ is elliptic.  Call the resulting family $\cf$.
Let $H$ be the pullback of a hyperplane to $\cf$, and let $Q^j_e$ be the
section given by $q^j_e$. 
\begin{tm}
Let $D = H - \sum m^j_e Q^j_e$.  Then
$$
\# Z_n(d,\vi)_{\sum m^j_e q^j_e} = S^2 - D^2 / 2.
$$
\end{tm}
\bpf
Let $Q$ be any section.  Let $S$ be the section given by $Q +
\pi^*(\oh(1)) - \sum m^j_e q^j_e$ in the Picard group of the generic fiber.  Then 
$$
\# Z_n (d, \vi)_{\sum m^j_e q^j_e} = S \cdot Q.
$$
(The sections $S$ and $Q$ intersect transversely from Subsubsection 
\ref{esubZ}.)  By the
previous proposition, as $S^2 = Q^2$,
\begin{eqnarray*}
S &=& Q + D + \left( \frac {S^2 - (Q+D)^2 }{ 2} \right) F \\ 
\text{so } S \cdot Q &=& \left(Q + D + \left( \frac {S^2 - (Q+D)^2 }{ 2} \right)F
\right) \cdot Q \\ 
&=& Q^2 + D \cdot Q + \frac {S^2 - Q^2 - D^2 }{ 2} - D \cdot Q \\ &=& S^2 -
D^2 / 2.
\end{eqnarray*}
\epf

To calculate 
$$
\#Z_n(d,\vi)_{\sum m^j_e q^j_e} = S^2 - (H - \sum m^j_e Q^j_e)^2 / 2,
$$
we need to calculate $H^2$, $H \cdot Q^j_e$,
and 
$Q^j_e \cdot Q^{j'}_{e'}$, and these correspond to simpler enumerative
problems. 

If $(e,j) \neq (e',j')$, $Q^j_e$ could intersect $Q^{j'}_{e'}$ in two ways.
If $e' + e \geq n$, the elliptic curve could pass through $\De^j_e \cap
\De^{j'}_{e'}$, which will happen $\# W (\ce') / d!$ times (where $d' = d$, 
$\vec{h'} = d \vep_{1,n-1}$, $\vec{i'} = \vi - \vep_e - \vep_{e'} +
\vep_{e+e'-n}$).  Or the curve could break into two intersecting
components, one 
rational (call it $R$) containing $Q^j_e$ and $Q^{j'}_{e'}$ (which will be blown down in
the construction of $\cf$), and the other (call it $E$) smooth elliptic.
This will happen  
$$
\sum_{\substack{{d(0) + d(1) = d }\\{ \vi(0) + \vi(1) = \vi}}}
(d(0) d(1))^{\de_{n,2}} 
\binom { \vi - \vep_e - \vep_{e'}} {\vi(1)} \left( \frac {\#
X(\ce(0)) }{ d(0)!} \right) \left( \frac 
{\# W(\ce(1)) }{ d(1)!} \right)
$$
times where $h_{n-1}(0) = d(0)$, $h_{n-1}(1) = d(1)$.  The factor of
$(d(0) d(1))^{\de_{n,2}}$ corresponds to the fact when $n=2$, $\pi(C) =
\pi(R \cup E)$ is a plane curve, and the point $R \cap E$ could map to any
node of the plane curve $\pi(R \cup E)$. 
Transversality in both cases is simple to check, and both possibilities are
of the right dimension. 
Thus
$$
Q^j_e \cdot Q^{j'}_{e'} = \frac {\# W(\ce') }{ d!} \quad \quad \quad \quad
\quad \quad
$$
$$
 + 
\sum_{ \substack{{ d(0) + d(1) = d }\\ { \vi(0) + \vi(1) = \vi}}} (d(0)
d(1))^{\de_{n,2}} 
\binom { \vi - \vep_e - \vep_{e'} }{ \vi(1)} \left( \frac {\#
X(\ce(0))}{ d(0)!} \right) \left(\frac {\# W(\ce(1)) }{d(1)!} \right). 
$$

To determine $H \cdot Q^j_e$, fix a general hyperplane
$h$ in $\proj^n$, and let $H$ be its pullback to the fibration $\cf$.  Then
$H$ is a 
multisection of the elliptic fibration.  The cycle $H$ could intersect $Q^j_e$
in two ways.  Either $\pi(q^j_e) \in h \cap \De^j_e$ --- which will happen
$\# W(\ce') / d!$ times with $(d', \vec{h'}, \vec{i'}) = (d, d
\vep_{1,n-1},
\vi - \vep_e + \vep_{e-1})$ --- or the curve breaks into two
pieces, one 
rational containing a point of $h$ and $p^j_e$, which will happen 
$$
\sum_{ \substack{{ d(0) + d(1) = d }\\{ \vi(0) + \vi(1) = \vi}}} (d(0)
d(1))^{\de_{n,2}} 
d(0) \binom { \vi - \vep_e }{ \vi(1)} \left(
\frac{ \# X(\ce(0))}{d(0)!} \right) \left(\frac {\# W(\ce(1)) }{ d(1)!}
\right) 
$$
times where $h_{1,n-1}(0) = d(0)$, $h_{1,n-1}(1) = d(1)$.  (The second $d(0)$
in the formula comes from the choice of point of $h$ on the degree $d(0)$
rational component.)  Thus
$$
H \cdot Q^j_e = 
\frac { \# W(\ce') }{  d! } \quad \quad \quad \quad \quad \quad 
$$
$$
 + 
\sum_{ \substack {{d(0) + d(1) = d }\\{ \vi(0) + \vi(1) = \vi}}}
(d(0)d(1))^{\de_{n,2}} d(0)
\binom { \vi - \vep_e }{\vi(1)} \left(
\frac {\# X(\ce(0))}{ d(0)!} \right) \left(\frac {\# W(\ce(1)) }{ d(1)!}
\right). 
$$

To determine $H^2$, fix a second general hyperplane $h'$ in $\proj^n$, and
let $H'$ be its pullback to $\cf$.  Once again, $H$ could intersect $H'$ in
two ways depending on if the curve passes through $h \cap h'$, or if the
curve breaks into two pieces.  Thus
$$
H^2 = \frac {\# W (\ce') }{ d! } \quad \quad \quad \quad \quad \quad
$$
$$
+ \sum_{ \substack{{ d(0) + d(1) = d }\\{ \vi(0) + \vi(1) = \vi}}}
(d(0)d(1))^{\de_{n,2}} 
d(0)^2 
\binom { \vi }{ \vi(1)} \left(
\frac {\# X(\ce(0))}{ d(0)!} \right) \left( \frac {\# W(\ce(1)) }{ d(1)!}
\right) 
$$
where $\vec{h'} = d \vep_{1,n-1}$, $\vec{i'} = \vi + \vep_{n-2}$, $\vh(0)
= d(0) \vep_{1,n-1}$, $\vh(1) = d(1) \vep_{1,n-1}$.

The self-intersection of a section $S^2$ ($= (Q^j_e)^2$) can be calculated
as follows.
Fix $e$ such that $i_e>0$.
We can calculate $H \cdot Q^1_e$, so if we can evaluate $(H - Q^1_e)
\cdot Q^1_e$ then we can find $S^2 = (Q^1_e)^2$.  Fix a general hyperplane $h$ containing $\De^1_e$, and let
$(H-Q^1_e)$ be the multisection that is the pullback
of $h$ to $\cf$, minus the section $Q^1_e$.  The cycle $(H-Q^1_e)$ intersects
$Q^1_e$ if the curve is tangent to $h$ along $\De^1_e$ or if the curve breaks into
two pieces, with $Q^1_e$ on the rational piece.  Thus
$$
(H-Q^1_e) \cdot Q^1_e = \frac { \# W(\ce') }{ (d-2)! }  \quad \quad \quad
\quad \quad \quad
$$
$$
+
\sum_{ \substack{{ d(0) + d(1) = d }\\{ \vi(0) + \vi(1) = \vi -
\vep_e}}} (d(0)d(1))^{\de_{n,2}} 
\binom { \vi - \vep_e }{  \vi(1)} \left( \frac { 
\# X(\ce(0)) }{ ( d(0)-2)!} \right) \left(\frac{\# W(\ce(1))}{  d(1)!} \right)
$$
where $\vec{h'} = (d-2) \vep_{1,n-1} + \vep_{2,e}$, $\vh(0) = ( d(0) - 1)
\vep_{1,n-1} + \vep_{1,e}$, $\vh(1) = d(1) \vep_{1,n-1}$.  The denominator
$(d(0)-2)!$ arises because we have a degree $d(0)$ (rational) curve passing
through an $e$-plane on $h$, and various incidence conditions $\vi(0)$.
The number of such curves with a choice of one of the  other intersections
of $C(0)$ with 
$h$ is $(d(0)-1) \# X(\ce(0)) / (d(0) - 1)!$.

As an example, consider the elliptic quartics in $\proj^2$ passing through
11 fixed points, including $q^1_0$, $q^2_0$, $q^3_0$, $q^4_0$.  How many
such two-nodal quartics have $\oh(1) = q^1_0 + \dots + q^4_0$ in the Picard
group of the normalization of the curve?  We construct the fibration $\cf$ over
the (normalized) variety of two-nodal plane quartics through 11 fixed
points.  We have sections $Q^1_0$, \dots, $Q^{11}_0$ and a multisection $H$.
If $i \neq j$, $Q^i_0 \cdot Q^j_0 = 3$, $H \cdot Q^j_0 = 30$, $H ^2 = 225 +
3 {\binom {11} 2} = 390$, $(H-Q^1_0) \cdot Q^1_0 = 185$, so
$$
S^2 = H \cdot Q^1_0  - (H - Q^1_0) \cdot Q^1_0 = -155.
$$
Let $D = H - Q^1_0 - \dots - Q^4_0$.  Then
\begin{eqnarray*}
D^2 &=& H^2 + 4 S^2 - 8 H \cdot Q^1_0 + 12 Q^1_0 \cdot Q^2_0 \\
&=& 390 + 4 (-155) - 8 (30) + 12(3) \\
&=& -434
\end{eqnarray*}
so the answer is $S^2 - D^2/2 = 62$.

To determine the enumerative geometry of quartic elliptic space curves
(see Subsubsection \ref{eqesc}), various $\# Z_2(d,\vi)_\cd$ were needed with $d=3$
and $d=4$.  When $d=3$, $i_0$ is necessarily 8, and the results are given
in the Table \ref{ez3}.  For convenience, we write $p^j = q^j_0$ for the
base points and $l^j = q^j_1$ for the marked points on lines.  These values
were independently confirmed by M. Roth ([Ro]).
When $d=4$, $i_0$ must be 11, and the results are given in the Table \ref{ez4}.
For convenience again, we write $p^j = q^j_0$ and $l^j = q^j_1$.
\begin{table}
\begin{center}
\begin{tabular}{c|c|c}
$i_1$ & $\cd$ & $\# Z_2(d,\vi)_\cd$ \\
\hline
0 & $p^1 + p^2 + p^3$ & 0 \\
1 & $p^1 + p^2 + l^1$ & 1 \\
2 & $p^1 + l^1 + l^2$ & 5 \\
3 & $l^1 + l^2 + l^3$ & 18 \\
0 & $p^1 + 2 p^2    $   & 1 \\
1 & $2p^1 + l^1     $   & 4 \\
1 & $p^1 + 2l^1     $   & 5 \\
2 & $l^1 + 2l^2     $   & 16 \\
0 & $3p^1           $     & 3 \\
1 & $3l^1           $     & 14 \\
0 & $p^1 + p^2 + p^3 + p^4 - p^5 $& 1 \\
1 & $p^1 + p^2 + p^3 + p^4 - l^1 $& 2 \\
1 & $p^1 + p^2 + p^3 + l^1 - p^4 $& 4 \\
2 & $p^1 + p^2 + p^3 + l^1 - l^2 $& 10 \\
2 & $p^1 + p^2 + l^1 + l^2 - p^3 $& 14 \\
3 & $p^1 + p^2 + l^1 + l^2 - l^3 $& 39 \\
3 & $p^1 + l^1 + l^2 + l^3 - p^2 $& 45 \\
4 & $p^1 + l^1 + l^2 + l^3 - l^4 $& 135 \\
4 & $l^1 + l^2 + l^3 + l^4 - p^1 $& 135 \\
5 & $l^1 + l^2 + l^3 + l^4 - l^5 $& 432 
\end{tabular}
\end{center}
\caption{Counting cubic elliptic plane curves with a divisorial condition}
\label{ez3}
\end{table}

\begin{table}
\begin{center}
\begin{tabular}{c|c|c}
$i_1$ & $\cd$ & $\# Z_2(d,h,\vi)_\cd$ \\
\hline
0 & $p^1 + p^2 + p^3 + p^4 $ & 62 \\
1 & $p^1 + p^2 + p^3 + l^1 $ & 464 \\
2 & $p^1 + p^2 + l^1 + l^2 $ & 2,522 \\
3 & $p^1 + l^1 + l^2 + l^3 $ & 11,960 \\
4 & $l^1 + l^2 + l^3 + l^4 $ & 52,160 
\end{tabular}
\end{center}
\caption{Counting quartic elliptic plane curves with a divisorial condition}
\label{ez4}
\end{table}

\subsection{Examples}

\subsubsection{Plane curves}
Type IIc
components in this case are never enumeratively meaningful, as the elliptic
curve $C(0)$ must map to the line $H$ with degree at least two.  The
recursive formulas we get are identical to the genus 1 recursive
formulas of Caporaso and Harris in [CH3].

\subsubsection{Cubic elliptic space curves}
The number of smooth cubic elliptic space curves through $j$ general points
and $12-2j$ general lines is 1500, 150, 14, and 1 for $j= 0$, 1, 2, and 3
respectively.  (The number is 0 for $j>3$ as cubic elliptic space curves
must lie in a plane.)  The degenerations involved in calculating the first
case appeared in Subsubsection \ref{ecubics}.  As the Chow ring of the space
of smooth elliptic cubics is not hard to calculate (see [H1], p. 36), these
results may be easily verified.

The number of cubic elliptics tangent to $H$, through $j$ general points
and $11-2j$ general lines is 4740, 498, 50, and 4 for $j=0$, 1, 2, and 3
respectively.  The number of cubic elliptics triply tangent to $H$ through
$j$ general points and $10-2j$ general lines is 2790, 306, 33, and 3 for
$j=0$, 1, 2, and 3 respectively. 

These numbers are needed for the next examples.

\begin{table}
\begin{center}
\begin{tabular}{c|r}
$j$ & \# quartics \\
\hline
0 & 52,832,040 \\
1 & 4,436,208 \\
2 & 385,656 \\
3 & 34,674 \\
4 & 3,220 \\
5 & 310 \\
6 & 32 \\
7 & 4 \\
8 & 1
\end{tabular}
\end{center}
\caption{Number of quartic elliptic space curves through $j$ general points
and $16-2j$ general lines}
\label{eqescnums}
\end{table}

\begin{table}
\begin{center}
\begin{tabular}{r|c|c|c|r}
&$(i_0,i_1,h_0,h_1)$ & point & line & \# curves\\
&  & degen. & degen. & \\
\hline
1 & (16,0,0,0) & & 10 & 52,832,040 \\
2 & (14,1,0,0) & 40 & 11 & 4,436,268 \\
3 & (12,2,0,0) & 41 & 12 & 385,656 \\
4 & (10,3,0,0) & 42 & 13 & 34,674 \\
5 & (8,4,0,0) & 43 & 14 & 3,220 \\
6 & (6,5,0,0) & 44 & 15 & 310 \\
7 & (4,6,0,0) & 45 & 16 & 32 \\
8 & (1,7,0,0) & 46 & 17 & 4 \\
9 & (0,8,0,0) & 47 & & 1 \\
10 & (15,0,1,0) & & 18, 40 & 52,832,040 \\
11 & (13,1,1,0) & 48 & 19, 41 & 4,436,268 \\
12 & (11,2,1,0) & 49 & 20, 42 & 385,656 \\
13 & (9,3,1,0) & 50 & 21, 43 & 34,674 \\
14 & (7,4,1,0) & 51 & 22, 44 & 3,220 \\
15 & (5,5,1,0) & 52 & 23, 45 & 310 \\
16 & (3,6,1,0) & 53 & 24, 46 & 32 \\
17 & (1,7,1,0) & 54 & 25, 47 & 4 \\
18 & (14,0,2,0) & & 26, 48+ & 48,395,772 \\
19 & (12,1,2,0) & 55+ & 27, 49+ & 4,050,612 \\
20 & (10,2,2,0) & 56+ & 28, 50+ & 350,982 \\
21 & (8,3,2,0) & 57+ & 29, 51+ &31,454 \\
22 & (6,4,2,0) & 58+ & 30, 52 & 2,910 \\
23 & (4,5,2,0) & 59 & 31, 53 & 278 \\
24 & (2,6,2,0) & 60 & 32, 54 & 28 \\
25 & (0,7,2,0) &  61 &  & 3 \\
26 & (13,0,3,0) & & 33, 55+ & 39,347,736 \\
27 & (11,1,3,0) & 62+ & 34, 56+ & 3,266,100 \\
28 & (9,2,3,0) & 63+ & 35, 57+ & 280,752 \\
29 & (7,3,3,0) & 64+ & 36, 58+ & 24,972 \\
30 & (5,4,3,0) & 65+ & 37, 59+ & 2,290 \\
31 & (3,5,3,0) & 66+ & 38, 60+ & 214 \\
32 & (1,6,3,0) & 67+ & 39, 61+ & 20 \\
33 & (12,0,4,0) & & 62+ & 23,962,326 \\
34 & (10,1,4,0) & + & 63+ & 1,939,857 
\end{tabular}
\end{center}
\caption{Quartic elliptic space curves with incidence conditions}
\label{eqesc1}
\end{table}

\begin{table}
\begin{center}
\begin{tabular}{r|c|c|c|r}
&$(i_0,i_1,h_0,h_1)$ & point & line & \# curves\\
&  & degen. & degen. & \\
\hline
35 & (8,2,4,0) & + & 64+ & 161,735 \\
36 & (6,3,4,0) & + & 65+ & 13,908 \\
37 & (4,4,4,0) & + & 66+ & 1,222 \\
38 & (2,5,4,0) & + & 67+ & 104 \\
39 & (0,6,4,0) & + & & 8 \\
40 & (14,0,0,1) & & 48+ & 4,436,268 \\
41 & (12,1,0,1) & 68 & 49 & 385,656 \\
42 & (10,2,0,1) & 69 & 50 & 34,674 \\
43 & (8,3,0,1) & 70 & 51 & 3,220 \\
44 & (6,4,0,1) & 71 & 52 & 310 \\
45 & (3,5,0,1) & 72 & 53 & 32 \\
46 & (2,6,0,1) & 73 & 54 & 4 \\
47 & (0,7,0,1) & 74 & & 1 \\
48 & (13,0,1,1) & & 55, 68+ & 4,436,268 \\
49 & (11,1,1,1) & 75+ & 56, 69+ & 385,656 \\
50 & (9,2,1,1) & 76+ & 57, 70+ & 34,674 \\
51 & (7,3,1,1) & 77+ & 58, 71+ & 3,220 \\
52 & (5,4,1,1) & 78+ & 59, 72 & 310 \\
53 & (3,5,1,1) & 79 & 60, 73 & 32 \\
54 & (1,6,1,1) & 80 & 61, 74 & 4 \\
55 & (12,0,2,1) & & 62, 75 & 4,028,112 \\
56 & (10,1,2,1) & 81+ & 63, 76+ & 349,032 \\
57 & (8,2,2,1) & 82+ & 64, 77+ & 28,340 \\
58 & (6,3,2,1) & 83+ & 65, 78+ & 2,901 \\
59 & (4,4,2,1) & 84+ & 66, 79+ & 278 \\
60 & (2,5,2,1) & 85+ & 67, 80+ & 28 \\
61 & (0,6,2,1) & 86+ & & 3 \\
62 & (11,0,3,1) & & 81+ & 2,849,436 \\
63 & (9,1,3,1) & + & 82+ & 243,507 \\
64 & (7,2,3,1) & + & 83+ & 21,310 \\
65 & (5,3,3,1) & + & 84+ & 1,909 \\
66 & (3,4,3,1) & + & 85+ & 172 \\
67 & (1,5,3,1) & + & 86+ & 14 \\
68 & (12,0,0,2) & & 75+ & 385,656 
\end{tabular}
\end{center}
\caption{Quartic elliptic space curves with incidence conditions, cont'd}
\end{table}

\begin{table}
\begin{center}
\begin{tabular}{r|c|c|c|r}
&$(i_0,i_1,h_0,h_1)$ & point & line & \# curves\\
&  & degen. & degen. & \\
\hline
69 & (10,1,0,2) & 87 & 76+ & 34,674 \\
70 & (8,2,0,2) & 88 & 77+ & 3,220 \\
71 & (6,3,0,2) & 89 & 78+ & 310 \\ 
72 & (4,4,0,2) & 90 & 79 & 32 \\
73 & (2,5,0,2) & 91 & 80 & 4 \\
74 & (0,6,0,2) & 92 & & 1 \\
75 & (11,0,1,2) & & 81, 87+ & 384,156 \\
76 & (9,1,1,2) & 93+ & 82, 88+ & 34,524 \\
77 & (7,2,1,2) & 94+ & 83, 89+ & 3,206 \\
78 & (5,3,1,2) & 95+ & 84, 90+ & 309 \\
79 & (3,4,1,2) & 96 & 85, 91+ & 32 \\
80 & (1,5,1,2) & 97 & 86, 92+ & 4 \\
81 & (10,0,2,2) & & 93+ & 312,348 \\
82 & (8,1,2,2) & + & 94+ & 28,340 \\
83 & (6,2,2,2) & + & 95+ & 2,612 \\
84 & (4,3,2,2) & + & 96+ & 246 \\
85 & (2,4,2,2) & + & 97+ & 24 \\
86 & (0,5,2,2) & + & & 2 \\
87 & (10,0,0,3) & & 93+ & 34,674 \\
88 & (8,1,0,3) & 98+ & 94+ & 3,220 \\
89 & (6,2,0,3) & 99+ & 95+ & 310 \\
90 & (4,3,0,3) & 100+ & 96 & 32 \\
91 & (2,4,0,3) & 101 & 97 & 4 \\
92 & (0,5,0,3) & 102 & & 1 \\
93 & (9,0,1,3) & & 98+ & 31,056 \\
94 & (7,1,1,3) & + & 99+ & 3,052 \\
95 & (5,2,1,3) & + & 100+ & 304 \\
96 & (3,3,1,3) & + & 101+ & 32 \\
97 & (1,4,1,3) & + & 102+ & 4 \\
98 & (8,0,0,4) & & + &2,519 \\
99 & (6,1,0,4) & + & + & 277 \\
100 & (4,2,0,4) & + & + & 31 \\
101 & (2,3,0,4) & + & + & 4 \\
102 & (0,4,0,4) & + & &  1 
\end{tabular}
\end{center}
\caption{Quartic elliptic space curves with incidence conditions, cont'd}
\label{eqesc3}
\end{table}

\subsubsection{Quartic elliptic space curves}
\label{eqesc}

The number of smooth quartic elliptic space curves through $j$ general
points and $16-2j$ general lines is given in the Table \ref{eqescnums}.
These numbers agree with those recently found by Getzler by means of genus
1 Gromov-Witten invariants (cf. [G3]).  


\begin{figure}
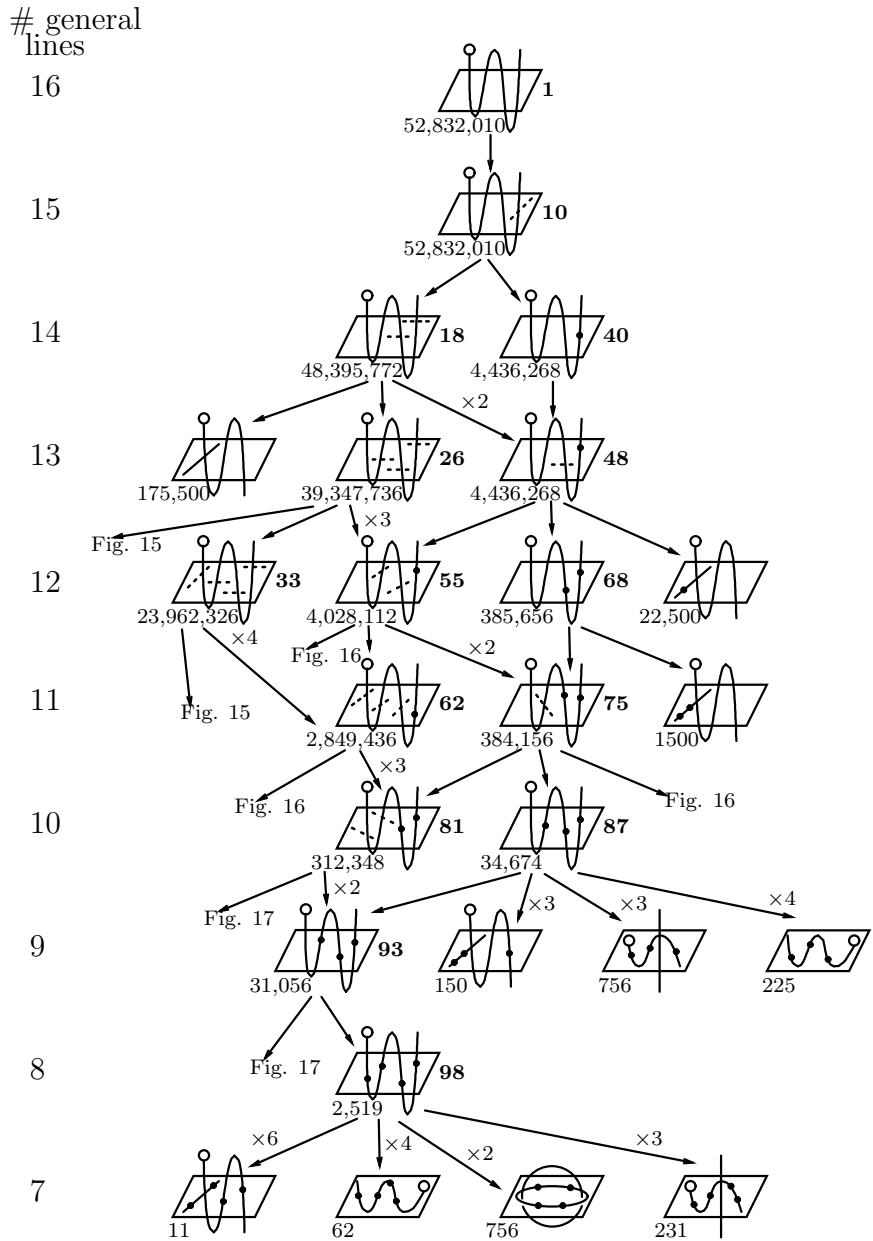

\begin{center}

	   \setlength{\unitlength}{.1\baseunit}
	    \input qesc1.tex 
\caption{Counting quartic elliptic space curves through 16 general lines}
\label{qesc1}
\end{center}
\end{figure}

The space of smooth quartic elliptic space curves is birational
to pencils in the space of space quadrics (as a quartic elliptic space
curve is the base locus of a unique pencil, and a general pencil defines a
smooth quartic elliptic).  By this means the last four numbers in Table
\ref{eqescnums} may be easily calculated.  D. Avritzer and I. Vainsencher used this
method (cf. [AV]) to calculate the top number, although they likely misprinted
their answer ([G5]).

\begin{figure}
\begin{center}

	   \setlength{\unitlength}{.1\baseunit}
	    \input qesc2.tex 
\caption{Additional degenerations of quartic elliptic space curves}
\label{qesc2}
\end{center}
\end{figure}

Other enumerative data can also be found.  For example, Tables \ref{eqesc1}
to \ref{eqesc3} 
give the number of smooth quartic elliptic space curves through $i_0$
general points and $i_1$ general lines, and $h_0$ general points and
$h_1$ general lines in $H$, with 
$$
2 i_0 + i_1 + 2 h_0 + h_1= 16.
$$
At each stage, the number may be computed by degenerating a point or a line
(assuming there is a point or line to degenerate).  Each row is labeled,
and the labels of the different degenerations that are also smooth
quartics are given in each case, and a ``+''
is added if there are other degenerations.  (This will help the reader to
follow through the degenerations.)  Keep in mind that these numbers
are not quite what the algorithm of this section produces; in the
algorithm, the intersections with $H$ are labeled, so the number computed
for $(i_0,i_1,h_0,h_1)$ will be $(4-h_0-h_1)!$ times the number in the
table.

These computations are not as difficult as one might think. For example, if
$i_0$ and $i_1$ are both positive, it is possible to degenerate a point and
then a line, or a line and then a point.  Both methods must yield the same
number, providing a means of double-checking.

\begin{figure}
\begin{center}

	   \setlength{\unitlength}{.1\baseunit}
	    \input qesc3.tex 
\caption{Additional degenerations of quartic elliptic space curves, cont'd}
\label{qesc3}
\end{center}
\end{figure}

\begin{figure}
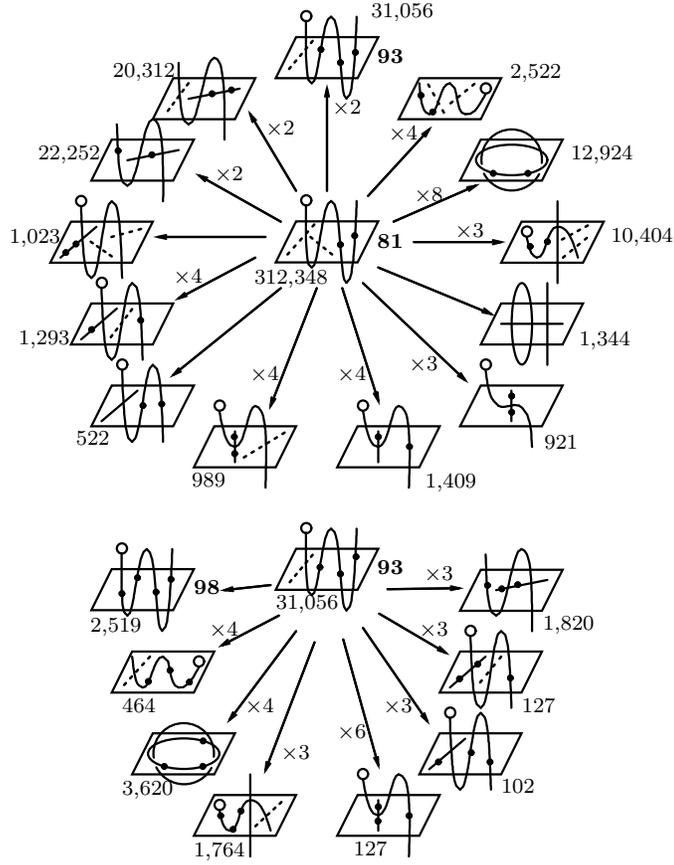

\begin{center}

	   \setlength{\unitlength}{.1\baseunit}
	    \input qesc4.tex 
\caption{Additional degenerations of quartic elliptic space curves, cont'd}
\label{qesc4}
\end{center}
\end{figure}

As an example, the degenerations used to compute the 52,832,040 quartic
space curves through 16 general lines are displayed in Figures \ref{qesc1}
to \ref{qesc4}, using the pictorial shorthand described earlier.  In
Figure \ref{qesc1}, degenerations involving nondegenerate
quartic elliptic space curves are given (as well as a few more).  The
remaining degenerations are given in Figures \ref{qesc2} to
\ref{qesc4}.  The boldfaced numbers
indicated the corresponding rows in Tables \ref{eqesc1} to
\ref{eqesc3}.

\subsection{Curves of higher genus}
\label{ehighgenus}
The genus 2 case seems potentially tractable.  An analog of Proposition
\ref{ecodim2} is needed, showing that in the space of stable maps of the
desired sort (with general assigned incidences and intersections with $H$),
no divisor representing maps with collapsed components of positive
genus is enumeratively meaningful.  New types of components arise,
including one in which $C(0)$ and $C(1)$ are both rational, and intersect
each other in 3 points (which will require the intersection theory of a
blow-up of $H^3$), and one in which $C(0)$, $C(1)$, and $C(2)$ are rational
and $C(0)$ intersects $C(i)$ ($i = 1,2$) in 2 points (which will require
the intersection theory of $(\Bl_\De H^2)^2$).  The main difficulty will
arise from components where $p_a(C(0)) = 2$, as the divisorial condition
analogous to that of Type IIc components is a codimension 2 condition, and
the calculation of $\# Z$ (when $p_a(C(0))=1$) using elliptic fibrations
now involves fibrations of abelian surfaces.

For genus greater than 2, the situation is more grave.  The map 
$$
\rho_A:  \mbar_{g,m}(\proj^n,d) \dashrightarrow \mbar_{g,m}(\proj^1,d)
$$
induced by projection from a general $(n-2)$-plane $A$ in $H$ is not
dominant, so the multiplicity calculations are no longer immediate from
the situation on $\cmbar_{g,m}(\proj^1,d)$.  Second, the divisorial
condition is even more complicated than for genus 2, and the other
computations (involving the intersection theory of products of repeated blow-ups of
powers of $H$ along various diagonals) will be horrendous.  It is also
awkward that the dimensions of these spaces may not be what one would naively
expect.  For example, the space of genus 3 quartic space curves is of
dimension 17 (as all genus 3 quartic space curves must lie in a plane, so
the dimension is $\dim 
\proj^{3*}$ plus the dimension of the space of plane quartics), not 16.
(This is because the normal bundle of the general such map has non-zero
$H^1$.)

\begin{figure}
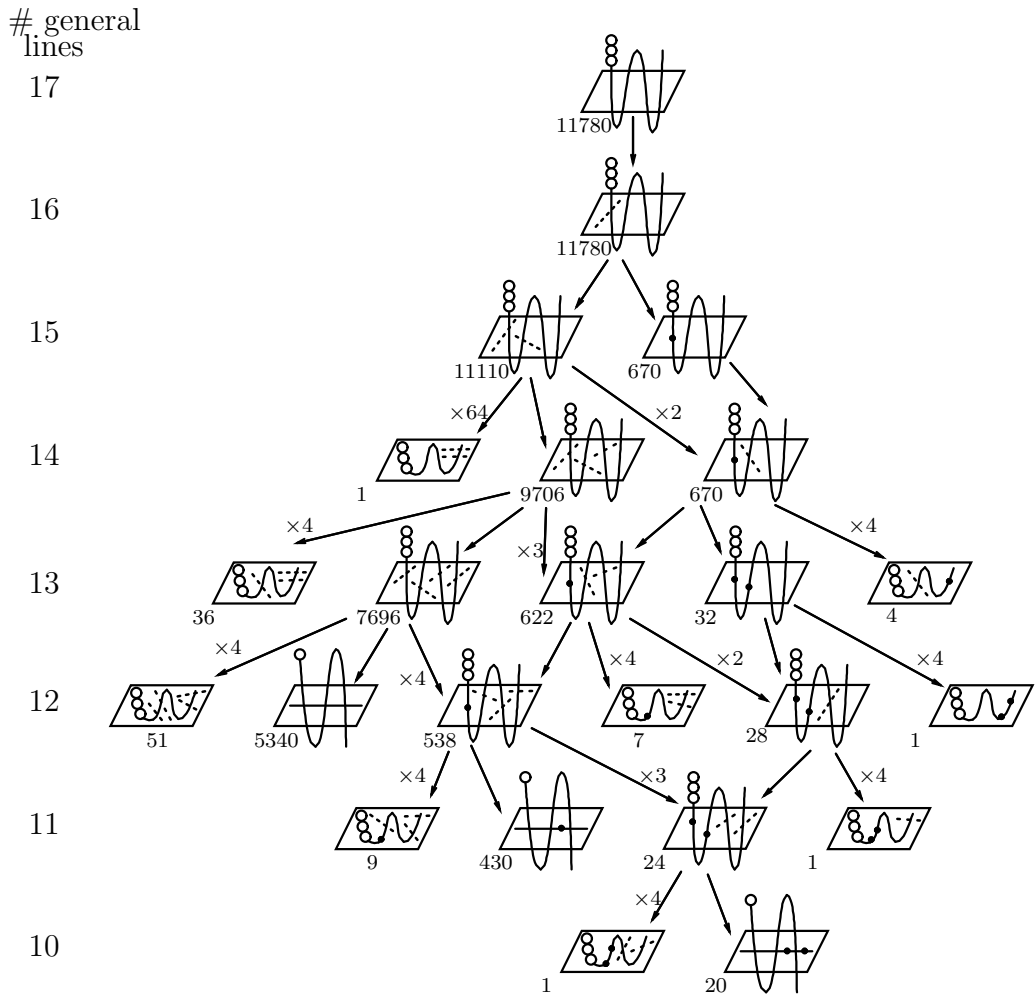

\begin{center}

	   \setlength{\unitlength}{.1\baseunit}
	    \input egenus3.tex 
\end{center}
\caption{Genus 3 quartic space curves.}
\label{egenus3}
\end{figure}

But all is not necessarily lost.  Even in the case of genus 3 quartic space
curves we can successfully follow through the degenerations (see
Figure \ref{egenus3}; genus 3 curves are indicated by 3 open circles).  The
unexpectedly high dimension of the space is compensated by unexpectedly
high dimensions of degenerations.  For example, one would naively expect
that requiring two curves to intersect in three points would impose three
conditions, but when one of the curves is cubic elliptic (and hence planar)
and the other is a line, the cost is only two conditions (as a line
intersecting the cubic at two points necessarily intersects it at a third).
This should work in general when $d$ is small enough (for fixed $g$) that
the curve must be planar, and the degenerations will look very much like
those of [CH3] and Figure \ref{egenus3}.  Perhaps when $d$ is small enough
that $C(0)$ must be genus 0 or 1 this analysis can still be carried
through.

\noindent
\address{Department of Mathematics \\ Princeton University \\ Fine Hall,
Washington Road \\ Princeton NJ 08544-1000}
\email{vakil@math.princeton.edu}
\end{document}

%% file: rcubics.tex
\begingroup\makeatletter\ifx\SetFigFont\undefined
\def\x#1#2#3#4#5#6#7\relax{\def\x{#1#2#3#4#5#6}}%
\expandafter\x\fmtname xxxxxx\relax \def\y{splain}%
\ifx\x\y   
\gdef\SetFigFont#1#2#3{%
  \ifnum #1<17\tiny\else \ifnum #1<20\small\else
  \ifnum #1<24\normalsize\else \ifnum #1<29\large\else
  \ifnum #1<34\Large\else \ifnum #1<41\LARGE\else
     \huge\fi\fi\fi\fi\fi\fi
  \csname #3\endcsname}%
\else
\gdef\SetFigFont#1#2#3{\begingroup
  \count@#1\relax \ifnum 25<\count@\count@25\fi
  \def\x{\endgroup\@setsize\SetFigFont{#2pt}}%
  \expandafter\x
    \csname \romannumeral\the\count@ pt\expandafter\endcsname
    \csname @\romannumeral\the\count@ pt\endcsname
  \csname #3\endcsname}%
\fi
\fi\endgroup
\begin{picture}(13332,14274)(0,-10)
\thicklines
\put(4748,7824){\ellipse{600}{1200}}
\path(4125,7815)(4725,7815)
\put(12698,6024){\ellipse{600}{1200}}
\path(12075,6015)(12675,6015)
\put(12698,4224){\ellipse{600}{1200}}
\path(12075,4215)(12675,4215)
\put(10148,2424){\ellipse{600}{1200}}
\path(9525,2415)(10125,2415)
\put(2948,624){\ellipse{600}{1200}}
\path(2325,615)(2925,615)
\put(3173,4224){\ellipse{600}{1200}}
\path(2550,4215)(3150,4215)
\put(5348,2424){\ellipse{600}{1200}}
\path(4725,2415)(5325,2415)
\put(3075,6053){\ellipse{900}{376}}
\path(3075,6090)(3075,6615)
\put(4650,4253){\ellipse{900}{376}}
\path(4650,4290)(4650,4815)
\put(6825,2453){\ellipse{900}{376}}
\path(6825,2490)(6825,3015)
\put(5250,653){\ellipse{900}{376}}
\path(5250,690)(5250,1215)
\put(6323,6024){\ellipse{600}{1200}}
\path(5700,6015)(6300,6015)
\put(7898,4224){\ellipse{600}{1200}}
\path(7275,4215)(7875,4215)
\path(8925,4215)(9900,4215)
\path(9225,4815)(9075,3615)
\path(9825,4815)(9600,3615)
\path(11175,2415)(12150,2415)
\path(11475,3015)(11325,1815)
\path(12075,3015)(11850,1815)
\path(945,6315)(645,5715)(1845,5715)
	(2145,6315)(945,6315)
\path(2595,6315)(2295,5715)(3495,5715)
	(3795,6315)(2595,6315)
\path(4170,6315)(3870,5715)(5070,5715)
	(5370,6315)(4170,6315)
\path(5745,6315)(5445,5715)(6645,5715)
	(6945,6315)(5745,6315)
\path(7320,6315)(7020,5715)(8220,5715)
	(8520,6315)(7320,6315)
\path(8970,6315)(8670,5715)(9870,5715)
	(10170,6315)(8970,6315)
\path(10545,6315)(10245,5715)(11445,5715)
	(11745,6315)(10545,6315)
\path(12120,6315)(11820,5715)(13020,5715)
	(13320,6315)(12120,6315)
\path(945,4515)(645,3915)(1845,3915)
	(2145,4515)(945,4515)
\path(2595,4515)(2295,3915)(3495,3915)
	(3795,4515)(2595,4515)
\path(4170,4515)(3870,3915)(5070,3915)
	(5370,4515)(4170,4515)
\path(5745,4515)(5445,3915)(6645,3915)
	(6945,4515)(5745,4515)
\path(7320,4515)(7020,3915)(8220,3915)
	(8520,4515)(7320,4515)
\path(8970,4515)(8670,3915)(9870,3915)
	(10170,4515)(8970,4515)
\path(10545,4515)(10245,3915)(11445,3915)
	(11745,4515)(10545,4515)
\path(12120,4515)(11820,3915)(13020,3915)
	(13320,4515)(12120,4515)
\path(1545,2715)(1245,2115)(2445,2115)
	(2745,2715)(1545,2715)
\path(3195,2715)(2895,2115)(4095,2115)
	(4395,2715)(3195,2715)
\path(4770,2715)(4470,2115)(5670,2115)
	(5970,2715)(4770,2715)
\path(6345,2715)(6045,2115)(7245,2115)
	(7545,2715)(6345,2715)
\path(7920,2715)(7620,2115)(8820,2115)
	(9120,2715)(7920,2715)
\path(9570,2715)(9270,2115)(10470,2115)
	(10770,2715)(9570,2715)
\path(11145,2715)(10845,2115)(12045,2115)
	(12345,2715)(11145,2715)
\path(6900,8115)(6600,7515)(7800,7515)
	(8100,8115)(6900,8115)
\path(9375,8115)(9075,7515)(10275,7515)
	(10575,8115)(9375,8115)
\path(4200,8115)(3900,7515)(5100,7515)
	(5400,8115)(4200,8115)
\path(5850,9915)(5550,9315)(6750,9315)
	(7050,9915)(5850,9915)
\path(8475,9915)(8175,9315)(9375,9315)
	(9675,9915)(8475,9915)
\path(7200,11715)(6900,11115)(8100,11115)
	(8400,11715)(7200,11715)
\path(7350,13515)(7050,12915)(8250,12915)
	(8550,13515)(7350,13515)
\path(4770,915)(4470,315)(5670,315)
	(5970,915)(4770,915)
\path(2370,915)(2070,315)(3270,315)
	(3570,915)(2370,915)
\path(7245,915)(6945,315)(8145,315)
	(8445,915)(7245,915)
\put(1425,6315){\ellipse{600}{600}}
\path(1350,6240)(1500,5790)
\path(1350,6240)(1500,5790)
\put(1500,4515){\ellipse{600}{600}}
\path(1425,4440)(1575,3990)
\path(1425,4440)(1575,3990)
\put(2100,2715){\ellipse{600}{600}}
\path(2025,2640)(2175,2190)
\path(2025,2640)(2175,2190)
\put(3600,2565){\blacken\ellipse{74}{74}}
\put(3600,2565){\ellipse{74}{74}}
\put(9555,7830){\blacken\ellipse{74}{74}}
\put(9555,7830){\ellipse{74}{74}}
\put(8085,585){\blacken\ellipse{74}{74}}
\put(8085,585){\ellipse{74}{74}}
\put(8625,9615){\blacken\ellipse{74}{74}}
\put(8625,9615){\ellipse{74}{74}}
\put(12225,6015){\blacken\ellipse{74}{74}}
\put(12225,6015){\ellipse{74}{74}}
\put(9105,6105){\blacken\ellipse{74}{74}}
\put(9105,6105){\ellipse{74}{74}}
\put(12525,4215){\blacken\ellipse{74}{74}}
\put(12525,4215){\ellipse{74}{74}}
\put(12225,4215){\blacken\ellipse{74}{74}}
\put(12225,4215){\ellipse{74}{74}}
\put(9375,4215){\blacken\ellipse{74}{74}}
\put(9375,4215){\ellipse{74}{74}}
\put(7425,4215){\blacken\ellipse{74}{74}}
\put(7425,4215){\ellipse{74}{74}}
\put(4275,4350){\blacken\ellipse{74}{74}}
\put(4275,4350){\ellipse{74}{74}}
\put(4875,2415){\blacken\ellipse{74}{74}}
\put(4875,2415){\ellipse{74}{74}}
\put(7050,2295){\blacken\ellipse{74}{74}}
\put(7050,2295){\ellipse{74}{74}}
\put(11700,2415){\blacken\ellipse{74}{74}}
\put(11700,2415){\ellipse{74}{74}}
\put(11250,2415){\blacken\ellipse{74}{74}}
\put(11250,2415){\ellipse{74}{74}}
\put(10050,2415){\blacken\ellipse{74}{74}}
\put(10050,2415){\ellipse{74}{74}}
\put(9675,2415){\blacken\ellipse{74}{74}}
\put(9675,2415){\ellipse{74}{74}}
\put(2475,615){\blacken\ellipse{74}{74}}
\put(2475,615){\ellipse{74}{74}}
\put(2850,615){\blacken\ellipse{74}{74}}
\put(2850,615){\ellipse{74}{74}}
\put(4875,765){\blacken\ellipse{74}{74}}
\put(4875,765){\ellipse{74}{74}}
\put(5175,465){\blacken\ellipse{74}{74}}
\put(5175,465){\ellipse{74}{74}}
\put(10703,6015){\blacken\ellipse{74}{74}}
\put(10703,6015){\ellipse{74}{74}}
\put(11310,6022){\blacken\ellipse{74}{74}}
\put(11310,6022){\ellipse{74}{74}}
\put(10688,4222){\blacken\ellipse{74}{74}}
\put(10688,4222){\ellipse{74}{74}}
\put(11295,4207){\blacken\ellipse{74}{74}}
\put(11295,4207){\ellipse{74}{74}}
\put(5798,4155){\blacken\ellipse{74}{74}}
\put(5798,4155){\ellipse{74}{74}}
\put(3465,4237){\blacken\ellipse{74}{74}}
\put(3465,4237){\ellipse{74}{74}}
\put(1440,4380){\blacken\ellipse{74}{74}}
\put(1440,4380){\ellipse{74}{74}}
\put(2040,2572){\blacken\ellipse{74}{74}}
\put(2040,2572){\ellipse{74}{74}}
\put(2153,2265){\blacken\ellipse{74}{74}}
\put(2153,2265){\ellipse{74}{74}}
\put(3330,2430){\blacken\ellipse{74}{74}}
\put(3330,2430){\ellipse{74}{74}}
\put(3938,2347){\blacken\ellipse{74}{74}}
\put(3938,2347){\ellipse{74}{74}}
\put(3240,690){\blacken\ellipse{74}{74}}
\put(3240,690){\ellipse{74}{74}}
\put(5640,2467){\blacken\ellipse{74}{74}}
\put(5640,2467){\ellipse{74}{74}}
\put(5625,750){\blacken\ellipse{74}{74}}
\put(5625,750){\ellipse{74}{74}}
\put(7965,2347){\blacken\ellipse{74}{74}}
\put(7965,2347){\ellipse{74}{74}}
\put(8228,2422){\blacken\ellipse{74}{74}}
\put(8228,2422){\ellipse{74}{74}}
\put(6570,2602){\blacken\ellipse{74}{74}}
\put(6570,2602){\ellipse{74}{74}}
\put(7635,727){\blacken\ellipse{74}{74}}
\put(7635,727){\ellipse{74}{74}}
\put(7320,495){\blacken\ellipse{74}{74}}
\put(7320,495){\ellipse{74}{74}}
\dottedline{135}(7500,11565)(7200,11265)
\dottedline{135}(6075,9840)(5775,9465)
\dottedline{135}(6150,9465)(6375,9615)
\dottedline{135}(7200,8040)(6825,7740)
\dottedline{135}(7200,7815)(7575,7665)
\dottedline{135}(7575,7890)(7950,8040)
\dottedline{135}(9750,7965)(9975,7665)
\dottedline{135}(6375,6165)(6750,6165)
\dottedline{135}(9225,6165)(9600,5865)
\dottedline{135}(9600,6165)(9975,6165)
\dottedline{135}(7875,4365)(8325,4440)
\dottedline{135}(10800,4365)(11175,4140)
\dottedline{135}(10200,2565)(10650,2640)
\path(7800,12690)(7800,12240)
\path(7800,12690)(7800,12240)
\path(7770.000,12360.000)(7800.000,12240.000)(7830.000,12360.000)
\path(7425,10740)(6750,10140)
\path(7425,10740)(6750,10140)
\path(6819.758,10242.146)(6750.000,10140.000)(6859.620,10197.301)
\path(7800,10815)(8400,10140)
\path(7800,10815)(8400,10140)
\path(8297.854,10209.758)(8400.000,10140.000)(8342.699,10249.620)
\path(9525,9090)(9750,8565)
\path(9525,9090)(9750,8565)
\path(9675.155,8663.480)(9750.000,8565.000)(9730.304,8687.115)
\path(6450,9015)(5475,8340)
\path(6450,9015)(5475,8340)
\path(5556.587,8432.971)(5475.000,8340.000)(5590.739,8383.639)
\path(6825,9015)(7050,8565)
\path(6825,9015)(7050,8565)
\path(6969.502,8658.915)(7050.000,8565.000)(7023.167,8685.748)
\path(6975,9165)(9225,8190)
\path(6975,9165)(9225,8190)
\path(9102.965,8210.186)(9225.000,8190.000)(9126.822,8265.240)
\path(10725,7890)(12825,6765)
\path(10725,7890)(12825,6765)
\path(12705.056,6795.222)(12825.000,6765.000)(12733.389,6848.111)
\path(10650,7590)(10950,6765)
\path(10650,7590)(10950,6765)
\path(10880.797,6867.523)(10950.000,6765.000)(10937.185,6888.028)
\path(10425,7365)(10200,6540)
\path(10425,7365)(10200,6540)
\path(10202.631,6663.665)(10200.000,6540.000)(10260.517,6647.878)
\path(6225,7440)(2175,6540)
\path(6225,7440)(2175,6540)
\path(2285.635,6595.317)(2175.000,6540.000)(2298.650,6536.746)
\path(6150,7215)(3900,6465)
\path(6150,7215)(3900,6465)
\path(4004.355,6531.408)(3900.000,6465.000)(4023.329,6474.487)
\path(7050,7215)(6825,6390)
\path(7050,7215)(6825,6390)
\path(6827.631,6513.665)(6825.000,6390.000)(6885.517,6497.878)
\path(7650,7215)(7800,6615)
\path(7650,7215)(7800,6615)
\path(7741.791,6724.141)(7800.000,6615.000)(7800.000,6738.693)
\path(8025,7665)(8850,6540)
\path(8025,7665)(8850,6540)
\path(8754.844,6619.028)(8850.000,6540.000)(8803.229,6654.510)
\path(11550,5415)(12375,4815)
\path(11550,5415)(12375,4815)
\path(12260.307,4861.319)(12375.000,4815.000)(12295.597,4909.843)
\path(11325,5415)(11325,4740)
\path(11325,5415)(11325,4740)
\path(11295.000,4860.000)(11325.000,4740.000)(11355.000,4860.000)
\path(2100,3390)(11700,3390)
\path(2100,3390)(11700,3390)
\path(11025,3390)(11025,3690)
\path(11025,3390)(11025,3690)
\path(2100,3390)(2100,3090)
\path(2100,3390)(2100,3090)
\path(2070.000,3210.000)(2100.000,3090.000)(2130.000,3210.000)
\path(4050,3390)(4050,2865)
\path(4050,3390)(4050,2865)
\path(4020.000,2985.000)(4050.000,2865.000)(4080.000,2985.000)
\path(5025,3390)(5025,2865)
\path(5025,3390)(5025,2865)
\path(4995.000,2985.000)(5025.000,2865.000)(5055.000,2985.000)
\path(6975,3390)(6975,2790)
\path(6975,3390)(6975,2790)
\path(6945.000,2910.000)(6975.000,2790.000)(7005.000,2910.000)
\path(8100,3390)(8100,2790)
\path(8100,3390)(8100,2790)
\path(8070.000,2910.000)(8100.000,2790.000)(8130.000,2910.000)
\path(9675,3390)(9675,2790)
\path(9675,3390)(9675,2790)
\path(9645.000,2910.000)(9675.000,2790.000)(9705.000,2910.000)
\path(11700,3390)(11700,2865)
\path(11700,3390)(11700,2865)
\path(11670.000,2985.000)(11700.000,2865.000)(11730.000,2985.000)
\path(4050,1665)(4575,1065)
\path(4050,1665)(4575,1065)
\path(4473.402,1135.554)(4575.000,1065.000)(4518.557,1175.064)
\path(4275,1740)(7200,1065)
\path(4275,1740)(7200,1065)
\path(7076.327,1062.751)(7200.000,1065.000)(7089.819,1121.215)
\path(975,5265)(10800,5265)
\path(975,5265)(10800,5265)
\path(9450,5565)(9450,5265)
\path(9450,5565)(9450,5265)
\path(975,5265)(975,4665)
\path(975,5265)(975,4665)
\path(945.000,4785.000)(975.000,4665.000)(1005.000,4785.000)
\path(2775,5265)(2775,4590)
\path(2775,5265)(2775,4590)
\path(2745.000,4710.000)(2775.000,4590.000)(2805.000,4710.000)
\path(4425,5265)(4425,4590)
\path(4425,5265)(4425,4590)
\path(4395.000,4710.000)(4425.000,4590.000)(4455.000,4710.000)
\path(6150,5265)(6150,4590)
\path(6150,5265)(6150,4590)
\path(6120.000,4710.000)(6150.000,4590.000)(6180.000,4710.000)
\path(7500,5265)(7500,4590)
\path(7500,5265)(7500,4590)
\path(7470.000,4710.000)(7500.000,4590.000)(7530.000,4710.000)
\path(9300,5265)(9300,4590)
\path(9300,5265)(9300,4590)
\path(9270.000,4710.000)(9300.000,4590.000)(9330.000,4710.000)
\path(10800,5265)(10800,4665)
\path(10800,5265)(10800,4665)
\path(10770.000,4785.000)(10800.000,4665.000)(10830.000,4785.000)
\path(6300,7065)(5400,6465)
\path(6300,7065)(5400,6465)
\path(5483.205,6556.526)(5400.000,6465.000)(5516.487,6506.603)
\path(3675,1740)(3375,1140)
\path(3401.833,1260.748)(3375.000,1140.000)(3455.498,1233.915)
\path(7275,6015)(8250,6015)
\path(7575,6615)(7425,5640)
\path(8175,6615)(7950,5415)
\path(7430,13816)	(7433.545,13774.251)
	(7437.063,13733.374)
	(7440.555,13693.362)
	(7444.022,13654.208)
	(7447.465,13615.906)
	(7450.885,13578.448)
	(7457.664,13506.039)
	(7464.369,13436.926)
	(7471.010,13371.054)
	(7477.598,13308.368)
	(7484.142,13248.812)
	(7490.653,13192.334)
	(7497.142,13138.876)
	(7503.618,13088.385)
	(7510.092,13040.805)
	(7516.574,12996.081)
	(7523.075,12954.159)
	(7529.605,12914.984)
	(7536.174,12878.500)
	(7549.470,12813.388)
	(7563.046,12758.383)
	(7576.985,12713.046)
	(7591.369,12676.938)
	(7655.000,12616.000)

\path(7655,12616)	(7694.506,12667.432)
	(7709.825,12726.470)
	(7716.592,12762.787)
	(7722.847,12803.053)
	(7728.652,12846.830)
	(7734.068,12893.680)
	(7739.157,12943.165)
	(7743.981,12994.849)
	(7748.602,13048.293)
	(7753.082,13103.059)
	(7757.482,13158.710)
	(7761.864,13214.807)
	(7766.290,13270.915)
	(7770.822,13326.593)
	(7775.522,13381.406)
	(7780.451,13434.915)
	(7785.671,13486.682)
	(7791.245,13536.271)
	(7797.233,13583.242)
	(7803.698,13627.158)
	(7810.702,13667.583)
	(7818.306,13704.077)
	(7835.563,13763.524)
	(7880.000,13816.000)

\path(7880,13816)	(7920.301,13801.276)
	(7958.043,13752.694)
	(7976.230,13714.539)
	(7994.107,13666.518)
	(8011.783,13608.166)
	(8029.370,13539.014)
	(8038.164,13500.242)
	(8046.976,13458.596)
	(8055.821,13414.016)
	(8064.712,13366.445)
	(8073.663,13315.824)
	(8082.687,13262.094)
	(8091.799,13205.198)
	(8101.012,13145.077)
	(8110.339,13081.672)
	(8119.796,13014.925)
	(8129.394,12944.779)
	(8139.148,12871.173)
	(8144.088,12833.056)
	(8149.072,12794.051)
	(8154.102,12754.153)
	(8159.180,12713.354)
	(8164.307,12671.646)
	(8169.484,12629.023)
	(8174.715,12585.477)
	(8180.000,12541.000)

\path(7314,12008)	(7317.545,11966.251)
	(7321.063,11925.374)
	(7324.555,11885.362)
	(7328.022,11846.208)
	(7331.465,11807.906)
	(7334.885,11770.448)
	(7341.664,11698.039)
	(7348.369,11628.926)
	(7355.010,11563.054)
	(7361.598,11500.368)
	(7368.142,11440.812)
	(7374.653,11384.334)
	(7381.142,11330.876)
	(7387.618,11280.385)
	(7394.092,11232.805)
	(7400.574,11188.081)
	(7407.075,11146.159)
	(7413.605,11106.984)
	(7420.174,11070.500)
	(7433.470,11005.388)
	(7447.046,10950.383)
	(7460.985,10905.046)
	(7475.369,10868.938)
	(7539.000,10808.000)

\path(7539,10808)	(7578.506,10859.432)
	(7593.825,10918.470)
	(7600.592,10954.787)
	(7606.847,10995.053)
	(7612.652,11038.830)
	(7618.068,11085.680)
	(7623.157,11135.165)
	(7627.981,11186.849)
	(7632.602,11240.293)
	(7637.082,11295.059)
	(7641.482,11350.710)
	(7645.864,11406.807)
	(7650.290,11462.915)
	(7654.822,11518.593)
	(7659.522,11573.406)
	(7664.451,11626.915)
	(7669.671,11678.682)
	(7675.245,11728.271)
	(7681.233,11775.242)
	(7687.698,11819.158)
	(7694.702,11859.583)
	(7702.306,11896.077)
	(7719.563,11955.524)
	(7764.000,12008.000)

\path(7764,12008)	(7804.301,11993.276)
	(7842.043,11944.694)
	(7860.230,11906.539)
	(7878.107,11858.518)
	(7895.783,11800.166)
	(7913.370,11731.014)
	(7922.164,11692.242)
	(7930.976,11650.596)
	(7939.821,11606.016)
	(7948.712,11558.445)
	(7957.663,11507.824)
	(7966.687,11454.094)
	(7975.799,11397.198)
	(7985.012,11337.077)
	(7994.339,11273.672)
	(8003.796,11206.925)
	(8013.394,11136.779)
	(8023.148,11063.173)
	(8028.088,11025.056)
	(8033.072,10986.051)
	(8038.102,10946.153)
	(8043.180,10905.354)
	(8048.307,10863.646)
	(8053.484,10821.023)
	(8058.715,10777.477)
	(8064.000,10733.000)

\path(5925,10215)	(5928.545,10173.251)
	(5932.063,10132.374)
	(5935.555,10092.362)
	(5939.022,10053.208)
	(5942.465,10014.906)
	(5945.885,9977.448)
	(5952.664,9905.039)
	(5959.369,9835.926)
	(5966.010,9770.054)
	(5972.598,9707.368)
	(5979.142,9647.812)
	(5985.653,9591.334)
	(5992.142,9537.876)
	(5998.618,9487.385)
	(6005.092,9439.805)
	(6011.574,9395.081)
	(6018.075,9353.159)
	(6024.605,9313.984)
	(6031.174,9277.500)
	(6044.470,9212.388)
	(6058.046,9157.383)
	(6071.985,9112.046)
	(6086.369,9075.938)
	(6150.000,9015.000)

\path(6150,9015)	(6189.506,9066.432)
	(6204.825,9125.470)
	(6211.592,9161.787)
	(6217.847,9202.053)
	(6223.652,9245.830)
	(6229.068,9292.680)
	(6234.157,9342.165)
	(6238.981,9393.849)
	(6243.602,9447.293)
	(6248.082,9502.059)
	(6252.482,9557.710)
	(6256.864,9613.807)
	(6261.290,9669.915)
	(6265.822,9725.593)
	(6270.522,9780.406)
	(6275.451,9833.915)
	(6280.671,9885.682)
	(6286.245,9935.271)
	(6292.233,9982.242)
	(6298.698,10026.158)
	(6305.702,10066.583)
	(6313.306,10103.077)
	(6330.563,10162.524)
	(6375.000,10215.000)

\path(6375,10215)	(6415.301,10200.276)
	(6453.043,10151.694)
	(6471.230,10113.539)
	(6489.107,10065.518)
	(6506.783,10007.166)
	(6524.370,9938.014)
	(6533.164,9899.242)
	(6541.976,9857.596)
	(6550.821,9813.016)
	(6559.712,9765.445)
	(6568.663,9714.824)
	(6577.687,9661.094)
	(6586.799,9604.198)
	(6596.012,9544.077)
	(6605.339,9480.672)
	(6614.796,9413.925)
	(6624.394,9343.779)
	(6634.148,9270.173)
	(6639.088,9232.056)
	(6644.072,9193.051)
	(6649.102,9153.153)
	(6654.180,9112.354)
	(6659.307,9070.646)
	(6664.484,9028.023)
	(6669.715,8984.477)
	(6675.000,8940.000)

\path(8568,10218)	(8571.545,10179.730)
	(8575.063,10142.259)
	(8582.022,10069.691)
	(8588.885,10000.244)
	(8595.664,9933.869)
	(8602.369,9870.515)
	(8609.010,9810.133)
	(8615.598,9752.670)
	(8622.142,9698.078)
	(8628.653,9646.306)
	(8635.142,9597.303)
	(8641.618,9551.019)
	(8648.092,9507.404)
	(8654.574,9466.408)
	(8661.075,9427.979)
	(8674.174,9358.625)
	(8687.470,9298.939)
	(8701.046,9248.518)
	(8729.369,9173.859)
	(8793.000,9118.000)

\path(8793,9118)	(8832.759,9165.186)
	(8848.354,9219.347)
	(8861.716,9289.601)
	(8867.706,9329.758)
	(8873.313,9372.734)
	(8878.595,9418.126)
	(8883.611,9465.533)
	(8888.419,9514.552)
	(8893.078,9564.782)
	(8897.646,9615.822)
	(8902.181,9667.269)
	(8906.743,9718.721)
	(8911.388,9769.778)
	(8916.177,9820.037)
	(8921.167,9869.096)
	(8926.416,9916.554)
	(8931.983,9962.009)
	(8937.927,10005.059)
	(8944.305,10045.302)
	(8958.601,10115.762)
	(8975.337,10170.175)
	(9018.000,10218.000)

\path(9018,10218)	(9059.731,10202.994)
	(9098.496,10154.210)
	(9117.040,10115.984)
	(9135.172,10067.914)
	(9153.002,10009.531)
	(9170.640,9940.369)
	(9179.421,9901.600)
	(9188.195,9859.961)
	(9196.976,9815.394)
	(9205.778,9767.840)
	(9214.614,9717.242)
	(9223.498,9663.540)
	(9232.443,9606.676)
	(9241.464,9546.593)
	(9250.575,9483.231)
	(9259.788,9416.532)
	(9269.118,9346.439)
	(9278.579,9272.891)
	(9283.362,9234.804)
	(9288.184,9195.832)
	(9293.044,9155.967)
	(9297.946,9115.203)
	(9302.891,9073.531)
	(9307.880,9030.945)
	(9312.916,8987.437)
	(9318.000,8943.000)

\path(10618,4884)	(10621.545,4845.730)
	(10625.063,4808.259)
	(10632.022,4735.691)
	(10638.885,4666.244)
	(10645.664,4599.869)
	(10652.369,4536.515)
	(10659.010,4476.133)
	(10665.598,4418.670)
	(10672.142,4364.078)
	(10678.653,4312.306)
	(10685.142,4263.303)
	(10691.618,4217.019)
	(10698.092,4173.404)
	(10704.574,4132.408)
	(10711.075,4093.979)
	(10724.174,4024.625)
	(10737.470,3964.939)
	(10751.046,3914.518)
	(10779.369,3839.859)
	(10843.000,3784.000)

\path(10843,3784)	(10882.759,3831.186)
	(10898.354,3885.347)
	(10911.716,3955.601)
	(10917.706,3995.758)
	(10923.313,4038.734)
	(10928.595,4084.126)
	(10933.611,4131.533)
	(10938.419,4180.552)
	(10943.078,4230.782)
	(10947.646,4281.822)
	(10952.181,4333.269)
	(10956.743,4384.721)
	(10961.388,4435.778)
	(10966.177,4486.037)
	(10971.167,4535.096)
	(10976.416,4582.554)
	(10981.983,4628.009)
	(10987.927,4671.059)
	(10994.305,4711.302)
	(11008.601,4781.762)
	(11025.337,4836.175)
	(11068.000,4884.000)

\path(11068,4884)	(11109.731,4868.994)
	(11148.496,4820.210)
	(11167.040,4781.984)
	(11185.172,4733.914)
	(11203.002,4675.531)
	(11220.640,4606.369)
	(11229.421,4567.600)
	(11238.195,4525.961)
	(11246.976,4481.394)
	(11255.778,4433.840)
	(11264.614,4383.242)
	(11273.498,4329.540)
	(11282.443,4272.676)
	(11291.464,4212.593)
	(11300.575,4149.231)
	(11309.788,4082.532)
	(11319.118,4012.439)
	(11328.579,3938.891)
	(11333.362,3900.804)
	(11338.184,3861.832)
	(11343.044,3821.967)
	(11347.946,3781.203)
	(11352.891,3739.531)
	(11357.880,3696.945)
	(11362.916,3653.437)
	(11368.000,3609.000)

\path(10650,6615)	(10653.545,6573.251)
	(10657.063,6532.374)
	(10660.555,6492.362)
	(10664.022,6453.208)
	(10667.465,6414.906)
	(10670.885,6377.448)
	(10677.664,6305.039)
	(10684.369,6235.926)
	(10691.010,6170.054)
	(10697.598,6107.368)
	(10704.142,6047.812)
	(10710.653,5991.334)
	(10717.142,5937.876)
	(10723.618,5887.385)
	(10730.092,5839.805)
	(10736.574,5795.081)
	(10743.075,5753.159)
	(10749.605,5713.984)
	(10756.174,5677.500)
	(10769.470,5612.388)
	(10783.046,5557.383)
	(10796.985,5512.046)
	(10811.369,5475.938)
	(10875.000,5415.000)

\path(10875,5415)	(10914.506,5466.432)
	(10929.825,5525.470)
	(10936.592,5561.787)
	(10942.847,5602.053)
	(10948.652,5645.830)
	(10954.068,5692.680)
	(10959.157,5742.165)
	(10963.981,5793.849)
	(10968.602,5847.293)
	(10973.082,5902.059)
	(10977.482,5957.710)
	(10981.864,6013.807)
	(10986.290,6069.915)
	(10990.822,6125.593)
	(10995.522,6180.406)
	(11000.451,6233.915)
	(11005.671,6285.682)
	(11011.245,6335.271)
	(11017.233,6382.242)
	(11023.698,6426.158)
	(11030.702,6466.583)
	(11038.306,6503.077)
	(11055.563,6562.524)
	(11100.000,6615.000)

\path(11100,6615)	(11140.301,6600.273)
	(11178.043,6551.688)
	(11196.230,6513.532)
	(11214.107,6465.511)
	(11231.783,6407.158)
	(11249.370,6338.006)
	(11258.164,6299.235)
	(11266.976,6257.588)
	(11275.821,6213.009)
	(11284.712,6165.438)
	(11293.663,6114.817)
	(11302.687,6061.088)
	(11311.799,6004.192)
	(11321.012,5944.071)
	(11330.339,5880.667)
	(11339.796,5813.921)
	(11349.394,5743.775)
	(11359.148,5670.170)
	(11364.088,5632.053)
	(11369.072,5593.049)
	(11374.102,5553.151)
	(11379.180,5512.352)
	(11384.307,5470.645)
	(11389.484,5428.022)
	(11394.715,5384.476)
	(11400.000,5340.000)

\path(7280,691)	(7285.392,649.893)
	(7290.997,612.218)
	(7303.055,546.703)
	(7316.591,493.541)
	(7332.020,451.816)
	(7370.225,399.022)
	(7421.000,381.000)

\path(7421,381)	(7467.654,396.553)
	(7502.768,441.464)
	(7531.304,506.209)
	(7544.657,543.042)
	(7558.225,581.260)
	(7572.630,619.673)
	(7588.492,657.091)
	(7627.069,724.176)
	(7678.917,772.987)
	(7749.000,794.000)

\path(7749,794)	(7816.204,788.162)
	(7883.121,766.438)
	(7946.989,732.088)
	(8005.050,688.370)
	(8054.542,638.543)
	(8092.704,585.867)
	(8116.777,533.599)
	(8124.000,485.000)

\path(8124,485)	(8082.064,417.438)
	(8037.908,391.757)
	(7983.000,381.000)

\path(7983,381)	(7912.177,393.699)
	(7854.464,425.938)
	(7814.268,468.207)
	(7796.000,511.000)

\path(7796,511)	(7802.968,568.470)
	(7847.643,630.394)
	(7886.142,665.088)
	(7936.496,703.371)
	(7999.512,746.066)
	(8036.021,769.327)
	(8076.000,794.000)

\path(7943,2503)	(7948.392,2461.893)
	(7953.997,2424.218)
	(7966.055,2358.703)
	(7979.591,2305.541)
	(7995.020,2263.816)
	(8033.225,2211.022)
	(8084.000,2193.000)

\path(8084,2193)	(8130.654,2208.553)
	(8165.768,2253.464)
	(8194.304,2318.209)
	(8207.657,2355.042)
	(8221.225,2393.260)
	(8235.630,2431.673)
	(8251.492,2469.091)
	(8290.069,2536.176)
	(8341.917,2584.987)
	(8412.000,2606.000)

\path(8412,2606)	(8479.204,2600.162)
	(8546.121,2578.438)
	(8609.989,2544.088)
	(8668.050,2500.370)
	(8717.542,2450.543)
	(8755.704,2397.867)
	(8779.777,2345.599)
	(8787.000,2297.000)

\path(8787,2297)	(8745.064,2229.438)
	(8700.908,2203.757)
	(8646.000,2193.000)

\path(8646,2193)	(8575.177,2205.699)
	(8517.464,2237.938)
	(8477.268,2280.207)
	(8459.000,2323.000)

\path(8459,2323)	(8465.968,2380.470)
	(8510.642,2442.394)
	(8549.142,2477.088)
	(8599.496,2515.371)
	(8662.512,2558.066)
	(8699.021,2581.327)
	(8739.000,2606.000)

\path(5768,4303)	(5773.392,4261.893)
	(5778.997,4224.218)
	(5791.055,4158.703)
	(5804.591,4105.541)
	(5820.020,4063.816)
	(5858.225,4011.022)
	(5909.000,3993.000)

\path(5909,3993)	(5955.654,4008.553)
	(5990.768,4053.464)
	(6019.304,4118.209)
	(6032.657,4155.042)
	(6046.225,4193.260)
	(6060.630,4231.673)
	(6076.492,4269.091)
	(6115.069,4336.176)
	(6166.917,4384.987)
	(6237.000,4406.000)

\path(6237,4406)	(6304.204,4400.162)
	(6371.121,4378.438)
	(6434.989,4344.088)
	(6493.050,4300.370)
	(6542.542,4250.543)
	(6580.704,4197.867)
	(6604.777,4145.599)
	(6612.000,4097.000)

\path(6612,4097)	(6570.064,4029.438)
	(6525.908,4003.757)
	(6471.000,3993.000)

\path(6471,3993)	(6400.177,4005.699)
	(6342.464,4037.938)
	(6302.268,4080.207)
	(6284.000,4123.000)

\path(6284,4123)	(6290.968,4180.470)
	(6335.643,4242.394)
	(6374.142,4277.088)
	(6424.496,4315.371)
	(6487.512,4358.066)
	(6524.021,4381.327)
	(6564.000,4406.000)

\path(4193,6103)	(4198.392,6061.893)
	(4203.997,6024.218)
	(4216.055,5958.703)
	(4229.591,5905.541)
	(4245.020,5863.816)
	(4283.225,5811.022)
	(4334.000,5793.000)

\path(4334,5793)	(4380.654,5808.553)
	(4415.768,5853.464)
	(4444.304,5918.209)
	(4457.657,5955.042)
	(4471.225,5993.260)
	(4485.630,6031.673)
	(4501.492,6069.091)
	(4540.069,6136.176)
	(4591.917,6184.987)
	(4662.000,6206.000)

\path(4662,6206)	(4729.204,6200.162)
	(4796.121,6178.438)
	(4859.989,6144.088)
	(4918.050,6100.370)
	(4967.542,6050.543)
	(5005.704,5997.867)
	(5029.777,5945.599)
	(5037.000,5897.000)

\path(5037,5897)	(4995.064,5829.438)
	(4950.908,5803.757)
	(4896.000,5793.000)

\path(4896,5793)	(4825.177,5805.699)
	(4767.464,5837.938)
	(4727.268,5880.207)
	(4709.000,5923.000)

\path(4709,5923)	(4715.968,5980.470)
	(4760.642,6042.394)
	(4799.142,6077.088)
	(4849.496,6115.371)
	(4912.512,6158.066)
	(4949.021,6181.327)
	(4989.000,6206.000)

\path(9052,6686)	(9055.545,6647.730)
	(9059.063,6610.259)
	(9066.022,6537.691)
	(9072.885,6468.244)
	(9079.664,6401.869)
	(9086.369,6338.515)
	(9093.010,6278.133)
	(9099.598,6220.670)
	(9106.142,6166.078)
	(9112.653,6114.306)
	(9119.142,6065.303)
	(9125.618,6019.019)
	(9132.092,5975.404)
	(9138.574,5934.408)
	(9145.075,5895.979)
	(9158.174,5826.625)
	(9171.470,5766.939)
	(9185.046,5716.518)
	(9213.369,5641.859)
	(9277.000,5586.000)

\path(9277,5586)	(9316.759,5633.186)
	(9332.354,5687.347)
	(9345.716,5757.601)
	(9351.706,5797.758)
	(9357.313,5840.734)
	(9362.595,5886.126)
	(9367.611,5933.533)
	(9372.419,5982.552)
	(9377.078,6032.782)
	(9381.646,6083.822)
	(9386.181,6135.269)
	(9390.743,6186.721)
	(9395.388,6237.778)
	(9400.177,6288.037)
	(9405.167,6337.096)
	(9410.416,6384.554)
	(9415.983,6430.009)
	(9421.927,6473.059)
	(9428.305,6513.302)
	(9442.601,6583.762)
	(9459.337,6638.175)
	(9502.000,6686.000)

\path(9502,6686)	(9543.731,6670.994)
	(9582.496,6622.210)
	(9601.040,6583.984)
	(9619.172,6535.914)
	(9637.002,6477.531)
	(9654.640,6408.369)
	(9663.421,6369.600)
	(9672.195,6327.961)
	(9680.976,6283.394)
	(9689.778,6235.840)
	(9698.614,6185.242)
	(9707.498,6131.540)
	(9716.443,6074.676)
	(9725.464,6014.593)
	(9734.575,5951.231)
	(9743.788,5884.532)
	(9753.118,5814.439)
	(9762.579,5740.891)
	(9767.362,5702.804)
	(9772.184,5663.832)
	(9777.044,5623.967)
	(9781.946,5583.203)
	(9786.891,5541.531)
	(9791.880,5498.945)
	(9796.916,5455.437)
	(9802.000,5411.000)

\path(9493,8409)	(9496.545,8370.730)
	(9500.063,8333.259)
	(9507.022,8260.691)
	(9513.885,8191.244)
	(9520.664,8124.869)
	(9527.369,8061.515)
	(9534.010,8001.133)
	(9540.598,7943.670)
	(9547.142,7889.078)
	(9553.653,7837.306)
	(9560.142,7788.303)
	(9566.618,7742.019)
	(9573.092,7698.404)
	(9579.574,7657.408)
	(9586.075,7618.979)
	(9599.174,7549.625)
	(9612.470,7489.939)
	(9626.046,7439.518)
	(9654.369,7364.859)
	(9718.000,7309.000)

\path(9718,7309)	(9757.759,7356.186)
	(9773.354,7410.347)
	(9786.716,7480.601)
	(9792.706,7520.758)
	(9798.313,7563.734)
	(9803.595,7609.126)
	(9808.611,7656.533)
	(9813.419,7705.552)
	(9818.078,7755.782)
	(9822.646,7806.822)
	(9827.181,7858.269)
	(9831.743,7909.721)
	(9836.388,7960.778)
	(9841.177,8011.037)
	(9846.167,8060.096)
	(9851.416,8107.554)
	(9856.983,8153.009)
	(9862.927,8196.059)
	(9869.305,8236.302)
	(9883.601,8306.762)
	(9900.337,8361.175)
	(9943.000,8409.000)

\path(9943,8409)	(9984.731,8393.994)
	(10023.496,8345.210)
	(10042.040,8306.984)
	(10060.172,8258.914)
	(10078.002,8200.531)
	(10095.640,8131.369)
	(10104.421,8092.600)
	(10113.195,8050.961)
	(10121.976,8006.394)
	(10130.778,7958.840)
	(10139.614,7908.242)
	(10148.498,7854.540)
	(10157.443,7797.676)
	(10166.464,7737.593)
	(10175.575,7674.231)
	(10184.788,7607.532)
	(10194.118,7537.439)
	(10203.579,7463.891)
	(10208.362,7425.804)
	(10213.184,7386.832)
	(10218.044,7346.967)
	(10222.946,7306.203)
	(10227.891,7264.531)
	(10232.880,7221.945)
	(10237.916,7178.437)
	(10243.000,7134.000)

\path(7027,8411)	(7030.545,8372.730)
	(7034.063,8335.259)
	(7041.022,8262.691)
	(7047.885,8193.244)
	(7054.664,8126.869)
	(7061.369,8063.515)
	(7068.010,8003.133)
	(7074.598,7945.670)
	(7081.142,7891.078)
	(7087.653,7839.306)
	(7094.142,7790.303)
	(7100.618,7744.019)
	(7107.092,7700.404)
	(7113.574,7659.408)
	(7120.075,7620.979)
	(7133.174,7551.625)
	(7146.470,7491.939)
	(7160.046,7441.518)
	(7188.369,7366.859)
	(7252.000,7311.000)

\path(7252,7311)	(7291.759,7358.186)
	(7307.354,7412.347)
	(7320.716,7482.601)
	(7326.706,7522.758)
	(7332.313,7565.734)
	(7337.595,7611.126)
	(7342.611,7658.533)
	(7347.419,7707.552)
	(7352.078,7757.782)
	(7356.646,7808.822)
	(7361.181,7860.269)
	(7365.743,7911.721)
	(7370.388,7962.778)
	(7375.177,8013.037)
	(7380.167,8062.096)
	(7385.416,8109.554)
	(7390.983,8155.009)
	(7396.927,8198.059)
	(7403.305,8238.302)
	(7417.601,8308.762)
	(7434.337,8363.175)
	(7477.000,8411.000)

\path(7477,8411)	(7518.731,8395.994)
	(7557.496,8347.210)
	(7576.040,8308.984)
	(7594.172,8260.914)
	(7612.002,8202.531)
	(7629.640,8133.369)
	(7638.421,8094.600)
	(7647.195,8052.961)
	(7655.976,8008.394)
	(7664.778,7960.840)
	(7673.614,7910.242)
	(7682.498,7856.540)
	(7691.443,7799.676)
	(7700.464,7739.593)
	(7709.575,7676.231)
	(7718.788,7609.532)
	(7728.118,7539.439)
	(7737.579,7465.891)
	(7742.362,7427.804)
	(7747.184,7388.832)
	(7752.044,7348.967)
	(7756.946,7308.203)
	(7761.891,7266.531)
	(7766.880,7223.945)
	(7771.916,7180.437)
	(7777.000,7136.000)

\path(3268,3009)	(3271.545,2970.730)
	(3275.063,2933.259)
	(3282.022,2860.691)
	(3288.885,2791.244)
	(3295.664,2724.869)
	(3302.369,2661.515)
	(3309.010,2601.133)
	(3315.598,2543.670)
	(3322.142,2489.078)
	(3328.653,2437.306)
	(3335.142,2388.303)
	(3341.618,2342.019)
	(3348.092,2298.404)
	(3354.574,2257.408)
	(3361.075,2218.979)
	(3374.174,2149.625)
	(3387.470,2089.939)
	(3401.046,2039.518)
	(3429.369,1964.859)
	(3493.000,1909.000)

\path(3493,1909)	(3532.759,1956.186)
	(3548.354,2010.347)
	(3561.716,2080.601)
	(3567.706,2120.758)
	(3573.313,2163.734)
	(3578.595,2209.126)
	(3583.611,2256.533)
	(3588.419,2305.552)
	(3593.078,2355.782)
	(3597.646,2406.822)
	(3602.181,2458.269)
	(3606.743,2509.721)
	(3611.388,2560.778)
	(3616.177,2611.037)
	(3621.167,2660.096)
	(3626.416,2707.554)
	(3631.983,2753.009)
	(3637.927,2796.059)
	(3644.305,2836.302)
	(3658.601,2906.762)
	(3675.337,2961.175)
	(3718.000,3009.000)

\path(3718,3009)	(3759.731,2993.994)
	(3798.496,2945.210)
	(3817.040,2906.984)
	(3835.172,2858.914)
	(3853.002,2800.531)
	(3870.640,2731.369)
	(3879.421,2692.600)
	(3888.195,2650.961)
	(3896.976,2606.394)
	(3905.778,2558.840)
	(3914.614,2508.242)
	(3923.498,2454.540)
	(3932.443,2397.676)
	(3941.464,2337.593)
	(3950.575,2274.231)
	(3959.788,2207.532)
	(3969.118,2137.439)
	(3978.579,2063.891)
	(3983.362,2025.804)
	(3988.184,1986.832)
	(3993.044,1946.967)
	(3997.946,1906.203)
	(4002.891,1864.531)
	(4007.880,1821.945)
	(4012.916,1778.437)
	(4018.000,1734.000)

\put(300,13140){\makebox(0,0)[lb]{\smash{{{\SetFigFont{12}{14.4}{rm}12}}}}}
\put(300,9540){\makebox(0,0)[lb]{\smash{{{\SetFigFont{12}{14.4}{rm}10}}}}}
\put(300,7740){\makebox(0,0)[lb]{\smash{{{\SetFigFont{12}{14.4}{rm}9}}}}}
\put(300,5865){\makebox(0,0)[lb]{\smash{{{\SetFigFont{12}{14.4}{rm}8}}}}}
\put(300,4140){\makebox(0,0)[lb]{\smash{{{\SetFigFont{12}{14.4}{rm}7}}}}}
\put(0,14115){\makebox(0,0)[lb]{\smash{{{\SetFigFont{12}{14.4}{rm}\# general}}}}}
\put(225,13740){\makebox(0,0)[lb]{\smash{{{\SetFigFont{12}{14.4}{rm}lines}}}}}
\put(300,2340){\makebox(0,0)[lb]{\smash{{{\SetFigFont{12}{14.4}{rm}6}}}}}
\put(300,465){\makebox(0,0)[lb]{\smash{{{\SetFigFont{12}{14.4}{rm}5}}}}}
\put(6600,12615){\makebox(0,0)[lb]{\smash{{{\SetFigFont{8}{9.6}{rm}80160}}}}}
\put(6600,10815){\makebox(0,0)[lb]{\smash{{{\SetFigFont{8}{9.6}{rm}80160}}}}}
\put(5175,9015){\makebox(0,0)[lb]{\smash{{{\SetFigFont{8}{9.6}{rm}70296}}}}}
\put(8100,9015){\makebox(0,0)[lb]{\smash{{{\SetFigFont{8}{9.6}{rm}9864}}}}}
\put(3825,7215){\makebox(0,0)[lb]{\smash{{{\SetFigFont{8}{9.6}{rm}4968}}}}}
\put(1125,5415){\makebox(0,0)[lb]{\smash{{{\SetFigFont{8}{9.6}{rm}2552}}}}}
\put(2850,5415){\makebox(0,0)[lb]{\smash{{{\SetFigFont{8}{9.6}{rm}364}}}}}
\put(4425,5415){\makebox(0,0)[lb]{\smash{{{\SetFigFont{8}{9.6}{rm}12}}}}}
\put(6975,5415){\makebox(0,0)[lb]{\smash{{{\SetFigFont{8}{9.6}{rm}2380}}}}}
\put(8175,1815){\makebox(0,0)[lb]{\smash{{{\SetFigFont{8}{9.6}{rm}12}}}}}
\put(9375,1815){\makebox(0,0)[lb]{\smash{{{\SetFigFont{8}{9.6}{rm}74}}}}}
\put(10875,1815){\makebox(0,0)[lb]{\smash{{{\SetFigFont{8}{9.6}{rm}40}}}}}
\put(6525,1815){\makebox(0,0)[lb]{\smash{{{\SetFigFont{8}{9.6}{rm}110}}}}}
\put(7500,15){\makebox(0,0)[lb]{\smash{{{\SetFigFont{8}{9.6}{rm}12}}}}}
\put(5025,15){\makebox(0,0)[lb]{\smash{{{\SetFigFont{8}{9.6}{rm}50}}}}}
\put(1800,1815){\makebox(0,0)[lb]{\smash{{{\SetFigFont{8}{9.6}{rm}58}}}}}
\put(2250,15){\makebox(0,0)[lb]{\smash{{{\SetFigFont{8}{9.6}{rm}18}}}}}
\put(4650,1815){\makebox(0,0)[lb]{\smash{{{\SetFigFont{8}{9.6}{rm}126}}}}}
\put(6000,3615){\makebox(0,0)[lb]{\smash{{{\SetFigFont{8}{9.6}{rm}12}}}}}
\put(7200,3615){\makebox(0,0)[lb]{\smash{{{\SetFigFont{8}{9.6}{rm}610}}}}}
\put(4425,3615){\makebox(0,0)[lb]{\smash{{{\SetFigFont{8}{9.6}{rm}210}}}}}
\put(2475,3615){\makebox(0,0)[lb]{\smash{{{\SetFigFont{8}{9.6}{rm}504}}}}}
\put(1200,3615){\makebox(0,0)[lb]{\smash{{{\SetFigFont{8}{9.6}{rm}522}}}}}
\put(12000,5415){\makebox(0,0)[lb]{\smash{{{\SetFigFont{8}{9.6}{rm}920}}}}}
\put(12075,3615){\makebox(0,0)[lb]{\smash{{{\SetFigFont{8}{9.6}{rm}92}}}}}
\put(10125,5415){\makebox(0,0)[lb]{\smash{{{\SetFigFont{8}{9.6}{rm}1312}}}}}
\put(8550,5415){\makebox(0,0)[lb]{\smash{{{\SetFigFont{8}{9.6}{rm}7632}}}}}
\put(5475,5415){\makebox(0,0)[lb]{\smash{{{\SetFigFont{8}{9.6}{rm}2808}}}}}
\put(8550,3615){\makebox(0,0)[lb]{\smash{{{\SetFigFont{8}{9.6}{rm}420}}}}}
\put(2925,1815){\makebox(0,0)[lb]{\smash{{{\SetFigFont{8}{9.6}{rm}190}}}}}
\put(3900,6615){\makebox(0,0)[lb]{\smash{{{\SetFigFont{8}{9.6}{rm}$\times 16$}}}}}
\put(8550,8490){\makebox(0,0)[lb]{\smash{{{\SetFigFont{8}{9.6}{rm}$\times 2$}}}}}
\put(8700,6840){\makebox(0,0)[lb]{\smash{{{\SetFigFont{8}{9.6}{rm}$\times 3$}}}}}
\put(1050,4890){\makebox(0,0)[lb]{\smash{{{\SetFigFont{8}{9.6}{rm}$\times 2$}}}}}
\put(4500,4890){\makebox(0,0)[lb]{\smash{{{\SetFigFont{8}{9.6}{rm}$\times 8$}}}}}
\put(6225,4890){\makebox(0,0)[lb]{\smash{{{\SetFigFont{8}{9.6}{rm}$\times 27$}}}}}
\put(7575,4890){\makebox(0,0)[lb]{\smash{{{\SetFigFont{8}{9.6}{rm}$\times 2$}}}}}
\put(8175,3090){\makebox(0,0)[lb]{\smash{{{\SetFigFont{8}{9.6}{rm}$\times 9$}}}}}
\put(7050,3090){\makebox(0,0)[lb]{\smash{{{\SetFigFont{8}{9.6}{rm}$\times 4$}}}}}
\put(5100,3090){\makebox(0,0)[lb]{\smash{{{\SetFigFont{8}{9.6}{rm}$\times 2$}}}}}
\put(2175,3090){\makebox(0,0)[lb]{\smash{{{\SetFigFont{8}{9.6}{rm}$\times 2$}}}}}
\put(4500,1215){\makebox(0,0)[lb]{\smash{{{\SetFigFont{8}{9.6}{rm}$\times 2$}}}}}
\put(6750,1215){\makebox(0,0)[lb]{\smash{{{\SetFigFont{8}{9.6}{rm}$\times 3$}}}}}
\put(6375,7215){\makebox(0,0)[lb]{\smash{{{\SetFigFont{8}{9.6}{rm}45600}}}}}
\put(8925,7215){\makebox(0,0)[lb]{\smash{{{\SetFigFont{8}{9.6}{rm}9864}}}}}
\put(10350,4890){\makebox(0,0)[lb]{\smash{{{\SetFigFont{8}{9.6}{rm}$\times 2$}}}}}
\put(300,11340){\makebox(0,0)[lb]{\smash{{{\SetFigFont{12}{14.4}{rm}11}}}}}
\put(2325,6690){\makebox(0,0)[lb]{\smash{{{\SetFigFont{8}{9.6}{rm}$\times 2$}}}}}
\put(5100,6540){\makebox(0,0)[lb]{\smash{{{\SetFigFont{8}{9.6}{rm}$\times 81$}}}}}
\put(6450,6690){\makebox(0,0)[lb]{\smash{{{\SetFigFont{8}{9.6}{rm}$\times 3$}}}}}
\put(10125,3615){\makebox(0,0)[lb]{\smash{{{\SetFigFont{8}{9.6}{rm}1220}}}}}
\put(3075,1440){\makebox(0,0)[lb]{\smash{{{\SetFigFont{8}{9.6}{rm}$\times 3$}}}}}
\end{picture}

%% file: iegC.tex
\begingroup\makeatletter\ifx\SetFigFont\undefined
\def\x#1#2#3#4#5#6#7\relax{\def\x{#1#2#3#4#5#6}}%
\expandafter\x\fmtname xxxxxx\relax \def\y{splain}%
\ifx\x\y   
\gdef\SetFigFont#1#2#3{%
  \ifnum #1<17\tiny\else \ifnum #1<20\small\else
  \ifnum #1<24\normalsize\else \ifnum #1<29\large\else
  \ifnum #1<34\Large\else \ifnum #1<41\LARGE\else
     \huge\fi\fi\fi\fi\fi\fi
  \csname #3\endcsname}%
\else
\gdef\SetFigFont#1#2#3{\begingroup
  \count@#1\relax \ifnum 25<\count@\count@25\fi
  \def\x{\endgroup\@setsize\SetFigFont{#2pt}}%
  \expandafter\x
    \csname \romannumeral\the\count@ pt\expandafter\endcsname
    \csname @\romannumeral\the\count@ pt\endcsname
  \csname #3\endcsname}%
\fi
\fi\endgroup
\begin{picture}(8274,7035)(0,-10)
\thicklines
\path(5037,2283)(5037,985)
\path(5007.000,1105.000)(5037.000,985.000)(5067.000,1105.000)
\put(5337,1608){\makebox(0,0)[lb]{\smash{{{\SetFigFont{12}{14.4}{rm}$\pi$}}}}}
\put(4062,2508){\blacken\ellipse{74}{74}}
\put(4062,2508){\ellipse{74}{74}}
\put(4062,3408){\blacken\ellipse{74}{74}}
\put(4062,3408){\ellipse{74}{74}}
\put(4062,558){\blacken\ellipse{74}{74}}
\put(4062,558){\ellipse{74}{74}}
\put(3957,4983){\blacken\ellipse{74}{74}}
\put(3957,4983){\ellipse{74}{74}}
\put(4047,4563){\blacken\ellipse{74}{74}}
\put(4047,4563){\ellipse{74}{74}}
\put(3949,6648){\blacken\ellipse{74}{74}}
\put(3949,6648){\ellipse{74}{74}}
\path(12,558)(8262,558)
\dottedline{135}(3762,3708)(4362,3708)(4362,7008)
	(3762,7008)(3762,3708)
\path(162,2208)	(206.121,2204.793)
	(249.854,2201.636)
	(293.202,2198.530)
	(336.166,2195.475)
	(378.750,2192.470)
	(420.955,2189.515)
	(462.783,2186.611)
	(504.238,2183.758)
	(545.322,2180.955)
	(586.036,2178.203)
	(626.384,2175.501)
	(666.368,2172.849)
	(705.989,2170.248)
	(745.251,2167.698)
	(784.156,2165.198)
	(822.705,2162.748)
	(860.902,2160.349)
	(898.749,2158.001)
	(973.402,2153.456)
	(1046.682,2149.112)
	(1118.609,2144.971)
	(1189.201,2141.031)
	(1258.478,2137.294)
	(1326.457,2133.759)
	(1393.159,2130.426)
	(1458.601,2127.294)
	(1522.803,2124.365)
	(1585.784,2121.638)
	(1647.562,2119.113)
	(1708.157,2116.790)
	(1767.587,2114.668)
	(1825.871,2112.749)
	(1883.029,2111.032)
	(1939.078,2109.517)
	(1994.038,2108.204)
	(2047.927,2107.093)
	(2100.766,2106.184)
	(2152.572,2105.477)
	(2203.364,2104.972)
	(2253.161,2104.669)
	(2301.983,2104.568)
	(2349.847,2104.669)
	(2396.773,2104.972)
	(2442.781,2105.477)
	(2487.887,2106.184)
	(2532.113,2107.093)
	(2575.475,2108.204)
	(2617.994,2109.517)
	(2659.689,2111.032)
	(2700.577,2112.749)
	(2740.678,2114.668)
	(2780.011,2116.790)
	(2818.595,2119.113)
	(2893.590,2124.365)
	(2965.815,2130.426)
	(3035.420,2137.294)
	(3102.557,2144.971)
	(3167.376,2153.456)
	(3230.029,2162.748)
	(3290.667,2172.849)
	(3349.441,2183.758)
	(3406.502,2195.475)
	(3462.000,2208.000)

\path(3462,2208)	(3533.242,2224.855)
	(3576.169,2235.427)
	(3622.550,2247.442)
	(3671.323,2260.910)
	(3721.427,2275.843)
	(3771.799,2292.250)
	(3821.378,2310.142)
	(3869.101,2329.530)
	(3913.906,2350.424)
	(3954.732,2372.835)
	(3990.517,2396.772)
	(4042.713,2449.269)
	(4062.000,2508.000)

\path(4062,2508)	(4042.713,2566.731)
	(3990.517,2619.228)
	(3954.732,2643.165)
	(3913.906,2665.576)
	(3869.101,2686.470)
	(3821.378,2705.858)
	(3771.799,2723.750)
	(3721.427,2740.157)
	(3671.323,2755.090)
	(3622.550,2768.558)
	(3576.169,2780.573)
	(3533.242,2791.145)
	(3462.000,2808.000)

\path(3462,2808)	(3399.992,2820.810)
	(3333.656,2831.141)
	(3263.278,2839.157)
	(3189.149,2845.025)
	(3150.767,2847.205)
	(3111.556,2848.909)
	(3071.551,2850.159)
	(3030.789,2850.975)
	(2989.305,2851.378)
	(2947.135,2851.388)
	(2904.317,2851.027)
	(2860.884,2850.314)
	(2816.875,2849.271)
	(2772.325,2847.917)
	(2727.269,2846.275)
	(2681.744,2844.364)
	(2635.787,2842.205)
	(2589.433,2839.818)
	(2542.718,2837.225)
	(2495.678,2834.446)
	(2448.349,2831.502)
	(2400.768,2828.413)
	(2352.971,2825.200)
	(2304.993,2821.884)
	(2256.871,2818.486)
	(2208.641,2815.025)
	(2160.339,2811.523)
	(2112.000,2808.000)
	(2063.661,2804.477)
	(2015.359,2800.975)
	(1967.129,2797.514)
	(1919.007,2794.116)
	(1871.029,2790.800)
	(1823.232,2787.587)
	(1775.651,2784.498)
	(1728.322,2781.554)
	(1681.282,2778.775)
	(1634.567,2776.182)
	(1588.213,2773.795)
	(1542.256,2771.636)
	(1496.731,2769.725)
	(1451.675,2768.083)
	(1407.125,2766.729)
	(1363.116,2765.686)
	(1319.683,2764.973)
	(1276.865,2764.612)
	(1234.695,2764.622)
	(1193.211,2765.025)
	(1152.449,2765.841)
	(1112.444,2767.091)
	(1073.233,2768.795)
	(1034.851,2770.975)
	(960.722,2776.843)
	(890.344,2784.859)
	(824.008,2795.190)
	(762.000,2808.000)

\path(762,2808)	(690.758,2824.855)
	(647.831,2835.427)
	(601.450,2847.442)
	(552.677,2860.910)
	(502.573,2875.843)
	(452.201,2892.250)
	(402.622,2910.142)
	(354.899,2929.530)
	(310.094,2950.424)
	(269.268,2972.835)
	(233.483,2996.772)
	(181.287,3049.269)
	(162.000,3108.000)

\path(162,3108)	(181.287,3166.731)
	(233.483,3219.228)
	(269.268,3243.165)
	(310.094,3265.576)
	(354.899,3286.470)
	(402.622,3305.858)
	(452.201,3323.750)
	(502.573,3340.157)
	(552.677,3355.090)
	(601.450,3368.558)
	(647.831,3380.573)
	(690.758,3391.145)
	(762.000,3408.000)

\path(762,3408)	(836.420,3424.439)
	(874.660,3432.068)
	(913.573,3439.312)
	(953.147,3446.177)
	(993.372,3452.670)
	(1034.236,3458.796)
	(1075.729,3464.563)
	(1117.839,3469.976)
	(1160.555,3475.042)
	(1203.866,3479.768)
	(1247.762,3484.158)
	(1292.230,3488.221)
	(1337.261,3491.961)
	(1382.844,3495.386)
	(1428.966,3498.501)
	(1475.617,3501.313)
	(1522.787,3503.829)
	(1570.464,3506.054)
	(1618.637,3507.995)
	(1667.294,3509.659)
	(1716.426,3511.051)
	(1766.021,3512.178)
	(1816.068,3513.046)
	(1866.555,3513.661)
	(1917.473,3514.031)
	(1968.809,3514.160)
	(2020.554,3514.056)
	(2072.695,3513.724)
	(2125.222,3513.172)
	(2178.124,3512.405)
	(2231.390,3511.430)
	(2285.008,3510.252)
	(2338.968,3508.879)
	(2393.259,3507.317)
	(2447.870,3505.571)
	(2502.789,3503.649)
	(2558.006,3501.556)
	(2613.509,3499.300)
	(2669.288,3496.885)
	(2725.332,3494.319)
	(2781.629,3491.607)
	(2838.168,3488.757)
	(2894.940,3485.774)
	(2951.931,3482.665)
	(3009.132,3479.436)
	(3066.532,3476.094)
	(3124.119,3472.644)
	(3181.882,3469.093)
	(3239.810,3465.447)
	(3297.893,3461.713)
	(3356.120,3457.897)
	(3414.478,3454.005)
	(3472.958,3450.044)
	(3531.548,3446.019)
	(3590.237,3441.938)
	(3649.015,3437.806)
	(3707.869,3433.630)
	(3766.790,3429.416)
	(3825.766,3425.171)
	(3884.786,3420.900)
	(3943.839,3416.611)
	(4002.914,3412.309)
	(4062.000,3408.000)
	(4121.086,3403.691)
	(4180.161,3399.389)
	(4239.214,3395.100)
	(4298.234,3390.829)
	(4357.210,3386.584)
	(4416.131,3382.370)
	(4474.985,3378.194)
	(4533.763,3374.062)
	(4592.452,3369.981)
	(4651.042,3365.956)
	(4709.522,3361.995)
	(4767.880,3358.103)
	(4826.107,3354.287)
	(4884.190,3350.553)
	(4942.118,3346.907)
	(4999.881,3343.356)
	(5057.468,3339.906)
	(5114.868,3336.564)
	(5172.069,3333.335)
	(5229.060,3330.226)
	(5285.832,3327.243)
	(5342.371,3324.393)
	(5398.668,3321.681)
	(5454.712,3319.115)
	(5510.491,3316.700)
	(5565.994,3314.444)
	(5621.211,3312.351)
	(5676.130,3310.429)
	(5730.741,3308.683)
	(5785.032,3307.121)
	(5838.992,3305.748)
	(5892.610,3304.570)
	(5945.876,3303.595)
	(5998.778,3302.828)
	(6051.305,3302.276)
	(6103.446,3301.944)
	(6155.191,3301.840)
	(6206.527,3301.969)
	(6257.445,3302.339)
	(6307.932,3302.954)
	(6357.979,3303.822)
	(6407.574,3304.949)
	(6456.706,3306.341)
	(6505.363,3308.005)
	(6553.536,3309.946)
	(6601.213,3312.171)
	(6648.383,3314.687)
	(6695.034,3317.499)
	(6741.156,3320.614)
	(6786.739,3324.039)
	(6831.770,3327.779)
	(6876.238,3331.842)
	(6920.134,3336.232)
	(6963.445,3340.958)
	(7006.161,3346.024)
	(7048.271,3351.437)
	(7089.764,3357.204)
	(7130.628,3363.330)
	(7170.853,3369.823)
	(7210.427,3376.688)
	(7249.340,3383.932)
	(7287.580,3391.561)
	(7362.000,3408.000)

\path(7362,3408)	(7433.242,3424.855)
	(7476.169,3435.427)
	(7522.550,3447.442)
	(7571.323,3460.910)
	(7621.427,3475.843)
	(7671.799,3492.250)
	(7721.378,3510.142)
	(7769.101,3529.530)
	(7813.906,3550.424)
	(7854.732,3572.835)
	(7890.517,3596.772)
	(7942.713,3649.269)
	(7962.000,3708.000)

\path(7962,3708)	(7942.713,3766.731)
	(7890.517,3819.228)
	(7854.732,3843.165)
	(7813.906,3865.576)
	(7769.101,3886.470)
	(7721.378,3905.858)
	(7671.799,3923.750)
	(7621.427,3940.157)
	(7571.323,3955.090)
	(7522.550,3968.558)
	(7476.169,3980.573)
	(7433.242,3991.145)
	(7362.000,4008.000)

\path(7362,4008)	(7303.115,4021.595)
	(7243.435,4034.757)
	(7182.922,4047.489)
	(7121.533,4059.789)
	(7059.230,4071.657)
	(6995.972,4083.094)
	(6931.717,4094.099)
	(6866.427,4104.672)
	(6800.061,4114.814)
	(6732.578,4124.525)
	(6663.938,4133.803)
	(6594.100,4142.651)
	(6523.025,4151.066)
	(6450.672,4159.050)
	(6377.000,4166.603)
	(6339.657,4170.217)
	(6301.970,4173.724)
	(6263.932,4177.122)
	(6225.540,4180.413)
	(6186.788,4183.596)
	(6147.672,4186.671)
	(6108.185,4189.638)
	(6068.323,4192.497)
	(6028.081,4195.248)
	(5987.455,4197.892)
	(5946.438,4200.427)
	(5905.026,4202.855)
	(5863.213,4205.175)
	(5820.996,4207.386)
	(5778.368,4209.490)
	(5735.325,4211.486)
	(5691.862,4213.374)
	(5647.973,4215.155)
	(5603.654,4216.827)
	(5558.899,4218.391)
	(5513.704,4219.848)
	(5468.063,4221.197)
	(5421.971,4222.437)
	(5375.424,4223.570)
	(5328.416,4224.595)
	(5280.943,4225.512)
	(5232.998,4226.322)
	(5184.578,4227.023)
	(5135.677,4227.616)
	(5086.290,4228.102)
	(5036.412,4228.479)
	(4986.038,4228.749)
	(4935.164,4228.911)
	(4883.782,4228.965)
	(4831.890,4228.911)
	(4779.482,4228.749)
	(4726.553,4228.479)
	(4673.097,4228.102)
	(4619.110,4227.616)
	(4564.586,4227.023)
	(4509.521,4226.322)
	(4453.910,4225.512)
	(4397.748,4224.595)
	(4341.028,4223.570)
	(4283.748,4222.437)
	(4225.900,4221.197)
	(4167.481,4219.848)
	(4108.486,4218.391)
	(4048.908,4216.827)
	(3988.744,4215.155)
	(3927.988,4213.374)
	(3866.636,4211.486)
	(3804.681,4209.490)
	(3742.119,4207.386)
	(3678.946,4205.175)
	(3615.155,4202.855)
	(3550.743,4200.427)
	(3485.703,4197.892)
	(3420.031,4195.248)
	(3353.722,4192.497)
	(3286.771,4189.638)
	(3219.172,4186.671)
	(3150.921,4183.596)
	(3082.013,4180.413)
	(3012.442,4177.122)
	(2942.204,4173.724)
	(2871.293,4170.217)
	(2799.705,4166.603)
	(2727.435,4162.880)
	(2654.476,4159.050)
	(2580.825,4155.112)
	(2506.477,4151.066)
	(2469.039,4149.003)
	(2431.425,4146.912)
	(2393.635,4144.795)
	(2355.666,4142.651)
	(2317.520,4140.479)
	(2279.194,4138.281)
	(2240.689,4136.056)
	(2202.004,4133.803)
	(2163.139,4131.524)
	(2124.091,4129.218)
	(2084.862,4126.885)
	(2045.451,4124.525)
	(2005.856,4122.137)
	(1966.077,4119.723)
	(1926.114,4117.282)
	(1885.965,4114.814)
	(1845.631,4112.319)
	(1805.110,4109.797)
	(1764.403,4107.248)
	(1723.507,4104.672)
	(1682.424,4102.069)
	(1641.151,4099.439)
	(1599.689,4096.783)
	(1558.037,4094.099)
	(1516.194,4091.388)
	(1474.159,4088.650)
	(1431.933,4085.885)
	(1389.513,4083.094)
	(1346.901,4080.275)
	(1304.094,4077.429)
	(1261.093,4074.557)
	(1217.896,4071.657)
	(1174.504,4068.730)
	(1130.915,4065.777)
	(1087.130,4062.796)
	(1043.146,4059.789)
	(998.964,4056.754)
	(954.583,4053.693)
	(910.002,4050.604)
	(865.221,4047.489)
	(820.240,4044.346)
	(775.056,4041.177)
	(729.671,4037.981)
	(684.082,4034.757)
	(638.291,4031.507)
	(592.295,4028.230)
	(546.094,4024.926)
	(499.689,4021.595)
	(453.077,4018.236)
	(406.259,4014.851)
	(359.233,4011.439)
	(312.000,4008.000)

\path(5637,3558)	(5692.548,3508.998)
	(5738.943,3465.039)
	(5776.844,3425.246)
	(5806.909,3388.740)
	(5846.170,3322.071)
	(5862.000,3258.000)

\path(5862,3258)	(5857.978,3219.530)
	(5841.617,3181.665)
	(5788.140,3107.096)
	(5734.093,3032.979)
	(5712.000,2958.000)

\path(5712,2958)	(5721.760,2912.177)
	(5747.670,2862.919)
	(5793.245,2804.951)
	(5862.000,2733.000)

\path(6912,3858)	(6847.475,3901.904)
	(6803.310,3939.195)
	(6762.000,4008.000)

\path(6762,4008)	(6759.479,4062.947)
	(6773.698,4119.919)
	(6799.325,4178.090)
	(6831.030,4236.634)
	(6863.481,4294.723)
	(6891.347,4351.531)
	(6909.297,4406.233)
	(6912.000,4458.000)

\path(6912,4458)	(6900.085,4501.511)
	(6873.458,4549.999)
	(6828.601,4608.737)
	(6762.000,4683.000)

\path(4062,3858)	(4074.944,3911.336)
	(4086.659,3961.038)
	(4097.172,4007.299)
	(4106.511,4050.311)
	(4114.704,4090.266)
	(4121.777,4127.356)
	(4132.676,4193.711)
	(4139.428,4250.915)
	(4142.253,4300.506)
	(4141.371,4344.022)
	(4137.000,4383.000)

\path(4137,4383)	(4110.462,4454.057)
	(4088.499,4493.834)
	(4063.916,4534.703)
	(4039.094,4575.359)
	(4016.413,4614.497)
	(3987.000,4683.000)

\path(3987,4683)	(3974.659,4751.212)
	(3969.811,4790.801)
	(3965.844,4833.278)
	(3962.759,4878.063)
	(3960.555,4924.580)
	(3959.233,4972.251)
	(3958.793,5020.500)
	(3959.233,5068.749)
	(3960.555,5116.420)
	(3962.759,5162.937)
	(3965.844,5207.722)
	(3969.811,5250.199)
	(3974.659,5289.788)
	(3987.000,5358.000)

\path(3987,5358)	(4010.255,5430.948)
	(4028.483,5474.014)
	(4051.991,5523.233)
	(4081.438,5579.922)
	(4117.482,5645.401)
	(4138.184,5681.848)
	(4160.783,5720.987)
	(4185.361,5762.983)
	(4212.000,5808.000)

\path(4137,5433)	(4116.776,5477.660)
	(4098.255,5519.346)
	(4081.382,5558.224)
	(4066.102,5594.457)
	(4040.101,5659.651)
	(4019.812,5716.245)
	(4004.798,5765.559)
	(3994.617,5808.910)
	(3987.000,5883.000)

\path(3987,5883)	(3994.617,5957.090)
	(4004.798,6000.441)
	(4019.812,6049.755)
	(4040.101,6106.349)
	(4066.102,6171.543)
	(4081.382,6207.776)
	(4098.255,6246.654)
	(4116.776,6288.340)
	(4137.000,6333.000)

\path(4137,6033)	(4107.425,6075.377)
	(4080.303,6115.084)
	(4033.088,6187.150)
	(3994.695,6250.517)
	(3964.466,6306.503)
	(3941.741,6356.425)
	(3925.861,6401.603)
	(3916.167,6443.356)
	(3912.000,6483.000)

\path(3912,6483)	(3917.660,6549.045)
	(3927.351,6586.248)
	(3942.203,6627.701)
	(3962.654,6674.505)
	(3989.144,6727.756)
	(4022.113,6788.555)
	(4062.000,6858.000)

\path(312,5958)	(372.392,5952.690)
	(431.716,5947.492)
	(489.986,5942.405)
	(547.215,5937.428)
	(603.415,5932.560)
	(658.599,5927.801)
	(712.781,5923.150)
	(765.973,5918.605)
	(818.187,5914.166)
	(869.437,5909.833)
	(919.736,5905.603)
	(969.097,5901.477)
	(1017.532,5897.453)
	(1065.054,5893.531)
	(1111.677,5889.709)
	(1157.412,5885.988)
	(1202.274,5882.365)
	(1246.274,5878.841)
	(1289.426,5875.414)
	(1331.742,5872.083)
	(1373.236,5868.848)
	(1413.920,5865.708)
	(1453.807,5862.662)
	(1492.911,5859.708)
	(1531.243,5856.847)
	(1568.817,5854.077)
	(1641.742,5848.808)
	(1711.789,5843.894)
	(1779.060,5839.328)
	(1843.659,5835.102)
	(1905.689,5831.212)
	(1965.252,5827.648)
	(2022.452,5824.406)
	(2077.391,5821.477)
	(2130.174,5818.855)
	(2180.901,5816.532)
	(2229.678,5814.503)
	(2276.606,5812.761)
	(2321.789,5811.297)
	(2365.329,5810.107)
	(2407.330,5809.181)
	(2447.895,5808.515)
	(2487.126,5808.101)
	(2562.000,5808.000)

\path(2562,5808)	(2627.642,5808.958)
	(2698.263,5811.134)
	(2735.337,5812.652)
	(2773.531,5814.444)
	(2812.805,5816.500)
	(2853.115,5818.808)
	(2894.422,5821.358)
	(2936.683,5824.141)
	(2979.857,5827.146)
	(3023.903,5830.362)
	(3068.778,5833.779)
	(3114.443,5837.387)
	(3160.855,5841.176)
	(3207.972,5845.135)
	(3255.754,5849.254)
	(3304.159,5853.522)
	(3353.145,5857.930)
	(3402.671,5862.466)
	(3452.695,5867.121)
	(3503.176,5871.885)
	(3554.073,5876.746)
	(3605.344,5881.695)
	(3656.948,5886.722)
	(3708.843,5891.815)
	(3760.987,5896.965)
	(3813.340,5902.161)
	(3865.860,5907.394)
	(3918.505,5912.652)
	(3971.234,5917.925)
	(4024.005,5923.204)
	(4076.777,5928.477)
	(4129.509,5933.735)
	(4182.159,5938.966)
	(4234.685,5944.162)
	(4287.047,5949.310)
	(4339.202,5954.402)
	(4391.109,5959.427)
	(4442.728,5964.374)
	(4494.015,5969.233)
	(4544.930,5973.994)
	(4595.432,5978.647)
	(4645.479,5983.181)
	(4695.029,5987.585)
	(4744.041,5991.850)
	(4792.474,5995.966)
	(4840.285,5999.921)
	(4887.435,6003.705)
	(4933.880,6007.309)
	(4979.581,6010.722)
	(5024.494,6013.933)
	(5068.579,6016.933)
	(5111.795,6019.711)
	(5154.100,6022.256)
	(5195.452,6024.558)
	(5235.810,6026.607)
	(5275.132,6028.393)
	(5313.378,6029.906)
	(5350.505,6031.134)
	(5421.238,6032.697)
	(5487.000,6033.000)

\path(5487,6033)	(5528.934,6032.609)
	(5572.137,6031.893)
	(5616.725,6030.841)
	(5662.815,6029.443)
	(5710.523,6027.689)
	(5759.967,6025.568)
	(5811.262,6023.070)
	(5864.527,6020.185)
	(5919.876,6016.902)
	(5977.429,6013.212)
	(6037.300,6009.103)
	(6099.606,6004.567)
	(6164.466,5999.591)
	(6231.994,5994.166)
	(6302.308,5988.282)
	(6375.525,5981.929)
	(6413.258,5978.573)
	(6451.761,5975.095)
	(6491.048,5971.495)
	(6531.133,5967.772)
	(6572.032,5963.923)
	(6613.758,5959.947)
	(6656.327,5955.844)
	(6699.753,5951.612)
	(6744.050,5947.250)
	(6789.234,5942.756)
	(6835.318,5938.130)
	(6882.318,5933.369)
	(6930.247,5928.473)
	(6979.121,5923.440)
	(7028.954,5918.269)
	(7079.761,5912.958)
	(7131.556,5907.507)
	(7184.354,5901.915)
	(7238.170,5896.179)
	(7293.017,5890.299)
	(7348.911,5884.273)
	(7405.867,5878.099)
	(7463.898,5871.778)
	(7523.020,5865.307)
	(7583.246,5858.685)
	(7644.593,5851.911)
	(7707.073,5844.984)
	(7770.703,5837.902)
	(7835.496,5830.664)
	(7901.466,5823.268)
	(7968.630,5815.714)
	(8037.000,5808.000)

\put(3912,3033){\makebox(0,0)[lb]{\smash{{{\SetFigFont{8}{9.6}{rm}$p^1_1$}}}}}
\put(4362,2358){\makebox(0,0)[lb]{\smash{{{\SetFigFont{8}{9.6}{rm}$p^1_2$}}}}}
\put(3237,6483){\makebox(0,0)[lb]{\smash{{{\SetFigFont{8}{9.6}{rm}$p^3_1$}}}}}
\put(3987,33){\makebox(0,0)[lb]{\smash{{{\SetFigFont{12}{14.4}{rm}$z$}}}}}
\put(7437,858){\makebox(0,0)[lb]{\smash{{{\SetFigFont{12}{14.4}{rm}$\proj^1$}}}}}
\put(4512,4983){\makebox(0,0)[lb]{\smash{{{\SetFigFont{8}{9.6}{rm}$q$}}}}}
\put(4512,4608){\makebox(0,0)[lb]{\smash{{{\SetFigFont{8}{9.6}{rm}$p^2_1$}}}}}
\end{picture}

%% file: ideg0.tex
\begingroup\makeatletter\ifx\SetFigFont\undefined
\def\x#1#2#3#4#5#6#7\relax{\def\x{#1#2#3#4#5#6}}%
\expandafter\x\fmtname xxxxxx\relax \def\y{splain}%
\ifx\x\y   
\gdef\SetFigFont#1#2#3{%
  \ifnum #1<17\tiny\else \ifnum #1<20\small\else
  \ifnum #1<24\normalsize\else \ifnum #1<29\large\else
  \ifnum #1<34\Large\else \ifnum #1<41\LARGE\else
     \huge\fi\fi\fi\fi\fi\fi
  \csname #3\endcsname}%
\else
\gdef\SetFigFont#1#2#3{\begingroup
  \count@#1\relax \ifnum 25<\count@\count@25\fi
  \def\x{\endgroup\@setsize\SetFigFont{#2pt}}%
  \expandafter\x
    \csname \romannumeral\the\count@ pt\expandafter\endcsname
    \csname @\romannumeral\the\count@ pt\endcsname
  \csname #3\endcsname}%
\fi
\fi\endgroup
\begin{picture}(10524,5439)(0,-10)
\thicklines
\put(4212.000,4437.000){\arc{150.000}{1.5708}{4.7124}}
\put(5712.000,4437.000){\arc{150.000}{4.7124}{7.8540}}
\put(4062.000,4437.000){\arc{150.000}{4.7124}{7.8540}}
\put(5862.000,4437.000){\arc{150.000}{1.5708}{4.7124}}
\put(3612,4512){\blacken\ellipse{74}{74}}
\put(3612,4512){\ellipse{74}{74}}
\put(3612,4362){\blacken\ellipse{74}{74}}
\put(3612,4362){\ellipse{74}{74}}
\put(5112,4512){\blacken\ellipse{74}{74}}
\put(5112,4512){\ellipse{74}{74}}
\path(4212,4512)(5712,4512)
\path(5712,4362)(4212,4362)
\path(5862,4512)(6912,4512)
\path(6912,4362)(5862,4362)
\path(4062,4362)(3312,4362)
\path(3312,4512)(4062,4512)
\put(912.000,2937.000){\arc{150.000}{1.5708}{4.7124}}
\put(2412.000,2937.000){\arc{150.000}{4.7124}{7.8540}}
\put(762.000,2937.000){\arc{150.000}{4.7124}{7.8540}}
\put(2562.000,2937.000){\arc{150.000}{1.5708}{4.7124}}
\put(1812.000,1137.000){\arc{150.000}{1.5708}{4.7124}}
\put(3312.000,1137.000){\arc{150.000}{4.7124}{7.8540}}
\put(1662.000,1137.000){\arc{150.000}{4.7124}{7.8540}}
\put(3462.000,1137.000){\arc{150.000}{1.5708}{4.7124}}
\put(7812.000,2937.000){\arc{150.000}{1.5708}{4.7124}}
\put(9312.000,2937.000){\arc{150.000}{4.7124}{7.8540}}
\put(7662.000,2937.000){\arc{150.000}{4.7124}{7.8540}}
\put(9462.000,2937.000){\arc{150.000}{1.5708}{4.7124}}
\put(312,2862){\blacken\ellipse{74}{74}}
\put(312,2862){\ellipse{74}{74}}
\put(1287,1212){\blacken\ellipse{74}{74}}
\put(1287,1212){\ellipse{74}{74}}
\put(312,3612){\blacken\ellipse{74}{74}}
\put(312,3612){\ellipse{74}{74}}
\put(312,3312){\blacken\ellipse{74}{74}}
\put(312,3312){\ellipse{74}{74}}
\put(1287,912){\blacken\ellipse{74}{74}}
\put(1287,912){\ellipse{74}{74}}
\put(1287,612){\blacken\ellipse{74}{74}}
\put(1287,612){\ellipse{74}{74}}
\put(7737,3612){\blacken\ellipse{74}{74}}
\put(7737,3612){\ellipse{74}{74}}
\put(7737,3312){\blacken\ellipse{74}{74}}
\put(7737,3312){\ellipse{74}{74}}
\put(6987,1812){\blacken\ellipse{74}{74}}
\put(6987,1812){\ellipse{74}{74}}
\put(6987,1512){\blacken\ellipse{74}{74}}
\put(6987,1512){\ellipse{74}{74}}
\put(6987,912){\blacken\ellipse{74}{74}}
\put(6987,912){\ellipse{74}{74}}
\put(7737,2712){\blacken\ellipse{74}{74}}
\put(7737,2712){\ellipse{74}{74}}
\dottedline{135}(3612,5412)(3612,3612)
\path(912,3012)(2412,3012)
\path(2412,2862)(912,2862)
\path(2562,3012)(3612,3012)
\path(3612,2862)(2562,2862)
\path(762,2862)(12,2862)
\path(12,3012)(762,3012)
\path(1812,1212)(3312,1212)
\path(3312,1062)(1812,1062)
\path(3462,1212)(4512,1212)
\path(4512,1062)(3462,1062)
\path(1662,1062)(912,1062)
\path(912,1212)(1662,1212)
\path(6987,2112)(6987,612)
\path(6012,1212)(9762,1212)
\path(6012,1062)(9762,1062)
\path(7812,3012)(9312,3012)
\path(9312,2862)(7812,2862)
\path(9462,3012)(10512,3012)
\path(10512,2862)(9462,2862)
\path(7662,2862)(6912,2862)
\path(6912,3012)(7662,3012)
\path(7737,2412)(7737,3912)
\path(312,3012)(312,3912)
\dottedline{135}(312,2862)(312,2112)
\dottedline{135}(1287,1212)(1287,2112)
\path(1287,1062)(1287,312)
\path(4512,3912)(3837,3162)
\path(3894.977,3271.264)(3837.000,3162.000)(3939.575,3231.126)
\path(4812,3912)(3987,1737)
\path(4001.509,1859.839)(3987.000,1737.000)(4057.608,1838.560)
\path(5112,3912)(6012,1737)
\path(5938.397,1836.411)(6012.000,1737.000)(5993.838,1859.353)
\path(5412,3912)(6687,3312)
\path(6565.648,3335.951)(6687.000,3312.000)(6591.196,3390.240)
\dottedline{135}(312,3912)(312,4212)
\dottedline{135}(1287,312)(1287,12)
\dottedline{135}(7737,3912)(7737,4212)
\dottedline{135}(7737,2412)(7737,2112)
\dottedline{135}(6987,2112)(6987,2412)
\dottedline{135}(6987,612)(6987,312)
\put(5037,4662){\makebox(0,0)[lb]{\smash{{{\SetFigFont{8}{9.6}{rm}$q$}}}}}
\put(612,3612){\makebox(0,0)[lb]{\smash{{{\SetFigFont{8}{9.6}{rm}$q$}}}}}
\put(537,3237){\makebox(0,0)[lb]{\smash{{{\SetFigFont{8}{9.6}{rm}$p^1_1$}}}}}
\put(1362,1362){\makebox(0,0)[lb]{\smash{{{\SetFigFont{8}{9.6}{rm}$p^1_1$}}}}}
\put(6387,1737){\makebox(0,0)[lb]{\smash{{{\SetFigFont{8}{9.6}{rm}$p^1_1$}}}}}
\put(7962,3237){\makebox(0,0)[lb]{\smash{{{\SetFigFont{8}{9.6}{rm}$p^2_1$}}}}}
\put(3762,4662){\makebox(0,0)[lb]{\smash{{{\SetFigFont{8}{9.6}{rm}$p^1_1$}}}}}
\put(5412,5262){\makebox(0,0)[lb]{\smash{{{\SetFigFont{12}{14.4}{rm}$X(2, h_1=2, 0)$}}}}}
\put(3762,4062){\makebox(0,0)[lb]{\smash{{{\SetFigFont{8}{9.6}{rm}$p^2_1$}}}}}
\put(462,2562){\makebox(0,0)[lb]{\smash{{{\SetFigFont{8}{9.6}{rm}$p^2_1$}}}}}
\put(1437,762){\makebox(0,0)[lb]{\smash{{{\SetFigFont{8}{9.6}{rm}$p^2_1$}}}}}
\put(987,462){\makebox(0,0)[lb]{\smash{{{\SetFigFont{8}{9.6}{rm}$q$}}}}}
\put(7212,1437){\makebox(0,0)[lb]{\smash{{{\SetFigFont{8}{9.6}{rm}$p^2_1$}}}}}
\put(7137,3537){\makebox(0,0)[lb]{\smash{{{\SetFigFont{8}{9.6}{rm}$p^1_1$}}}}}
\put(7962,2487){\makebox(0,0)[lb]{\smash{{{\SetFigFont{8}{9.6}{rm}$q$}}}}}
\put(7287,687){\makebox(0,0)[lb]{\smash{{{\SetFigFont{8}{9.6}{rm}$q$}}}}}
\end{picture}

%% file: inap.tex
\begingroup\makeatletter\ifx\SetFigFont\undefined
\def\x#1#2#3#4#5#6#7\relax{\def\x{#1#2#3#4#5#6}}%
\expandafter\x\fmtname xxxxxx\relax \def\y{splain}%
\ifx\x\y   
\gdef\SetFigFont#1#2#3{%
  \ifnum #1<17\tiny\else \ifnum #1<20\small\else
  \ifnum #1<24\normalsize\else \ifnum #1<29\large\else
  \ifnum #1<34\Large\else \ifnum #1<41\LARGE\else
     \huge\fi\fi\fi\fi\fi\fi
  \csname #3\endcsname}%
\else
\gdef\SetFigFont#1#2#3{\begingroup
  \count@#1\relax \ifnum 25<\count@\count@25\fi
  \def\x{\endgroup\@setsize\SetFigFont{#2pt}}%
  \expandafter\x
    \csname \romannumeral\the\count@ pt\expandafter\endcsname
    \csname @\romannumeral\the\count@ pt\endcsname
  \csname #3\endcsname}%
\fi
\fi\endgroup
\begin{picture}(2124,6660)(0,-10)
\thicklines
\put(1212.000,3174.000){\arc{150.000}{4.7124}{7.8540}}
\put(1362.000,3174.000){\arc{150.000}{1.5708}{4.7124}}
\put(1212.000,6174.000){\arc{150.000}{4.7124}{7.8540}}
\put(1362.000,6174.000){\arc{150.000}{1.5708}{4.7124}}
\put(1287,4599){\blacken\ellipse{74}{74}}
\put(1287,4599){\ellipse{74}{74}}
\put(612,4599){\blacken\ellipse{74}{74}}
\put(612,4599){\ellipse{74}{74}}
\put(612,3249){\blacken\ellipse{74}{74}}
\put(612,3249){\ellipse{74}{74}}
\put(612,3099){\blacken\ellipse{74}{74}}
\put(612,3099){\ellipse{74}{74}}
\put(1287,3174){\blacken\ellipse{74}{74}}
\put(1287,3174){\ellipse{74}{74}}
\put(612,324){\blacken\ellipse{74}{74}}
\put(612,324){\ellipse{74}{74}}
\put(1287,324){\blacken\ellipse{74}{74}}
\put(1287,324){\ellipse{74}{74}}
\put(612,6249){\blacken\ellipse{74}{74}}
\put(612,6249){\ellipse{74}{74}}
\put(612,6099){\blacken\ellipse{74}{74}}
\put(612,6099){\ellipse{74}{74}}
\put(1287,6174){\blacken\ellipse{74}{74}}
\put(1287,6174){\ellipse{74}{74}}
\path(1062,2424)(1062,924)
\path(1032.000,1044.000)(1062.000,924.000)(1092.000,1044.000)
\path(12,4599)(2112,4599)
\path(1212,3249)(12,3249)
\path(1212,3099)(12,3099)
\path(1362,3249)(2112,3249)
\path(1362,3099)(2112,3099)
\path(12,324)(2112,324)
\path(1362,6249)(2112,6249)
\path(1362,6099)(2112,6099)
\path(12,6249)(1212,6249)
\path(1212,6099)(12,6099)
\put(1212,1599){\makebox(0,0)[lb]{\smash{{{\SetFigFont{8}{9.6}{rm}$\pi$}}}}}
\put(462,5649){\makebox(0,0)[lb]{\smash{{{\SetFigFont{8}{9.6}{rm}$x_4$}}}}}
\put(1212,4899){\makebox(0,0)[lb]{\smash{{{\SetFigFont{8}{9.6}{rm}$y_2$}}}}}
\put(462,4899){\makebox(0,0)[lb]{\smash{{{\SetFigFont{8}{9.6}{rm}$x_2$}}}}}
\put(1212,3549){\makebox(0,0)[lb]{\smash{{{\SetFigFont{8}{9.6}{rm}$y_3$}}}}}
\put(462,3549){\makebox(0,0)[lb]{\smash{{{\SetFigFont{8}{9.6}{rm}$x_3$}}}}}
\put(462,2649){\makebox(0,0)[lb]{\smash{{{\SetFigFont{8}{9.6}{rm}$x_5$}}}}}
\put(537,24){\makebox(0,0)[lb]{\smash{{{\SetFigFont{8}{9.6}{rm}$x$}}}}}
\put(462,6549){\makebox(0,0)[lb]{\smash{{{\SetFigFont{8}{9.6}{rm}$x_1$}}}}}
\put(1212,6549){\makebox(0,0)[lb]{\smash{{{\SetFigFont{8}{9.6}{rm}$y_1$}}}}}
\put(1212,24){\makebox(0,0)[lb]{\smash{{{\SetFigFont{8}{9.6}{rm}$y$}}}}}
\end{picture}

%% file: r2lines.tex
\begingroup\makeatletter\ifx\SetFigFont\undefined
\def\x#1#2#3#4#5#6#7\relax{\def\x{#1#2#3#4#5#6}}%
\expandafter\x\fmtname xxxxxx\relax \def\y{splain}%
\ifx\x\y   
\gdef\SetFigFont#1#2#3{%
  \ifnum #1<17\tiny\else \ifnum #1<20\small\else
  \ifnum #1<24\normalsize\else \ifnum #1<29\large\else
  \ifnum #1<34\Large\else \ifnum #1<41\LARGE\else
     \huge\fi\fi\fi\fi\fi\fi
  \csname #3\endcsname}%
\else
\gdef\SetFigFont#1#2#3{\begingroup
  \count@#1\relax \ifnum 25<\count@\count@25\fi
  \def\x{\endgroup\@setsize\SetFigFont{#2pt}}%
  \expandafter\x
    \csname \romannumeral\the\count@ pt\expandafter\endcsname
    \csname @\romannumeral\the\count@ pt\endcsname
  \csname #3\endcsname}%
\fi
\fi\endgroup
\begin{picture}(6024,2559)(0,-10)
\thicklines
\path(3612,387)(4212,1587)(6012,1587)
	(5412,387)(3612,387)
\path(5487,2337)(4437,12)
\dottedline{135}(4362,762)(5787,1437)
\dottedline{135}(4662,1287)(5262,612)
\dottedline{135}(4587,1812)(6012,1887)
\dottedline{135}(4662,2337)(6012,2037)
\put(4062,1737){\makebox(0,0)[lb]{\smash{{{\SetFigFont{12}{14.4}{rm}$L_4$}}}}}
\put(4212,1287){\makebox(0,0)[lb]{\smash{{{\SetFigFont{12}{14.4}{rm}$L_1$}}}}}
\put(3912,537){\makebox(0,0)[lb]{\smash{{{\SetFigFont{12}{14.4}{rm}$L_2$}}}}}
\put(5487,2412){\makebox(0,0)[lb]{\smash{{{\SetFigFont{12}{14.4}{rm}$l$}}}}}
\put(4137,2337){\makebox(0,0)[lb]{\smash{{{\SetFigFont{12}{14.4}{rm}$L_3$}}}}}
\path(12,387)(612,1587)(2412,1587)
	(1812,387)(12,387)
\path(1962,1062)(462,1137)
\dottedline{135}(912,2262)(912,87)
\dottedline{135}(1737,2262)(1437,87)
\dottedline{135}(1287,837)(1287,1437)
\dottedline{135}(612,837)(762,1437)
\put(762,2337){\makebox(0,0)[lb]{\smash{{{\SetFigFont{12}{14.4}{rm}$L_3$}}}}}
\put(1512,2337){\makebox(0,0)[lb]{\smash{{{\SetFigFont{12}{14.4}{rm}$L_4$}}}}}
\put(387,537){\makebox(0,0)[lb]{\smash{{{\SetFigFont{12}{14.4}{rm}$L_1$}}}}}
\put(1062,537){\makebox(0,0)[lb]{\smash{{{\SetFigFont{12}{14.4}{rm}$L_2$}}}}}
\put(1812,1137){\makebox(0,0)[lb]{\smash{{{\SetFigFont{12}{14.4}{rm}$l$}}}}}
\end{picture}

%% file: r92conics.tex
\begingroup\makeatletter\ifx\SetFigFont\undefined
\def\x#1#2#3#4#5#6#7\relax{\def\x{#1#2#3#4#5#6}}%
\expandafter\x\fmtname xxxxxx\relax \def\y{splain}%
\ifx\x\y   
\gdef\SetFigFont#1#2#3{%
  \ifnum #1<17\tiny\else \ifnum #1<20\small\else
  \ifnum #1<24\normalsize\else \ifnum #1<29\large\else
  \ifnum #1<34\Large\else \ifnum #1<41\LARGE\else
     \huge\fi\fi\fi\fi\fi\fi
  \csname #3\endcsname}%
\else
\gdef\SetFigFont#1#2#3{\begingroup
  \count@#1\relax \ifnum 25<\count@\count@25\fi
  \def\x{\endgroup\@setsize\SetFigFont{#2pt}}%
  \expandafter\x
    \csname \romannumeral\the\count@ pt\expandafter\endcsname
    \csname @\romannumeral\the\count@ pt\endcsname
  \csname #3\endcsname}%
\fi
\fi\endgroup
\begin{picture}(7032,10671)(0,-10)
\thicklines
\path(6825,4362)(5775,4062)
\path(6450,4812)(6150,3612)
\path(2250,2562)(1200,2262)
\path(1875,3012)(1575,1812)
\path(3225,762)(2175,462)
\path(2850,1212)(2550,12)
\put(3450,9612){\ellipse{600}{1200}}
\put(3450,7812){\ellipse{600}{1200}}
\put(2252,6030){\ellipse{600}{1200}}
\put(4727,6030){\ellipse{600}{1200}}
\put(1800,4212){\ellipse{600}{1200}}
\put(3973,2411){\ellipse{600}{1200}}
\put(5175,612){\ellipse{900}{450}}
\put(6300,2412){\ellipse{900}{450}}
\put(3975,4212){\ellipse{900}{450}}
\put(4275,2487){\blacken\ellipse{74}{74}}
\put(4275,2487){\ellipse{74}{74}}
\put(1950,6012){\blacken\ellipse{74}{74}}
\put(1950,6012){\ellipse{74}{74}}
\put(1500,4212){\blacken\ellipse{74}{74}}
\put(1500,4212){\ellipse{74}{74}}
\put(2025,2487){\blacken\ellipse{74}{74}}
\put(2025,2487){\ellipse{74}{74}}
\put(3675,2337){\blacken\ellipse{74}{74}}
\put(3675,2337){\ellipse{74}{74}}
\put(5850,2412){\blacken\ellipse{74}{74}}
\put(5850,2412){\ellipse{74}{74}}
\put(3000,687){\blacken\ellipse{74}{74}}
\put(3000,687){\ellipse{74}{74}}
\put(2347,514){\blacken\ellipse{74}{74}}
\put(2347,514){\ellipse{74}{74}}
\put(5535,754){\blacken\ellipse{74}{74}}
\put(5535,754){\ellipse{74}{74}}
\put(4778,507){\blacken\ellipse{74}{74}}
\put(4778,507){\ellipse{74}{74}}
\path(3000,9912)(2700,9312)(3900,9312)
	(4200,9912)(3000,9912)
\path(2970,8127)(2670,7527)(3870,7527)
	(4170,8127)(2970,8127)
\path(1320,2712)(1020,2112)(2220,2112)
	(2520,2712)(1320,2712)
\path(3495,2712)(3195,2112)(4395,2112)
	(4695,2712)(3495,2712)
\path(5820,2712)(5520,2112)(6720,2112)
	(7020,2712)(5820,2712)
\path(1320,4512)(1020,3912)(2220,3912)
	(2520,4512)(1320,4512)
\path(3495,4512)(3195,3912)(4395,3912)
	(4695,4512)(3495,4512)
\path(5820,4512)(5520,3912)(6720,3912)
	(7020,4512)(5820,4512)
\path(1770,6312)(1470,5712)(2670,5712)
	(2970,6312)(1770,6312)
\path(4245,6312)(3945,5712)(5145,5712)
	(5445,6312)(4245,6312)
\path(2220,912)(1920,312)(3120,312)
	(3420,912)(2220,912)
\path(4695,912)(4395,312)(5595,312)
	(5895,912)(4695,912)
\path(3450,8937)(3450,8562)
\path(3420.000,8682.000)(3450.000,8562.000)(3480.000,8682.000)
\path(3150,7362)(2550,6687)
\path(2607.301,6796.620)(2550.000,6687.000)(2652.146,6756.758)
\path(3750,7362)(4500,6687)
\path(4390.736,6744.977)(4500.000,6687.000)(4430.874,6789.575)
\path(2175,5337)(2025,4887)
\path(2034.487,5010.329)(2025.000,4887.000)(2091.408,4991.355)
\path(4350,5562)(2625,4737)
\path(2720.312,4815.839)(2625.000,4737.000)(2746.200,4761.711)
\path(4650,5337)(4425,4737)
\path(4439.045,4859.893)(4425.000,4737.000)(4495.225,4838.826)
\path(5100,5562)(6150,4812)
\path(6034.915,4857.337)(6150.000,4812.000)(6069.789,4906.161)
\path(2175,3762)(2175,2937)
\path(2145.000,3057.000)(2175.000,2937.000)(2205.000,3057.000)
\path(2325,3762)(3600,3087)
\path(3479.909,3116.633)(3600.000,3087.000)(3507.982,3169.660)
\path(2550,3912)(5475,2787)
\path(5352.229,2802.077)(5475.000,2787.000)(5373.768,2858.078)
\path(3750,1812)(3300,1137)
\path(3341.603,1253.487)(3300.000,1137.000)(3391.526,1220.205)
\path(4350,1962)(4800,1212)
\path(4712.536,1299.464)(4800.000,1212.000)(4763.985,1330.334)
\dottedline{135}(3000,8037)(3450,7662)
\dottedline{135}(4200,6012)(4725,5862)
\dottedline{135}(4725,6087)(5325,6237)
\dottedline{135}(1725,4062)(2325,4362)
\put(300,9537){\makebox(0,0)[lb]{\smash{{{\SetFigFont{12}{14.4}{rm}8}}}}}
\put(300,7662){\makebox(0,0)[lb]{\smash{{{\SetFigFont{12}{14.4}{rm}7}}}}}
\put(300,5937){\makebox(0,0)[lb]{\smash{{{\SetFigFont{12}{14.4}{rm}6}}}}}
\put(300,4137){\makebox(0,0)[lb]{\smash{{{\SetFigFont{12}{14.4}{rm}5}}}}}
\put(300,2262){\makebox(0,0)[lb]{\smash{{{\SetFigFont{12}{14.4}{rm}4}}}}}
\put(300,537){\makebox(0,0)[lb]{\smash{{{\SetFigFont{12}{14.4}{rm}3}}}}}
\put(0,10512){\makebox(0,0)[lb]{\smash{{{\SetFigFont{12}{14.4}{rm}\# general}}}}}
\put(225,10137){\makebox(0,0)[lb]{\smash{{{\SetFigFont{12}{14.4}{rm}lines}}}}}
\put(4350,9537){\makebox(0,0)[lb]{\smash{{{\SetFigFont{8}{9.6}{rm}92}}}}}
\put(4275,7662){\makebox(0,0)[lb]{\smash{{{\SetFigFont{8}{9.6}{rm}92}}}}}
\put(1650,6462){\makebox(0,0)[lb]{\smash{{{\SetFigFont{8}{9.6}{rm}18}}}}}
\put(5175,6462){\makebox(0,0)[lb]{\smash{{{\SetFigFont{8}{9.6}{rm}74}}}}}
\put(2925,5112){\makebox(0,0)[lb]{\smash{{{\SetFigFont{8}{9.6}{rm}$\times 2$}}}}}
\put(4575,4887){\makebox(0,0)[lb]{\smash{{{\SetFigFont{8}{9.6}{rm}$\times 8$}}}}}
\put(6750,4737){\makebox(0,0)[lb]{\smash{{{\SetFigFont{8}{9.6}{rm}30}}}}}
\put(4725,3912){\makebox(0,0)[lb]{\smash{{{\SetFigFont{8}{9.6}{rm}1}}}}}
\put(1050,4587){\makebox(0,0)[lb]{\smash{{{\SetFigFont{8}{9.6}{rm}18}}}}}
\put(4350,3237){\makebox(0,0)[lb]{\smash{{{\SetFigFont{8}{9.6}{rm}$\times 4$}}}}}
\put(4650,1587){\makebox(0,0)[lb]{\smash{{{\SetFigFont{8}{9.6}{rm}$\times 2$}}}}}
\put(1125,2787){\makebox(0,0)[lb]{\smash{{{\SetFigFont{8}{9.6}{rm}10}}}}}
\put(4650,2187){\makebox(0,0)[lb]{\smash{{{\SetFigFont{8}{9.6}{rm}4}}}}}
\put(6450,2862){\makebox(0,0)[lb]{\smash{{{\SetFigFont{8}{9.6}{rm}1}}}}}
\put(2250,1062){\makebox(0,0)[lb]{\smash{{{\SetFigFont{8}{9.6}{rm}2}}}}}
\put(5475,1062){\makebox(0,0)[lb]{\smash{{{\SetFigFont{8}{9.6}{rm}1}}}}}
\end{picture}

%% file: r92conics2.tex
\begingroup\makeatletter\ifx\SetFigFont\undefined
\def\x#1#2#3#4#5#6#7\relax{\def\x{#1#2#3#4#5#6}}%
\expandafter\x\fmtname xxxxxx\relax \def\y{splain}%
\ifx\x\y   
\gdef\SetFigFont#1#2#3{%
  \ifnum #1<17\tiny\else \ifnum #1<20\small\else
  \ifnum #1<24\normalsize\else \ifnum #1<29\large\else
  \ifnum #1<34\Large\else \ifnum #1<41\LARGE\else
     \huge\fi\fi\fi\fi\fi\fi
  \csname #3\endcsname}%
\else
\gdef\SetFigFont#1#2#3{\begingroup
  \count@#1\relax \ifnum 25<\count@\count@25\fi
  \def\x{\endgroup\@setsize\SetFigFont{#2pt}}%
  \expandafter\x
    \csname \romannumeral\the\count@ pt\expandafter\endcsname
    \csname @\romannumeral\the\count@ pt\endcsname
  \csname #3\endcsname}%
\fi
\fi\endgroup
\begin{picture}(4626,3466)(0,-10)
\thicklines
\put(3450,2407){\ellipse{600}{1200}}
\put(3450,607){\ellipse{600}{1200}}
\path(3000,2707)(2700,2107)(3900,2107)
	(4200,2707)(3000,2707)
\path(2970,922)(2670,322)(3870,322)
	(4170,922)(2970,922)
\path(3450,1732)(3450,1357)
\path(3420.000,1477.000)(3450.000,1357.000)(3480.000,1477.000)
\dottedline{135}(3000,832)(3450,457)
\put(4350,2332){\makebox(0,0)[lb]{\smash{{{\SetFigFont{8}{14.4}{rm}184}}}}}
\put(4275,457){\makebox(0,0)[lb]{\smash{{{\SetFigFont{8}{14.4}{rm}92}}}}}
\put(300,2332){\makebox(0,0)[lb]{\smash{{{\SetFigFont{12}{14.4}{rm}8}}}}}
\put(300,457){\makebox(0,0)[lb]{\smash{{{\SetFigFont{12}{14.4}{rm}7}}}}}
\put(0,3307){\makebox(0,0)[lb]{\smash{{{\SetFigFont{12}{14.4}{rm}\# general}}}}}
\put(225,2932){\makebox(0,0)[lb]{\smash{{{\SetFigFont{12}{14.4}{rm}lines}}}}}
\put(3600,1507){\makebox(0,0)[lb]{\smash{{{\SetFigFont{8}{14.4}{rm}$\times 2$}}}}}
\end{picture}

%% file: rtype2eg.tex
\begingroup\makeatletter\ifx\SetFigFont\undefined
\def\x#1#2#3#4#5#6#7\relax{\def\x{#1#2#3#4#5#6}}%
\expandafter\x\fmtname xxxxxx\relax \def\y{splain}%
\ifx\x\y   
\gdef\SetFigFont#1#2#3{%
  \ifnum #1<17\tiny\else \ifnum #1<20\small\else
  \ifnum #1<24\normalsize\else \ifnum #1<29\large\else
  \ifnum #1<34\Large\else \ifnum #1<41\LARGE\else
     \huge\fi\fi\fi\fi\fi\fi
  \csname #3\endcsname}%
\else
\gdef\SetFigFont#1#2#3{\begingroup
  \count@#1\relax \ifnum 25<\count@\count@25\fi
  \def\x{\endgroup\@setsize\SetFigFont{#2pt}}%
  \expandafter\x
    \csname \romannumeral\the\count@ pt\expandafter\endcsname
    \csname @\romannumeral\the\count@ pt\endcsname
  \csname #3\endcsname}%
\fi
\fi\endgroup
\begin{picture}(10824,7409)(0,-10)
\thicklines
\path(2472,7287)(1812,12)
\path(1212,4287)(12,687)(9612,687)
	(10812,4287)(1212,4287)
\path(1587,2637)	(1626.783,2616.118)
	(1666.384,2595.371)
	(1705.803,2574.761)
	(1745.041,2554.286)
	(1784.098,2533.948)
	(1822.975,2513.745)
	(1861.672,2493.677)
	(1900.190,2473.745)
	(1938.529,2453.949)
	(1976.690,2434.287)
	(2014.673,2414.761)
	(2052.479,2395.370)
	(2090.108,2376.113)
	(2127.561,2356.992)
	(2201.940,2319.153)
	(2275.620,2281.852)
	(2348.605,2245.088)
	(2420.899,2208.860)
	(2492.507,2173.168)
	(2563.432,2138.011)
	(2633.678,2103.387)
	(2703.249,2069.295)
	(2772.149,2035.736)
	(2840.383,2002.708)
	(2907.954,1970.210)
	(2974.866,1938.242)
	(3041.124,1906.802)
	(3106.731,1875.890)
	(3171.692,1845.504)
	(3236.009,1815.645)
	(3299.689,1786.311)
	(3362.734,1757.500)
	(3425.148,1729.214)
	(3486.936,1701.449)
	(3548.101,1674.207)
	(3608.648,1647.485)
	(3668.580,1621.283)
	(3727.903,1595.601)
	(3786.618,1570.436)
	(3844.732,1545.789)
	(3902.247,1521.659)
	(3959.167,1498.044)
	(4015.498,1474.944)
	(4071.242,1452.357)
	(4126.404,1430.284)
	(4180.988,1408.724)
	(4234.998,1387.674)
	(4288.437,1367.135)
	(4341.310,1347.106)
	(4393.622,1327.585)
	(4445.375,1308.573)
	(4496.574,1290.067)
	(4547.223,1272.068)
	(4597.326,1254.574)
	(4646.887,1237.584)
	(4695.910,1221.099)
	(4744.399,1205.116)
	(4792.358,1189.634)
	(4839.791,1174.654)
	(4886.702,1160.175)
	(4933.096,1146.194)
	(4978.976,1132.712)
	(5024.345,1119.728)
	(5069.209,1107.240)
	(5113.572,1095.249)
	(5157.436,1083.752)
	(5200.807,1072.750)
	(5243.688,1062.241)
	(5286.084,1052.224)
	(5327.998,1042.699)
	(5369.434,1033.665)
	(5410.396,1025.121)
	(5450.889,1017.066)
	(5490.917,1009.499)
	(5530.483,1002.419)
	(5569.591,995.826)
	(5608.246,989.718)
	(5646.452,984.095)
	(5684.212,978.956)
	(5758.412,970.127)
	(5830.878,963.223)
	(5901.643,958.237)
	(5970.739,955.163)
	(6038.198,953.993)
	(6104.053,954.721)
	(6168.335,957.341)
	(6231.076,961.844)
	(6292.309,968.225)
	(6352.067,976.476)
	(6410.380,986.590)
	(6467.282,998.561)
	(6522.805,1012.382)
	(6576.980,1028.046)
	(6629.841,1045.546)
	(6681.418,1064.875)
	(6731.745,1086.026)
	(6780.853,1108.993)
	(6828.775,1133.769)
	(6875.543,1160.346)
	(6921.189,1188.718)
	(6965.745,1218.877)
	(7009.244,1250.819)
	(7051.717,1284.534)
	(7093.197,1320.016)
	(7133.716,1357.260)
	(7173.306,1396.256)
	(7212.000,1437.000)

\path(7212,1437)	(7244.772,1475.085)
	(7274.262,1514.570)
	(7300.433,1555.569)
	(7323.253,1598.193)
	(7342.687,1642.558)
	(7358.701,1688.775)
	(7371.260,1736.959)
	(7380.329,1787.222)
	(7385.875,1839.679)
	(7387.864,1894.441)
	(7386.260,1951.623)
	(7381.029,2011.338)
	(7372.138,2073.699)
	(7359.552,2138.819)
	(7343.237,2206.812)
	(7323.157,2277.791)
	(7299.280,2351.870)
	(7285.906,2390.106)
	(7271.570,2429.160)
	(7256.267,2469.046)
	(7239.994,2509.777)
	(7222.744,2551.368)
	(7204.516,2593.833)
	(7185.303,2637.186)
	(7165.103,2681.442)
	(7143.910,2726.613)
	(7121.720,2772.716)
	(7098.529,2819.763)
	(7074.333,2867.769)
	(7049.127,2916.748)
	(7022.907,2966.715)
	(6995.669,3017.682)
	(6967.409,3069.666)
	(6938.122,3122.679)
	(6907.804,3176.736)
	(6876.450,3231.851)
	(6844.057,3288.038)
	(6810.620,3345.312)
	(6776.134,3403.687)
	(6740.596,3463.175)
	(6704.002,3523.793)
	(6666.346,3585.554)
	(6627.625,3648.472)
	(6587.834,3712.562)
	(6546.969,3777.837)
	(6505.026,3844.311)
	(6462.000,3912.000)

\path(5487,7347)	(5497.335,7277.966)
	(5507.627,7209.653)
	(5517.877,7142.057)
	(5528.086,7075.175)
	(5538.256,7009.006)
	(5548.386,6943.545)
	(5558.478,6878.791)
	(5568.532,6814.740)
	(5578.551,6751.389)
	(5588.533,6688.737)
	(5598.482,6626.779)
	(5608.397,6565.513)
	(5618.279,6504.936)
	(5628.129,6445.046)
	(5637.949,6385.840)
	(5647.739,6327.314)
	(5657.500,6269.467)
	(5667.233,6212.294)
	(5676.938,6155.794)
	(5686.618,6099.963)
	(5696.272,6044.799)
	(5705.902,5990.299)
	(5715.508,5936.460)
	(5725.092,5883.279)
	(5734.655,5830.754)
	(5744.197,5778.881)
	(5753.719,5727.658)
	(5763.222,5677.082)
	(5772.707,5627.150)
	(5782.176,5577.859)
	(5791.628,5529.207)
	(5801.065,5481.190)
	(5810.488,5433.806)
	(5819.898,5387.052)
	(5829.296,5340.926)
	(5838.682,5295.423)
	(5848.057,5250.542)
	(5857.423,5206.280)
	(5866.780,5162.634)
	(5876.130,5119.601)
	(5885.473,5077.178)
	(5894.809,5035.363)
	(5904.141,4994.152)
	(5913.469,4953.543)
	(5922.794,4913.533)
	(5932.116,4874.119)
	(5941.438,4835.299)
	(5950.758,4797.069)
	(5969.403,4722.369)
	(5988.057,4649.997)
	(6006.728,4579.931)
	(6025.423,4512.148)
	(6044.151,4446.625)
	(6062.917,4383.340)
	(6081.730,4322.271)
	(6100.597,4263.394)
	(6119.526,4206.687)
	(6138.523,4152.128)
	(6157.596,4099.693)
	(6176.753,4049.362)
	(6196.000,4001.110)
	(6215.346,3954.915)
	(6234.798,3910.755)
	(6254.362,3868.607)
	(6274.047,3828.449)
	(6293.859,3790.257)
	(6333.897,3719.685)
	(6374.535,3656.710)
	(6415.831,3601.151)
	(6457.845,3552.828)
	(6500.636,3511.560)
	(6544.262,3477.167)
	(6588.784,3449.468)
	(6634.259,3428.282)
	(6680.748,3413.429)
	(6728.308,3404.729)
	(6777.000,3402.000)

\path(6777,3402)	(6841.034,3413.742)
	(6901.689,3447.459)
	(6959.145,3501.676)
	(7013.580,3574.919)
	(7039.721,3618.215)
	(7065.173,3665.713)
	(7089.960,3717.231)
	(7114.103,3772.583)
	(7137.625,3831.585)
	(7160.548,3894.053)
	(7182.894,3959.802)
	(7204.687,4028.649)
	(7225.948,4100.408)
	(7246.699,4174.896)
	(7256.891,4213.106)
	(7266.963,4251.928)
	(7276.920,4291.340)
	(7286.763,4331.319)
	(7296.495,4371.842)
	(7306.120,4412.886)
	(7315.639,4454.428)
	(7325.057,4496.444)
	(7334.375,4538.912)
	(7343.596,4581.808)
	(7352.724,4625.110)
	(7361.760,4668.795)
	(7370.708,4712.839)
	(7379.571,4757.219)
	(7388.351,4801.913)
	(7397.051,4846.897)
	(7405.674,4892.149)
	(7414.223,4937.644)
	(7422.700,4983.361)
	(7431.109,5029.276)
	(7439.451,5075.366)
	(7447.731,5121.608)
	(7455.950,5167.979)
	(7464.112,5214.457)
	(7472.219,5261.017)
	(7480.274,5307.636)
	(7488.279,5354.293)
	(7496.239,5400.964)
	(7504.155,5447.625)
	(7512.030,5494.254)
	(7519.867,5540.828)
	(7527.669,5587.323)
	(7535.438,5633.717)
	(7543.178,5679.986)
	(7550.891,5726.108)
	(7558.580,5772.059)
	(7566.248,5817.816)
	(7573.897,5863.357)
	(7581.531,5908.659)
	(7589.152,5953.697)
	(7596.763,5998.450)
	(7604.366,6042.893)
	(7611.965,6087.005)
	(7619.563,6130.762)
	(7627.161,6174.141)
	(7634.763,6217.119)
	(7642.373,6259.673)
	(7649.991,6301.780)
	(7657.622,6343.416)
	(7665.268,6384.560)
	(7672.932,6425.187)
	(7680.616,6465.275)
	(7688.324,6504.800)
	(7696.058,6543.740)
	(7703.822,6582.072)
	(7711.617,6619.772)
	(7727.313,6693.186)
	(7743.171,6763.797)
	(7759.212,6831.421)
	(7775.458,6895.874)
	(7791.932,6956.972)
	(7808.657,7014.530)
	(7825.654,7068.363)
	(7842.946,7118.288)
	(7860.555,7164.119)
	(7878.504,7205.674)
	(7915.509,7275.213)
	(7954.141,7325.430)
	(7994.579,7354.851)
	(8037.000,7362.000)

\path(8037,7362)	(8093.143,7350.890)
	(8146.833,7329.159)
	(8198.104,7296.476)
	(8246.991,7252.506)
	(8293.527,7196.916)
	(8337.746,7129.374)
	(8358.999,7091.017)
	(8379.685,7049.547)
	(8399.809,7004.922)
	(8419.375,6957.101)
	(8438.389,6906.042)
	(8456.853,6851.703)
	(8474.773,6794.044)
	(8492.152,6733.022)
	(8508.995,6668.595)
	(8525.306,6600.722)
	(8541.090,6529.362)
	(8556.351,6454.472)
	(8563.786,6415.691)
	(8571.093,6376.012)
	(8578.270,6335.429)
	(8585.320,6293.939)
	(8592.242,6251.534)
	(8599.036,6208.211)
	(8605.705,6163.965)
	(8612.247,6118.788)
	(8618.664,6072.678)
	(8624.956,6025.628)
	(8631.124,5977.633)
	(8637.167,5928.689)
	(8643.088,5878.789)
	(8648.886,5827.929)
	(8654.561,5776.103)
	(8660.115,5723.306)
	(8665.548,5669.534)
	(8670.860,5614.780)
	(8676.052,5559.040)
	(8681.124,5502.308)
	(8686.078,5444.580)
	(8690.913,5385.849)
	(8695.630,5326.111)
	(8700.229,5265.361)
	(8704.712,5203.594)
	(8709.079,5140.803)
	(8713.329,5076.985)
	(8717.465,5012.133)
	(8721.485,4946.243)
	(8725.392,4879.309)
	(8729.185,4811.326)
	(8732.865,4742.290)
	(8736.432,4672.194)
	(8739.887,4601.033)
	(8743.231,4528.803)
	(8746.464,4455.499)
	(8749.586,4381.114)
	(8751.106,4343.515)
	(8752.598,4305.644)
	(8754.064,4267.500)
	(8755.501,4229.083)
	(8756.912,4190.392)
	(8758.296,4151.427)
	(8759.652,4112.186)
	(8760.982,4072.670)
	(8762.284,4032.877)
	(8763.560,3992.806)
	(8764.809,3952.458)
	(8766.031,3911.832)
	(8767.227,3870.926)
	(8768.395,3829.741)
	(8769.538,3788.275)
	(8770.654,3746.528)
	(8771.743,3704.500)
	(8772.806,3662.189)
	(8773.843,3619.595)
	(8774.854,3576.717)
	(8775.839,3533.556)
	(8776.797,3490.109)
	(8777.730,3446.376)
	(8778.636,3402.357)
	(8779.517,3358.052)
	(8780.372,3313.458)
	(8781.201,3268.577)
	(8782.005,3223.406)
	(8782.783,3177.946)
	(8783.535,3132.196)
	(8784.262,3086.155)
	(8784.964,3039.823)
	(8785.640,2993.198)
	(8786.291,2946.281)
	(8786.917,2899.070)
	(8787.517,2851.566)
	(8788.093,2803.766)
	(8788.644,2755.671)
	(8789.169,2707.280)
	(8789.670,2658.592)
	(8790.146,2609.607)
	(8790.598,2560.325)
	(8791.024,2510.743)
	(8791.426,2460.862)
	(8791.804,2410.681)
	(8792.157,2360.200)
	(8792.486,2309.417)
	(8792.790,2258.333)
	(8793.070,2206.946)
	(8793.326,2155.256)
	(8793.558,2103.262)
	(8793.766,2050.963)
	(8793.950,1998.359)
	(8794.110,1945.450)
	(8794.246,1892.234)
	(8794.358,1838.711)
	(8794.447,1784.880)
	(8794.512,1730.741)
	(8794.553,1676.293)
	(8794.571,1621.535)
	(8794.565,1566.467)
	(8794.536,1511.088)
	(8794.484,1455.397)
	(8794.408,1399.394)
	(8794.310,1343.078)
	(8794.188,1286.448)
	(8794.043,1229.504)
	(8793.875,1172.246)
	(8793.684,1114.671)
	(8793.470,1056.781)
	(8793.234,998.574)
	(8792.975,940.049)
	(8792.693,881.206)
	(8792.389,822.044)
	(8792.062,762.563)
	(8791.713,702.761)
	(8791.341,642.639)
	(8790.947,582.196)
	(8790.531,521.430)
	(8790.092,460.342)
	(8789.632,398.931)
	(8789.149,337.195)
	(8788.645,275.135)
	(8788.118,212.749)
	(8787.570,150.038)
	(8787.000,87.000)

\put(5517,1602){\makebox(0,0)[lb]{\smash{{{\SetFigFont{12}{14.4}{rm}$C(0)$}}}}}
\put(2637,6477){\makebox(0,0)[lb]{\smash{{{\SetFigFont{12}{14.4}{rm}$C(1)$}}}}}
\put(8802,6102){\makebox(0,0)[lb]{\smash{{{\SetFigFont{12}{14.4}{rm}$C(2)$}}}}}
\end{picture}

%% file: rspecialize.tex
\begingroup\makeatletter\ifx\SetFigFont\undefined
\def\x#1#2#3#4#5#6#7\relax{\def\x{#1#2#3#4#5#6}}%
\expandafter\x\fmtname xxxxxx\relax \def\y{splain}%
\ifx\x\y   
\gdef\SetFigFont#1#2#3{%
  \ifnum #1<17\tiny\else \ifnum #1<20\small\else
  \ifnum #1<24\normalsize\else \ifnum #1<29\large\else
  \ifnum #1<34\Large\else \ifnum #1<41\LARGE\else
     \huge\fi\fi\fi\fi\fi\fi
  \csname #3\endcsname}%
\else
\gdef\SetFigFont#1#2#3{\begingroup
  \count@#1\relax \ifnum 25<\count@\count@25\fi
  \def\x{\endgroup\@setsize\SetFigFont{#2pt}}%
  \expandafter\x
    \csname \romannumeral\the\count@ pt\expandafter\endcsname
    \csname @\romannumeral\the\count@ pt\endcsname
  \csname #3\endcsname}%
\fi
\fi\endgroup
\begin{picture}(6938,5166)(0,-10)
\thicklines
\path(2712,1539)(3912,2439)
\path(2712,3339)(3912,3039)
\dottedline{135}(2712,39)(3912,939)(3912,5139)
	(2712,4239)(2712,39)
\dottedline{135}(12,1539)(1212,2439)(6312,2439)
	(5112,1539)(12,1539)
\path(4512,3939)(3312,3339)
\path(3405.915,3419.498)(3312.000,3339.000)(3432.748,3365.833)
\path(4212,939)(3612,1839)
\path(3703.526,1755.795)(3612.000,1839.000)(3653.603,1722.513)
\put(912,1839){\makebox(0,0)[lb]{\smash{{{\SetFigFont{12}{14.4}{rm}$H$}}}}}
\put(3312,4239){\makebox(0,0)[lb]{\smash{{{\SetFigFont{12}{14.4}{rm}$Q$}}}}}
\put(4512,3939){\makebox(0,0)[lb]{\smash{{{\SetFigFont{12}{14.4}{rm}$H' \cap Q$}}}}}
\put(4212,639){\makebox(0,0)[lb]{\smash{{{\SetFigFont{12}{14.4}{rm}$H \cap Q$}}}}}
\put(4662,3339){\makebox(0,0)[lb]{\smash{{{\SetFigFont{12}{14.4}{rm}(the general $(E-1)$-plane)}}}}}
\put(4362,39){\makebox(0,0)[lb]{\smash{{{\SetFigFont{12}{14.4}{rm}(the specialized $(E-1)$-plane)}}}}}
\end{picture}

%% file: rYeg.tex
\begingroup\makeatletter\ifx\SetFigFont\undefined
\def\x#1#2#3#4#5#6#7\relax{\def\x{#1#2#3#4#5#6}}%
\expandafter\x\fmtname xxxxxx\relax \def\y{splain}%
\ifx\x\y   
\gdef\SetFigFont#1#2#3{%
  \ifnum #1<17\tiny\else \ifnum #1<20\small\else
  \ifnum #1<24\normalsize\else \ifnum #1<29\large\else
  \ifnum #1<34\Large\else \ifnum #1<41\LARGE\else
     \huge\fi\fi\fi\fi\fi\fi
  \csname #3\endcsname}%
\else
\gdef\SetFigFont#1#2#3{\begingroup
  \count@#1\relax \ifnum 25<\count@\count@25\fi
  \def\x{\endgroup\@setsize\SetFigFont{#2pt}}%
  \expandafter\x
    \csname \romannumeral\the\count@ pt\expandafter\endcsname
    \csname @\romannumeral\the\count@ pt\endcsname
  \csname #3\endcsname}%
\fi
\fi\endgroup
\begin{picture}(9549,9429)(0,-10)
\thicklines
\path(1437,2697)(7437,3747)
\path(1812,4497)(12,1497)(7512,1497)
	(9312,4497)(1812,4497)
\path(8112,9402)(4407,12)
\dottedline{135}(2037,3672)(2112,1947)
\dottedline{135}(2487,3972)(3987,2097)
\dottedline{135}(4587,4122)(4662,2322)
\dottedline{135}(5112,5097)(4512,4797)(6612,4797)
	(7212,5097)(5112,5097)
\dottedline{135}(5712,5997)(5112,5697)(7212,5697)
	(7812,5997)(5712,5997)
\dottedline{135}(6312,6897)(5712,6597)(7812,6597)
	(8412,6897)(6312,6897)
\dottedline{135}(6912,7797)(6312,7497)(8412,7497)
	(9012,7797)(6912,7797)
\dottedline{135}(7437,8697)(6837,8397)(8937,8397)
	(9537,8697)(7437,8697)
\put(7812,3897){\makebox(0,0)[lb]{\smash{{{\SetFigFont{12}{14.4}{rm}$L_0$}}}}}
\put(7212,1797){\makebox(0,0)[lb]{\smash{{{\SetFigFont{12}{14.4}{rm}$H$}}}}}
\put(4872,297){\makebox(0,0)[lb]{\smash{{{\SetFigFont{12}{14.4}{rm}$L_1$}}}}}
\end{picture}

%% file: e1500cubics.tex
\begingroup\makeatletter\ifx\SetFigFont\undefined
\def\x#1#2#3#4#5#6#7\relax{\def\x{#1#2#3#4#5#6}}%
\expandafter\x\fmtname xxxxxx\relax \def\y{splain}%
\ifx\x\y   
\gdef\SetFigFont#1#2#3{%
  \ifnum #1<17\tiny\else \ifnum #1<20\small\else
  \ifnum #1<24\normalsize\else \ifnum #1<29\large\else
  \ifnum #1<34\Large\else \ifnum #1<41\LARGE\else
     \huge\fi\fi\fi\fi\fi\fi
  \csname #3\endcsname}%
\else
\gdef\SetFigFont#1#2#3{\begingroup
  \count@#1\relax \ifnum 25<\count@\count@25\fi
  \def\x{\endgroup\@setsize\SetFigFont{#2pt}}%
  \expandafter\x
    \csname \romannumeral\the\count@ pt\expandafter\endcsname
    \csname @\romannumeral\the\count@ pt\endcsname
  \csname #3\endcsname}%
\fi
\fi\endgroup
\begin{picture}(13212,12474)(0,-10)
\thicklines
\put(4665,4200){\ellipse{660}{1200}}
\path(4185,4515)(3885,3915)(5085,3915)
	(5385,4515)(4185,4515)
\path(4185,4215)(5115,4215)
\put(4680,2415){\ellipse{660}{1200}}
\path(4200,2730)(3900,2130)(5100,2130)
	(5400,2730)(4200,2730)
\path(4200,2430)(5130,2430)
\put(11280,615){\ellipse{660}{1200}}
\path(10800,930)(10500,330)(11700,330)
	(12000,930)(10800,930)
\path(10800,630)(11730,630)
\put(8625,7815){\blacken\ellipse{74}{74}}
\put(8625,7815){\ellipse{74}{74}}
\put(12150,4215){\blacken\ellipse{74}{74}}
\put(12150,4215){\ellipse{74}{74}}
\put(9750,2265){\blacken\ellipse{74}{74}}
\put(9750,2265){\ellipse{74}{74}}
\put(4650,2415){\blacken\ellipse{74}{74}}
\put(4650,2415){\ellipse{74}{74}}
\put(7050,2415){\blacken\ellipse{74}{74}}
\put(7050,2415){\ellipse{74}{74}}
\put(7650,615){\blacken\ellipse{74}{74}}
\put(7650,615){\ellipse{74}{74}}
\put(7415,11950){\ellipse{150}{150}}
\put(7325,10110){\ellipse{150}{150}}
\put(9750,4065){\blacken\ellipse{74}{74}}
\put(9750,4065){\ellipse{74}{74}}
\put(5925,8330){\ellipse{150}{150}}
\put(8565,8340){\ellipse{150}{150}}
\put(7035,6540){\ellipse{150}{150}}
\put(9485,6540){\ellipse{150}{150}}
\put(7025,4700){\ellipse{150}{150}}
\put(9435,4750){\ellipse{150}{150}}
\put(9435,2930){\ellipse{150}{150}}
\put(12150,2415){\blacken\ellipse{74}{74}}
\put(12150,2415){\ellipse{74}{74}}
\put(12750,4145){\ellipse{150}{150}}
\put(12810,2330){\ellipse{150}{150}}
\put(8325,545){\ellipse{150}{150}}
\put(7740,2315){\ellipse{150}{150}}
\put(4995,5902){\ellipse{150}{150}}
\put(2260,4112){\ellipse{150}{150}}
\put(7080,4252){\blacken\ellipse{74}{74}}
\put(7080,4252){\ellipse{74}{74}}
\put(9473,4380){\blacken\ellipse{74}{74}}
\put(9473,4380){\ellipse{74}{74}}
\put(9548,6022){\blacken\ellipse{74}{74}}
\put(9548,6022){\ellipse{74}{74}}
\put(9465,2572){\blacken\ellipse{74}{74}}
\put(9465,2572){\ellipse{74}{74}}
\put(12593,2527){\blacken\ellipse{74}{74}}
\put(12593,2527){\ellipse{74}{74}}
\put(8100,735){\blacken\ellipse{74}{74}}
\put(8100,735){\ellipse{74}{74}}
\put(11093,622){\blacken\ellipse{74}{74}}
\put(11093,622){\ellipse{74}{74}}
\put(11408,637){\blacken\ellipse{74}{74}}
\put(11408,637){\ellipse{74}{74}}
\dottedline{135}(7500,9765)(7200,9465)
\dottedline{135}(6075,8040)(5775,7665)
\dottedline{135}(6150,7665)(6375,7815)
\dottedline{135}(9750,6165)(9975,5865)
\path(7800,10890)(7800,10440)
\path(7800,10890)(7800,10440)
\path(7770.000,10560.000)(7800.000,10440.000)(7830.000,10560.000)
\path(7425,8940)(6750,8340)
\path(7425,8940)(6750,8340)
\path(6819.758,8442.146)(6750.000,8340.000)(6859.620,8397.301)
\path(7800,9015)(8400,8340)
\path(7800,9015)(8400,8340)
\path(8297.854,8409.758)(8400.000,8340.000)(8342.699,8449.620)
\path(9525,7290)(9750,6765)
\path(9525,7290)(9750,6765)
\path(9675.155,6863.480)(9750.000,6765.000)(9730.304,6887.115)
\path(6450,7215)(5475,6540)
\path(6450,7215)(5475,6540)
\path(5556.587,6632.971)(5475.000,6540.000)(5590.739,6583.639)
\path(6825,7215)(7050,6765)
\path(6825,7215)(7050,6765)
\path(6969.502,6858.915)(7050.000,6765.000)(7023.167,6885.748)
\path(6975,7365)(9225,6390)
\path(6975,7365)(9225,6390)
\path(9102.965,6410.186)(9225.000,6390.000)(9126.822,6465.240)
\path(9375,6315)(9075,5715)(10275,5715)
	(10575,6315)(9375,6315)
\path(5850,8115)(5550,7515)(6750,7515)
	(7050,8115)(5850,8115)
\path(8475,8115)(8175,7515)(9375,7515)
	(9675,8115)(8475,8115)
\path(7200,9915)(6900,9315)(8100,9315)
	(8400,9915)(7200,9915)
\path(7350,11715)(7050,11115)(8250,11115)
	(8550,11715)(7350,11715)
\dottedline{135}(7185,4218)(7560,4068)
\path(6885,4518)(6585,3918)(7785,3918)
	(8085,4518)(6885,4518)
\path(9300,2718)(9000,2118)(10200,2118)
	(10500,2718)(9300,2718)
\path(6450,5340)(2400,4665)
\path(6450,5340)(2400,4665)
\path(2513.435,4714.320)(2400.000,4665.000)(2523.299,4655.136)
\path(6750,5340)(5625,4590)
\path(6750,5340)(5625,4590)
\path(5708.205,4681.526)(5625.000,4590.000)(5741.487,4631.603)
\path(7050,5415)(7275,4740)
\path(7050,5415)(7275,4740)
\path(7208.592,4844.355)(7275.000,4740.000)(7265.513,4863.329)
\path(9000,5265)(7875,4665)
\path(9000,5265)(7875,4665)
\path(7966.765,4747.941)(7875.000,4665.000)(7995.000,4695.000)
\path(9675,5415)(9675,4740)
\path(9675,5415)(9675,4740)
\path(9645.000,4860.000)(9675.000,4740.000)(9705.000,4860.000)
\path(10425,5565)(12225,4665)
\path(10425,5565)(12225,4665)
\path(12104.252,4691.833)(12225.000,4665.000)(12131.085,4745.498)
\path(6375,3540)(5100,2865)
\path(6375,3540)(5100,2865)
\path(5192.018,2947.660)(5100.000,2865.000)(5220.091,2894.633)
\path(7050,3615)(6975,2865)
\path(7050,3615)(6975,2865)
\path(6957.089,2987.390)(6975.000,2865.000)(7016.792,2981.419)
\path(7575,3615)(9300,2865)
\path(7575,3615)(9300,2865)
\path(9177.990,2885.335)(9300.000,2865.000)(9201.913,2940.359)
\path(10050,3690)(10050,3015)
\path(10050,3690)(10050,3015)
\path(10020.000,3135.000)(10050.000,3015.000)(10080.000,3135.000)
\path(10275,3765)(11850,2940)
\path(10275,3765)(11850,2940)
\path(11729.780,2969.106)(11850.000,2940.000)(11757.620,3022.256)
\path(9225,1740)(8400,1065)
\path(9225,1740)(8400,1065)
\path(8473.878,1164.207)(8400.000,1065.000)(8511.872,1117.770)
\path(10050,1815)(10725,990)
\path(10050,1815)(10725,990)
\path(10625.793,1063.878)(10725.000,990.000)(10672.230,1101.872)
\dottedline{135}(7560,4293)(7935,4443)
\dottedline{135}(1275,3990)(1800,4440)
\dottedline{135}(1950,4440)(1500,3990)
\dottedline{135}(1875,3990)(2475,4440)
\dottedline{135}(7125,2640)(7350,2190)
\dottedline{135}(7493,2197)(7815,2670)
\dottedline{135}(9975,2265)(10275,2640)
\dottedline{135}(7703,840)(8175,390)
\dottedline{135}(7200,6240)(6825,5940)
\dottedline{135}(7200,6015)(7575,5865)
\dottedline{135}(7575,6090)(7950,6240)
\path(6900,6315)(6600,5715)(7800,5715)
	(8100,6315)(6900,6315)
\path(9300,4518)(9000,3918)(10200,3918)
	(10500,4518)(9300,4518)
\path(4200,6315)(3900,5715)(5100,5715)
	(5400,6315)(4200,6315)
\path(12000,4515)(11700,3915)(12900,3915)
	(13200,4515)(12000,4515)
\path(12000,2715)(11700,2115)(12900,2115)
	(13200,2715)(12000,2715)
\path(7500,915)(7200,315)(8400,315)
	(8700,915)(7500,915)
\path(6900,2715)(6600,2115)(7800,2115)
	(8100,2715)(6900,2715)
\path(1485,4515)(1185,3915)(2385,3915)
	(2685,4515)(1485,4515)
\path(7425,11865)	(7432.821,11792.648)
	(7440.496,11723.330)
	(7448.035,11656.996)
	(7455.449,11593.600)
	(7462.748,11533.092)
	(7469.944,11475.426)
	(7477.046,11420.552)
	(7484.065,11368.424)
	(7491.011,11318.993)
	(7497.896,11272.211)
	(7504.729,11228.029)
	(7511.521,11186.401)
	(7518.283,11147.277)
	(7525.025,11110.611)
	(7538.493,11044.456)
	(7552.007,10987.553)
	(7565.652,10939.517)
	(7593.674,10868.510)
	(7655.000,10816.000)

\path(7655,10816)	(7697.385,10868.018)
	(7713.549,10927.227)
	(7720.610,10963.604)
	(7727.079,11003.913)
	(7733.024,11047.719)
	(7738.513,11094.584)
	(7743.612,11144.072)
	(7748.389,11195.746)
	(7752.912,11249.169)
	(7757.247,11303.906)
	(7761.462,11359.519)
	(7765.625,11415.573)
	(7769.802,11471.629)
	(7774.062,11527.252)
	(7778.470,11582.006)
	(7783.095,11635.453)
	(7788.004,11687.157)
	(7793.265,11736.681)
	(7798.944,11783.590)
	(7805.109,11827.445)
	(7811.827,11867.812)
	(7819.166,11904.252)
	(7835.974,11963.608)
	(7880.000,12016.000)

\path(7880,12016)	(7920.301,12001.276)
	(7958.043,11952.694)
	(7976.230,11914.539)
	(7994.107,11866.518)
	(8011.783,11808.166)
	(8029.370,11739.014)
	(8038.164,11700.242)
	(8046.976,11658.596)
	(8055.821,11614.016)
	(8064.712,11566.445)
	(8073.663,11515.824)
	(8082.687,11462.094)
	(8091.799,11405.198)
	(8101.012,11345.077)
	(8110.339,11281.672)
	(8119.796,11214.925)
	(8129.394,11144.779)
	(8139.148,11071.173)
	(8144.088,11033.056)
	(8149.072,10994.051)
	(8154.102,10954.153)
	(8159.180,10913.354)
	(8164.307,10871.646)
	(8169.484,10829.023)
	(8174.715,10785.477)
	(8180.000,10741.000)

\path(7325,10030)	(7332.122,9959.550)
	(7339.117,9892.052)
	(7345.997,9827.461)
	(7352.770,9765.728)
	(7359.446,9706.807)
	(7366.036,9650.652)
	(7372.549,9597.215)
	(7378.994,9546.450)
	(7385.382,9498.310)
	(7391.723,9452.749)
	(7398.025,9409.719)
	(7404.300,9369.173)
	(7416.805,9295.349)
	(7429.315,9230.903)
	(7441.909,9175.459)
	(7454.665,9128.644)
	(7480.978,9059.403)
	(7539.000,9008.000)

\path(7539,9008)	(7580.498,9059.828)
	(7596.402,9118.982)
	(7603.372,9155.340)
	(7609.775,9195.635)
	(7615.677,9239.431)
	(7621.143,9286.291)
	(7626.239,9335.778)
	(7631.031,9387.455)
	(7635.584,9440.886)
	(7639.964,9495.632)
	(7644.236,9551.257)
	(7648.466,9607.325)
	(7652.720,9663.398)
	(7657.063,9719.039)
	(7661.562,9773.812)
	(7666.281,9827.279)
	(7671.286,9879.003)
	(7676.642,9928.549)
	(7682.417,9975.477)
	(7688.674,10019.353)
	(7695.480,10059.737)
	(7702.901,10096.195)
	(7719.847,10155.581)
	(7764.000,10208.000)

\path(7764,10208)	(7804.301,10193.276)
	(7842.043,10144.694)
	(7860.230,10106.539)
	(7878.107,10058.518)
	(7895.783,10000.166)
	(7913.370,9931.014)
	(7922.164,9892.242)
	(7930.976,9850.596)
	(7939.821,9806.016)
	(7948.712,9758.445)
	(7957.663,9707.824)
	(7966.687,9654.094)
	(7975.799,9597.198)
	(7985.012,9537.077)
	(7994.339,9473.672)
	(8003.796,9406.925)
	(8013.394,9336.779)
	(8023.148,9263.173)
	(8028.088,9225.056)
	(8033.072,9186.051)
	(8038.102,9146.153)
	(8043.180,9105.354)
	(8048.307,9063.646)
	(8053.484,9021.023)
	(8058.715,8977.477)
	(8064.000,8933.000)

\path(5945,8250)	(5951.634,8178.697)
	(5958.159,8110.381)
	(5964.585,8045.004)
	(5970.921,7982.521)
	(5977.177,7922.882)
	(5983.361,7866.042)
	(5989.485,7811.952)
	(5995.556,7760.565)
	(6001.584,7711.833)
	(6007.579,7665.710)
	(6013.550,7622.148)
	(6019.506,7581.099)
	(6031.412,7506.352)
	(6043.371,7441.091)
	(6055.460,7384.937)
	(6067.753,7337.510)
	(6093.251,7267.322)
	(6150.000,7215.000)

\path(6150,7215)	(6190.470,7266.619)
	(6206.072,7325.711)
	(6212.938,7362.047)
	(6219.265,7402.327)
	(6225.117,7446.113)
	(6230.557,7492.968)
	(6235.649,7542.454)
	(6240.458,7594.135)
	(6245.046,7647.572)
	(6249.477,7702.329)
	(6253.815,7757.968)
	(6258.124,7814.051)
	(6262.467,7870.142)
	(6266.907,7925.803)
	(6271.509,7980.597)
	(6276.337,8034.086)
	(6281.453,8085.834)
	(6286.921,8135.402)
	(6292.806,8182.353)
	(6299.171,8226.250)
	(6306.079,8266.656)
	(6313.594,8303.133)
	(6330.700,8362.551)
	(6375.000,8415.000)

\path(6375,8415)	(6415.301,8400.276)
	(6453.043,8351.694)
	(6471.230,8313.539)
	(6489.107,8265.518)
	(6506.783,8207.166)
	(6524.370,8138.014)
	(6533.164,8099.242)
	(6541.976,8057.596)
	(6550.821,8013.016)
	(6559.712,7965.445)
	(6568.663,7914.824)
	(6577.687,7861.094)
	(6586.799,7804.198)
	(6596.012,7744.077)
	(6605.339,7680.672)
	(6614.796,7613.925)
	(6624.394,7543.779)
	(6634.148,7470.173)
	(6639.088,7432.056)
	(6644.072,7393.051)
	(6649.102,7353.153)
	(6654.180,7312.354)
	(6659.307,7270.646)
	(6664.484,7228.023)
	(6669.715,7184.477)
	(6675.000,7140.000)

\path(8575,8270)	(8582.261,8204.359)
	(8589.393,8141.469)
	(8596.406,8081.287)
	(8603.311,8023.770)
	(8610.116,7968.872)
	(8616.833,7916.552)
	(8623.472,7866.766)
	(8630.041,7819.469)
	(8636.551,7774.618)
	(8643.013,7732.170)
	(8649.436,7692.082)
	(8655.829,7654.309)
	(8668.570,7585.535)
	(8681.315,7525.500)
	(8694.144,7473.856)
	(8707.136,7430.254)
	(8733.932,7365.781)
	(8793.000,7318.000)

\path(8793,7318)	(8834.780,7365.625)
	(8850.969,7419.915)
	(8864.686,7490.246)
	(8870.775,7530.425)
	(8876.432,7573.412)
	(8881.722,7618.806)
	(8886.705,7666.205)
	(8891.444,7715.209)
	(8896.001,7765.418)
	(8900.440,7816.429)
	(8904.821,7867.842)
	(8909.208,7919.257)
	(8913.662,7970.272)
	(8918.246,8020.487)
	(8923.023,8069.500)
	(8928.054,8116.910)
	(8933.401,8162.317)
	(8939.128,8205.320)
	(8945.295,8245.517)
	(8959.204,8315.893)
	(8975.626,8370.237)
	(9018.000,8418.000)

\path(9018,8418)	(9059.731,8402.994)
	(9098.496,8354.210)
	(9117.040,8315.984)
	(9135.172,8267.914)
	(9153.002,8209.531)
	(9170.640,8140.369)
	(9179.421,8101.600)
	(9188.195,8059.961)
	(9196.976,8015.394)
	(9205.778,7967.840)
	(9214.614,7917.242)
	(9223.498,7863.540)
	(9232.443,7806.676)
	(9241.464,7746.593)
	(9250.575,7683.231)
	(9259.788,7616.532)
	(9269.118,7546.439)
	(9278.579,7472.891)
	(9283.362,7434.804)
	(9288.184,7395.832)
	(9293.044,7355.967)
	(9297.946,7315.203)
	(9302.891,7273.531)
	(9307.880,7230.945)
	(9312.916,7187.437)
	(9318.000,7143.000)

\path(9485,6460)	(9493.003,6394.363)
	(9500.853,6331.479)
	(9508.560,6271.304)
	(9516.135,6213.795)
	(9523.589,6158.907)
	(9530.933,6106.599)
	(9538.176,6056.825)
	(9545.331,6009.542)
	(9552.407,5964.707)
	(9559.415,5922.277)
	(9566.365,5882.208)
	(9573.269,5844.455)
	(9586.981,5775.729)
	(9600.634,5715.749)
	(9614.314,5664.167)
	(9628.106,5620.636)
	(9656.370,5556.331)
	(9718.000,5509.000)

\path(9718,5509)	(9761.090,5556.935)
	(9777.662,5611.316)
	(9791.610,5681.702)
	(9797.762,5721.896)
	(9803.453,5764.891)
	(9808.747,5810.285)
	(9813.708,5857.680)
	(9818.403,5906.673)
	(9822.895,5956.866)
	(9827.250,6007.858)
	(9831.531,6059.247)
	(9835.805,6110.635)
	(9840.135,6161.621)
	(9844.587,6211.804)
	(9849.225,6260.784)
	(9854.114,6308.161)
	(9859.320,6353.535)
	(9864.905,6396.504)
	(9870.937,6436.669)
	(9884.595,6506.986)
	(9900.813,6561.282)
	(9943.000,6609.000)

\path(9943,6609)	(9984.731,6593.994)
	(10023.496,6545.210)
	(10042.040,6506.984)
	(10060.172,6458.914)
	(10078.002,6400.531)
	(10095.640,6331.369)
	(10104.421,6292.600)
	(10113.195,6250.961)
	(10121.976,6206.394)
	(10130.778,6158.840)
	(10139.614,6108.242)
	(10148.498,6054.540)
	(10157.443,5997.676)
	(10166.464,5937.593)
	(10175.575,5874.231)
	(10184.788,5807.532)
	(10194.118,5737.439)
	(10203.579,5663.891)
	(10208.362,5625.804)
	(10213.184,5586.832)
	(10218.044,5546.967)
	(10222.946,5506.203)
	(10227.891,5464.531)
	(10232.880,5421.945)
	(10237.916,5378.437)
	(10243.000,5334.000)

\path(7035,4630)	(7041.604,4566.873)
	(7048.097,4506.392)
	(7054.488,4448.513)
	(7060.786,4393.196)
	(7067.000,4340.399)
	(7073.140,4290.079)
	(7079.215,4242.194)
	(7085.235,4196.704)
	(7091.207,4153.565)
	(7097.143,4112.735)
	(7108.939,4037.839)
	(7120.697,3971.678)
	(7132.490,3913.917)
	(7144.393,3864.222)
	(7156.480,3822.257)
	(7181.500,3760.172)
	(7237.000,3714.000)

\path(7237,3714)	(7278.086,3761.467)
	(7294.070,3815.710)
	(7307.665,3886.014)
	(7313.720,3926.186)
	(7319.360,3969.168)
	(7324.647,4014.561)
	(7329.641,4061.963)
	(7334.404,4110.973)
	(7338.996,4161.189)
	(7343.479,4212.211)
	(7347.914,4263.636)
	(7352.360,4315.065)
	(7356.881,4366.095)
	(7361.535,4416.325)
	(7366.385,4465.355)
	(7371.491,4512.782)
	(7376.914,4558.206)
	(7382.715,4601.226)
	(7388.955,4641.440)
	(7402.997,4711.846)
	(7419.527,4766.215)
	(7462.000,4814.000)

\path(7462,4814)	(7503.731,4798.994)
	(7542.496,4750.210)
	(7561.040,4711.984)
	(7579.172,4663.914)
	(7597.002,4605.531)
	(7614.640,4536.369)
	(7623.421,4497.600)
	(7632.195,4455.961)
	(7640.976,4411.394)
	(7649.778,4363.840)
	(7658.614,4313.242)
	(7667.498,4259.540)
	(7676.443,4202.676)
	(7685.464,4142.593)
	(7694.575,4079.231)
	(7703.788,4012.532)
	(7713.118,3942.439)
	(7722.579,3868.891)
	(7727.362,3830.804)
	(7732.184,3791.832)
	(7737.044,3751.967)
	(7741.946,3711.203)
	(7746.891,3669.531)
	(7751.880,3626.945)
	(7756.916,3583.437)
	(7762.000,3539.000)

\path(9435,2860)	(9442.234,2794.771)
	(9449.339,2732.276)
	(9456.325,2672.473)
	(9463.203,2615.316)
	(9469.982,2560.764)
	(9476.673,2508.772)
	(9483.284,2459.299)
	(9489.827,2412.299)
	(9496.311,2367.730)
	(9502.745,2325.550)
	(9509.141,2285.713)
	(9515.507,2248.178)
	(9528.192,2179.836)
	(9540.879,2120.180)
	(9553.648,2068.862)
	(9566.578,2025.535)
	(9593.241,1961.471)
	(9652.000,1914.000)

\path(9652,1914)	(9693.815,1961.631)
	(9710.013,2015.922)
	(9723.737,2086.255)
	(9729.827,2126.434)
	(9735.486,2169.421)
	(9740.775,2214.815)
	(9745.757,2262.214)
	(9750.495,2311.218)
	(9755.051,2361.426)
	(9759.487,2412.437)
	(9763.866,2463.850)
	(9768.250,2515.264)
	(9772.701,2566.279)
	(9777.282,2616.493)
	(9782.054,2665.505)
	(9787.081,2712.915)
	(9792.425,2758.321)
	(9798.148,2801.323)
	(9804.312,2841.520)
	(9818.215,2911.895)
	(9834.631,2966.238)
	(9877.000,3014.000)

\path(9877,3014)	(9918.731,2998.994)
	(9957.496,2950.210)
	(9976.040,2911.984)
	(9994.172,2863.914)
	(10012.002,2805.531)
	(10029.640,2736.369)
	(10038.421,2697.600)
	(10047.195,2655.961)
	(10055.976,2611.394)
	(10064.778,2563.840)
	(10073.614,2513.242)
	(10082.498,2459.540)
	(10091.443,2402.676)
	(10100.464,2342.593)
	(10109.575,2279.231)
	(10118.788,2212.532)
	(10128.118,2142.439)
	(10137.579,2068.891)
	(10142.362,2030.804)
	(10147.184,1991.832)
	(10152.044,1951.967)
	(10156.946,1911.203)
	(10161.891,1869.531)
	(10166.880,1826.945)
	(10171.916,1783.437)
	(10177.000,1739.000)

\path(7035,6470)	(7042.186,6403.887)
	(7049.246,6340.545)
	(7056.190,6279.930)
	(7063.029,6221.998)
	(7069.773,6166.705)
	(7076.431,6114.007)
	(7083.013,6063.861)
	(7089.529,6016.222)
	(7095.989,5971.047)
	(7102.404,5928.291)
	(7108.782,5887.912)
	(7115.134,5849.864)
	(7127.800,5780.588)
	(7140.481,5720.114)
	(7153.257,5668.089)
	(7166.206,5624.162)
	(7192.943,5559.199)
	(7252.000,5511.000)

\path(7252,5511)	(7293.559,5558.573)
	(7309.683,5612.848)
	(7323.362,5683.170)
	(7329.439,5723.347)
	(7335.091,5766.332)
	(7340.380,5811.726)
	(7345.366,5859.126)
	(7350.113,5908.132)
	(7354.682,5958.343)
	(7359.134,6009.358)
	(7363.533,6060.775)
	(7367.938,6112.194)
	(7372.414,6163.214)
	(7377.020,6213.434)
	(7381.820,6262.452)
	(7386.874,6309.868)
	(7392.246,6355.281)
	(7397.996,6398.289)
	(7404.187,6438.492)
	(7418.138,6508.878)
	(7434.594,6563.230)
	(7477.000,6611.000)

\path(7477,6611)	(7518.731,6595.994)
	(7557.496,6547.210)
	(7576.040,6508.984)
	(7594.172,6460.914)
	(7612.002,6402.531)
	(7629.640,6333.369)
	(7638.421,6294.600)
	(7647.195,6252.961)
	(7655.976,6208.394)
	(7664.778,6160.840)
	(7673.614,6110.242)
	(7682.498,6056.540)
	(7691.443,5999.676)
	(7700.464,5939.593)
	(7709.575,5876.231)
	(7718.788,5809.532)
	(7728.118,5739.439)
	(7737.579,5665.891)
	(7742.362,5627.804)
	(7747.184,5588.832)
	(7752.044,5548.967)
	(7756.946,5508.203)
	(7761.891,5466.531)
	(7766.880,5423.945)
	(7771.916,5380.437)
	(7777.000,5336.000)

\path(9435,4690)	(9442.124,4622.731)
	(9449.126,4558.281)
	(9456.017,4496.605)
	(9462.806,4437.659)
	(9469.504,4381.398)
	(9476.120,4327.777)
	(9482.664,4276.751)
	(9489.146,4228.277)
	(9495.577,4182.309)
	(9501.965,4138.802)
	(9508.322,4097.712)
	(9514.656,4058.995)
	(9527.298,3988.497)
	(9539.971,3926.953)
	(9552.754,3874.003)
	(9565.727,3829.291)
	(9592.560,3763.152)
	(9652.000,3714.000)

\path(9652,3714)	(9693.235,3761.501)
	(9709.263,3815.755)
	(9722.885,3886.065)
	(9728.947,3926.238)
	(9734.590,3969.222)
	(9739.878,4014.615)
	(9744.870,4062.016)
	(9749.627,4111.025)
	(9754.212,4161.239)
	(9758.686,4212.258)
	(9763.109,4263.681)
	(9767.543,4315.107)
	(9772.049,4366.133)
	(9776.688,4416.360)
	(9781.522,4465.386)
	(9786.612,4512.810)
	(9792.018,4558.231)
	(9797.804,4601.247)
	(9804.028,4641.457)
	(9818.041,4711.856)
	(9834.548,4766.220)
	(9877.000,4814.000)

\path(9877,4814)	(9918.731,4798.994)
	(9957.496,4750.210)
	(9976.040,4711.984)
	(9994.172,4663.914)
	(10012.002,4605.531)
	(10029.640,4536.369)
	(10038.421,4497.600)
	(10047.195,4455.961)
	(10055.976,4411.394)
	(10064.778,4363.840)
	(10073.614,4313.242)
	(10082.498,4259.540)
	(10091.443,4202.676)
	(10100.464,4142.593)
	(10109.575,4079.231)
	(10118.788,4012.532)
	(10128.118,3942.439)
	(10137.579,3868.891)
	(10142.362,3830.804)
	(10147.184,3791.832)
	(10152.044,3751.967)
	(10156.946,3711.203)
	(10161.891,3669.531)
	(10166.880,3626.945)
	(10171.916,3583.437)
	(10177.000,3539.000)

\spline(4245,6210)
(4455,5835)(4800,6210)(4950,5960)
\spline(12045,4410)
(12255,4035)(12600,4410)(12720,4205)
\spline(12045,2610)
(12255,2235)(12600,2610)(12780,2390)
\spline(7545,810)
(7755,435)(8100,810)(8280,590)
\spline(6945,2610)
(7155,2235)(7500,2610)(7695,2360)
\spline(1530,4410)
(1740,4035)(2085,4410)(2235,4175)
\put(300,11340){\makebox(0,0)[lb]{\smash{{{\SetFigFont{12}{14.4}{rm}12}}}}}
\put(300,9465){\makebox(0,0)[lb]{\smash{{{\SetFigFont{12}{14.4}{rm}11}}}}}
\put(300,7740){\makebox(0,0)[lb]{\smash{{{\SetFigFont{12}{14.4}{rm}10}}}}}
\put(300,5940){\makebox(0,0)[lb]{\smash{{{\SetFigFont{12}{14.4}{rm}9}}}}}
\put(300,4065){\makebox(0,0)[lb]{\smash{{{\SetFigFont{12}{14.4}{rm}8}}}}}
\put(300,2340){\makebox(0,0)[lb]{\smash{{{\SetFigFont{12}{14.4}{rm}7}}}}}
\put(0,12315){\makebox(0,0)[lb]{\smash{{{\SetFigFont{12}{14.4}{rm}\# general}}}}}
\put(225,11940){\makebox(0,0)[lb]{\smash{{{\SetFigFont{12}{14.4}{rm}lines}}}}}
\put(6600,10815){\makebox(0,0)[lb]{\smash{{{\SetFigFont{8}{9.6}{rm}1500}}}}}
\put(6600,9015){\makebox(0,0)[lb]{\smash{{{\SetFigFont{8}{9.6}{rm}1500}}}}}
\put(5175,7215){\makebox(0,0)[lb]{\smash{{{\SetFigFont{8}{9.6}{rm}1350}}}}}
\put(8100,7215){\makebox(0,0)[lb]{\smash{{{\SetFigFont{8}{9.6}{rm}150}}}}}
\put(3825,5415){\makebox(0,0)[lb]{\smash{{{\SetFigFont{8}{9.6}{rm}1}}}}}
\put(8550,6690){\makebox(0,0)[lb]{\smash{{{\SetFigFont{8}{9.6}{rm}$\times 2$}}}}}
\put(8925,5415){\makebox(0,0)[lb]{\smash{{{\SetFigFont{8}{9.6}{rm}150}}}}}
\put(5175,6765){\makebox(0,0)[lb]{\smash{{{\SetFigFont{8}{9.6}{rm}$\times 27$}}}}}
\put(300,540){\makebox(0,0)[lb]{\smash{{{\SetFigFont{12}{14.4}{rm}6}}}}}
\put(6375,5415){\makebox(0,0)[lb]{\smash{{{\SetFigFont{8}{9.6}{rm}1023}}}}}
\put(1500,3615){\makebox(0,0)[lb]{\smash{{{\SetFigFont{8}{9.6}{rm}18}}}}}
\put(3825,3615){\makebox(0,0)[lb]{\smash{{{\SetFigFont{8}{9.6}{rm}588}}}}}
\put(6525,3615){\makebox(0,0)[lb]{\smash{{{\SetFigFont{8}{9.6}{rm}127}}}}}
\put(12300,3615){\makebox(0,0)[lb]{\smash{{{\SetFigFont{8}{9.6}{rm}1}}}}}
\put(12300,1815){\makebox(0,0)[lb]{\smash{{{\SetFigFont{8}{9.6}{rm}1}}}}}
\put(8925,1815){\makebox(0,0)[lb]{\smash{{{\SetFigFont{8}{9.6}{rm}11}}}}}
\put(3825,1815){\makebox(0,0)[lb]{\smash{{{\SetFigFont{8}{9.6}{rm}90}}}}}
\put(7200,1815){\makebox(0,0)[lb]{\smash{{{\SetFigFont{8}{9.6}{rm}5}}}}}
\put(7800,15){\makebox(0,0)[lb]{\smash{{{\SetFigFont{8}{9.6}{rm}1}}}}}
\put(10500,15){\makebox(0,0)[lb]{\smash{{{\SetFigFont{8}{9.6}{rm}8}}}}}
\put(8100,2415){\makebox(0,0)[lb]{\smash{{{\SetFigFont{12}{14.4}{rm}*}}}}}
\put(8700,615){\makebox(0,0)[lb]{\smash{{{\SetFigFont{12}{14.4}{rm}*}}}}}
\put(2700,4215){\makebox(0,0)[lb]{\smash{{{\SetFigFont{12}{14.4}{rm}*}}}}}
\put(2850,4965){\makebox(0,0)[lb]{\smash{{{\SetFigFont{8}{9.6}{rm}$\times 3$}}}}}
\put(7200,5115){\makebox(0,0)[lb]{\smash{{{\SetFigFont{8}{9.6}{rm}$\times 3$}}}}}
\put(11550,5115){\makebox(0,0)[lb]{\smash{{{\SetFigFont{8}{9.6}{rm}$\times 9$}}}}}
\put(7125,3090){\makebox(0,0)[lb]{\smash{{{\SetFigFont{8}{9.6}{rm}$\times 3$}}}}}
\put(8925,3090){\makebox(0,0)[lb]{\smash{{{\SetFigFont{8}{9.6}{rm}$\times 2$}}}}}
\put(11550,3165){\makebox(0,0)[lb]{\smash{{{\SetFigFont{8}{9.6}{rm}$\times 3$}}}}}
\put(8925,1215){\makebox(0,0)[lb]{\smash{{{\SetFigFont{8}{9.6}{rm}$\times 3$}}}}}
\put(8925,3615){\makebox(0,0)[lb]{\smash{{{\SetFigFont{8}{9.6}{rm}14}}}}}
\end{picture}

%% file: ewxy.tex
\begingroup\makeatletter\ifx\SetFigFont\undefined
\def\x#1#2#3#4#5#6#7\relax{\def\x{#1#2#3#4#5#6}}%
\expandafter\x\fmtname xxxxxx\relax \def\y{splain}%
\ifx\x\y   
\gdef\SetFigFont#1#2#3{%
  \ifnum #1<17\tiny\else \ifnum #1<20\small\else
  \ifnum #1<24\normalsize\else \ifnum #1<29\large\else
  \ifnum #1<34\Large\else \ifnum #1<41\LARGE\else
     \huge\fi\fi\fi\fi\fi\fi
  \csname #3\endcsname}%
\else
\gdef\SetFigFont#1#2#3{\begingroup
  \count@#1\relax \ifnum 25<\count@\count@25\fi
  \def\x{\endgroup\@setsize\SetFigFont{#2pt}}%
  \expandafter\x
    \csname \romannumeral\the\count@ pt\expandafter\endcsname
    \csname @\romannumeral\the\count@ pt\endcsname
  \csname #3\endcsname}%
\fi
\fi\endgroup
\begin{picture}(9475,2696)(0,-10)
\thicklines
\put(105,2598){\ellipse{150}{150}}
\put(2880,2598){\ellipse{150}{150}}
\put(1530,2598){\blacken\ellipse{74}{74}}
\put(1530,2598){\ellipse{74}{74}}
\put(3937,2611){\blacken\ellipse{74}{74}}
\put(3937,2611){\ellipse{74}{74}}
\put(5310,2611){\blacken\ellipse{74}{74}}
\put(5310,2611){\ellipse{74}{74}}
\put(9430,2611){\blacken\ellipse{74}{74}}
\put(9430,2611){\ellipse{74}{74}}
\put(8881,2611){\blacken\ellipse{74}{74}}
\put(8881,2611){\ellipse{74}{74}}
\put(5860,2611){\blacken\ellipse{74}{74}}
\put(5860,2611){\ellipse{74}{74}}
\put(6409,2611){\blacken\ellipse{74}{74}}
\put(6409,2611){\ellipse{74}{74}}
\put(6958,2611){\blacken\ellipse{74}{74}}
\put(6958,2611){\ellipse{74}{74}}
\put(7782,2611){\blacken\ellipse{74}{74}}
\put(7782,2611){\ellipse{74}{74}}
\put(8331,2611){\blacken\ellipse{74}{74}}
\put(8331,2611){\ellipse{74}{74}}
\put(3630,1038){\blacken\ellipse{74}{74}}
\put(3630,1038){\ellipse{74}{74}}
\put(8580,963){\ellipse{150}{150}}
\put(4455,2598){\blacken\ellipse{74}{74}}
\put(4455,2598){\ellipse{74}{74}}
\put(6105,1038){\blacken\ellipse{74}{74}}
\put(6105,1038){\ellipse{74}{74}}
\put(3405,2613){\blacken\ellipse{74}{74}}
\put(3405,2613){\ellipse{74}{74}}
\path(6958,2611)(6105,1098)
\path(2880,2523)(3630,1098)(3405,2598)
\path(3930,2598)(3630,1098)(4455,2598)
\path(5860,2611)(6105,1098)(6409,2611)
\path(7782,2611)(8580,1098)(8331,2611)
\path(8881,2611)(8580,1098)(9430,2611)
\spline(5310,2611)
(5585,1787)(6105,1098)
\spline(5310,2611)
(5722,2062)(6105,1098)
\put(0,63){\makebox(0,0)[lb]{\smash{{{\SetFigFont{12}{14.4}{rm}$W$}}}}}
\put(1350,63){\makebox(0,0)[lb]{\smash{{{\SetFigFont{12}{14.4}{rm}$X$}}}}}
\put(3540,33){\makebox(0,0)[lb]{\smash{{{\SetFigFont{12}{14.4}{rm}$Y^a$}}}}}
\put(5880,63){\makebox(0,0)[lb]{\smash{{{\SetFigFont{12}{14.4}{rm}$Y^b$}}}}}
\put(8430,63){\makebox(0,0)[lb]{\smash{{{\SetFigFont{12}{14.4}{rm}$Y^c$}}}}}
\put(4110,933){\makebox(0,0)[lb]{\smash{{{\SetFigFont{12}{14.4}{rm}$H$}}}}}
\put(6540,963){\makebox(0,0)[lb]{\smash{{{\SetFigFont{12}{14.4}{rm}$H$}}}}}
\put(8985,978){\makebox(0,0)[lb]{\smash{{{\SetFigFont{12}{14.4}{rm}$H$}}}}}
\end{picture}

%% file: eYb.tex
\begingroup\makeatletter\ifx\SetFigFont\undefined
\def\x#1#2#3#4#5#6#7\relax{\def\x{#1#2#3#4#5#6}}%
\expandafter\x\fmtname xxxxxx\relax \def\y{splain}%
\ifx\x\y   
\gdef\SetFigFont#1#2#3{%
  \ifnum #1<17\tiny\else \ifnum #1<20\small\else
  \ifnum #1<24\normalsize\else \ifnum #1<29\large\else
  \ifnum #1<34\Large\else \ifnum #1<41\LARGE\else
     \huge\fi\fi\fi\fi\fi\fi
  \csname #3\endcsname}%
\else
\gdef\SetFigFont#1#2#3{\begingroup
  \count@#1\relax \ifnum 25<\count@\count@25\fi
  \def\x{\endgroup\@setsize\SetFigFont{#2pt}}%
  \expandafter\x
    \csname \romannumeral\the\count@ pt\expandafter\endcsname
    \csname @\romannumeral\the\count@ pt\endcsname
  \csname #3\endcsname}%
\fi
\fi\endgroup
\begin{picture}(20262,7529)(0,-10)
\thicklines
\put(5700,6458){\ellipse{1200}{2100}}
\path(4950,6908)(4350,6008)(6450,6008)
	(7050,6908)(4950,6908)
\put(4200,6383){\makebox(0,0)[lb]{\smash{{{\SetFigFont{12}{14.4}{rm}$\#$}}}}}
\path(8850,6908)(8250,6008)(10350,6008)
	(10950,6908)(8850,6908)
\put(8100,6383){\makebox(0,0)[lb]{\smash{{{\SetFigFont{12}{14.4}{rm}$\#$}}}}}
\path(750,6608)(1050,6308)
\path(750,6608)(1050,6308)
\path(1050,6608)(750,6308)
\path(1050,6608)(750,6308)
\path(4950,6608)(5250,6308)
\path(4950,6608)(5250,6308)
\path(5250,6608)(4950,6308)
\path(5250,6608)(4950,6308)
\put(1500,6458){\ellipse{1200}{2100}}
\path(750,6908)(150,6008)(2250,6008)
	(2850,6908)(750,6908)
\put(0,6383){\makebox(0,0)[lb]{\smash{{{\SetFigFont{12}{14.4}{rm}$\#$}}}}}
\path(900,3908)(1200,3608)
\path(900,3908)(1200,3608)
\path(1200,3908)(900,3608)
\path(1200,3908)(900,3608)
\path(900,1208)(1200,908)
\path(900,1208)(1200,908)
\path(1200,1208)(900,908)
\path(1200,1208)(900,908)
\put(6450,1058){\ellipse{1200}{2100}}
\path(5700,1508)(5100,608)(7200,608)
	(7800,1508)(5700,1508)
\put(4950,983){\makebox(0,0)[lb]{\smash{{{\SetFigFont{12}{14.4}{rm}$\#$}}}}}
\path(14400,1508)(13800,608)(15900,608)
	(16500,1508)(14400,1508)
\put(13650,983){\makebox(0,0)[lb]{\smash{{{\SetFigFont{12}{14.4}{rm}$\#$}}}}}
\path(5700,1208)(6000,908)
\path(5700,1208)(6000,908)
\path(6000,1208)(5700,908)
\path(6000,1208)(5700,908)
\path(10350,1208)(10650,908)
\path(10350,1208)(10650,908)
\path(10650,1208)(10350,908)
\path(10650,1208)(10350,908)
\put(10528,1584){\ellipse{1124}{1050}}
\path(9600,1508)(9000,608)(11100,608)
	(11700,1508)(9600,1508)
\put(8850,983){\makebox(0,0)[lb]{\smash{{{\SetFigFont{12}{14.4}{rm}$\#$}}}}}
\put(6187.500,1058.000){\arc{3975.000}{2.5850}{3.6982}}
\put(10312.500,1058.000){\arc{3975.000}{5.7266}{6.8398}}
\put(14625,1058){\blacken\ellipse{150}{150}}
\put(14625,1058){\ellipse{150}{150}}
\put(15525,1058){\blacken\ellipse{150}{150}}
\put(15525,1058){\ellipse{150}{150}}
\path(14250,1058)(15900,1058)
\path(14250,1058)(15900,1058)
\dottedline{225}(10500,1808)(11400,1808)
\dottedline{225}(10500,758)(10500,1358)
\dottedline{225}(6413,735)(7448,1223)
\dottedline{225}(5550,1058)(6300,1058)
\dottedline{225}(6450,1808)(7200,1808)
\put(8100,983){\makebox(0,0)[lb]{\smash{{{\SetFigFont{12}{14.4}{rm}$-$}}}}}
\put(7380,1733){\makebox(0,0)[lb]{\smash{{{\SetFigFont{8}{9.6}{rm}6}}}}}
\put(11595,1748){\makebox(0,0)[lb]{\smash{{{\SetFigFont{8}{9.6}{rm}6}}}}}
\put(12675,983){\makebox(0,0)[lb]{\smash{{{\SetFigFont{12}{14.4}{rm}$\times$}}}}}
\put(16050,308){\makebox(0,0)[lb]{\smash{{{\SetFigFont{12}{14.4}{rm}1}}}}}
\put(11250,308){\makebox(0,0)[lb]{\smash{{{\SetFigFont{12}{14.4}{rm}58}}}}}
\put(7350,308){\makebox(0,0)[lb]{\smash{{{\SetFigFont{12}{14.4}{rm}74}}}}}
\put(10350,3758){\ellipse{1200}{2100}}
\path(9600,4208)(9000,3308)(11100,3308)
	(11700,4208)(9600,4208)
\put(8850,3683){\makebox(0,0)[lb]{\smash{{{\SetFigFont{12}{14.4}{rm}$\#$}}}}}
\put(6450,3758){\ellipse{1200}{2100}}
\path(5700,4208)(5100,3308)(7200,3308)
	(7800,4208)(5700,4208)
\put(4950,3683){\makebox(0,0)[lb]{\smash{{{\SetFigFont{12}{14.4}{rm}$\#$}}}}}
\path(5700,3908)(6000,3608)
\path(5700,3908)(6000,3608)
\path(6000,3908)(5700,3608)
\path(6000,3908)(5700,3608)
\path(9600,3908)(9900,3608)
\path(9600,3908)(9900,3608)
\path(9900,3908)(9600,3608)
\path(9900,3908)(9600,3608)
\path(14250,3908)(14550,3608)
\path(14250,3908)(14550,3608)
\path(14550,3908)(14250,3608)
\path(14550,3908)(14250,3608)
\put(14356,4292){\ellipse{1124}{1050}}
\path(13425,4208)(12825,3308)(14925,3308)
	(15525,4208)(13425,4208)
\put(12675,3683){\makebox(0,0)[lb]{\smash{{{\SetFigFont{12}{14.4}{rm}$\#$}}}}}
\path(18150,4208)(17550,3308)(19650,3308)
	(20250,4208)(18150,4208)
\put(17400,3683){\makebox(0,0)[lb]{\smash{{{\SetFigFont{12}{14.4}{rm}$\#$}}}}}
\put(14287.500,3758.000){\arc{3975.000}{5.7266}{6.8398}}
\put(18300,3758){\blacken\ellipse{150}{150}}
\put(18300,3758){\ellipse{150}{150}}
\put(19275,3758){\blacken\ellipse{150}{150}}
\put(19275,3758){\ellipse{150}{150}}
\path(18000,3758)(19650,3758)
\path(18000,3758)(19650,3758)
\put(19800,3008){\makebox(0,0)[lb]{\smash{{{\SetFigFont{12}{14.4}{rm}1}}}}}
\put(16575,3683){\makebox(0,0)[lb]{\smash{{{\SetFigFont{12}{14.4}{rm}$\times$}}}}}
\put(6187.500,3758.000){\arc{3975.000}{2.5850}{3.6982}}
\put(9000,6458){\blacken\ellipse{150}{150}}
\put(9000,6458){\ellipse{150}{150}}
\put(9975,6458){\blacken\ellipse{150}{150}}
\put(9975,6458){\ellipse{150}{150}}
\put(1576,3769){\ellipse{1050}{2100}}
\put(1500,1058){\ellipse{900}{2100}}
\path(600,6458)(2250,6458)
\path(600,6458)(2250,6458)
\path(8700,6458)(10350,6458)
\path(8700,6458)(10350,6458)
\path(600,3758)(2250,3758)
\path(600,3758)(2250,3758)
\path(600,1058)(2250,1058)
\path(600,1058)(2250,1058)
\dottedline{225}(450,383)(1020,1965)
\dottedline{225}(1500,1793)(2415,1823)
\dottedline{225}(2070,1650)(2145,248)
\dottedline{225}(495,3173)(1005,4883)
\dottedline{225}(1500,4508)(2325,4508)
\dottedline{225}(1500,7208)(2250,7208)
\dottedline{225}(5700,7208)(6450,7208)
\path(750,4208)(150,3308)(2250,3308)
	(2850,4208)(750,4208)
\path(750,1508)(150,608)(2250,608)
	(2850,1508)(750,1508)
\dottedline{225}(6450,4508)(7200,4508)
\dottedline{225}(5550,3758)(6300,3758)
\dottedline{225}(10350,4508)(11100,4508)
\dottedline{225}(10500,3758)(11100,3758)
\dottedline{225}(14400,4508)(15150,4508)
\put(3450,6383){\makebox(0,0)[lb]{\smash{{{\SetFigFont{12}{14.4}{rm}$=$}}}}}
\put(3450,3683){\makebox(0,0)[lb]{\smash{{{\SetFigFont{12}{14.4}{rm}$=$}}}}}
\put(3450,983){\makebox(0,0)[lb]{\smash{{{\SetFigFont{12}{14.4}{rm}$=$}}}}}
\put(7350,6383){\makebox(0,0)[lb]{\smash{{{\SetFigFont{12}{14.4}{rm}$\times$}}}}}
\put(6630,7118){\makebox(0,0)[lb]{\smash{{{\SetFigFont{8}{9.6}{rm}8}}}}}
\put(2400,7103){\makebox(0,0)[lb]{\smash{{{\SetFigFont{8}{9.6}{rm}8}}}}}
\put(2415,4433){\makebox(0,0)[lb]{\smash{{{\SetFigFont{8}{9.6}{rm}7}}}}}
\put(2625,1718){\makebox(0,0)[lb]{\smash{{{\SetFigFont{8}{9.6}{rm}6}}}}}
\put(0,3683){\makebox(0,0)[lb]{\smash{{{\SetFigFont{12}{14.4}{rm}$\#$}}}}}
\put(0,983){\makebox(0,0)[lb]{\smash{{{\SetFigFont{12}{14.4}{rm}$\#$}}}}}
\put(2400,5708){\makebox(0,0)[lb]{\smash{{{\SetFigFont{12}{14.4}{rm}184}}}}}
\put(6600,5708){\makebox(0,0)[lb]{\smash{{{\SetFigFont{12}{14.4}{rm}184}}}}}
\put(10500,5708){\makebox(0,0)[lb]{\smash{{{\SetFigFont{12}{14.4}{rm}1}}}}}
\put(2400,3008){\makebox(0,0)[lb]{\smash{{{\SetFigFont{12}{14.4}{rm}68}}}}}
\put(2400,308){\makebox(0,0)[lb]{\smash{{{\SetFigFont{12}{14.4}{rm}16}}}}}
\put(8100,3683){\makebox(0,0)[lb]{\smash{{{\SetFigFont{12}{14.4}{rm}$+$}}}}}
\put(7380,4433){\makebox(0,0)[lb]{\smash{{{\SetFigFont{8}{9.6}{rm}7}}}}}
\put(11250,4448){\makebox(0,0)[lb]{\smash{{{\SetFigFont{8}{9.6}{rm}7}}}}}
\put(15360,4418){\makebox(0,0)[lb]{\smash{{{\SetFigFont{8}{9.6}{rm}7}}}}}
\put(12000,3683){\makebox(0,0)[lb]{\smash{{{\SetFigFont{12}{14.4}{rm}$-$}}}}}
\put(7350,3008){\makebox(0,0)[lb]{\smash{{{\SetFigFont{12}{14.4}{rm}92}}}}}
\put(11250,3008){\makebox(0,0)[lb]{\smash{{{\SetFigFont{12}{14.4}{rm}92}}}}}
\put(15000,3008){\makebox(0,0)[lb]{\smash{{{\SetFigFont{12}{14.4}{rm}116}}}}}
\end{picture}

%% file: qesc2.tex
\begingroup\makeatletter\ifx\SetFigFont\undefined
\def\x#1#2#3#4#5#6#7\relax{\def\x{#1#2#3#4#5#6}}%
\expandafter\x\fmtname xxxxxx\relax \def\y{splain}%
\ifx\x\y   
\gdef\SetFigFont#1#2#3{%
  \ifnum #1<17\tiny\else \ifnum #1<20\small\else
  \ifnum #1<24\normalsize\else \ifnum #1<29\large\else
  \ifnum #1<34\Large\else \ifnum #1<41\LARGE\else
     \huge\fi\fi\fi\fi\fi\fi
  \csname #3\endcsname}%
\else
\gdef\SetFigFont#1#2#3{\begingroup
  \count@#1\relax \ifnum 25<\count@\count@25\fi
  \def\x{\endgroup\@setsize\SetFigFont{#2pt}}%
  \expandafter\x
    \csname \romannumeral\the\count@ pt\expandafter\endcsname
    \csname @\romannumeral\the\count@ pt\endcsname
  \csname #3\endcsname}%
\fi
\fi\endgroup
\begin{picture}(9762,13922)(0,-10)
\thicklines
\path(8250,13218)(7950,12618)(9150,12618)
	(9450,13218)(8250,13218)
\path(8475,12918)(9225,13068)
\path(8325,13443)	(8326.383,13381.937)
	(8327.911,13323.410)
	(8329.591,13267.380)
	(8331.428,13213.804)
	(8333.431,13162.642)
	(8335.605,13113.851)
	(8337.959,13067.392)
	(8340.498,13023.223)
	(8343.229,12981.303)
	(8346.161,12941.589)
	(8352.649,12868.621)
	(8360.017,12803.987)
	(8368.320,12747.360)
	(8377.613,12698.409)
	(8387.951,12656.804)
	(8411.982,12594.317)
	(8475.000,12543.000)

\path(8475,12543)	(8528.589,12581.213)
	(8553.547,12628.770)
	(8577.510,12691.295)
	(8600.649,12765.792)
	(8611.964,12806.593)
	(8623.137,12849.263)
	(8634.190,12893.426)
	(8645.145,12938.709)
	(8656.024,12984.737)
	(8666.846,13031.134)
	(8677.635,13077.526)
	(8688.411,13123.538)
	(8699.196,13168.797)
	(8710.012,13212.925)
	(8720.880,13255.551)
	(8731.821,13296.297)
	(8754.010,13370.655)
	(8776.750,13433.003)
	(8800.214,13480.341)
	(8850.000,13518.000)

\path(8850,13518)	(8891.364,13498.528)
	(8925.387,13448.794)
	(8939.784,13411.479)
	(8952.510,13365.281)
	(8963.622,13309.760)
	(8973.173,13244.475)
	(8977.380,13208.034)
	(8981.218,13168.988)
	(8984.693,13127.281)
	(8987.813,13082.859)
	(8990.584,13035.667)
	(8993.013,12985.649)
	(8995.106,12932.752)
	(8996.872,12876.919)
	(8998.315,12818.096)
	(8999.445,12756.228)
	(9000.267,12691.260)
	(9000.787,12623.138)
	(9001.014,12551.805)
	(9000.954,12477.208)
	(9000.819,12438.668)
	(9000.614,12399.291)
	(9000.341,12359.071)
	(9000.000,12318.000)

\path(7650,11718)(7350,11118)(8550,11118)
	(8850,11718)(7650,11718)
\path(7950,11193)(7950,11643)
\path(7725,11943)	(7737.918,11872.728)
	(7750.604,11808.347)
	(7763.140,11749.666)
	(7775.610,11696.492)
	(7788.094,11648.634)
	(7800.676,11605.898)
	(7826.464,11535.026)
	(7853.631,11482.339)
	(7882.836,11446.297)
	(7950.000,11418.000)

\path(7950,11418)	(8003.119,11443.212)
	(8047.518,11512.365)
	(8067.113,11558.775)
	(8085.326,11610.599)
	(8102.423,11665.979)
	(8118.671,11723.059)
	(8134.336,11779.980)
	(8149.683,11834.886)
	(8164.978,11885.920)
	(8180.489,11931.224)
	(8213.218,11997.214)
	(8250.000,12018.000)

\path(8250,12018)	(8298.586,11995.342)
	(8337.768,11943.332)
	(8353.970,11905.221)
	(8367.986,11858.454)
	(8379.871,11802.591)
	(8389.680,11737.193)
	(8393.823,11700.780)
	(8397.468,11661.819)
	(8400.621,11620.255)
	(8403.289,11576.032)
	(8405.479,11529.096)
	(8407.199,11479.391)
	(8408.454,11426.863)
	(8409.252,11371.457)
	(8409.600,11313.118)
	(8409.504,11251.790)
	(8408.972,11187.420)
	(8408.009,11119.951)
	(8406.624,11049.330)
	(8404.823,10975.501)
	(8403.769,10937.367)
	(8402.613,10898.409)
	(8401.356,10858.623)
	(8400.000,10818.000)

\put(6900,10518){\ellipse{150}{150}}
\path(6750,10218)(6450,9618)(7650,9618)
	(7950,10218)(6750,10218)
\path(7125,10143)(6600,9693)
\path(6900,10443)	(6902.400,10371.173)
	(6904.878,10302.344)
	(6907.443,10236.465)
	(6910.100,10173.487)
	(6912.856,10113.363)
	(6915.719,10056.044)
	(6918.696,10001.482)
	(6921.792,9949.630)
	(6925.016,9900.439)
	(6928.373,9853.862)
	(6931.871,9809.849)
	(6935.517,9768.353)
	(6939.317,9729.327)
	(6943.278,9692.721)
	(6951.713,9626.580)
	(6960.875,9569.546)
	(6970.819,9521.234)
	(6993.277,9449.240)
	(7050.000,9393.000)

\path(7050,9393)	(7094.458,9439.489)
	(7114.439,9494.688)
	(7133.229,9566.704)
	(7142.250,9607.967)
	(7151.061,9652.173)
	(7159.689,9698.901)
	(7168.164,9747.731)
	(7176.514,9798.242)
	(7184.770,9850.013)
	(7192.959,9902.623)
	(7201.110,9955.654)
	(7209.253,10008.683)
	(7217.416,10061.290)
	(7225.628,10113.055)
	(7233.918,10163.557)
	(7242.314,10212.376)
	(7250.847,10259.091)
	(7259.545,10303.281)
	(7268.436,10344.526)
	(7286.914,10416.499)
	(7306.513,10471.646)
	(7350.000,10518.000)

\path(7350,10518)	(7408.951,10457.580)
	(7420.951,10423.249)
	(7431.965,10380.403)
	(7442.049,10328.630)
	(7451.257,10267.519)
	(7459.645,10196.657)
	(7463.549,10157.441)
	(7467.268,10115.632)
	(7470.809,10071.180)
	(7474.179,10024.033)
	(7477.386,9974.140)
	(7480.435,9921.448)
	(7483.335,9865.906)
	(7486.090,9807.464)
	(7488.710,9746.068)
	(7491.200,9681.669)
	(7493.567,9614.214)
	(7495.818,9543.652)
	(7497.960,9469.931)
	(7498.992,9431.870)
	(7500.000,9393.000)

\put(1793,12024){\ellipse{150}{150}}
\path(1650,11723)(1350,11123)(2550,11123)
	(2850,11723)(1650,11723)
\path(1793,11960)	(1792.281,11897.481)
	(1791.735,11837.566)
	(1791.368,11780.214)
	(1791.182,11725.383)
	(1791.182,11673.029)
	(1791.371,11623.113)
	(1791.753,11575.591)
	(1792.332,11530.422)
	(1793.111,11487.563)
	(1794.095,11446.973)
	(1796.690,11372.431)
	(1800.146,11306.460)
	(1804.495,11248.724)
	(1809.766,11198.887)
	(1815.988,11156.614)
	(1831.410,11093.414)
	(1875.000,11043.000)

\path(1875,11043)	(1911.574,11082.346)
	(1928.524,11129.964)
	(1944.835,11192.300)
	(1960.699,11266.417)
	(1968.522,11306.977)
	(1976.304,11349.382)
	(1984.070,11393.266)
	(1991.842,11438.260)
	(1999.645,11483.999)
	(2007.503,11530.116)
	(2015.438,11576.244)
	(2023.475,11622.015)
	(2031.638,11667.064)
	(2039.950,11711.023)
	(2048.436,11753.526)
	(2057.118,11794.205)
	(2075.168,11868.625)
	(2094.292,11931.349)
	(2114.678,11979.443)
	(2160.000,12020.000)

\path(2160,12020)	(2210.182,11966.103)
	(2229.785,11906.135)
	(2238.424,11869.426)
	(2246.369,11828.816)
	(2253.680,11784.741)
	(2260.416,11737.640)
	(2266.638,11687.949)
	(2272.407,11636.104)
	(2277.781,11582.544)
	(2282.822,11527.704)
	(2287.589,11472.022)
	(2292.142,11415.935)
	(2296.543,11359.880)
	(2300.849,11304.293)
	(2305.123,11249.612)
	(2309.423,11196.274)
	(2313.810,11144.715)
	(2318.344,11095.373)
	(2323.086,11048.685)
	(2328.095,11005.087)
	(2333.431,10965.017)
	(2339.155,10928.911)
	(2352.005,10870.340)
	(2385.000,10820.000)

\path(2385,10820)	(2435.879,10882.259)
	(2446.702,10918.561)
	(2456.875,10964.035)
	(2466.455,11019.123)
	(2475.495,11084.263)
	(2479.830,11120.739)
	(2484.051,11159.894)
	(2488.165,11201.782)
	(2492.178,11246.457)
	(2496.098,11293.975)
	(2499.930,11344.391)
	(2503.683,11397.759)
	(2507.363,11454.134)
	(2510.977,11513.572)
	(2514.532,11576.128)
	(2518.034,11641.855)
	(2521.491,11710.810)
	(2524.909,11783.047)
	(2526.606,11820.414)
	(2528.295,11858.621)
	(2529.979,11897.677)
	(2531.657,11937.587)
	(2533.330,11978.359)
	(2535.000,12020.000)

\put(900,10818){\makebox(0,0)[lb]{\smash{{{\SetFigFont{8}{9.6}{rm}4,028,112}}}}}
\put(2850,11343){\makebox(0,0)[lb]{\smash{{{\SetFigFont{8}{9.6}{rm}{\bf 55}}}}}}
\put(3525,10068){\ellipse{150}{150}}
\path(2550,10218)(2250,9618)(3450,9618)
	(3750,10218)(2550,10218)
\path(2550,10143)	(2555.246,10084.334)
	(2560.786,10030.509)
	(2566.675,9981.360)
	(2572.968,9936.723)
	(2579.720,9896.433)
	(2586.986,9860.324)
	(2603.280,9799.995)
	(2622.289,9754.416)
	(2644.452,9722.270)
	(2700.000,9693.000)

\path(2700,9693)	(2743.968,9707.735)
	(2780.565,9757.298)
	(2812.552,9830.205)
	(2827.680,9871.824)
	(2842.691,9914.974)
	(2857.931,9958.218)
	(2873.744,10000.121)
	(2908.472,10074.163)
	(2949.637,10125.617)
	(3000.000,10143.000)

\path(3000,10143)	(3071.177,10072.692)
	(3089.997,10003.659)
	(3097.497,9964.756)
	(3104.591,9924.484)
	(3111.933,9884.009)
	(3120.174,9844.498)
	(3141.960,9773.034)
	(3175.163,9719.424)
	(3225.000,9693.000)

\path(3225,9693)	(3297.053,9699.015)
	(3366.383,9746.333)
	(3402.224,9787.677)
	(3440.021,9841.984)
	(3480.654,9910.132)
	(3502.308,9949.671)
	(3525.000,9993.000)

\path(8550,6618)(8250,6018)(9450,6018)
	(9750,6618)(8550,6618)
\path(8880,6543)(9293,6063)
\path(8625,6918)	(8653.568,6844.418)
	(8680.811,6777.168)
	(8706.867,6716.060)
	(8731.872,6660.900)
	(8755.964,6611.497)
	(8779.281,6567.658)
	(8824.136,6495.904)
	(8867.537,6444.099)
	(8910.583,6410.705)
	(8954.371,6394.185)
	(9000.000,6393.000)

\path(9000,6393)	(9063.216,6460.811)
	(9084.063,6524.456)
	(9100.046,6596.025)
	(9112.899,6666.591)
	(9124.353,6727.226)
	(9150.000,6783.000)

\path(9150,6783)	(9195.315,6719.335)
	(9201.923,6685.370)
	(9206.644,6643.387)
	(9209.476,6592.990)
	(9210.420,6533.780)
	(9209.476,6465.358)
	(9208.296,6427.568)
	(9206.644,6387.325)
	(9204.520,6344.580)
	(9201.923,6299.284)
	(9198.855,6251.385)
	(9195.315,6200.835)
	(9191.303,6147.583)
	(9186.818,6091.581)
	(9181.862,6032.777)
	(9176.434,5971.122)
	(9170.533,5906.567)
	(9164.161,5839.062)
	(9157.316,5768.556)
	(9150.000,5695.000)

\path(8235,5118)(7935,4518)(9135,4518)
	(9435,5118)(8235,5118)
\path(8460,4818)(9210,4968)
\path(8310,5343)	(8311.383,5281.937)
	(8312.911,5223.410)
	(8314.591,5167.380)
	(8316.428,5113.804)
	(8318.431,5062.642)
	(8320.605,5013.851)
	(8322.959,4967.392)
	(8325.498,4923.223)
	(8328.229,4881.303)
	(8331.161,4841.589)
	(8337.649,4768.621)
	(8345.017,4703.987)
	(8353.320,4647.360)
	(8362.613,4598.409)
	(8372.951,4556.804)
	(8396.982,4494.317)
	(8460.000,4443.000)

\path(8460,4443)	(8513.589,4481.213)
	(8538.547,4528.770)
	(8562.510,4591.295)
	(8585.649,4665.792)
	(8596.964,4706.593)
	(8608.137,4749.263)
	(8619.190,4793.426)
	(8630.145,4838.709)
	(8641.024,4884.737)
	(8651.846,4931.134)
	(8662.635,4977.526)
	(8673.411,5023.538)
	(8684.196,5068.797)
	(8695.012,5112.925)
	(8705.880,5155.551)
	(8716.821,5196.297)
	(8739.010,5270.655)
	(8761.750,5333.003)
	(8785.214,5380.341)
	(8835.000,5418.000)

\path(8835,5418)	(8876.364,5398.528)
	(8910.387,5348.794)
	(8924.784,5311.479)
	(8937.510,5265.281)
	(8948.622,5209.760)
	(8958.173,5144.475)
	(8962.380,5108.034)
	(8966.218,5068.988)
	(8969.693,5027.281)
	(8972.813,4982.859)
	(8975.584,4935.667)
	(8978.013,4885.649)
	(8980.106,4832.752)
	(8981.872,4776.919)
	(8983.315,4718.096)
	(8984.445,4656.228)
	(8985.267,4591.260)
	(8985.787,4523.137)
	(8986.014,4451.805)
	(8985.954,4377.208)
	(8985.819,4338.668)
	(8985.614,4299.291)
	(8985.341,4259.071)
	(8985.000,4218.000)

\put(2025,4968){\ellipse{150}{150}}
\path(1050,5118)(750,4518)(1950,4518)
	(2250,5118)(1050,5118)
\path(1050,5043)	(1055.246,4984.334)
	(1060.786,4930.509)
	(1066.675,4881.360)
	(1072.968,4836.723)
	(1079.720,4796.433)
	(1086.986,4760.324)
	(1103.280,4699.995)
	(1122.289,4654.416)
	(1144.452,4622.270)
	(1200.000,4593.000)

\path(1200,4593)	(1243.968,4607.735)
	(1280.565,4657.298)
	(1312.552,4730.205)
	(1327.680,4771.824)
	(1342.691,4814.974)
	(1357.931,4858.218)
	(1373.744,4900.121)
	(1408.472,4974.163)
	(1449.637,5025.617)
	(1500.000,5043.000)

\path(1500,5043)	(1571.177,4972.692)
	(1589.997,4903.659)
	(1597.497,4864.756)
	(1604.591,4824.484)
	(1611.933,4784.009)
	(1620.174,4744.498)
	(1641.960,4673.034)
	(1675.163,4619.424)
	(1725.000,4593.000)

\path(1725,4593)	(1797.053,4599.015)
	(1866.383,4646.333)
	(1902.224,4687.677)
	(1940.021,4741.984)
	(1980.654,4810.132)
	(2002.308,4849.671)
	(2025.000,4893.000)

\put(893,6924){\ellipse{150}{150}}
\path(750,6623)(450,6023)(1650,6023)
	(1950,6623)(750,6623)
\path(893,6860)	(892.281,6797.481)
	(891.735,6737.566)
	(891.368,6680.214)
	(891.182,6625.383)
	(891.182,6573.029)
	(891.371,6523.113)
	(891.753,6475.591)
	(892.332,6430.422)
	(893.111,6387.563)
	(894.095,6346.973)
	(896.690,6272.431)
	(900.146,6206.460)
	(904.495,6148.724)
	(909.766,6098.887)
	(915.988,6056.614)
	(931.410,5993.414)
	(975.000,5943.000)

\path(975,5943)	(1011.574,5982.346)
	(1028.524,6029.964)
	(1044.835,6092.300)
	(1060.699,6166.417)
	(1068.522,6206.977)
	(1076.304,6249.382)
	(1084.070,6293.266)
	(1091.842,6338.260)
	(1099.645,6383.999)
	(1107.503,6430.116)
	(1115.438,6476.244)
	(1123.475,6522.015)
	(1131.638,6567.064)
	(1139.950,6611.023)
	(1148.436,6653.526)
	(1157.118,6694.205)
	(1175.168,6768.625)
	(1194.292,6831.349)
	(1214.678,6879.443)
	(1260.000,6920.000)

\path(1260,6920)	(1310.182,6866.103)
	(1329.785,6806.135)
	(1338.424,6769.426)
	(1346.369,6728.816)
	(1353.680,6684.741)
	(1360.416,6637.640)
	(1366.638,6587.949)
	(1372.407,6536.104)
	(1377.781,6482.544)
	(1382.822,6427.704)
	(1387.589,6372.022)
	(1392.142,6315.935)
	(1396.543,6259.880)
	(1400.849,6204.293)
	(1405.123,6149.612)
	(1409.423,6096.274)
	(1413.810,6044.715)
	(1418.344,5995.373)
	(1423.086,5948.685)
	(1428.095,5905.087)
	(1433.431,5865.017)
	(1439.155,5828.911)
	(1452.005,5770.340)
	(1485.000,5720.000)

\path(1485,5720)	(1535.879,5782.259)
	(1546.702,5818.561)
	(1556.875,5864.035)
	(1566.455,5919.123)
	(1575.495,5984.263)
	(1579.830,6020.739)
	(1584.051,6059.894)
	(1588.165,6101.782)
	(1592.178,6146.457)
	(1596.098,6193.975)
	(1599.930,6244.391)
	(1603.683,6297.759)
	(1607.363,6354.134)
	(1610.977,6413.572)
	(1614.532,6476.128)
	(1618.034,6541.855)
	(1621.491,6610.810)
	(1624.909,6683.047)
	(1626.606,6720.414)
	(1628.295,6758.621)
	(1629.979,6797.677)
	(1631.657,6837.587)
	(1633.330,6878.359)
	(1635.000,6920.000)

\put(7800,3918){\ellipse{150}{150}}
\path(7650,3618)(7350,3018)(8550,3018)
	(8850,3618)(7650,3618)
\path(8025,3543)(7500,3093)
\path(7800,3843)	(7802.400,3771.173)
	(7804.878,3702.344)
	(7807.443,3636.465)
	(7810.100,3573.487)
	(7812.856,3513.363)
	(7815.719,3456.044)
	(7818.696,3401.482)
	(7821.792,3349.630)
	(7825.016,3300.439)
	(7828.373,3253.862)
	(7831.871,3209.849)
	(7835.517,3168.353)
	(7839.317,3129.327)
	(7843.278,3092.721)
	(7851.713,3026.580)
	(7860.875,2969.546)
	(7870.819,2921.234)
	(7893.277,2849.240)
	(7950.000,2793.000)

\path(7950,2793)	(7994.458,2839.489)
	(8014.439,2894.688)
	(8033.229,2966.704)
	(8042.250,3007.967)
	(8051.061,3052.173)
	(8059.689,3098.901)
	(8068.164,3147.731)
	(8076.514,3198.242)
	(8084.770,3250.013)
	(8092.959,3302.623)
	(8101.110,3355.654)
	(8109.253,3408.683)
	(8117.416,3461.290)
	(8125.628,3513.055)
	(8133.918,3563.557)
	(8142.314,3612.376)
	(8150.847,3659.091)
	(8159.545,3703.281)
	(8168.436,3744.526)
	(8186.914,3816.499)
	(8206.513,3871.646)
	(8250.000,3918.000)

\path(8250,3918)	(8308.951,3857.580)
	(8320.951,3823.249)
	(8331.965,3780.403)
	(8342.049,3728.630)
	(8351.257,3667.519)
	(8359.645,3596.657)
	(8363.549,3557.441)
	(8367.268,3515.632)
	(8370.809,3471.180)
	(8374.179,3424.033)
	(8377.386,3374.140)
	(8380.435,3321.448)
	(8383.335,3265.906)
	(8386.090,3207.464)
	(8388.710,3146.068)
	(8391.200,3081.669)
	(8393.567,3014.214)
	(8395.818,2943.652)
	(8397.960,2869.931)
	(8398.992,2831.870)
	(8400.000,2793.000)

\put(2100,1968){\ellipse{150}{150}}
\path(2025,2118)(1725,1518)(2925,1518)
	(3225,2118)(2025,2118)
\path(2100,1893)	(2108.041,1853.311)
	(2116.003,1816.927)
	(2131.912,1753.633)
	(2148.166,1702.239)
	(2165.205,1661.865)
	(2203.396,1610.664)
	(2250.000,1593.000)

\path(2250,1593)	(2298.905,1609.672)
	(2336.535,1659.397)
	(2366.899,1731.305)
	(2380.608,1772.178)
	(2394.004,1814.520)
	(2407.587,1856.970)
	(2421.858,1898.170)
	(2454.470,1971.382)
	(2495.848,2023.283)
	(2550.000,2043.000)

\path(2550,2043)	(2594.268,2037.507)
	(2634.844,2020.148)
	(2672.607,1989.606)
	(2708.438,1944.562)
	(2743.213,1883.698)
	(2760.480,1846.922)
	(2777.812,1805.695)
	(2795.321,1759.855)
	(2813.115,1709.235)
	(2831.305,1653.672)
	(2850.000,1593.000)

\path(2550,2418)(2550,1218)
\path(6750,2118)(6450,1518)(7650,1518)
	(7950,2118)(6750,2118)
\path(7050,1593)(7050,2043)
\path(6825,2343)	(6837.918,2272.728)
	(6850.604,2208.347)
	(6863.140,2149.666)
	(6875.610,2096.492)
	(6888.094,2048.634)
	(6900.676,2005.898)
	(6926.464,1935.026)
	(6953.631,1882.339)
	(6982.836,1846.297)
	(7050.000,1818.000)

\path(7050,1818)	(7103.119,1843.212)
	(7147.518,1912.365)
	(7167.113,1958.775)
	(7185.326,2010.599)
	(7202.423,2065.979)
	(7218.671,2123.059)
	(7234.336,2179.980)
	(7249.683,2234.886)
	(7264.978,2285.920)
	(7280.489,2331.224)
	(7313.218,2397.214)
	(7350.000,2418.000)

\path(7350,2418)	(7398.586,2395.342)
	(7437.768,2343.332)
	(7453.970,2305.221)
	(7467.986,2258.454)
	(7479.871,2202.591)
	(7489.680,2137.193)
	(7493.823,2100.780)
	(7497.468,2061.819)
	(7500.621,2020.255)
	(7503.289,1976.032)
	(7505.479,1929.096)
	(7507.199,1879.391)
	(7508.454,1826.863)
	(7509.252,1771.457)
	(7509.600,1713.118)
	(7509.504,1651.790)
	(7508.972,1587.420)
	(7508.009,1519.951)
	(7506.624,1449.330)
	(7504.823,1375.501)
	(7503.769,1337.367)
	(7502.613,1298.409)
	(7501.356,1258.623)
	(7500.000,1218.000)

\put(3525,618){\ellipse{300}{1200}}
\path(3225,918)(2925,318)(4125,318)
	(4425,918)(3225,918)
\path(3225,618)(4125,618)
\path(3900,1218)(3900,18)
\path(5175,918)(4875,318)(6075,318)
	(6375,918)(5175,918)
\path(5250,1218)	(5252.997,1153.522)
	(5256.839,1093.666)
	(5261.597,1038.214)
	(5267.337,986.944)
	(5274.131,939.637)
	(5282.045,896.074)
	(5291.149,856.035)
	(5301.513,819.300)
	(5326.291,754.863)
	(5356.931,701.006)
	(5393.984,655.971)
	(5438.000,618.000)

\path(5438,618)	(5477.262,604.077)
	(5523.080,605.381)
	(5567.608,612.995)
	(5603.000,618.000)

\path(5603,618)	(5667.762,631.790)
	(5738.000,625.000)

\path(5738,625)	(5784.347,586.789)
	(5823.034,541.374)
	(5854.611,486.976)
	(5879.624,421.816)
	(5889.840,384.645)
	(5898.621,344.117)
	(5906.036,300.010)
	(5912.151,252.100)
	(5917.037,200.167)
	(5920.762,143.987)
	(5923.393,83.339)
	(5925.000,18.000)

\path(5625,843)(5625,393)
\put(4725,8868){\ellipse{150}{150}}
\path(4650,9018)(4350,8418)(5550,8418)
	(5850,9018)(4650,9018)
\path(4725,8793)	(4733.041,8753.311)
	(4741.003,8716.927)
	(4756.912,8653.633)
	(4773.166,8602.239)
	(4790.205,8561.865)
	(4828.396,8510.664)
	(4875.000,8493.000)

\path(4875,8493)	(4923.905,8509.672)
	(4961.535,8559.398)
	(4991.899,8631.305)
	(5005.608,8672.178)
	(5019.004,8714.520)
	(5032.587,8756.970)
	(5046.858,8798.170)
	(5079.470,8871.382)
	(5120.848,8923.283)
	(5175.000,8943.000)

\path(5175,8943)	(5219.268,8937.507)
	(5259.844,8920.148)
	(5297.607,8889.606)
	(5333.438,8844.562)
	(5368.213,8783.698)
	(5385.480,8746.922)
	(5402.812,8705.695)
	(5420.321,8659.855)
	(5438.115,8609.235)
	(5456.305,8553.672)
	(5475.000,8493.000)

\path(5175,9318)(5175,8118)
\put(4793,13824){\ellipse{150}{150}}
\path(4650,13523)(4350,12923)(5550,12923)
	(5850,13523)(4650,13523)
\path(4793,13760)	(4792.281,13697.481)
	(4791.735,13637.566)
	(4791.368,13580.214)
	(4791.182,13525.383)
	(4791.182,13473.029)
	(4791.371,13423.113)
	(4791.753,13375.591)
	(4792.332,13330.422)
	(4793.111,13287.563)
	(4794.095,13246.973)
	(4796.690,13172.431)
	(4800.146,13106.460)
	(4804.495,13048.724)
	(4809.766,12998.887)
	(4815.988,12956.614)
	(4831.410,12893.414)
	(4875.000,12843.000)

\path(4875,12843)	(4911.574,12882.346)
	(4928.524,12929.964)
	(4944.835,12992.300)
	(4960.699,13066.417)
	(4968.522,13106.977)
	(4976.304,13149.382)
	(4984.070,13193.266)
	(4991.842,13238.260)
	(4999.645,13283.999)
	(5007.503,13330.116)
	(5015.438,13376.244)
	(5023.475,13422.015)
	(5031.638,13467.064)
	(5039.950,13511.023)
	(5048.436,13553.526)
	(5057.118,13594.205)
	(5075.168,13668.625)
	(5094.292,13731.349)
	(5114.678,13779.443)
	(5160.000,13820.000)

\path(5160,13820)	(5210.182,13766.103)
	(5229.785,13706.135)
	(5238.424,13669.426)
	(5246.369,13628.816)
	(5253.680,13584.741)
	(5260.416,13537.640)
	(5266.638,13487.949)
	(5272.407,13436.104)
	(5277.781,13382.544)
	(5282.822,13327.704)
	(5287.589,13272.022)
	(5292.142,13215.935)
	(5296.543,13159.880)
	(5300.849,13104.293)
	(5305.123,13049.612)
	(5309.423,12996.274)
	(5313.810,12944.715)
	(5318.344,12895.373)
	(5323.086,12848.685)
	(5328.095,12805.087)
	(5333.431,12765.017)
	(5339.155,12728.911)
	(5352.005,12670.340)
	(5385.000,12620.000)

\path(5385,12620)	(5435.879,12682.259)
	(5446.702,12718.561)
	(5456.875,12764.035)
	(5466.455,12819.123)
	(5475.495,12884.263)
	(5479.830,12920.739)
	(5484.051,12959.894)
	(5488.165,13001.782)
	(5492.178,13046.457)
	(5496.098,13093.975)
	(5499.930,13144.391)
	(5503.683,13197.759)
	(5507.363,13254.134)
	(5510.977,13313.572)
	(5514.532,13376.128)
	(5518.034,13441.855)
	(5521.491,13510.810)
	(5524.909,13583.047)
	(5526.606,13620.414)
	(5528.295,13658.621)
	(5529.979,13697.677)
	(5531.657,13737.587)
	(5533.330,13778.359)
	(5535.000,13820.000)

\put(1193,13524){\ellipse{150}{150}}
\path(1050,13223)(750,12623)(1950,12623)
	(2250,13223)(1050,13223)
\path(1193,13460)	(1192.281,13397.481)
	(1191.735,13337.566)
	(1191.368,13280.214)
	(1191.182,13225.383)
	(1191.182,13173.029)
	(1191.371,13123.113)
	(1191.753,13075.591)
	(1192.332,13030.422)
	(1193.111,12987.563)
	(1194.095,12946.973)
	(1196.690,12872.431)
	(1200.146,12806.460)
	(1204.495,12748.724)
	(1209.766,12698.887)
	(1215.988,12656.614)
	(1231.410,12593.414)
	(1275.000,12543.000)

\path(1275,12543)	(1311.574,12582.346)
	(1328.524,12629.964)
	(1344.835,12692.300)
	(1360.699,12766.417)
	(1368.522,12806.977)
	(1376.304,12849.382)
	(1384.070,12893.266)
	(1391.842,12938.260)
	(1399.645,12983.999)
	(1407.503,13030.116)
	(1415.438,13076.244)
	(1423.475,13122.015)
	(1431.638,13167.064)
	(1439.950,13211.023)
	(1448.436,13253.526)
	(1457.118,13294.205)
	(1475.168,13368.625)
	(1494.292,13431.349)
	(1514.678,13479.443)
	(1560.000,13520.000)

\path(1560,13520)	(1610.182,13466.103)
	(1629.785,13406.135)
	(1638.424,13369.426)
	(1646.369,13328.816)
	(1653.680,13284.741)
	(1660.416,13237.640)
	(1666.638,13187.949)
	(1672.407,13136.104)
	(1677.781,13082.544)
	(1682.822,13027.704)
	(1687.589,12972.022)
	(1692.142,12915.935)
	(1696.543,12859.880)
	(1700.849,12804.293)
	(1705.123,12749.612)
	(1709.423,12696.274)
	(1713.810,12644.715)
	(1718.344,12595.373)
	(1723.086,12548.685)
	(1728.095,12505.087)
	(1733.431,12465.017)
	(1739.155,12428.911)
	(1752.005,12370.340)
	(1785.000,12320.000)

\path(1785,12320)	(1835.879,12382.259)
	(1846.702,12418.561)
	(1856.875,12464.035)
	(1866.455,12519.123)
	(1875.495,12584.263)
	(1879.830,12620.739)
	(1884.051,12659.894)
	(1888.165,12701.782)
	(1892.178,12746.457)
	(1896.098,12793.975)
	(1899.930,12844.391)
	(1903.683,12897.759)
	(1907.363,12954.134)
	(1910.977,13013.572)
	(1914.532,13076.128)
	(1918.034,13141.855)
	(1921.491,13210.810)
	(1924.909,13283.047)
	(1926.606,13320.414)
	(1928.295,13358.621)
	(1929.979,13397.677)
	(1931.657,13437.587)
	(1933.330,13478.359)
	(1935.000,13520.000)

\dottedline{135}(1290,12858)(1545,13128)
\dottedline{135}(1815,13068)(1575,12738)
\dottedline{135}(1800,12783)(2085,13083)
\dottedline{135}(1275,13143)(900,12693)
\put(4493,6024){\ellipse{150}{150}}
\path(4350,5723)(4050,5123)(5250,5123)
	(5550,5723)(4350,5723)
\path(4493,5960)	(4492.281,5897.481)
	(4491.735,5837.566)
	(4491.368,5780.214)
	(4491.182,5725.383)
	(4491.182,5673.029)
	(4491.371,5623.113)
	(4491.753,5575.591)
	(4492.332,5530.422)
	(4493.111,5487.563)
	(4494.095,5446.973)
	(4496.690,5372.431)
	(4500.146,5306.460)
	(4504.495,5248.724)
	(4509.766,5198.887)
	(4515.988,5156.614)
	(4531.410,5093.414)
	(4575.000,5043.000)

\path(4575,5043)	(4611.574,5082.346)
	(4628.524,5129.964)
	(4644.835,5192.300)
	(4660.699,5266.417)
	(4668.522,5306.977)
	(4676.304,5349.382)
	(4684.070,5393.266)
	(4691.842,5438.260)
	(4699.645,5483.999)
	(4707.503,5530.116)
	(4715.438,5576.244)
	(4723.475,5622.015)
	(4731.638,5667.064)
	(4739.950,5711.023)
	(4748.436,5753.526)
	(4757.118,5794.205)
	(4775.168,5868.625)
	(4794.292,5931.349)
	(4814.678,5979.443)
	(4860.000,6020.000)

\path(4860,6020)	(4910.182,5966.103)
	(4929.785,5906.135)
	(4938.424,5869.426)
	(4946.369,5828.816)
	(4953.680,5784.741)
	(4960.416,5737.640)
	(4966.638,5687.949)
	(4972.407,5636.104)
	(4977.781,5582.544)
	(4982.822,5527.704)
	(4987.589,5472.022)
	(4992.142,5415.935)
	(4996.543,5359.880)
	(5000.849,5304.293)
	(5005.123,5249.612)
	(5009.423,5196.274)
	(5013.810,5144.715)
	(5018.344,5095.373)
	(5023.086,5048.685)
	(5028.095,5005.087)
	(5033.431,4965.017)
	(5039.155,4928.911)
	(5052.005,4870.340)
	(5085.000,4820.000)

\path(5085,4820)	(5135.879,4882.259)
	(5146.702,4918.561)
	(5156.875,4964.035)
	(5166.455,5019.123)
	(5175.495,5084.263)
	(5179.830,5120.739)
	(5184.051,5159.894)
	(5188.165,5201.782)
	(5192.178,5246.457)
	(5196.098,5293.975)
	(5199.930,5344.391)
	(5203.683,5397.759)
	(5207.363,5454.134)
	(5210.977,5513.572)
	(5214.532,5576.128)
	(5218.034,5641.855)
	(5221.491,5710.810)
	(5224.909,5783.047)
	(5226.606,5820.414)
	(5228.295,5858.621)
	(5229.979,5897.677)
	(5231.657,5937.587)
	(5233.330,5978.359)
	(5235.000,6020.000)

\dottedline{135}(4590,5358)(4845,5628)
\dottedline{135}(5115,5568)(4875,5238)
\dottedline{135}(5100,5283)(5385,5583)
\dottedline{135}(4575,5643)(4200,5193)
\put(1802.018,3322.349){\arc{891.445}{3.0126}{6.3941}}
\put(1805.047,3326.276){\arc{909.963}{0.4829}{2.6784}}
\put(1600,6381){\blacken\ellipse{74}{74}}
\put(1600,6381){\ellipse{74}{74}}
\put(5250,1293){\ellipse{150}{150}}
\put(6825,2418){\ellipse{150}{150}}
\put(7725,12018){\ellipse{150}{150}}
\put(2295,11388){\blacken\ellipse{74}{74}}
\put(2295,11388){\ellipse{74}{74}}
\put(1800,3318){\ellipse{1050}{450}}
\dottedline{135}(1890,11643)(1605,11373)
\dottedline{135}(1890,11358)(2160,11508)
\dottedline{135}(2700,10083)(2445,9798)
\dottedline{135}(2745,9933)(2925,9708)
\dottedline{135}(2970,9918)(3345,10098)
\dottedline{135}(4875,8808)(4575,8538)
\dottedline{135}(5265,8748)(5580,8928)
\dottedline{135}(5265,8553)(5640,8808)
\dottedline{135}(7035,9798)(7380,10098)
\path(4125,12543)(3075,11793)
\path(3155.211,11887.161)(3075.000,11793.000)(3190.085,11838.337)
\path(4575,12243)(3675,10668)
\path(3708.489,10787.073)(3675.000,10668.000)(3760.584,10757.305)
\path(5100,12093)(5100,9618)
\path(5070.000,9738.000)(5100.000,9618.000)(5130.000,9738.000)
\path(5550,12318)(6525,10443)
\path(6443.021,10535.625)(6525.000,10443.000)(6496.254,10563.307)
\path(5775,12693)(7350,11718)
\path(7232.178,11755.655)(7350.000,11718.000)(7263.759,11806.671)
\path(6225,13068)(7875,13068)
\path(7755.000,13038.000)(7875.000,13068.000)(7755.000,13098.000)
\path(4050,12993)(2475,13143)
\path(2597.304,13161.488)(2475.000,13143.000)(2591.615,13101.758)
\dottedline{135}(995,6563)(665,6223)
\dottedline{135}(985,6363)(1285,6113)
\dottedline{135}(1225,6353)(1495,6533)
\dottedline{135}(1195,5013)(985,4793)
\dottedline{135}(1215,4903)(975,4623)
\dottedline{135}(1415,4803)(1725,5023)
\dottedline{135}(1465,4613)(1845,4873)
\dottedline{135}(2255,1943)(1955,1673)
\dottedline{135}(2335,1863)(2065,1613)
\dottedline{135}(2635,1743)(2855,2023)
\dottedline{135}(2685,1593)(3005,1963)
\dottedline{135}(7255,1773)(7755,1993)
\dottedline{135}(7965,3213)(8265,3453)
\dottedline{135}(8245,3173)(8675,3503)
\dottedline{135}(8455,5033)(8165,4653)
\path(3525,5493)(2400,6168)
\path(3525,5493)(2400,6168)
\path(2518.334,6131.985)(2400.000,6168.000)(2487.464,6080.536)
\path(3300,5043)(2400,4818)
\path(3300,5043)(2400,4818)
\path(2509.141,4876.209)(2400.000,4818.000)(2523.693,4818.000)
\path(3975,4593)(2700,3543)
\path(3975,4593)(2700,3543)
\path(2773.560,3642.443)(2700.000,3543.000)(2811.703,3596.127)
\path(4500,4443)(3375,2343)
\path(4500,4443)(3375,2343)
\path(3405.222,2462.944)(3375.000,2343.000)(3458.111,2434.611)
\path(6150,5418)(7875,6018)
\path(6150,5418)(7875,6018)
\path(7771.516,5950.243)(7875.000,6018.000)(7751.805,6006.912)
\path(6000,5118)(7650,4743)
\path(6000,5118)(7650,4743)
\path(7526.335,4740.341)(7650.000,4743.000)(7539.633,4798.849)
\path(5700,4668)(7125,3468)
\path(5700,4668)(7125,3468)
\path(7013.887,3522.349)(7125.000,3468.000)(7052.535,3568.244)
\path(5400,4443)(6450,2418)
\path(5400,4443)(6450,2418)
\path(6368.129,2510.721)(6450.000,2418.000)(6421.395,2538.340)
\path(4800,4368)(4125,1293)
\path(4800,4368)(4125,1293)
\path(4121.427,1416.642)(4125.000,1293.000)(4180.031,1403.777)
\path(5025,4368)(5475,1368)
\path(5025,4368)(5475,1368)
\path(5427.531,1482.222)(5475.000,1368.000)(5486.867,1491.123)
\path(1350,3618)(1050,3018)(2250,3018)
	(2550,3618)(1350,3618)
\dottedline{135}(5400,13443)(5700,13443)
\dottedline{135}(5400,12993)(5100,12993)
\dottedline{135}(5175,13293)(4875,13293)
\put(2250,9318){\makebox(0,0)[lb]{\smash{{{\SetFigFont{8}{9.6}{rm}14,400}}}}}
\put(3075,12168){\makebox(0,0)[lb]{\smash{{{\SetFigFont{8}{9.6}{rm}$\times 3$}}}}}
\put(5175,10668){\makebox(0,0)[lb]{\smash{{{\SetFigFont{8}{9.6}{rm}$\times 3$}}}}}
\put(6225,11268){\makebox(0,0)[lb]{\smash{{{\SetFigFont{8}{9.6}{rm}$\times 3$}}}}}
\put(6450,12318){\makebox(0,0)[lb]{\smash{{{\SetFigFont{8}{9.6}{rm}$\times 2$}}}}}
\put(5550,5343){\makebox(0,0)[lb]{\smash{{{\SetFigFont{8}{9.6}{rm}{\bf 33}}}}}}
\put(1950,6243){\makebox(0,0)[lb]{\smash{{{\SetFigFont{8}{9.6}{rm}{\bf 62}}}}}}
\put(2250,12843){\makebox(0,0)[lb]{\smash{{{\SetFigFont{8}{9.6}{rm}{\bf 33}}}}}}
\put(0,5718){\makebox(0,0)[lb]{\smash{{{\SetFigFont{8}{9.6}{rm}2,849,436}}}}}
\put(1050,4218){\makebox(0,0)[lb]{\smash{{{\SetFigFont{8}{9.6}{rm}52,160}}}}}
\put(4650,18){\makebox(0,0)[lb]{\smash{{{\SetFigFont{8}{9.6}{rm}81,345}}}}}
\put(6300,1218){\makebox(0,0)[lb]{\smash{{{\SetFigFont{8}{9.6}{rm}79,475}}}}}
\put(6975,2718){\makebox(0,0)[lb]{\smash{{{\SetFigFont{8}{9.6}{rm}71,115}}}}}
\put(7650,4218){\makebox(0,0)[lb]{\smash{{{\SetFigFont{8}{9.6}{rm}1,114,864}}}}}
\put(2850,6018){\makebox(0,0)[lb]{\smash{{{\SetFigFont{8}{9.6}{rm}$\times 4$}}}}}
\put(2400,4968){\makebox(0,0)[lb]{\smash{{{\SetFigFont{8}{9.6}{rm}$\times 4$}}}}}
\put(5475,1893){\makebox(0,0)[lb]{\smash{{{\SetFigFont{8}{9.6}{rm}$\times 3$}}}}}
\put(6225,3093){\makebox(0,0)[lb]{\smash{{{\SetFigFont{8}{9.6}{rm}$\times 8$}}}}}
\put(6525,3993){\makebox(0,0)[lb]{\smash{{{\SetFigFont{8}{9.6}{rm}$\times 6$}}}}}
\put(6900,4968){\makebox(0,0)[lb]{\smash{{{\SetFigFont{8}{9.6}{rm}$\times 4$}}}}}
\put(6750,5868){\makebox(0,0)[lb]{\smash{{{\SetFigFont{8}{9.6}{rm}$\times 2$}}}}}
\put(3225,2943){\makebox(0,0)[lb]{\smash{{{\SetFigFont{8}{9.6}{rm}$\times 3$}}}}}
\put(2325,18){\makebox(0,0)[lb]{\smash{{{\SetFigFont{8}{9.6}{rm}253,440}}}}}
\put(2625,3993){\makebox(0,0)[lb]{\smash{{{\SetFigFont{8}{9.6}{rm}$\times 32$}}}}}
\put(675,2718){\makebox(0,0)[lb]{\smash{{{\SetFigFont{8}{9.6}{rm}98,340}}}}}
\put(225,12318){\makebox(0,0)[lb]{\smash{{{\SetFigFont{8}{9.6}{rm}23,962,326}}}}}
\put(5775,13143){\makebox(0,0)[lb]{\smash{{{\SetFigFont{8}{9.6}{rm}{\bf 26}}}}}}
\put(3375,11118){\makebox(0,0)[lb]{\smash{{{\SetFigFont{8}{9.6}{rm}$\times 4$}}}}}
\put(3825,12618){\makebox(0,0)[lb]{\smash{{{\SetFigFont{8}{9.6}{rm}39,347,736}}}}}
\put(7650,12318){\makebox(0,0)[lb]{\smash{{{\SetFigFont{8}{9.6}{rm}2,404,944}}}}}
\put(7275,10818){\makebox(0,0)[lb]{\smash{{{\SetFigFont{8}{9.6}{rm}138,060}}}}}
\put(6000,9318){\makebox(0,0)[lb]{\smash{{{\SetFigFont{8}{9.6}{rm}125,100}}}}}
\put(4200,8118){\makebox(0,0)[lb]{\smash{{{\SetFigFont{8}{9.6}{rm}62,370}}}}}
\put(7800,5718){\makebox(0,0)[lb]{\smash{{{\SetFigFont{8}{9.6}{rm}1,240,648}}}}}
\put(3525,4818){\makebox(0,0)[lb]{\smash{{{\SetFigFont{8}{9.6}{rm}23,962,326}}}}}
\put(1425,1218){\makebox(0,0)[lb]{\smash{{{\SetFigFont{8}{9.6}{rm}236,115}}}}}
\end{picture}

%% file: qesc4.tex
\begingroup\makeatletter\ifx\SetFigFont\undefined
\def\x#1#2#3#4#5#6#7\relax{\def\x{#1#2#3#4#5#6}}%
\expandafter\x\fmtname xxxxxx\relax \def\y{splain}%
\ifx\x\y   
\gdef\SetFigFont#1#2#3{%
  \ifnum #1<17\tiny\else \ifnum #1<20\small\else
  \ifnum #1<24\normalsize\else \ifnum #1<29\large\else
  \ifnum #1<34\Large\else \ifnum #1<41\LARGE\else
     \huge\fi\fi\fi\fi\fi\fi
  \csname #3\endcsname}%
\else
\gdef\SetFigFont#1#2#3{\begingroup
  \count@#1\relax \ifnum 25<\count@\count@25\fi
  \def\x{\endgroup\@setsize\SetFigFont{#2pt}}%
  \expandafter\x
    \csname \romannumeral\the\count@ pt\expandafter\endcsname
    \csname @\romannumeral\the\count@ pt\endcsname
  \csname #3\endcsname}%
\fi
\fi\endgroup
\begin{picture}(9143,12422)(0,-10)
\thicklines
\put(4343,9624){\ellipse{150}{150}}
\path(4200,9323)(3900,8723)(5100,8723)
	(5400,9323)(4200,9323)
\path(4343,9560)	(4342.281,9497.481)
	(4341.735,9437.566)
	(4341.368,9380.214)
	(4341.182,9325.383)
	(4341.182,9273.029)
	(4341.371,9223.113)
	(4341.753,9175.591)
	(4342.332,9130.422)
	(4343.111,9087.563)
	(4344.095,9046.973)
	(4346.690,8972.431)
	(4350.146,8906.460)
	(4354.495,8848.724)
	(4359.766,8798.887)
	(4365.988,8756.614)
	(4381.410,8693.414)
	(4425.000,8643.000)

\path(4425,8643)	(4461.574,8682.346)
	(4478.524,8729.964)
	(4494.835,8792.300)
	(4510.699,8866.417)
	(4518.522,8906.977)
	(4526.304,8949.382)
	(4534.070,8993.266)
	(4541.842,9038.260)
	(4549.645,9083.999)
	(4557.503,9130.116)
	(4565.438,9176.244)
	(4573.475,9222.015)
	(4581.638,9267.064)
	(4589.950,9311.023)
	(4598.436,9353.526)
	(4607.118,9394.205)
	(4625.168,9468.625)
	(4644.292,9531.349)
	(4664.678,9579.443)
	(4710.000,9620.000)

\path(4710,9620)	(4760.182,9566.103)
	(4779.785,9506.135)
	(4788.424,9469.426)
	(4796.369,9428.816)
	(4803.680,9384.741)
	(4810.416,9337.640)
	(4816.638,9287.949)
	(4822.407,9236.104)
	(4827.781,9182.544)
	(4832.822,9127.704)
	(4837.589,9072.022)
	(4842.142,9015.935)
	(4846.543,8959.880)
	(4850.849,8904.293)
	(4855.123,8849.612)
	(4859.423,8796.274)
	(4863.810,8744.715)
	(4868.344,8695.373)
	(4873.086,8648.685)
	(4878.095,8605.087)
	(4883.431,8565.017)
	(4889.155,8528.911)
	(4902.005,8470.340)
	(4935.000,8420.000)

\path(4935,8420)	(4985.879,8482.259)
	(4996.702,8518.561)
	(5006.875,8564.035)
	(5016.455,8619.123)
	(5025.495,8684.263)
	(5029.830,8720.739)
	(5034.051,8759.894)
	(5038.165,8801.782)
	(5042.178,8846.457)
	(5046.098,8893.975)
	(5049.930,8944.391)
	(5053.683,8997.759)
	(5057.363,9054.134)
	(5060.977,9113.572)
	(5064.532,9176.128)
	(5068.034,9241.855)
	(5071.491,9310.810)
	(5074.909,9383.047)
	(5076.606,9420.414)
	(5078.295,9458.621)
	(5079.979,9497.677)
	(5081.657,9537.587)
	(5083.330,9578.359)
	(5085.000,9620.000)

\put(4343,12324){\ellipse{150}{150}}
\path(4200,12023)(3900,11423)(5100,11423)
	(5400,12023)(4200,12023)
\path(4343,12260)	(4342.281,12197.481)
	(4341.735,12137.566)
	(4341.368,12080.214)
	(4341.182,12025.383)
	(4341.182,11973.029)
	(4341.371,11923.113)
	(4341.753,11875.591)
	(4342.332,11830.422)
	(4343.111,11787.563)
	(4344.095,11746.973)
	(4346.690,11672.431)
	(4350.146,11606.460)
	(4354.495,11548.724)
	(4359.766,11498.887)
	(4365.988,11456.614)
	(4381.410,11393.414)
	(4425.000,11343.000)

\path(4425,11343)	(4461.574,11382.346)
	(4478.524,11429.964)
	(4494.835,11492.300)
	(4510.699,11566.417)
	(4518.522,11606.977)
	(4526.304,11649.382)
	(4534.070,11693.266)
	(4541.842,11738.260)
	(4549.645,11783.999)
	(4557.503,11830.116)
	(4565.438,11876.244)
	(4573.475,11922.015)
	(4581.638,11967.064)
	(4589.950,12011.023)
	(4598.436,12053.526)
	(4607.118,12094.205)
	(4625.168,12168.625)
	(4644.292,12231.349)
	(4664.678,12279.443)
	(4710.000,12320.000)

\path(4710,12320)	(4760.182,12266.103)
	(4779.785,12206.135)
	(4788.424,12169.426)
	(4796.369,12128.816)
	(4803.680,12084.741)
	(4810.416,12037.640)
	(4816.638,11987.949)
	(4822.407,11936.104)
	(4827.781,11882.544)
	(4832.822,11827.704)
	(4837.589,11772.022)
	(4842.142,11715.935)
	(4846.543,11659.880)
	(4850.849,11604.293)
	(4855.123,11549.612)
	(4859.423,11496.274)
	(4863.810,11444.715)
	(4868.344,11395.373)
	(4873.086,11348.685)
	(4878.095,11305.087)
	(4883.431,11265.017)
	(4889.155,11228.911)
	(4902.005,11170.340)
	(4935.000,11120.000)

\path(4935,11120)	(4985.879,11182.259)
	(4996.702,11218.561)
	(5006.875,11264.035)
	(5016.455,11319.123)
	(5025.495,11384.263)
	(5029.830,11420.739)
	(5034.051,11459.894)
	(5038.165,11501.782)
	(5042.178,11546.457)
	(5046.098,11593.975)
	(5049.930,11644.391)
	(5053.683,11697.759)
	(5057.363,11754.134)
	(5060.977,11813.572)
	(5064.532,11876.128)
	(5068.034,11941.855)
	(5071.491,12010.810)
	(5074.909,12083.047)
	(5076.606,12120.414)
	(5078.295,12158.621)
	(5079.979,12197.677)
	(5081.657,12237.587)
	(5083.330,12278.359)
	(5085.000,12320.000)

\put(6975,11268){\ellipse{150}{150}}
\path(6000,11418)(5700,10818)(6900,10818)
	(7200,11418)(6000,11418)
\path(6000,11343)	(6005.246,11284.334)
	(6010.786,11230.509)
	(6016.675,11181.360)
	(6022.968,11136.723)
	(6029.720,11096.433)
	(6036.986,11060.324)
	(6053.280,10999.995)
	(6072.289,10954.416)
	(6094.452,10922.270)
	(6150.000,10893.000)

\path(6150,10893)	(6193.968,10907.735)
	(6230.565,10957.298)
	(6262.552,11030.205)
	(6277.680,11071.824)
	(6292.691,11114.974)
	(6307.931,11158.218)
	(6323.744,11200.121)
	(6358.472,11274.163)
	(6399.637,11325.617)
	(6450.000,11343.000)

\path(6450,11343)	(6521.177,11272.692)
	(6539.997,11203.659)
	(6547.497,11164.756)
	(6554.591,11124.484)
	(6561.933,11084.009)
	(6570.174,11044.498)
	(6591.960,10973.034)
	(6625.163,10919.424)
	(6675.000,10893.000)

\path(6675,10893)	(6747.052,10899.015)
	(6816.382,10946.333)
	(6852.224,10987.677)
	(6890.021,11041.984)
	(6930.654,11110.132)
	(6952.308,11149.671)
	(6975.000,11193.000)

\put(3075,6618){\ellipse{150}{150}}
\path(3000,6318)(2700,5718)(3900,5718)
	(4200,6318)(3000,6318)
\path(3300,5793)(3300,6243)
\path(3075,6543)	(3087.918,6472.728)
	(3100.604,6408.347)
	(3113.140,6349.666)
	(3125.610,6296.492)
	(3138.094,6248.634)
	(3150.676,6205.898)
	(3176.464,6135.026)
	(3203.631,6082.339)
	(3232.836,6046.297)
	(3300.000,6018.000)

\path(3300,6018)	(3353.119,6043.212)
	(3397.518,6112.365)
	(3417.113,6158.775)
	(3435.326,6210.599)
	(3452.423,6265.979)
	(3468.671,6323.059)
	(3484.336,6379.980)
	(3499.683,6434.886)
	(3514.978,6485.920)
	(3530.489,6531.224)
	(3563.218,6597.214)
	(3600.000,6618.000)

\path(3600,6618)	(3648.586,6595.342)
	(3687.768,6543.332)
	(3703.970,6505.221)
	(3717.986,6458.454)
	(3729.871,6402.591)
	(3739.680,6337.193)
	(3743.823,6300.780)
	(3747.468,6261.819)
	(3750.621,6220.255)
	(3753.289,6176.032)
	(3755.479,6129.096)
	(3757.199,6079.391)
	(3758.454,6026.863)
	(3759.252,5971.457)
	(3759.600,5913.118)
	(3759.504,5851.790)
	(3758.972,5787.420)
	(3758.009,5719.951)
	(3756.624,5649.330)
	(3754.823,5575.501)
	(3753.769,5537.367)
	(3752.613,5498.409)
	(3751.356,5458.623)
	(3750.000,5418.000)

\put(5175,6618){\ellipse{150}{150}}
\path(5100,6318)(4800,5718)(6000,5718)
	(6300,6318)(5100,6318)
\path(5400,5793)(5400,6243)
\path(5175,6543)	(5187.918,6472.728)
	(5200.604,6408.347)
	(5213.140,6349.666)
	(5225.610,6296.492)
	(5238.094,6248.634)
	(5250.676,6205.898)
	(5276.464,6135.026)
	(5303.631,6082.339)
	(5332.836,6046.297)
	(5400.000,6018.000)

\path(5400,6018)	(5453.119,6043.212)
	(5497.518,6112.365)
	(5517.113,6158.775)
	(5535.326,6210.599)
	(5552.423,6265.979)
	(5568.671,6323.059)
	(5584.336,6379.980)
	(5599.683,6434.886)
	(5614.978,6485.920)
	(5630.489,6531.224)
	(5663.218,6597.214)
	(5700.000,6618.000)

\path(5700,6618)	(5748.586,6595.342)
	(5787.768,6543.332)
	(5803.970,6505.221)
	(5817.986,6458.454)
	(5829.871,6402.591)
	(5839.680,6337.193)
	(5843.823,6300.780)
	(5847.468,6261.819)
	(5850.621,6220.255)
	(5853.289,6176.032)
	(5855.479,6129.096)
	(5857.199,6079.391)
	(5858.454,6026.863)
	(5859.252,5971.457)
	(5859.600,5913.118)
	(5859.504,5851.790)
	(5858.972,5787.420)
	(5858.009,5719.951)
	(5856.624,5649.330)
	(5854.823,5575.501)
	(5853.769,5537.367)
	(5852.613,5498.409)
	(5851.356,5458.623)
	(5850.000,5418.000)

\put(1350,8418){\ellipse{150}{150}}
\path(1200,8118)(900,7518)(2100,7518)
	(2400,8118)(1200,8118)
\path(1575,8043)(1050,7593)
\path(1350,8343)	(1352.400,8271.173)
	(1354.878,8202.344)
	(1357.443,8136.465)
	(1360.100,8073.487)
	(1362.856,8013.363)
	(1365.719,7956.044)
	(1368.696,7901.482)
	(1371.792,7849.630)
	(1375.016,7800.439)
	(1378.373,7753.862)
	(1381.871,7709.849)
	(1385.517,7668.353)
	(1389.317,7629.327)
	(1393.278,7592.721)
	(1401.712,7526.580)
	(1410.875,7469.546)
	(1420.819,7421.234)
	(1443.277,7349.240)
	(1500.000,7293.000)

\path(1500,7293)	(1544.458,7339.489)
	(1564.439,7394.688)
	(1583.229,7466.704)
	(1592.250,7507.967)
	(1601.061,7552.173)
	(1609.689,7598.901)
	(1618.164,7647.731)
	(1626.514,7698.242)
	(1634.770,7750.013)
	(1642.959,7802.623)
	(1651.110,7855.654)
	(1659.253,7908.683)
	(1667.416,7961.290)
	(1675.628,8013.055)
	(1683.918,8063.557)
	(1692.314,8112.376)
	(1700.847,8159.091)
	(1709.545,8203.281)
	(1718.436,8244.526)
	(1736.914,8316.499)
	(1756.513,8371.646)
	(1800.000,8418.000)

\path(1800,8418)	(1858.951,8357.580)
	(1870.951,8323.249)
	(1881.965,8280.403)
	(1892.049,8228.630)
	(1901.258,8167.519)
	(1909.645,8096.657)
	(1913.549,8057.441)
	(1917.268,8015.632)
	(1920.809,7971.180)
	(1924.179,7924.033)
	(1927.386,7874.140)
	(1930.435,7821.448)
	(1933.335,7765.906)
	(1936.090,7707.464)
	(1938.710,7646.068)
	(1941.200,7581.669)
	(1943.567,7514.214)
	(1945.818,7443.652)
	(1947.960,7369.931)
	(1948.992,7331.870)
	(1950.000,7293.000)

\put(1050,9618){\ellipse{150}{150}}
\path(900,9318)(600,8718)(1800,8718)
	(2100,9318)(900,9318)
\path(1275,9243)(750,8793)
\path(1050,9543)	(1052.400,9471.173)
	(1054.878,9402.344)
	(1057.443,9336.465)
	(1060.100,9273.487)
	(1062.856,9213.363)
	(1065.719,9156.044)
	(1068.696,9101.482)
	(1071.792,9049.630)
	(1075.016,9000.439)
	(1078.373,8953.862)
	(1081.871,8909.849)
	(1085.517,8868.353)
	(1089.317,8829.327)
	(1093.278,8792.721)
	(1101.712,8726.580)
	(1110.875,8669.546)
	(1120.819,8621.234)
	(1143.277,8549.240)
	(1200.000,8493.000)

\path(1200,8493)	(1244.458,8539.489)
	(1264.439,8594.688)
	(1283.229,8666.704)
	(1292.250,8707.967)
	(1301.061,8752.173)
	(1309.689,8798.901)
	(1318.164,8847.731)
	(1326.514,8898.242)
	(1334.770,8950.013)
	(1342.959,9002.623)
	(1351.110,9055.654)
	(1359.253,9108.683)
	(1367.416,9161.290)
	(1375.628,9213.055)
	(1383.918,9263.557)
	(1392.314,9312.376)
	(1400.847,9359.091)
	(1409.545,9403.281)
	(1418.436,9444.526)
	(1436.914,9516.499)
	(1456.513,9571.646)
	(1500.000,9618.000)

\path(1500,9618)	(1558.951,9557.580)
	(1570.951,9523.249)
	(1581.965,9480.403)
	(1592.049,9428.630)
	(1601.258,9367.519)
	(1609.645,9296.657)
	(1613.549,9257.441)
	(1617.268,9215.632)
	(1620.809,9171.180)
	(1624.179,9124.033)
	(1627.386,9074.140)
	(1630.435,9021.448)
	(1633.335,8965.906)
	(1636.090,8907.464)
	(1638.710,8846.068)
	(1641.200,8781.669)
	(1643.567,8714.214)
	(1645.818,8643.652)
	(1647.960,8569.931)
	(1648.992,8531.870)
	(1650.000,8493.000)

\path(1500,10518)(1200,9918)(2400,9918)
	(2700,10518)(1500,10518)
\path(1725,10218)(2475,10368)
\path(1575,10743)	(1576.383,10681.937)
	(1577.911,10623.410)
	(1579.591,10567.380)
	(1581.428,10513.804)
	(1583.431,10462.642)
	(1585.605,10413.851)
	(1587.959,10367.392)
	(1590.498,10323.223)
	(1593.229,10281.303)
	(1596.161,10241.589)
	(1602.649,10168.621)
	(1610.017,10103.987)
	(1618.320,10047.360)
	(1627.613,9998.409)
	(1637.951,9956.804)
	(1661.982,9894.317)
	(1725.000,9843.000)

\path(1725,9843)	(1778.589,9881.213)
	(1803.547,9928.770)
	(1827.510,9991.295)
	(1850.649,10065.792)
	(1861.964,10106.593)
	(1873.137,10149.263)
	(1884.190,10193.426)
	(1895.145,10238.709)
	(1906.024,10284.737)
	(1916.846,10331.134)
	(1927.635,10377.526)
	(1938.411,10423.538)
	(1949.196,10468.797)
	(1960.012,10512.925)
	(1970.880,10555.551)
	(1981.821,10596.297)
	(2004.010,10670.655)
	(2026.750,10733.003)
	(2050.214,10780.341)
	(2100.000,10818.000)

\path(2100,10818)	(2141.364,10798.528)
	(2175.387,10748.794)
	(2189.784,10711.479)
	(2202.510,10665.281)
	(2213.622,10609.760)
	(2223.172,10544.475)
	(2227.380,10508.034)
	(2231.218,10468.988)
	(2234.693,10427.281)
	(2237.813,10382.859)
	(2240.584,10335.667)
	(2243.013,10285.649)
	(2245.106,10232.752)
	(2246.872,10176.919)
	(2248.315,10118.096)
	(2249.445,10056.228)
	(2250.267,9991.260)
	(2250.787,9923.138)
	(2251.014,9851.805)
	(2250.954,9777.208)
	(2250.819,9738.668)
	(2250.614,9699.291)
	(2250.341,9659.071)
	(2250.000,9618.000)

\path(2400,11418)(2100,10818)(3300,10818)
	(3600,11418)(2400,11418)
\path(2625,11118)(3375,11268)
\path(2475,11643)	(2476.383,11581.937)
	(2477.911,11523.410)
	(2479.591,11467.380)
	(2481.428,11413.804)
	(2483.431,11362.642)
	(2485.605,11313.851)
	(2487.959,11267.392)
	(2490.498,11223.223)
	(2493.229,11181.303)
	(2496.161,11141.589)
	(2502.649,11068.621)
	(2510.017,11003.987)
	(2518.320,10947.360)
	(2527.613,10898.409)
	(2537.951,10856.804)
	(2561.982,10794.317)
	(2625.000,10743.000)

\path(2625,10743)	(2678.589,10781.213)
	(2703.547,10828.770)
	(2727.510,10891.295)
	(2750.649,10965.792)
	(2761.964,11006.593)
	(2773.137,11049.263)
	(2784.190,11093.426)
	(2795.145,11138.709)
	(2806.024,11184.737)
	(2816.846,11231.134)
	(2827.635,11277.526)
	(2838.411,11323.538)
	(2849.196,11368.797)
	(2860.012,11412.925)
	(2870.880,11455.551)
	(2881.821,11496.297)
	(2904.010,11570.655)
	(2926.750,11633.003)
	(2950.214,11680.341)
	(3000.000,11718.000)

\path(3000,11718)	(3041.364,11698.528)
	(3075.387,11648.794)
	(3089.784,11611.479)
	(3102.510,11565.281)
	(3113.622,11509.760)
	(3123.172,11444.475)
	(3127.380,11408.034)
	(3131.218,11368.988)
	(3134.693,11327.281)
	(3137.813,11282.859)
	(3140.584,11235.667)
	(3143.013,11185.649)
	(3145.106,11132.752)
	(3146.872,11076.919)
	(3148.315,11018.096)
	(3149.445,10956.228)
	(3150.267,10891.260)
	(3150.787,10823.138)
	(3151.014,10751.805)
	(3150.954,10677.208)
	(3150.819,10638.668)
	(3150.614,10599.291)
	(3150.341,10559.071)
	(3150.000,10518.000)

\put(7500,7818){\ellipse{300}{1200}}
\path(7200,8118)(6900,7518)(8100,7518)
	(8400,8118)(7200,8118)
\path(7200,7818)(8100,7818)
\path(7875,8418)(7875,7218)
\put(4343,4824){\ellipse{150}{150}}
\path(4200,4523)(3900,3923)(5100,3923)
	(5400,4523)(4200,4523)
\path(4343,4760)	(4342.281,4697.481)
	(4341.735,4637.566)
	(4341.368,4580.214)
	(4341.182,4525.383)
	(4341.182,4473.029)
	(4341.371,4423.113)
	(4341.753,4375.591)
	(4342.332,4330.422)
	(4343.111,4287.563)
	(4344.095,4246.973)
	(4346.690,4172.431)
	(4350.146,4106.460)
	(4354.495,4048.724)
	(4359.766,3998.887)
	(4365.988,3956.614)
	(4381.410,3893.414)
	(4425.000,3843.000)

\path(4425,3843)	(4461.574,3882.346)
	(4478.524,3929.964)
	(4494.835,3992.300)
	(4510.699,4066.417)
	(4518.522,4106.977)
	(4526.304,4149.382)
	(4534.070,4193.266)
	(4541.842,4238.260)
	(4549.645,4283.999)
	(4557.503,4330.116)
	(4565.438,4376.244)
	(4573.475,4422.015)
	(4581.638,4467.064)
	(4589.950,4511.023)
	(4598.436,4553.526)
	(4607.118,4594.205)
	(4625.168,4668.625)
	(4644.292,4731.349)
	(4664.678,4779.443)
	(4710.000,4820.000)

\path(4710,4820)	(4760.182,4766.103)
	(4779.785,4706.135)
	(4788.424,4669.426)
	(4796.369,4628.816)
	(4803.680,4584.741)
	(4810.416,4537.640)
	(4816.638,4487.949)
	(4822.407,4436.104)
	(4827.781,4382.544)
	(4832.822,4327.704)
	(4837.589,4272.022)
	(4842.142,4215.935)
	(4846.543,4159.880)
	(4850.849,4104.293)
	(4855.123,4049.612)
	(4859.423,3996.274)
	(4863.810,3944.715)
	(4868.344,3895.373)
	(4873.086,3848.685)
	(4878.095,3805.087)
	(4883.431,3765.017)
	(4889.155,3728.911)
	(4902.005,3670.340)
	(4935.000,3620.000)

\path(4935,3620)	(4985.879,3682.259)
	(4996.702,3718.561)
	(5006.875,3764.035)
	(5016.455,3819.123)
	(5025.495,3884.263)
	(5029.830,3920.739)
	(5034.051,3959.894)
	(5038.165,4001.782)
	(5042.178,4046.457)
	(5046.098,4093.975)
	(5049.930,4144.391)
	(5053.683,4197.759)
	(5057.363,4254.134)
	(5060.977,4313.572)
	(5064.532,4376.128)
	(5068.034,4441.855)
	(5071.491,4510.810)
	(5074.909,4583.047)
	(5076.606,4620.414)
	(5078.295,4658.621)
	(5079.979,4697.677)
	(5081.657,4737.587)
	(5083.330,4778.359)
	(5085.000,4820.000)

\put(1643,4524){\ellipse{150}{150}}
\path(1500,4223)(1200,3623)(2400,3623)
	(2700,4223)(1500,4223)
\path(1643,4460)	(1642.281,4397.481)
	(1641.735,4337.566)
	(1641.368,4280.214)
	(1641.182,4225.383)
	(1641.182,4173.029)
	(1641.371,4123.113)
	(1641.753,4075.591)
	(1642.332,4030.422)
	(1643.111,3987.563)
	(1644.095,3946.973)
	(1646.690,3872.431)
	(1650.146,3806.460)
	(1654.495,3748.724)
	(1659.766,3698.887)
	(1665.988,3656.614)
	(1681.410,3593.414)
	(1725.000,3543.000)

\path(1725,3543)	(1761.574,3582.346)
	(1778.524,3629.964)
	(1794.835,3692.300)
	(1810.699,3766.417)
	(1818.522,3806.977)
	(1826.304,3849.382)
	(1834.070,3893.266)
	(1841.842,3938.260)
	(1849.645,3983.999)
	(1857.503,4030.116)
	(1865.438,4076.244)
	(1873.475,4122.015)
	(1881.638,4167.064)
	(1889.950,4211.023)
	(1898.436,4253.526)
	(1907.118,4294.205)
	(1925.168,4368.625)
	(1944.292,4431.349)
	(1964.678,4479.443)
	(2010.000,4520.000)

\path(2010,4520)	(2060.182,4466.103)
	(2079.785,4406.135)
	(2088.424,4369.426)
	(2096.369,4328.816)
	(2103.680,4284.741)
	(2110.416,4237.640)
	(2116.638,4187.949)
	(2122.407,4136.104)
	(2127.781,4082.544)
	(2132.822,4027.704)
	(2137.589,3972.022)
	(2142.142,3915.935)
	(2146.543,3859.880)
	(2150.849,3804.293)
	(2155.123,3749.612)
	(2159.423,3696.274)
	(2163.810,3644.715)
	(2168.344,3595.373)
	(2173.086,3548.685)
	(2178.095,3505.087)
	(2183.431,3465.017)
	(2189.155,3428.911)
	(2202.005,3370.340)
	(2235.000,3320.000)

\path(2235,3320)	(2285.879,3382.259)
	(2296.702,3418.561)
	(2306.875,3464.035)
	(2316.455,3519.123)
	(2325.495,3584.263)
	(2329.830,3620.739)
	(2334.051,3659.894)
	(2338.165,3701.782)
	(2342.178,3746.457)
	(2346.098,3793.975)
	(2349.930,3844.391)
	(2353.683,3897.759)
	(2357.363,3954.134)
	(2360.977,4013.572)
	(2364.532,4076.128)
	(2368.034,4141.855)
	(2371.491,4210.810)
	(2374.909,4283.047)
	(2376.606,4320.414)
	(2378.295,4358.621)
	(2379.979,4397.677)
	(2381.657,4437.587)
	(2383.330,4478.359)
	(2385.000,4520.000)

\put(2775,2868){\ellipse{150}{150}}
\path(1800,3018)(1500,2418)(2700,2418)
	(3000,3018)(1800,3018)
\path(1800,2943)	(1805.246,2884.334)
	(1810.786,2830.509)
	(1816.675,2781.360)
	(1822.968,2736.723)
	(1829.720,2696.433)
	(1836.986,2660.324)
	(1853.280,2599.995)
	(1872.289,2554.416)
	(1894.452,2522.270)
	(1950.000,2493.000)

\path(1950,2493)	(1993.968,2507.735)
	(2030.565,2557.298)
	(2062.552,2630.205)
	(2077.680,2671.824)
	(2092.691,2714.974)
	(2107.931,2758.218)
	(2123.744,2800.121)
	(2158.472,2874.163)
	(2199.637,2925.617)
	(2250.000,2943.000)

\path(2250,2943)	(2321.177,2872.692)
	(2339.997,2803.659)
	(2347.497,2764.756)
	(2354.591,2724.484)
	(2361.933,2684.009)
	(2370.174,2644.498)
	(2391.960,2573.034)
	(2425.163,2519.424)
	(2475.000,2493.000)

\path(2475,2493)	(2547.053,2499.015)
	(2616.383,2546.333)
	(2652.224,2587.677)
	(2690.021,2641.984)
	(2730.654,2710.132)
	(2752.308,2749.671)
	(2775.000,2793.000)

\put(5175,1218){\ellipse{150}{150}}
\path(5100,918)(4800,318)(6000,318)
	(6300,918)(5100,918)
\path(5400,393)(5400,843)
\path(5175,1143)	(5187.918,1072.728)
	(5200.604,1008.347)
	(5213.140,949.666)
	(5225.610,896.492)
	(5238.094,848.634)
	(5250.676,805.898)
	(5276.464,735.026)
	(5303.631,682.339)
	(5332.836,646.297)
	(5400.000,618.000)

\path(5400,618)	(5453.119,643.212)
	(5497.518,712.365)
	(5517.113,758.775)
	(5535.326,810.599)
	(5552.423,865.979)
	(5568.671,923.059)
	(5584.336,979.980)
	(5599.683,1034.886)
	(5614.978,1085.920)
	(5630.489,1131.224)
	(5663.218,1197.214)
	(5700.000,1218.000)

\path(5700,1218)	(5748.586,1195.342)
	(5787.768,1143.332)
	(5803.970,1105.221)
	(5817.986,1058.454)
	(5829.871,1002.591)
	(5839.680,937.193)
	(5843.823,900.780)
	(5847.468,861.819)
	(5850.621,820.255)
	(5853.289,776.032)
	(5855.479,729.096)
	(5857.199,679.391)
	(5858.454,626.863)
	(5859.252,571.457)
	(5859.600,513.118)
	(5859.504,451.790)
	(5858.972,387.420)
	(5858.009,319.951)
	(5856.624,249.330)
	(5854.823,175.501)
	(5853.769,137.367)
	(5852.613,98.409)
	(5851.356,58.623)
	(5850.000,18.000)

\put(1650,7218){\ellipse{150}{150}}
\path(1500,6918)(1200,6318)(2400,6318)
	(2700,6918)(1500,6918)
\path(1875,6843)(1350,6393)
\path(1650,7143)	(1652.400,7071.173)
	(1654.878,7002.344)
	(1657.443,6936.465)
	(1660.100,6873.487)
	(1662.856,6813.363)
	(1665.719,6756.044)
	(1668.696,6701.482)
	(1671.792,6649.630)
	(1675.016,6600.439)
	(1678.373,6553.862)
	(1681.871,6509.849)
	(1685.517,6468.353)
	(1689.317,6429.327)
	(1693.278,6392.721)
	(1701.712,6326.580)
	(1710.875,6269.546)
	(1720.819,6221.234)
	(1743.277,6149.240)
	(1800.000,6093.000)

\path(1800,6093)	(1844.458,6139.489)
	(1864.439,6194.688)
	(1883.229,6266.704)
	(1892.250,6307.967)
	(1901.061,6352.173)
	(1909.689,6398.901)
	(1918.164,6447.731)
	(1926.514,6498.242)
	(1934.770,6550.013)
	(1942.959,6602.623)
	(1951.110,6655.654)
	(1959.253,6708.683)
	(1967.416,6761.290)
	(1975.628,6813.055)
	(1983.918,6863.557)
	(1992.314,6912.376)
	(2000.847,6959.091)
	(2009.545,7003.281)
	(2018.436,7044.526)
	(2036.914,7116.499)
	(2056.513,7171.646)
	(2100.000,7218.000)

\path(2100,7218)	(2158.951,7157.580)
	(2170.951,7123.249)
	(2181.965,7080.403)
	(2192.049,7028.630)
	(2201.258,6967.519)
	(2209.645,6896.657)
	(2213.549,6857.441)
	(2217.268,6815.632)
	(2220.809,6771.180)
	(2224.179,6724.033)
	(2227.386,6674.140)
	(2230.435,6621.448)
	(2233.335,6565.906)
	(2236.090,6507.464)
	(2238.710,6446.068)
	(2241.200,6381.669)
	(2243.567,6314.214)
	(2245.818,6243.652)
	(2247.960,6169.931)
	(2248.992,6131.870)
	(2250.000,6093.000)

\put(6450,2118){\ellipse{150}{150}}
\path(6300,1818)(6000,1218)(7200,1218)
	(7500,1818)(6300,1818)
\path(6675,1743)(6150,1293)
\path(6450,2043)	(6452.400,1971.173)
	(6454.878,1902.344)
	(6457.443,1836.465)
	(6460.100,1773.487)
	(6462.856,1713.363)
	(6465.719,1656.044)
	(6468.696,1601.482)
	(6471.792,1549.630)
	(6475.016,1500.439)
	(6478.373,1453.862)
	(6481.871,1409.849)
	(6485.517,1368.353)
	(6489.317,1329.327)
	(6493.278,1292.721)
	(6501.713,1226.580)
	(6510.875,1169.546)
	(6520.819,1121.234)
	(6543.277,1049.240)
	(6600.000,993.000)

\path(6600,993)	(6644.458,1039.489)
	(6664.439,1094.688)
	(6683.229,1166.704)
	(6692.250,1207.967)
	(6701.061,1252.173)
	(6709.689,1298.901)
	(6718.164,1347.731)
	(6726.514,1398.242)
	(6734.770,1450.013)
	(6742.959,1502.623)
	(6751.110,1555.654)
	(6759.253,1608.683)
	(6767.416,1661.290)
	(6775.628,1713.055)
	(6783.918,1763.557)
	(6792.314,1812.376)
	(6800.847,1859.091)
	(6809.545,1903.281)
	(6818.436,1944.526)
	(6836.914,2016.499)
	(6856.513,2071.646)
	(6900.000,2118.000)

\path(6900,2118)	(6958.951,2057.580)
	(6970.951,2023.249)
	(6981.965,1980.403)
	(6992.049,1928.630)
	(7001.257,1867.519)
	(7009.645,1796.657)
	(7013.549,1757.441)
	(7017.268,1715.632)
	(7020.809,1671.180)
	(7024.179,1624.033)
	(7027.386,1574.140)
	(7030.435,1521.448)
	(7033.335,1465.906)
	(7036.090,1407.464)
	(7038.710,1346.068)
	(7041.200,1281.669)
	(7043.567,1214.214)
	(7045.818,1143.652)
	(7047.960,1069.931)
	(7048.992,1031.870)
	(7050.000,993.000)

\put(6750,3318){\ellipse{150}{150}}
\path(6600,3018)(6300,2418)(7500,2418)
	(7800,3018)(6600,3018)
\path(6975,2943)(6450,2493)
\path(6750,3243)	(6752.400,3171.173)
	(6754.878,3102.344)
	(6757.443,3036.465)
	(6760.100,2973.487)
	(6762.856,2913.363)
	(6765.719,2856.044)
	(6768.696,2801.482)
	(6771.792,2749.630)
	(6775.016,2700.439)
	(6778.373,2653.862)
	(6781.871,2609.849)
	(6785.517,2568.353)
	(6789.317,2529.327)
	(6793.278,2492.721)
	(6801.713,2426.580)
	(6810.875,2369.546)
	(6820.819,2321.234)
	(6843.277,2249.240)
	(6900.000,2193.000)

\path(6900,2193)	(6944.458,2239.489)
	(6964.439,2294.688)
	(6983.229,2366.704)
	(6992.250,2407.967)
	(7001.061,2452.173)
	(7009.689,2498.901)
	(7018.164,2547.731)
	(7026.514,2598.242)
	(7034.770,2650.013)
	(7042.959,2702.623)
	(7051.110,2755.654)
	(7059.253,2808.683)
	(7067.416,2861.290)
	(7075.628,2913.055)
	(7083.918,2963.557)
	(7092.314,3012.376)
	(7100.847,3059.091)
	(7109.545,3103.281)
	(7118.436,3144.526)
	(7136.914,3216.499)
	(7156.513,3271.646)
	(7200.000,3318.000)

\path(7200,3318)	(7258.951,3257.580)
	(7270.951,3223.249)
	(7281.965,3180.403)
	(7292.049,3128.630)
	(7301.257,3067.519)
	(7309.645,2996.657)
	(7313.549,2957.441)
	(7317.268,2915.632)
	(7320.809,2871.180)
	(7324.179,2824.033)
	(7327.386,2774.140)
	(7330.435,2721.448)
	(7333.335,2665.906)
	(7336.090,2607.464)
	(7338.710,2546.068)
	(7341.200,2481.669)
	(7343.567,2414.214)
	(7345.818,2343.652)
	(7347.960,2269.931)
	(7348.992,2231.870)
	(7350.000,2193.000)

\path(6900,4218)(6600,3618)(7800,3618)
	(8100,4218)(6900,4218)
\path(7125,3918)(7875,4068)
\path(6975,4443)	(6976.383,4381.937)
	(6977.911,4323.410)
	(6979.591,4267.380)
	(6981.428,4213.804)
	(6983.431,4162.642)
	(6985.605,4113.851)
	(6987.959,4067.392)
	(6990.498,4023.223)
	(6993.229,3981.303)
	(6996.161,3941.589)
	(7002.649,3868.621)
	(7010.017,3803.987)
	(7018.320,3747.360)
	(7027.613,3698.409)
	(7037.951,3656.804)
	(7061.982,3594.317)
	(7125.000,3543.000)

\path(7125,3543)	(7178.589,3581.213)
	(7203.547,3628.770)
	(7227.510,3691.295)
	(7250.649,3765.792)
	(7261.964,3806.593)
	(7273.137,3849.263)
	(7284.190,3893.426)
	(7295.145,3938.709)
	(7306.024,3984.737)
	(7316.846,4031.134)
	(7327.635,4077.526)
	(7338.411,4123.538)
	(7349.196,4168.797)
	(7360.012,4212.925)
	(7370.880,4255.551)
	(7381.821,4296.297)
	(7404.010,4370.655)
	(7426.750,4433.003)
	(7450.214,4480.341)
	(7500.000,4518.000)

\path(7500,4518)	(7541.364,4498.528)
	(7575.387,4448.794)
	(7589.784,4411.479)
	(7602.510,4365.281)
	(7613.622,4309.760)
	(7623.172,4244.475)
	(7627.380,4208.034)
	(7631.218,4168.988)
	(7634.693,4127.281)
	(7637.813,4082.859)
	(7640.584,4035.667)
	(7643.013,3985.649)
	(7645.106,3932.752)
	(7646.872,3876.919)
	(7648.315,3818.096)
	(7649.445,3756.228)
	(7650.267,3691.260)
	(7650.787,3623.137)
	(7651.014,3551.805)
	(7650.954,3477.208)
	(7650.819,3438.668)
	(7650.614,3399.291)
	(7650.341,3359.071)
	(7650.000,3318.000)

\put(7575,9168){\ellipse{150}{150}}
\path(7500,9318)(7200,8718)(8400,8718)
	(8700,9318)(7500,9318)
\path(7575,9093)	(7583.041,9053.311)
	(7591.003,9016.927)
	(7606.912,8953.633)
	(7623.166,8902.239)
	(7640.205,8861.865)
	(7678.396,8810.664)
	(7725.000,8793.000)

\path(7725,8793)	(7773.905,8809.672)
	(7811.535,8859.398)
	(7841.899,8931.305)
	(7855.608,8972.178)
	(7869.004,9014.520)
	(7882.587,9056.970)
	(7896.858,9098.170)
	(7929.470,9171.382)
	(7970.848,9223.283)
	(8025.000,9243.000)

\path(8025,9243)	(8069.268,9237.507)
	(8109.844,9220.148)
	(8147.607,9189.606)
	(8183.438,9144.562)
	(8218.213,9083.698)
	(8235.480,9046.922)
	(8252.812,9005.695)
	(8270.321,8959.855)
	(8288.115,8909.235)
	(8306.305,8853.672)
	(8325.000,8793.000)

\path(8025,9618)(8025,8418)
\put(7349.544,10221.904){\arc{892.289}{2.9979}{6.3930}}
\put(7351.851,10225.671){\arc{916.319}{0.4951}{2.6638}}
\put(7350,10218){\ellipse{1050}{450}}
\path(6900,10518)(6600,9918)(7800,9918)
	(8100,10518)(6900,10518)
\put(3075,768){\ellipse{150}{150}}
\path(3000,918)(2700,318)(3900,318)
	(4200,918)(3000,918)
\path(3075,693)	(3083.041,653.311)
	(3091.003,616.927)
	(3106.912,553.633)
	(3123.166,502.239)
	(3140.205,461.865)
	(3178.396,410.664)
	(3225.000,393.000)

\path(3225,393)	(3273.905,409.672)
	(3311.535,459.397)
	(3341.899,531.305)
	(3355.608,572.178)
	(3369.004,614.520)
	(3382.587,656.970)
	(3396.858,698.170)
	(3429.470,771.382)
	(3470.848,823.283)
	(3525.000,843.000)

\path(3525,843)	(3569.268,837.507)
	(3609.844,820.148)
	(3647.607,789.606)
	(3683.438,744.562)
	(3718.213,683.698)
	(3735.480,646.922)
	(3752.812,605.695)
	(3770.321,559.855)
	(3788.115,509.235)
	(3806.305,453.672)
	(3825.000,393.000)

\path(3525,1218)(3525,18)
\put(2549.544,1521.903){\arc{892.289}{2.9979}{6.3930}}
\put(2551.851,1525.671){\arc{916.319}{0.4951}{2.6638}}
\put(2550,1518){\ellipse{1050}{450}}
\path(2100,1818)(1800,1218)(3000,1218)
	(3300,1818)(2100,1818)
\put(2835,1705){\blacken\ellipse{74}{74}}
\put(2835,1705){\ellipse{74}{74}}
\put(2265,1330){\blacken\ellipse{74}{74}}
\put(2265,1330){\ellipse{74}{74}}
\put(2827,1332){\blacken\ellipse{74}{74}}
\put(2827,1332){\ellipse{74}{74}}
\put(6975,7218){\ellipse{150}{150}}
\put(4845,8980){\blacken\ellipse{74}{74}}
\put(4845,8980){\ellipse{74}{74}}
\put(5055,9108){\blacken\ellipse{74}{74}}
\put(5055,9108){\ellipse{74}{74}}
\put(5063,11793){\blacken\ellipse{74}{74}}
\put(5063,11793){\ellipse{74}{74}}
\put(4845,11643){\blacken\ellipse{74}{74}}
\put(4845,11643){\ellipse{74}{74}}
\put(4560,11838){\blacken\ellipse{74}{74}}
\put(4560,11838){\ellipse{74}{74}}
\put(6023,11163){\blacken\ellipse{74}{74}}
\put(6023,11163){\ellipse{74}{74}}
\put(6195,10923){\blacken\ellipse{74}{74}}
\put(6195,10923){\ellipse{74}{74}}
\put(7620,8950){\blacken\ellipse{74}{74}}
\put(7620,8950){\ellipse{74}{74}}
\put(7875,9018){\blacken\ellipse{74}{74}}
\put(7875,9018){\ellipse{74}{74}}
\put(7350,6760){\blacken\ellipse{74}{74}}
\put(7350,6760){\ellipse{74}{74}}
\put(7350,6520){\blacken\ellipse{74}{74}}
\put(7350,6520){\ellipse{74}{74}}
\put(5397,6162){\blacken\ellipse{74}{74}}
\put(5397,6162){\ellipse{74}{74}}
\put(5855,6019){\blacken\ellipse{74}{74}}
\put(5855,6019){\ellipse{74}{74}}
\put(3290,6162){\blacken\ellipse{74}{74}}
\put(3290,6162){\ellipse{74}{74}}
\put(3293,5913){\blacken\ellipse{74}{74}}
\put(3293,5913){\ellipse{74}{74}}
\put(2228,6640){\blacken\ellipse{74}{74}}
\put(2228,6640){\ellipse{74}{74}}
\put(1200,7735){\blacken\ellipse{74}{74}}
\put(1200,7735){\ellipse{74}{74}}
\put(1920,7870){\blacken\ellipse{74}{74}}
\put(1920,7870){\ellipse{74}{74}}
\put(1950,6618){\blacken\ellipse{74}{74}}
\put(1950,6618){\ellipse{74}{74}}
\put(975,8988){\blacken\ellipse{74}{74}}
\put(975,8988){\ellipse{74}{74}}
\put(825,8860){\blacken\ellipse{74}{74}}
\put(825,8860){\ellipse{74}{74}}
\put(1590,10353){\blacken\ellipse{74}{74}}
\put(1590,10353){\ellipse{74}{74}}
\put(2085,10293){\blacken\ellipse{74}{74}}
\put(2085,10293){\ellipse{74}{74}}
\put(2963,11186){\blacken\ellipse{74}{74}}
\put(2963,11186){\ellipse{74}{74}}
\put(3248,11246){\blacken\ellipse{74}{74}}
\put(3248,11246){\ellipse{74}{74}}
\put(5070,4405){\blacken\ellipse{74}{74}}
\put(5070,4405){\ellipse{74}{74}}
\put(4845,4150){\blacken\ellipse{74}{74}}
\put(4845,4150){\ellipse{74}{74}}
\put(4545,4248){\blacken\ellipse{74}{74}}
\put(4545,4248){\ellipse{74}{74}}
\put(2363,4075){\blacken\ellipse{74}{74}}
\put(2363,4075){\ellipse{74}{74}}
\put(2145,3888){\blacken\ellipse{74}{74}}
\put(2145,3888){\ellipse{74}{74}}
\put(1875,4098){\blacken\ellipse{74}{74}}
\put(1875,4098){\ellipse{74}{74}}
\put(1643,3865){\blacken\ellipse{74}{74}}
\put(1643,3865){\ellipse{74}{74}}
\put(2033,2575){\blacken\ellipse{74}{74}}
\put(2033,2575){\ellipse{74}{74}}
\put(2348,2740){\blacken\ellipse{74}{74}}
\put(2348,2740){\ellipse{74}{74}}
\put(2633,2568){\blacken\ellipse{74}{74}}
\put(2633,2568){\ellipse{74}{74}}
\put(3090,610){\blacken\ellipse{74}{74}}
\put(3090,610){\ellipse{74}{74}}
\put(3270,408){\blacken\ellipse{74}{74}}
\put(3270,408){\ellipse{74}{74}}
\put(3390,678){\blacken\ellipse{74}{74}}
\put(3390,678){\ellipse{74}{74}}
\put(5850,610){\blacken\ellipse{74}{74}}
\put(5850,610){\ellipse{74}{74}}
\put(5393,528){\blacken\ellipse{74}{74}}
\put(5393,528){\ellipse{74}{74}}
\put(5400,738){\blacken\ellipse{74}{74}}
\put(5400,738){\ellipse{74}{74}}
\put(6278,1398){\blacken\ellipse{74}{74}}
\put(6278,1398){\ellipse{74}{74}}
\put(6758,1533){\blacken\ellipse{74}{74}}
\put(6758,1533){\ellipse{74}{74}}
\put(7035,1398){\blacken\ellipse{74}{74}}
\put(7035,1398){\ellipse{74}{74}}
\put(6585,2613){\blacken\ellipse{74}{74}}
\put(6585,2613){\ellipse{74}{74}}
\put(6855,2830){\blacken\ellipse{74}{74}}
\put(6855,2830){\ellipse{74}{74}}
\put(7335,2628){\blacken\ellipse{74}{74}}
\put(7335,2628){\ellipse{74}{74}}
\put(6998,3985){\blacken\ellipse{74}{74}}
\put(6998,3985){\ellipse{74}{74}}
\put(7440,3993){\blacken\ellipse{74}{74}}
\put(7440,3993){\ellipse{74}{74}}
\put(7208,3933){\blacken\ellipse{74}{74}}
\put(7208,3933){\ellipse{74}{74}}
\put(7080,10023){\blacken\ellipse{74}{74}}
\put(7080,10023){\ellipse{74}{74}}
\put(7590,10023){\blacken\ellipse{74}{74}}
\put(7590,10023){\ellipse{74}{74}}
\path(7350,6843)(7350,6393)
\path(6900,6918)(6600,6318)(7800,6318)
	(8100,6918)(6900,6918)
\dottedline{135}(4463,9273)(4073,8838)
\dottedline{135}(4410,9093)(4725,8838)
\dottedline{135}(4440,11965)(4095,11553)
\dottedline{135}(6135,11328)(6360,10915)
\dottedline{135}(6458,11050)(6810,11335)
\dottedline{135}(8093,9055)(8468,9258)
\dottedline{135}(8108,8793)(8453,9078)
\dottedline{135}(3435,5860)(4028,6228)
\dottedline{135}(1485,7623)(1800,8035)
\dottedline{135}(1193,9025)(1500,8808)
\dottedline{135}(1493,9093)(1928,9213)
\dottedline{135}(2625,11343)(2318,10938)
\dottedline{135}(4440,4435)(4103,4038)
\dottedline{135}(2055,2943)(1680,2530)
\dottedline{135}(3600,460)(3983,813)
\dottedline{135}(6893,2568)(7200,2943)
\path(3945,3543)(3105,3108)
\path(3945,3543)(3105,3108)
\path(3197.764,3189.822)(3105.000,3108.000)(3225.355,3136.543)
\path(4200,3318)(3225,2013)
\path(4200,3318)(3225,2013)
\path(3272.790,2127.088)(3225.000,2013.000)(3320.856,2091.177)
\path(4470,3153)(3750,1278)
\path(4470,3153)(3750,1278)
\path(3765.011,1400.779)(3750.000,1278.000)(3821.024,1379.270)
\path(4890,3198)(5400,1353)
\path(4890,3198)(5400,1353)
\path(5339.113,1460.670)(5400.000,1353.000)(5396.944,1476.655)
\path(5175,3363)(6120,1998)
\path(5175,3363)(6120,1998)
\path(6027.029,2079.587)(6120.000,1998.000)(6076.361,2113.739)
\path(5415,3573)(6390,3003)
\path(5415,3573)(6390,3003)
\path(6271.263,3037.665)(6390.000,3003.000)(6301.545,3089.462)
\path(5535,3933)(6585,3948)
\path(5535,3933)(6585,3948)
\path(6465.441,3916.289)(6585.000,3948.000)(6464.584,3976.283)
\path(4650,9768)(4650,11268)
\path(4650,9768)(4650,11268)
\path(4680.000,11148.000)(4650.000,11268.000)(4620.000,11148.000)
\path(5250,9693)(6150,10668)
\path(5250,9693)(6150,10668)
\path(6090.650,10559.475)(6150.000,10668.000)(6046.562,10600.172)
\path(5625,9318)(6825,9843)
\path(5625,9318)(6825,9843)
\path(6727.086,9767.417)(6825.000,9843.000)(6703.037,9822.386)
\path(5925,9018)(7200,9018)
\path(5925,9018)(7200,9018)
\path(7080.000,8988.000)(7200.000,9018.000)(7080.000,9048.000)
\path(5400,8643)(7050,7968)
\path(5400,8643)(7050,7968)
\path(6927.575,7985.670)(7050.000,7968.000)(6950.293,8041.202)
\path(5175,8418)(6675,6993)
\path(5175,8418)(6675,6993)
\path(6567.338,7053.900)(6675.000,6993.000)(6608.663,7097.400)
\path(4875,8343)(5400,6693)
\path(4875,8343)(5400,6693)
\path(5335.028,6798.255)(5400.000,6693.000)(5392.203,6816.447)
\path(4500,8343)(3825,6618)
\path(4500,8343)(3825,6618)
\path(3840.791,6740.681)(3825.000,6618.000)(3896.665,6718.817)
\path(3975,8343)(2400,7068)
\path(3975,8343)(2400,7068)
\path(2474.393,7166.821)(2400.000,7068.000)(2512.145,7120.186)
\path(3600,8793)(2475,8193)
\path(3600,8793)(2475,8193)
\path(2566.765,8275.941)(2475.000,8193.000)(2595.000,8223.000)
\path(3750,9093)(2175,9093)
\path(3750,9093)(2175,9093)
\path(2295.000,9123.000)(2175.000,9093.000)(2295.000,9063.000)
\path(3975,9318)(2775,9993)
\path(3975,9318)(2775,9993)
\path(2894.297,9960.316)(2775.000,9993.000)(2864.881,9908.021)
\path(4200,9768)(3525,10893)
\path(4200,9768)(3525,10893)
\path(3612.464,10805.536)(3525.000,10893.000)(3561.015,10774.666)
\path(3825,3993)(3150,3918)
\path(3825,3993)(3150,3918)
\path(3265.953,3961.068)(3150.000,3918.000)(3272.579,3901.435)
\path(6975,7143)	(6980.625,7086.103)
	(6986.745,7033.314)
	(6993.429,6984.441)
	(7000.747,6939.293)
	(7008.766,6897.676)
	(7017.557,6859.399)
	(7037.725,6792.094)
	(7061.803,6735.839)
	(7090.341,6689.098)
	(7163.000,6618.000)

\path(7163,6618)	(7202.165,6604.714)
	(7247.994,6605.948)
	(7292.576,6613.208)
	(7328.000,6618.000)

\path(7328,6618)	(7392.762,6631.790)
	(7463.000,6625.000)

\path(7463,6625)	(7509.347,6586.789)
	(7548.034,6541.374)
	(7579.611,6486.976)
	(7604.624,6421.816)
	(7614.840,6384.645)
	(7623.621,6344.117)
	(7631.036,6300.010)
	(7637.151,6252.100)
	(7642.037,6200.167)
	(7645.762,6143.987)
	(7648.393,6083.339)
	(7650.000,6018.000)

\put(3585,8403){\makebox(0,0)[lb]{\smash{{{\SetFigFont{8}{9.6}{rm}312,348}}}}}
\put(5385,8928){\makebox(0,0)[lb]{\smash{{{\SetFigFont{8}{9.6}{rm}{\bf 81}}}}}}
\put(5385,11643){\makebox(0,0)[lb]{\smash{{{\SetFigFont{8}{9.6}{rm}{\bf 93}}}}}}
\put(7335,11418){\makebox(0,0)[lb]{\smash{{{\SetFigFont{8}{9.6}{rm}2,522}}}}}
\put(8220,10233){\makebox(0,0)[lb]{\smash{{{\SetFigFont{8}{9.6}{rm}12,924}}}}}
\put(8805,9018){\makebox(0,0)[lb]{\smash{{{\SetFigFont{8}{9.6}{rm}10,404}}}}}
\put(8340,7503){\makebox(0,0)[lb]{\smash{{{\SetFigFont{8}{9.6}{rm}1,344}}}}}
\put(7830,5988){\makebox(0,0)[lb]{\smash{{{\SetFigFont{8}{9.6}{rm}921}}}}}
\put(6090,5418){\makebox(0,0)[lb]{\smash{{{\SetFigFont{8}{9.6}{rm}1,409}}}}}
\put(2670,5433){\makebox(0,0)[lb]{\smash{{{\SetFigFont{8}{9.6}{rm}989}}}}}
\put(975,6048){\makebox(0,0)[lb]{\smash{{{\SetFigFont{8}{9.6}{rm}522}}}}}
\put(0,9003){\makebox(0,0)[lb]{\smash{{{\SetFigFont{8}{9.6}{rm}1,023}}}}}
\put(5385,4158){\makebox(0,0)[lb]{\smash{{{\SetFigFont{8}{9.6}{rm}{\bf 93}}}}}}
\put(2700,3858){\makebox(0,0)[lb]{\smash{{{\SetFigFont{8}{9.6}{rm}{\bf 98}}}}}}
\put(3900,3633){\makebox(0,0)[lb]{\smash{{{\SetFigFont{8}{9.6}{rm}31,056}}}}}
\put(1185,3303){\makebox(0,0)[lb]{\smash{{{\SetFigFont{8}{9.6}{rm}2,519}}}}}
\put(1650,2118){\makebox(0,0)[lb]{\smash{{{\SetFigFont{8}{9.6}{rm}464}}}}}
\put(2700,18){\makebox(0,0)[lb]{\smash{{{\SetFigFont{8}{9.6}{rm}1,764}}}}}
\put(5040,33){\makebox(0,0)[lb]{\smash{{{\SetFigFont{8}{9.6}{rm}127}}}}}
\put(7200,948){\makebox(0,0)[lb]{\smash{{{\SetFigFont{8}{9.6}{rm}102}}}}}
\put(7500,2118){\makebox(0,0)[lb]{\smash{{{\SetFigFont{8}{9.6}{rm}127}}}}}
\put(7815,3333){\makebox(0,0)[lb]{\smash{{{\SetFigFont{8}{9.6}{rm}1,820}}}}}
\put(6525,9093){\makebox(0,0)[lb]{\smash{{{\SetFigFont{8}{9.6}{rm}$\times 3$}}}}}
\put(5550,10518){\makebox(0,0)[lb]{\smash{{{\SetFigFont{8}{9.6}{rm}$\times 4$}}}}}
\put(4725,10893){\makebox(0,0)[lb]{\smash{{{\SetFigFont{8}{9.6}{rm}$\times 2$}}}}}
\put(3525,6993){\makebox(0,0)[lb]{\smash{{{\SetFigFont{8}{9.6}{rm}$\times 4$}}}}}
\put(2400,8418){\makebox(0,0)[lb]{\smash{{{\SetFigFont{8}{9.6}{rm}$\times 4$}}}}}
\put(3000,9918){\makebox(0,0)[lb]{\smash{{{\SetFigFont{8}{9.6}{rm}$\times 2$}}}}}
\put(3750,10593){\makebox(0,0)[lb]{\smash{{{\SetFigFont{8}{9.6}{rm}$\times 2$}}}}}
\put(3975,1443){\makebox(0,0)[lb]{\smash{{{\SetFigFont{8}{9.6}{rm}$\times 3$}}}}}
\put(6000,3243){\makebox(0,0)[lb]{\smash{{{\SetFigFont{8}{9.6}{rm}$\times 3$}}}}}
\put(6075,4068){\makebox(0,0)[lb]{\smash{{{\SetFigFont{8}{9.6}{rm}$\times 3$}}}}}
\put(5295,12318){\makebox(0,0)[lb]{\smash{{{\SetFigFont{8}{9.6}{rm}31,056}}}}}
\put(5925,9618){\makebox(0,0)[lb]{\smash{{{\SetFigFont{8}{9.6}{rm}$\times 8$}}}}}
\put(4800,6993){\makebox(0,0)[lb]{\smash{{{\SetFigFont{8}{9.6}{rm}$\times 4$}}}}}
\put(1515,11433){\makebox(0,0)[lb]{\smash{{{\SetFigFont{8}{9.6}{rm}20,312}}}}}
\put(2925,3243){\makebox(0,0)[lb]{\smash{{{\SetFigFont{8}{9.6}{rm}$\times 4$}}}}}
\put(1650,918){\makebox(0,0)[lb]{\smash{{{\SetFigFont{8}{9.6}{rm}3,620}}}}}
\put(4800,1743){\makebox(0,0)[lb]{\smash{{{\SetFigFont{8}{9.6}{rm}$\times 6$}}}}}
\put(420,10233){\makebox(0,0)[lb]{\smash{{{\SetFigFont{8}{9.6}{rm}22,252}}}}}
\put(135,7518){\makebox(0,0)[lb]{\smash{{{\SetFigFont{8}{9.6}{rm}1,293}}}}}
\put(5850,7143){\makebox(0,0)[lb]{\smash{{{\SetFigFont{8}{9.6}{rm}$\times 3$}}}}}
\put(5475,2118){\makebox(0,0)[lb]{\smash{{{\SetFigFont{8}{9.6}{rm}$\times 3$}}}}}
\put(3450,2118){\makebox(0,0)[lb]{\smash{{{\SetFigFont{8}{9.6}{rm}$\times 4$}}}}}
\end{picture}

%% file: egenus3.tex
\begingroup\makeatletter\ifx\SetFigFont\undefined
\def\x#1#2#3#4#5#6#7\relax{\def\x{#1#2#3#4#5#6}}%
\expandafter\x\fmtname xxxxxx\relax \def\y{splain}%
\ifx\x\y   
\gdef\SetFigFont#1#2#3{%
  \ifnum #1<17\tiny\else \ifnum #1<20\small\else
  \ifnum #1<24\normalsize\else \ifnum #1<29\large\else
  \ifnum #1<34\Large\else \ifnum #1<41\LARGE\else
     \huge\fi\fi\fi\fi\fi\fi
  \csname #3\endcsname}%
\else
\gdef\SetFigFont#1#2#3{\begingroup
  \count@#1\relax \ifnum 25<\count@\count@25\fi
  \def\x{\endgroup\@setsize\SetFigFont{#2pt}}%
  \expandafter\x
    \csname \romannumeral\the\count@ pt\expandafter\endcsname
    \csname @\romannumeral\the\count@ pt\endcsname
  \csname #3\endcsname}%
\fi
\fi\endgroup
\begin{picture}(15012,14274)(0,-10)
\thicklines
\put(10650,6465){\ellipse{150}{150}}
\put(10650,6615){\ellipse{150}{150}}
\put(10650,6765){\ellipse{150}{150}}
\path(10650,6390)	(10649.401,6328.543)
	(10648.957,6269.649)
	(10648.672,6213.276)
	(10648.550,6159.383)
	(10648.593,6107.929)
	(10648.806,6058.873)
	(10649.191,6012.172)
	(10649.753,5967.788)
	(10650.494,5925.677)
	(10651.418,5885.799)
	(10653.829,5812.576)
	(10657.014,5747.791)
	(10660.999,5691.112)
	(10665.812,5642.212)
	(10671.481,5600.760)
	(10685.495,5538.881)
	(10725.000,5490.000)

\path(10725,5490)	(10763.407,5535.771)
	(10781.068,5590.672)
	(10797.998,5662.434)
	(10806.257,5703.588)
	(10814.417,5747.696)
	(10822.505,5794.339)
	(10830.550,5843.095)
	(10838.578,5893.546)
	(10846.617,5945.270)
	(10854.696,5997.847)
	(10862.842,6050.858)
	(10871.084,6103.881)
	(10879.448,6156.497)
	(10887.962,6208.285)
	(10896.655,6258.825)
	(10905.554,6307.697)
	(10914.687,6354.481)
	(10924.082,6398.756)
	(10933.766,6440.102)
	(10954.115,6512.326)
	(10975.955,6567.792)
	(11025.000,6615.000)

\path(11025,6615)	(11073.026,6561.706)
	(11091.996,6501.915)
	(11100.416,6465.267)
	(11108.201,6424.702)
	(11115.406,6380.657)
	(11122.088,6333.571)
	(11128.303,6283.881)
	(11134.107,6232.027)
	(11139.555,6178.446)
	(11144.703,6123.576)
	(11149.608,6067.856)
	(11154.326,6011.722)
	(11158.913,5955.615)
	(11163.423,5899.971)
	(11167.915,5845.229)
	(11172.443,5791.828)
	(11177.063,5740.204)
	(11181.832,5690.796)
	(11186.805,5644.043)
	(11192.039,5600.382)
	(11197.589,5560.252)
	(11203.511,5524.091)
	(11216.697,5465.426)
	(11250.000,5415.000)

\path(11250,5415)	(11300.879,5477.259)
	(11311.702,5513.561)
	(11321.875,5559.035)
	(11331.455,5614.123)
	(11340.495,5679.263)
	(11344.830,5715.739)
	(11349.051,5754.894)
	(11353.165,5796.782)
	(11357.178,5841.457)
	(11361.098,5888.975)
	(11364.930,5939.391)
	(11368.683,5992.759)
	(11372.363,6049.134)
	(11375.977,6108.572)
	(11379.532,6171.128)
	(11383.034,6236.855)
	(11386.491,6305.810)
	(11389.909,6378.047)
	(11391.606,6415.414)
	(11393.295,6453.621)
	(11394.979,6492.677)
	(11396.657,6532.587)
	(11398.330,6573.359)
	(11400.000,6615.000)

\path(10515,6318)(10215,5718)(11415,5718)
	(11715,6318)(10515,6318)
\put(8235,6465){\ellipse{150}{150}}
\put(8235,6615){\ellipse{150}{150}}
\put(8235,6765){\ellipse{150}{150}}
\path(8235,6390)	(8234.401,6328.543)
	(8233.957,6269.649)
	(8233.672,6213.276)
	(8233.550,6159.383)
	(8233.593,6107.929)
	(8233.806,6058.873)
	(8234.191,6012.172)
	(8234.753,5967.788)
	(8235.494,5925.677)
	(8236.418,5885.799)
	(8238.829,5812.576)
	(8242.014,5747.791)
	(8245.999,5691.112)
	(8250.812,5642.212)
	(8256.481,5600.760)
	(8270.495,5538.881)
	(8310.000,5490.000)

\path(8310,5490)	(8348.407,5535.771)
	(8366.068,5590.672)
	(8382.998,5662.434)
	(8391.257,5703.588)
	(8399.417,5747.696)
	(8407.505,5794.339)
	(8415.550,5843.095)
	(8423.578,5893.546)
	(8431.617,5945.270)
	(8439.696,5997.847)
	(8447.842,6050.858)
	(8456.084,6103.881)
	(8464.448,6156.497)
	(8472.962,6208.285)
	(8481.655,6258.825)
	(8490.554,6307.697)
	(8499.687,6354.481)
	(8509.082,6398.756)
	(8518.766,6440.102)
	(8539.115,6512.326)
	(8560.955,6567.792)
	(8610.000,6615.000)

\path(8610,6615)	(8658.026,6561.706)
	(8676.996,6501.915)
	(8685.416,6465.267)
	(8693.201,6424.702)
	(8700.406,6380.657)
	(8707.088,6333.571)
	(8713.303,6283.881)
	(8719.107,6232.027)
	(8724.555,6178.446)
	(8729.703,6123.576)
	(8734.608,6067.856)
	(8739.326,6011.722)
	(8743.913,5955.615)
	(8748.423,5899.971)
	(8752.915,5845.229)
	(8757.443,5791.828)
	(8762.063,5740.204)
	(8766.832,5690.796)
	(8771.805,5644.043)
	(8777.039,5600.382)
	(8782.589,5560.252)
	(8788.511,5524.091)
	(8801.697,5465.426)
	(8835.000,5415.000)

\path(8835,5415)	(8885.879,5477.259)
	(8896.702,5513.561)
	(8906.875,5559.035)
	(8916.455,5614.123)
	(8925.495,5679.263)
	(8929.830,5715.739)
	(8934.051,5754.894)
	(8938.165,5796.782)
	(8942.178,5841.457)
	(8946.098,5888.975)
	(8949.930,5939.391)
	(8953.683,5992.759)
	(8957.363,6049.134)
	(8960.977,6108.572)
	(8964.532,6171.128)
	(8968.034,6236.855)
	(8971.491,6305.810)
	(8974.909,6378.047)
	(8976.606,6415.414)
	(8978.295,6453.621)
	(8979.979,6492.677)
	(8981.657,6532.587)
	(8983.330,6573.359)
	(8985.000,6615.000)

\path(8100,6318)(7800,5718)(9000,5718)
	(9300,6318)(8100,6318)
\put(8235,8265){\ellipse{150}{150}}
\put(8235,8415){\ellipse{150}{150}}
\put(8235,8565){\ellipse{150}{150}}
\path(8235,8190)	(8234.401,8128.543)
	(8233.957,8069.649)
	(8233.672,8013.276)
	(8233.550,7959.383)
	(8233.593,7907.929)
	(8233.806,7858.873)
	(8234.191,7812.172)
	(8234.753,7767.788)
	(8235.494,7725.677)
	(8236.418,7685.799)
	(8238.829,7612.576)
	(8242.014,7547.791)
	(8245.999,7491.112)
	(8250.812,7442.212)
	(8256.481,7400.760)
	(8270.495,7338.881)
	(8310.000,7290.000)

\path(8310,7290)	(8348.407,7335.771)
	(8366.068,7390.672)
	(8382.998,7462.434)
	(8391.257,7503.588)
	(8399.417,7547.696)
	(8407.505,7594.339)
	(8415.550,7643.095)
	(8423.578,7693.546)
	(8431.617,7745.270)
	(8439.696,7797.847)
	(8447.842,7850.858)
	(8456.084,7903.881)
	(8464.448,7956.497)
	(8472.962,8008.285)
	(8481.655,8058.825)
	(8490.554,8107.697)
	(8499.687,8154.481)
	(8509.082,8198.756)
	(8518.766,8240.102)
	(8539.115,8312.326)
	(8560.955,8367.792)
	(8610.000,8415.000)

\path(8610,8415)	(8658.026,8361.706)
	(8676.996,8301.915)
	(8685.416,8265.267)
	(8693.201,8224.702)
	(8700.406,8180.657)
	(8707.088,8133.571)
	(8713.303,8083.881)
	(8719.107,8032.027)
	(8724.555,7978.446)
	(8729.703,7923.576)
	(8734.608,7867.856)
	(8739.326,7811.722)
	(8743.913,7755.615)
	(8748.423,7699.971)
	(8752.915,7645.229)
	(8757.443,7591.828)
	(8762.063,7540.204)
	(8766.832,7490.796)
	(8771.805,7444.043)
	(8777.039,7400.382)
	(8782.589,7360.252)
	(8788.511,7324.091)
	(8801.697,7265.426)
	(8835.000,7215.000)

\path(8835,7215)	(8885.879,7277.259)
	(8896.702,7313.561)
	(8906.875,7359.035)
	(8916.455,7414.123)
	(8925.495,7479.263)
	(8929.830,7515.739)
	(8934.051,7554.894)
	(8938.165,7596.782)
	(8942.178,7641.457)
	(8946.098,7688.975)
	(8949.930,7739.391)
	(8953.683,7792.759)
	(8957.363,7849.134)
	(8960.977,7908.572)
	(8964.532,7971.128)
	(8968.034,8036.855)
	(8971.491,8105.810)
	(8974.909,8178.047)
	(8976.606,8215.414)
	(8978.295,8253.621)
	(8979.979,8292.677)
	(8981.657,8332.587)
	(8983.330,8373.359)
	(8985.000,8415.000)

\path(8100,8118)(7800,7518)(9000,7518)
	(9300,8118)(8100,8118)
\put(10635,8265){\ellipse{150}{150}}
\put(10635,8415){\ellipse{150}{150}}
\put(10635,8565){\ellipse{150}{150}}
\path(10635,8190)	(10634.401,8128.543)
	(10633.957,8069.649)
	(10633.672,8013.276)
	(10633.550,7959.383)
	(10633.593,7907.929)
	(10633.806,7858.873)
	(10634.191,7812.172)
	(10634.753,7767.788)
	(10635.494,7725.677)
	(10636.418,7685.799)
	(10638.829,7612.576)
	(10642.014,7547.791)
	(10645.999,7491.112)
	(10650.812,7442.212)
	(10656.481,7400.760)
	(10670.495,7338.881)
	(10710.000,7290.000)

\path(10710,7290)	(10748.407,7335.771)
	(10766.068,7390.672)
	(10782.998,7462.434)
	(10791.257,7503.588)
	(10799.417,7547.696)
	(10807.505,7594.339)
	(10815.550,7643.095)
	(10823.578,7693.546)
	(10831.617,7745.270)
	(10839.696,7797.847)
	(10847.842,7850.858)
	(10856.084,7903.881)
	(10864.448,7956.497)
	(10872.962,8008.285)
	(10881.655,8058.825)
	(10890.554,8107.697)
	(10899.687,8154.481)
	(10909.082,8198.756)
	(10918.766,8240.102)
	(10939.115,8312.326)
	(10960.955,8367.792)
	(11010.000,8415.000)

\path(11010,8415)	(11058.026,8361.706)
	(11076.996,8301.915)
	(11085.416,8265.267)
	(11093.201,8224.702)
	(11100.406,8180.657)
	(11107.088,8133.571)
	(11113.303,8083.881)
	(11119.107,8032.027)
	(11124.555,7978.446)
	(11129.703,7923.576)
	(11134.608,7867.856)
	(11139.326,7811.722)
	(11143.913,7755.615)
	(11148.423,7699.971)
	(11152.915,7645.229)
	(11157.443,7591.828)
	(11162.063,7540.204)
	(11166.832,7490.796)
	(11171.805,7444.043)
	(11177.039,7400.382)
	(11182.589,7360.252)
	(11188.511,7324.091)
	(11201.697,7265.426)
	(11235.000,7215.000)

\path(11235,7215)	(11285.879,7277.259)
	(11296.702,7313.561)
	(11306.875,7359.035)
	(11316.455,7414.123)
	(11325.495,7479.263)
	(11329.830,7515.739)
	(11334.051,7554.894)
	(11338.165,7596.782)
	(11342.178,7641.457)
	(11346.098,7688.975)
	(11349.930,7739.391)
	(11353.683,7792.759)
	(11357.363,7849.134)
	(11360.977,7908.572)
	(11364.532,7971.128)
	(11368.034,8036.855)
	(11371.491,8105.810)
	(11374.909,8178.047)
	(11376.606,8215.414)
	(11378.295,8253.621)
	(11379.979,8292.677)
	(11381.657,8332.587)
	(11383.330,8373.359)
	(11385.000,8415.000)

\path(10500,8118)(10200,7518)(11400,7518)
	(11700,8118)(10500,8118)
\put(9735,10065){\ellipse{150}{150}}
\put(9735,10215){\ellipse{150}{150}}
\put(9735,10365){\ellipse{150}{150}}
\path(9735,9990)	(9734.401,9928.543)
	(9733.957,9869.649)
	(9733.672,9813.276)
	(9733.550,9759.383)
	(9733.593,9707.929)
	(9733.806,9658.873)
	(9734.191,9612.172)
	(9734.753,9567.787)
	(9735.494,9525.677)
	(9736.418,9485.799)
	(9738.829,9412.576)
	(9742.014,9347.791)
	(9745.999,9291.112)
	(9750.812,9242.212)
	(9756.481,9200.760)
	(9770.495,9138.881)
	(9810.000,9090.000)

\path(9810,9090)	(9848.407,9135.771)
	(9866.068,9190.672)
	(9882.998,9262.434)
	(9891.257,9303.588)
	(9899.417,9347.696)
	(9907.505,9394.339)
	(9915.550,9443.095)
	(9923.578,9493.546)
	(9931.617,9545.270)
	(9939.696,9597.847)
	(9947.842,9650.858)
	(9956.084,9703.881)
	(9964.448,9756.497)
	(9972.962,9808.285)
	(9981.655,9858.825)
	(9990.554,9907.697)
	(9999.687,9954.481)
	(10009.082,9998.756)
	(10018.766,10040.102)
	(10039.115,10112.326)
	(10060.955,10167.792)
	(10110.000,10215.000)

\path(10110,10215)	(10158.026,10161.706)
	(10176.996,10101.915)
	(10185.416,10065.267)
	(10193.201,10024.702)
	(10200.406,9980.657)
	(10207.088,9933.571)
	(10213.303,9883.881)
	(10219.107,9832.027)
	(10224.555,9778.446)
	(10229.703,9723.576)
	(10234.608,9667.856)
	(10239.326,9611.722)
	(10243.913,9555.615)
	(10248.423,9499.971)
	(10252.915,9445.229)
	(10257.443,9391.828)
	(10262.063,9340.204)
	(10266.832,9290.796)
	(10271.805,9244.043)
	(10277.039,9200.382)
	(10282.589,9160.252)
	(10288.511,9124.091)
	(10301.697,9065.426)
	(10335.000,9015.000)

\path(10335,9015)	(10385.879,9077.259)
	(10396.702,9113.561)
	(10406.875,9159.035)
	(10416.455,9214.123)
	(10425.495,9279.263)
	(10429.830,9315.739)
	(10434.051,9354.894)
	(10438.165,9396.782)
	(10442.178,9441.457)
	(10446.098,9488.975)
	(10449.930,9539.391)
	(10453.683,9592.759)
	(10457.363,9649.134)
	(10460.977,9708.572)
	(10464.532,9771.128)
	(10468.034,9836.855)
	(10471.491,9905.810)
	(10474.909,9978.047)
	(10476.606,10015.414)
	(10478.295,10053.621)
	(10479.979,10092.677)
	(10481.657,10132.587)
	(10483.330,10173.359)
	(10485.000,10215.000)

\path(9600,9918)(9300,9318)(10500,9318)
	(10800,9918)(9600,9918)
\put(8835,11865){\ellipse{150}{150}}
\put(8835,12015){\ellipse{150}{150}}
\put(8835,12165){\ellipse{150}{150}}
\path(8835,11790)	(8834.401,11728.543)
	(8833.957,11669.649)
	(8833.672,11613.276)
	(8833.550,11559.383)
	(8833.593,11507.929)
	(8833.806,11458.873)
	(8834.191,11412.172)
	(8834.753,11367.787)
	(8835.494,11325.677)
	(8836.418,11285.799)
	(8838.829,11212.576)
	(8842.014,11147.791)
	(8845.999,11091.112)
	(8850.812,11042.212)
	(8856.481,11000.760)
	(8870.495,10938.881)
	(8910.000,10890.000)

\path(8910,10890)	(8948.407,10935.771)
	(8966.068,10990.672)
	(8982.998,11062.434)
	(8991.257,11103.588)
	(8999.417,11147.696)
	(9007.505,11194.339)
	(9015.550,11243.095)
	(9023.578,11293.546)
	(9031.617,11345.270)
	(9039.696,11397.847)
	(9047.842,11450.858)
	(9056.084,11503.881)
	(9064.448,11556.497)
	(9072.962,11608.285)
	(9081.655,11658.825)
	(9090.554,11707.697)
	(9099.687,11754.481)
	(9109.082,11798.756)
	(9118.766,11840.102)
	(9139.115,11912.326)
	(9160.955,11967.792)
	(9210.000,12015.000)

\path(9210,12015)	(9258.026,11961.706)
	(9276.996,11901.915)
	(9285.416,11865.267)
	(9293.201,11824.702)
	(9300.406,11780.657)
	(9307.088,11733.571)
	(9313.303,11683.881)
	(9319.107,11632.027)
	(9324.555,11578.446)
	(9329.703,11523.576)
	(9334.608,11467.856)
	(9339.326,11411.722)
	(9343.913,11355.615)
	(9348.423,11299.971)
	(9352.915,11245.229)
	(9357.443,11191.828)
	(9362.063,11140.204)
	(9366.832,11090.796)
	(9371.805,11044.043)
	(9377.039,11000.382)
	(9382.589,10960.252)
	(9388.511,10924.091)
	(9401.697,10865.426)
	(9435.000,10815.000)

\path(9435,10815)	(9485.879,10877.259)
	(9496.702,10913.561)
	(9506.875,10959.035)
	(9516.455,11014.123)
	(9525.495,11079.263)
	(9529.830,11115.739)
	(9534.051,11154.894)
	(9538.165,11196.782)
	(9542.178,11241.457)
	(9546.098,11288.975)
	(9549.930,11339.391)
	(9553.683,11392.759)
	(9557.363,11449.134)
	(9560.977,11508.572)
	(9564.532,11571.128)
	(9568.034,11636.855)
	(9571.491,11705.810)
	(9574.909,11778.047)
	(9576.606,11815.414)
	(9578.295,11853.621)
	(9579.979,11892.677)
	(9581.657,11932.587)
	(9583.330,11973.359)
	(9585.000,12015.000)

\path(8700,11718)(8400,11118)(9600,11118)
	(9900,11718)(8700,11718)
\put(8835,13665){\ellipse{150}{150}}
\put(8835,13815){\ellipse{150}{150}}
\put(8835,13965){\ellipse{150}{150}}
\path(8835,13590)	(8834.401,13528.543)
	(8833.957,13469.649)
	(8833.672,13413.276)
	(8833.550,13359.383)
	(8833.593,13307.929)
	(8833.806,13258.873)
	(8834.191,13212.172)
	(8834.753,13167.787)
	(8835.494,13125.677)
	(8836.418,13085.799)
	(8838.829,13012.576)
	(8842.014,12947.791)
	(8845.999,12891.112)
	(8850.812,12842.212)
	(8856.481,12800.760)
	(8870.495,12738.881)
	(8910.000,12690.000)

\path(8910,12690)	(8948.407,12735.771)
	(8966.068,12790.672)
	(8982.998,12862.434)
	(8991.257,12903.588)
	(8999.417,12947.696)
	(9007.505,12994.339)
	(9015.550,13043.095)
	(9023.578,13093.546)
	(9031.617,13145.270)
	(9039.696,13197.847)
	(9047.842,13250.858)
	(9056.084,13303.881)
	(9064.448,13356.497)
	(9072.962,13408.285)
	(9081.655,13458.825)
	(9090.554,13507.697)
	(9099.687,13554.481)
	(9109.082,13598.756)
	(9118.766,13640.102)
	(9139.115,13712.326)
	(9160.955,13767.792)
	(9210.000,13815.000)

\path(9210,13815)	(9258.026,13761.706)
	(9276.996,13701.915)
	(9285.416,13665.267)
	(9293.201,13624.702)
	(9300.406,13580.657)
	(9307.088,13533.571)
	(9313.303,13483.881)
	(9319.107,13432.027)
	(9324.555,13378.446)
	(9329.703,13323.576)
	(9334.608,13267.856)
	(9339.326,13211.722)
	(9343.913,13155.615)
	(9348.423,13099.971)
	(9352.915,13045.229)
	(9357.443,12991.828)
	(9362.063,12940.204)
	(9366.832,12890.796)
	(9371.805,12844.043)
	(9377.039,12800.382)
	(9382.589,12760.252)
	(9388.511,12724.091)
	(9401.697,12665.426)
	(9435.000,12615.000)

\path(9435,12615)	(9485.879,12677.259)
	(9496.702,12713.561)
	(9506.875,12759.035)
	(9516.455,12814.123)
	(9525.495,12879.263)
	(9529.830,12915.739)
	(9534.051,12954.894)
	(9538.165,12996.782)
	(9542.178,13041.457)
	(9546.098,13088.975)
	(9549.930,13139.391)
	(9553.683,13192.759)
	(9557.363,13249.134)
	(9560.977,13308.572)
	(9564.532,13371.128)
	(9568.034,13436.855)
	(9571.491,13505.810)
	(9574.909,13578.047)
	(9576.606,13615.414)
	(9578.295,13653.621)
	(9579.979,13692.677)
	(9581.657,13732.587)
	(9583.330,13773.359)
	(9585.000,13815.000)

\path(8700,13518)(8400,12918)(9600,12918)
	(9900,13518)(8700,13518)
\put(5835,6465){\ellipse{150}{150}}
\put(5835,6615){\ellipse{150}{150}}
\put(5835,6765){\ellipse{150}{150}}
\path(5835,6390)	(5834.401,6328.543)
	(5833.957,6269.649)
	(5833.672,6213.276)
	(5833.550,6159.383)
	(5833.593,6107.929)
	(5833.806,6058.873)
	(5834.191,6012.172)
	(5834.753,5967.788)
	(5835.494,5925.677)
	(5836.418,5885.799)
	(5838.829,5812.576)
	(5842.014,5747.791)
	(5845.999,5691.112)
	(5850.812,5642.212)
	(5856.481,5600.760)
	(5870.495,5538.881)
	(5910.000,5490.000)

\path(5910,5490)	(5948.407,5535.771)
	(5966.068,5590.672)
	(5982.997,5662.434)
	(5991.257,5703.588)
	(5999.417,5747.696)
	(6007.505,5794.339)
	(6015.550,5843.095)
	(6023.578,5893.546)
	(6031.617,5945.270)
	(6039.696,5997.847)
	(6047.842,6050.858)
	(6056.084,6103.881)
	(6064.448,6156.497)
	(6072.962,6208.285)
	(6081.655,6258.825)
	(6090.554,6307.697)
	(6099.687,6354.481)
	(6109.082,6398.756)
	(6118.766,6440.102)
	(6139.115,6512.326)
	(6160.955,6567.792)
	(6210.000,6615.000)

\path(6210,6615)	(6258.026,6561.706)
	(6276.996,6501.915)
	(6285.416,6465.267)
	(6293.201,6424.702)
	(6300.406,6380.657)
	(6307.088,6333.571)
	(6313.303,6283.881)
	(6319.107,6232.027)
	(6324.555,6178.446)
	(6329.703,6123.576)
	(6334.608,6067.856)
	(6339.326,6011.722)
	(6343.913,5955.615)
	(6348.423,5899.971)
	(6352.915,5845.229)
	(6357.443,5791.828)
	(6362.063,5740.204)
	(6366.832,5690.796)
	(6371.805,5644.043)
	(6377.039,5600.382)
	(6382.589,5560.252)
	(6388.511,5524.091)
	(6401.697,5465.426)
	(6435.000,5415.000)

\path(6435,5415)	(6485.879,5477.259)
	(6496.702,5513.561)
	(6506.875,5559.035)
	(6516.455,5614.123)
	(6525.495,5679.263)
	(6529.830,5715.739)
	(6534.051,5754.894)
	(6538.165,5796.782)
	(6542.178,5841.457)
	(6546.098,5888.975)
	(6549.930,5939.391)
	(6553.683,5992.759)
	(6557.363,6049.134)
	(6560.977,6108.572)
	(6564.532,6171.128)
	(6568.034,6236.855)
	(6571.491,6305.810)
	(6574.909,6378.047)
	(6576.606,6415.414)
	(6578.295,6453.621)
	(6579.979,6492.677)
	(6581.657,6532.587)
	(6583.330,6573.359)
	(6585.000,6615.000)

\path(5700,6318)(5400,5718)(6600,5718)
	(6900,6318)(5700,6318)
\put(11535,4665){\ellipse{150}{150}}
\put(11535,4815){\ellipse{150}{150}}
\put(11535,4965){\ellipse{150}{150}}
\path(11535,4590)	(11534.401,4528.543)
	(11533.957,4469.649)
	(11533.672,4413.276)
	(11533.550,4359.383)
	(11533.593,4307.929)
	(11533.806,4258.873)
	(11534.191,4212.172)
	(11534.753,4167.788)
	(11535.494,4125.677)
	(11536.418,4085.799)
	(11538.829,4012.576)
	(11542.014,3947.791)
	(11545.999,3891.112)
	(11550.812,3842.212)
	(11556.481,3800.760)
	(11570.495,3738.881)
	(11610.000,3690.000)

\path(11610,3690)	(11648.407,3735.771)
	(11666.068,3790.672)
	(11682.998,3862.434)
	(11691.257,3903.588)
	(11699.417,3947.696)
	(11707.505,3994.339)
	(11715.550,4043.095)
	(11723.578,4093.546)
	(11731.617,4145.270)
	(11739.696,4197.847)
	(11747.842,4250.858)
	(11756.084,4303.881)
	(11764.448,4356.497)
	(11772.962,4408.285)
	(11781.655,4458.825)
	(11790.554,4507.697)
	(11799.687,4554.481)
	(11809.082,4598.756)
	(11818.766,4640.102)
	(11839.115,4712.326)
	(11860.955,4767.792)
	(11910.000,4815.000)

\path(11910,4815)	(11958.026,4761.706)
	(11976.996,4701.915)
	(11985.416,4665.267)
	(11993.201,4624.702)
	(12000.406,4580.657)
	(12007.088,4533.571)
	(12013.303,4483.881)
	(12019.107,4432.027)
	(12024.555,4378.446)
	(12029.703,4323.576)
	(12034.608,4267.856)
	(12039.326,4211.722)
	(12043.913,4155.615)
	(12048.423,4099.971)
	(12052.915,4045.229)
	(12057.443,3991.828)
	(12062.063,3940.204)
	(12066.832,3890.796)
	(12071.805,3844.043)
	(12077.039,3800.382)
	(12082.589,3760.252)
	(12088.511,3724.091)
	(12101.697,3665.426)
	(12135.000,3615.000)

\path(12135,3615)	(12185.879,3677.259)
	(12196.702,3713.561)
	(12206.875,3759.035)
	(12216.455,3814.123)
	(12225.495,3879.263)
	(12229.830,3915.739)
	(12234.051,3954.894)
	(12238.165,3996.782)
	(12242.178,4041.457)
	(12246.098,4088.975)
	(12249.930,4139.391)
	(12253.683,4192.759)
	(12257.363,4249.134)
	(12260.977,4308.572)
	(12264.532,4371.128)
	(12268.034,4436.855)
	(12271.491,4505.810)
	(12274.909,4578.047)
	(12276.606,4615.414)
	(12278.295,4653.621)
	(12279.979,4692.677)
	(12281.657,4732.587)
	(12283.330,4773.359)
	(12285.000,4815.000)

\path(11400,4518)(11100,3918)(12300,3918)
	(12600,4518)(11400,4518)
\put(6735,4665){\ellipse{150}{150}}
\put(6735,4815){\ellipse{150}{150}}
\put(6735,4965){\ellipse{150}{150}}
\path(6735,4590)	(6734.401,4528.543)
	(6733.957,4469.649)
	(6733.672,4413.276)
	(6733.550,4359.383)
	(6733.593,4307.929)
	(6733.806,4258.873)
	(6734.191,4212.172)
	(6734.753,4167.788)
	(6735.494,4125.677)
	(6736.418,4085.799)
	(6738.829,4012.576)
	(6742.014,3947.791)
	(6745.999,3891.112)
	(6750.812,3842.212)
	(6756.481,3800.760)
	(6770.495,3738.881)
	(6810.000,3690.000)

\path(6810,3690)	(6848.407,3735.771)
	(6866.068,3790.672)
	(6882.997,3862.434)
	(6891.257,3903.588)
	(6899.417,3947.696)
	(6907.505,3994.339)
	(6915.550,4043.095)
	(6923.578,4093.546)
	(6931.617,4145.270)
	(6939.696,4197.847)
	(6947.842,4250.858)
	(6956.084,4303.881)
	(6964.448,4356.497)
	(6972.962,4408.285)
	(6981.655,4458.825)
	(6990.554,4507.697)
	(6999.687,4554.481)
	(7009.082,4598.756)
	(7018.766,4640.102)
	(7039.115,4712.326)
	(7060.955,4767.792)
	(7110.000,4815.000)

\path(7110,4815)	(7158.026,4761.706)
	(7176.996,4701.915)
	(7185.416,4665.267)
	(7193.201,4624.702)
	(7200.406,4580.657)
	(7207.088,4533.571)
	(7213.303,4483.881)
	(7219.107,4432.027)
	(7224.555,4378.446)
	(7229.703,4323.576)
	(7234.608,4267.856)
	(7239.326,4211.722)
	(7243.913,4155.615)
	(7248.423,4099.971)
	(7252.915,4045.229)
	(7257.443,3991.828)
	(7262.063,3940.204)
	(7266.832,3890.796)
	(7271.805,3844.043)
	(7277.039,3800.382)
	(7282.589,3760.252)
	(7288.511,3724.091)
	(7301.697,3665.426)
	(7335.000,3615.000)

\path(7335,3615)	(7385.879,3677.259)
	(7396.702,3713.561)
	(7406.875,3759.035)
	(7416.455,3814.123)
	(7425.495,3879.263)
	(7429.830,3915.739)
	(7434.051,3954.894)
	(7438.165,3996.782)
	(7442.178,4041.457)
	(7446.098,4088.975)
	(7449.930,4139.391)
	(7453.683,4192.759)
	(7457.363,4249.134)
	(7460.977,4308.572)
	(7464.532,4371.128)
	(7468.034,4436.855)
	(7471.491,4505.810)
	(7474.909,4578.047)
	(7476.606,4615.414)
	(7478.295,4653.621)
	(7479.979,4692.677)
	(7481.657,4732.587)
	(7483.330,4773.359)
	(7485.000,4815.000)

\path(6600,4518)(6300,3918)(7500,3918)
	(7800,4518)(6600,4518)
\put(10035,2865){\ellipse{150}{150}}
\put(10035,3015){\ellipse{150}{150}}
\put(10035,3165){\ellipse{150}{150}}
\path(10035,2790)	(10034.401,2728.543)
	(10033.957,2669.649)
	(10033.672,2613.276)
	(10033.550,2559.383)
	(10033.593,2507.929)
	(10033.806,2458.873)
	(10034.191,2412.172)
	(10034.753,2367.787)
	(10035.494,2325.677)
	(10036.418,2285.799)
	(10038.829,2212.576)
	(10042.014,2147.791)
	(10045.999,2091.112)
	(10050.812,2042.212)
	(10056.481,2000.760)
	(10070.495,1938.881)
	(10110.000,1890.000)

\path(10110,1890)	(10148.407,1935.771)
	(10166.068,1990.672)
	(10182.998,2062.434)
	(10191.257,2103.588)
	(10199.417,2147.696)
	(10207.505,2194.339)
	(10215.550,2243.095)
	(10223.578,2293.546)
	(10231.617,2345.270)
	(10239.696,2397.847)
	(10247.842,2450.858)
	(10256.084,2503.881)
	(10264.448,2556.497)
	(10272.962,2608.285)
	(10281.655,2658.825)
	(10290.554,2707.697)
	(10299.687,2754.481)
	(10309.082,2798.756)
	(10318.766,2840.102)
	(10339.115,2912.326)
	(10360.955,2967.792)
	(10410.000,3015.000)

\path(10410,3015)	(10458.026,2961.706)
	(10476.996,2901.915)
	(10485.416,2865.267)
	(10493.201,2824.702)
	(10500.406,2780.657)
	(10507.088,2733.571)
	(10513.303,2683.881)
	(10519.107,2632.027)
	(10524.555,2578.446)
	(10529.703,2523.576)
	(10534.608,2467.856)
	(10539.326,2411.722)
	(10543.913,2355.615)
	(10548.423,2299.971)
	(10552.915,2245.229)
	(10557.443,2191.828)
	(10562.063,2140.204)
	(10566.832,2090.796)
	(10571.805,2044.043)
	(10577.039,2000.382)
	(10582.589,1960.252)
	(10588.511,1924.091)
	(10601.697,1865.426)
	(10635.000,1815.000)

\path(10635,1815)	(10685.879,1877.259)
	(10696.702,1913.561)
	(10706.875,1959.035)
	(10716.455,2014.123)
	(10725.495,2079.263)
	(10729.830,2115.739)
	(10734.051,2154.894)
	(10738.165,2196.782)
	(10742.178,2241.457)
	(10746.098,2288.975)
	(10749.930,2339.391)
	(10753.683,2392.759)
	(10757.363,2449.134)
	(10760.977,2508.572)
	(10764.532,2571.128)
	(10768.034,2636.855)
	(10771.491,2705.810)
	(10774.909,2778.047)
	(10776.606,2815.414)
	(10778.295,2853.621)
	(10779.979,2892.677)
	(10781.657,2932.587)
	(10783.330,2973.359)
	(10785.000,3015.000)

\path(9900,2718)(9600,2118)(10800,2118)
	(11100,2718)(9900,2718)
\put(5790,7845){\ellipse{150}{150}}
\put(5835,7695){\ellipse{150}{150}}
\put(5775,8002){\ellipse{150}{150}}
\path(5700,8115)(5400,7515)(6600,7515)
	(6900,8115)(5700,8115)
\path(5865,7635)	(5919.728,7613.340)
	(5961.835,7602.646)
	(6023.000,7612.000)

\path(6023,7612)	(6068.241,7652.889)
	(6100.001,7712.860)
	(6122.165,7784.296)
	(6130.863,7821.936)
	(6138.618,7859.586)
	(6153.243,7931.115)
	(6169.925,7991.270)
	(6225.000,8047.000)

\path(6225,8047)	(6276.596,7980.370)
	(6291.016,7913.203)
	(6296.919,7875.326)
	(6302.539,7836.174)
	(6308.313,7796.938)
	(6314.682,7758.810)
	(6330.962,7690.641)
	(6390.000,7620.000)

\path(6390,7620)	(6452.798,7636.287)
	(6512.246,7703.945)
	(6542.582,7760.056)
	(6574.321,7832.631)
	(6590.951,7875.469)
	(6608.212,7922.876)
	(6626.197,7975.003)
	(6645.000,8032.000)

\put(3390,6045){\ellipse{150}{150}}
\put(3435,5895){\ellipse{150}{150}}
\put(3375,6202){\ellipse{150}{150}}
\path(3300,6315)(3000,5715)(4200,5715)
	(4500,6315)(3300,6315)
\path(3465,5835)	(3519.728,5813.340)
	(3561.835,5802.646)
	(3623.000,5812.000)

\path(3623,5812)	(3668.241,5852.889)
	(3700.001,5912.860)
	(3722.165,5984.296)
	(3730.863,6021.936)
	(3738.617,6059.586)
	(3753.243,6131.115)
	(3769.925,6191.270)
	(3825.000,6247.000)

\path(3825,6247)	(3876.596,6180.370)
	(3891.016,6113.203)
	(3896.919,6075.326)
	(3902.539,6036.174)
	(3908.313,5996.938)
	(3914.682,5958.810)
	(3930.962,5890.641)
	(3990.000,5820.000)

\path(3990,5820)	(4052.798,5836.287)
	(4112.246,5903.945)
	(4142.582,5960.056)
	(4174.321,6032.631)
	(4190.951,6075.469)
	(4208.212,6122.876)
	(4226.197,6175.003)
	(4245.000,6232.000)

\put(12990,6045){\ellipse{150}{150}}
\put(13035,5895){\ellipse{150}{150}}
\put(12975,6202){\ellipse{150}{150}}
\path(12900,6315)(12600,5715)(13800,5715)
	(14100,6315)(12900,6315)
\path(13065,5835)	(13119.728,5813.340)
	(13161.835,5802.646)
	(13223.000,5812.000)

\path(13223,5812)	(13268.241,5852.889)
	(13300.001,5912.860)
	(13322.165,5984.296)
	(13330.863,6021.936)
	(13338.618,6059.586)
	(13353.243,6131.115)
	(13369.925,6191.270)
	(13425.000,6247.000)

\path(13425,6247)	(13476.596,6180.370)
	(13491.016,6113.203)
	(13496.919,6075.326)
	(13502.539,6036.174)
	(13508.313,5996.938)
	(13514.682,5958.810)
	(13530.962,5890.641)
	(13590.000,5820.000)

\path(13590,5820)	(13652.798,5836.287)
	(13712.246,5903.945)
	(13742.582,5960.056)
	(13774.321,6032.631)
	(13790.951,6075.469)
	(13808.212,6122.876)
	(13826.197,6175.003)
	(13845.000,6232.000)

\put(9090,4245){\ellipse{150}{150}}
\put(9135,4095){\ellipse{150}{150}}
\put(9075,4402){\ellipse{150}{150}}
\path(9000,4515)(8700,3915)(9900,3915)
	(10200,4515)(9000,4515)
\path(9165,4035)	(9219.728,4013.340)
	(9261.835,4002.646)
	(9323.000,4012.000)

\path(9323,4012)	(9368.241,4052.889)
	(9400.001,4112.860)
	(9422.165,4184.296)
	(9430.863,4221.936)
	(9438.618,4259.586)
	(9453.243,4331.115)
	(9469.925,4391.270)
	(9525.000,4447.000)

\path(9525,4447)	(9576.596,4380.370)
	(9591.016,4313.203)
	(9596.919,4275.326)
	(9602.539,4236.174)
	(9608.313,4196.938)
	(9614.682,4158.810)
	(9630.962,4090.641)
	(9690.000,4020.000)

\path(9690,4020)	(9752.798,4036.287)
	(9812.246,4103.945)
	(9842.582,4160.056)
	(9874.321,4232.631)
	(9890.951,4275.469)
	(9908.212,4322.876)
	(9926.197,4375.003)
	(9945.000,4432.000)

\put(1890,4245){\ellipse{150}{150}}
\put(1935,4095){\ellipse{150}{150}}
\put(1875,4402){\ellipse{150}{150}}
\path(1800,4515)(1500,3915)(2700,3915)
	(3000,4515)(1800,4515)
\path(1965,4035)	(2019.728,4013.340)
	(2061.835,4002.646)
	(2123.000,4012.000)

\path(2123,4012)	(2168.241,4052.889)
	(2200.001,4112.860)
	(2222.165,4184.296)
	(2230.863,4221.936)
	(2238.617,4259.586)
	(2253.243,4331.115)
	(2269.925,4391.270)
	(2325.000,4447.000)

\path(2325,4447)	(2376.596,4380.370)
	(2391.016,4313.203)
	(2396.919,4275.326)
	(2402.539,4236.174)
	(2408.313,4196.938)
	(2414.682,4158.810)
	(2430.962,4090.641)
	(2490.000,4020.000)

\path(2490,4020)	(2552.798,4036.287)
	(2612.246,4103.945)
	(2642.582,4160.056)
	(2674.321,4232.631)
	(2690.951,4275.469)
	(2708.212,4322.876)
	(2726.197,4375.003)
	(2745.000,4432.000)

\put(13890,4245){\ellipse{150}{150}}
\put(13935,4095){\ellipse{150}{150}}
\put(13875,4402){\ellipse{150}{150}}
\path(13800,4515)(13500,3915)(14700,3915)
	(15000,4515)(13800,4515)
\path(13965,4035)	(14019.728,4013.340)
	(14061.835,4002.646)
	(14123.000,4012.000)

\path(14123,4012)	(14168.241,4052.889)
	(14200.001,4112.860)
	(14222.165,4184.296)
	(14230.863,4221.936)
	(14238.618,4259.586)
	(14253.243,4331.115)
	(14269.925,4391.270)
	(14325.000,4447.000)

\path(14325,4447)	(14376.596,4380.370)
	(14391.016,4313.203)
	(14396.919,4275.326)
	(14402.539,4236.174)
	(14408.313,4196.938)
	(14414.682,4158.810)
	(14430.962,4090.641)
	(14490.000,4020.000)

\path(14490,4020)	(14552.798,4036.287)
	(14612.246,4103.945)
	(14642.582,4160.056)
	(14674.321,4232.631)
	(14690.951,4275.469)
	(14708.212,4322.876)
	(14726.197,4375.003)
	(14745.000,4432.000)

\put(5190,2445){\ellipse{150}{150}}
\put(5235,2295){\ellipse{150}{150}}
\put(5175,2602){\ellipse{150}{150}}
\path(5100,2715)(4800,2115)(6000,2115)
	(6300,2715)(5100,2715)
\path(5265,2235)	(5319.728,2213.340)
	(5361.835,2202.646)
	(5423.000,2212.000)

\path(5423,2212)	(5468.241,2252.889)
	(5500.001,2312.860)
	(5522.165,2384.296)
	(5530.863,2421.936)
	(5538.618,2459.586)
	(5553.243,2531.115)
	(5569.925,2591.270)
	(5625.000,2647.000)

\path(5625,2647)	(5676.596,2580.370)
	(5691.016,2513.203)
	(5696.919,2475.326)
	(5702.539,2436.174)
	(5708.313,2396.938)
	(5714.682,2358.810)
	(5730.962,2290.641)
	(5790.000,2220.000)

\path(5790,2220)	(5852.798,2236.287)
	(5912.246,2303.945)
	(5942.582,2360.056)
	(5974.321,2432.631)
	(5990.951,2475.469)
	(6008.212,2522.876)
	(6026.197,2575.003)
	(6045.000,2632.000)

\put(8490,645){\ellipse{150}{150}}
\put(8535,495){\ellipse{150}{150}}
\put(8475,802){\ellipse{150}{150}}
\path(8400,915)(8100,315)(9300,315)
	(9600,915)(8400,915)
\path(8565,435)	(8619.728,413.340)
	(8661.835,402.646)
	(8723.000,412.000)

\path(8723,412)	(8768.241,452.889)
	(8800.001,512.860)
	(8822.165,584.296)
	(8830.863,621.936)
	(8838.618,659.586)
	(8853.243,731.115)
	(8869.925,791.270)
	(8925.000,847.000)

\path(8925,847)	(8976.596,780.370)
	(8991.016,713.203)
	(8996.919,675.326)
	(9002.539,636.174)
	(9008.313,596.938)
	(9014.682,558.810)
	(9030.962,490.641)
	(9090.000,420.000)

\path(9090,420)	(9152.798,436.287)
	(9212.246,503.945)
	(9242.582,560.056)
	(9274.321,632.631)
	(9290.951,675.469)
	(9308.212,722.876)
	(9326.197,775.003)
	(9345.000,832.000)

\put(12390,2445){\ellipse{150}{150}}
\put(12435,2295){\ellipse{150}{150}}
\put(12375,2602){\ellipse{150}{150}}
\path(12300,2715)(12000,2115)(13200,2115)
	(13500,2715)(12300,2715)
\path(12465,2235)	(12519.728,2213.340)
	(12561.835,2202.646)
	(12623.000,2212.000)

\path(12623,2212)	(12668.241,2252.889)
	(12700.001,2312.860)
	(12722.165,2384.296)
	(12730.863,2421.936)
	(12738.618,2459.586)
	(12753.243,2531.115)
	(12769.925,2591.270)
	(12825.000,2647.000)

\path(12825,2647)	(12876.596,2580.370)
	(12891.016,2513.203)
	(12896.919,2475.326)
	(12902.539,2436.174)
	(12908.313,2396.938)
	(12914.682,2358.810)
	(12930.962,2290.641)
	(12990.000,2220.000)

\path(12990,2220)	(13052.798,2236.287)
	(13112.246,2303.945)
	(13142.582,2360.056)
	(13174.321,2432.631)
	(13190.951,2475.469)
	(13208.212,2522.876)
	(13226.197,2575.003)
	(13245.000,2632.000)

\put(4275,4965){\ellipse{150}{150}}
\path(4200,4515)(3900,3915)(5100,3915)
	(5400,4515)(4200,4515)
\path(4125,4215)(5175,4215)
\path(4275,4890)	(4278.148,4845.719)
	(4281.281,4802.361)
	(4284.401,4759.921)
	(4287.508,4718.390)
	(4290.604,4677.761)
	(4293.691,4638.027)
	(4296.769,4599.180)
	(4299.839,4561.214)
	(4305.964,4487.894)
	(4312.075,4418.008)
	(4318.184,4351.497)
	(4324.299,4288.303)
	(4330.432,4228.367)
	(4336.592,4171.632)
	(4342.791,4118.039)
	(4349.038,4067.529)
	(4355.344,4020.045)
	(4361.719,3975.527)
	(4368.173,3933.918)
	(4374.716,3895.159)
	(4388.113,3825.957)
	(4401.992,3767.456)
	(4416.436,3719.188)
	(4431.525,3680.685)
	(4500.000,3615.000)

\path(4500,3615)	(4551.454,3675.775)
	(4573.500,3746.126)
	(4583.749,3789.448)
	(4593.540,3837.500)
	(4602.913,3889.755)
	(4611.911,3945.685)
	(4620.576,4004.764)
	(4628.950,4066.465)
	(4637.074,4130.262)
	(4644.991,4195.628)
	(4652.743,4262.036)
	(4660.372,4328.959)
	(4667.921,4395.870)
	(4675.430,4462.244)
	(4682.942,4527.553)
	(4690.499,4591.270)
	(4698.143,4652.869)
	(4705.916,4711.823)
	(4713.861,4767.605)
	(4722.019,4819.688)
	(4730.432,4867.546)
	(4739.142,4910.652)
	(4757.622,4980.502)
	(4800.000,5040.000)

\path(4800,5040)	(4831.959,5021.060)
	(4859.266,4964.987)
	(4871.312,4921.722)
	(4882.359,4867.609)
	(4892.462,4802.125)
	(4897.178,4764.956)
	(4901.677,4724.749)
	(4905.969,4681.438)
	(4910.059,4634.959)
	(4913.954,4585.245)
	(4917.662,4532.233)
	(4921.188,4475.856)
	(4924.540,4416.049)
	(4927.726,4352.747)
	(4930.750,4285.886)
	(4933.622,4215.399)
	(4936.346,4141.221)
	(4937.656,4102.728)
	(4938.931,4063.288)
	(4940.173,4022.892)
	(4941.383,3981.533)
	(4942.562,3939.203)
	(4943.709,3895.892)
	(4944.827,3851.594)
	(4945.916,3806.300)
	(4946.977,3760.002)
	(4948.011,3712.691)
	(4949.018,3664.360)
	(4950.000,3615.000)

\put(7575,3165){\ellipse{150}{150}}
\path(7500,2715)(7200,2115)(8400,2115)
	(8700,2715)(7500,2715)
\path(7425,2415)(8475,2415)
\path(7575,3090)	(7578.148,3045.719)
	(7581.281,3002.361)
	(7584.401,2959.921)
	(7587.508,2918.390)
	(7590.604,2877.761)
	(7593.691,2838.027)
	(7596.769,2799.180)
	(7599.839,2761.214)
	(7605.964,2687.894)
	(7612.075,2618.008)
	(7618.184,2551.497)
	(7624.299,2488.303)
	(7630.432,2428.367)
	(7636.592,2371.632)
	(7642.791,2318.039)
	(7649.038,2267.529)
	(7655.344,2220.045)
	(7661.719,2175.527)
	(7668.173,2133.918)
	(7674.716,2095.159)
	(7688.113,2025.957)
	(7701.992,1967.456)
	(7716.436,1919.188)
	(7731.525,1880.685)
	(7800.000,1815.000)

\path(7800,1815)	(7851.454,1875.775)
	(7873.500,1946.126)
	(7883.749,1989.448)
	(7893.540,2037.500)
	(7902.913,2089.755)
	(7911.911,2145.685)
	(7920.576,2204.764)
	(7928.950,2266.465)
	(7937.074,2330.262)
	(7944.991,2395.628)
	(7952.743,2462.036)
	(7960.372,2528.959)
	(7967.921,2595.870)
	(7975.430,2662.244)
	(7982.942,2727.553)
	(7990.499,2791.270)
	(7998.143,2852.869)
	(8005.916,2911.823)
	(8013.861,2967.605)
	(8022.019,3019.688)
	(8030.432,3067.546)
	(8039.142,3110.652)
	(8057.622,3180.502)
	(8100.000,3240.000)

\path(8100,3240)	(8131.959,3221.060)
	(8159.266,3164.987)
	(8171.312,3121.722)
	(8182.359,3067.609)
	(8192.462,3002.125)
	(8197.178,2964.956)
	(8201.677,2924.749)
	(8205.969,2881.438)
	(8210.059,2834.959)
	(8213.954,2785.245)
	(8217.662,2732.233)
	(8221.188,2675.856)
	(8224.540,2616.049)
	(8227.726,2552.747)
	(8230.750,2485.886)
	(8233.622,2415.399)
	(8236.346,2341.221)
	(8237.656,2302.728)
	(8238.931,2263.288)
	(8240.173,2222.892)
	(8241.383,2181.533)
	(8242.562,2139.203)
	(8243.709,2095.892)
	(8244.827,2051.594)
	(8245.916,2006.300)
	(8246.977,1960.002)
	(8248.011,1912.691)
	(8249.018,1864.360)
	(8250.000,1815.000)

\put(10875,1365){\ellipse{150}{150}}
\path(10800,915)(10500,315)(11700,315)
	(12000,915)(10800,915)
\path(10725,615)(11775,615)
\path(10875,1290)	(10878.148,1245.719)
	(10881.281,1202.361)
	(10884.401,1159.921)
	(10887.508,1118.390)
	(10890.604,1077.761)
	(10893.691,1038.027)
	(10896.769,999.180)
	(10899.839,961.214)
	(10905.964,887.894)
	(10912.075,818.008)
	(10918.184,751.497)
	(10924.299,688.303)
	(10930.432,628.367)
	(10936.592,571.632)
	(10942.791,518.039)
	(10949.038,467.529)
	(10955.344,420.045)
	(10961.719,375.527)
	(10968.173,333.918)
	(10974.716,295.159)
	(10988.113,225.957)
	(11001.992,167.456)
	(11016.436,119.188)
	(11031.525,80.685)
	(11100.000,15.000)

\path(11100,15)	(11151.454,75.775)
	(11173.500,146.126)
	(11183.749,189.448)
	(11193.540,237.500)
	(11202.913,289.755)
	(11211.911,345.685)
	(11220.576,404.764)
	(11228.950,466.465)
	(11237.074,530.262)
	(11244.991,595.628)
	(11252.743,662.036)
	(11260.372,728.959)
	(11267.921,795.870)
	(11275.430,862.244)
	(11282.942,927.553)
	(11290.499,991.270)
	(11298.143,1052.869)
	(11305.916,1111.823)
	(11313.861,1167.605)
	(11322.019,1219.688)
	(11330.432,1267.546)
	(11339.142,1310.652)
	(11357.622,1380.502)
	(11400.000,1440.000)

\path(11400,1440)	(11431.959,1421.060)
	(11459.266,1364.987)
	(11471.312,1321.722)
	(11482.359,1267.609)
	(11492.462,1202.125)
	(11497.178,1164.956)
	(11501.677,1124.749)
	(11505.969,1081.438)
	(11510.059,1034.959)
	(11513.954,985.245)
	(11517.662,932.233)
	(11521.188,875.856)
	(11524.540,816.049)
	(11527.726,752.747)
	(11530.750,685.886)
	(11533.622,615.399)
	(11536.346,541.221)
	(11537.656,502.728)
	(11538.931,463.288)
	(11540.173,422.892)
	(11541.383,381.533)
	(11542.562,339.203)
	(11543.709,295.892)
	(11544.827,251.594)
	(11545.916,206.300)
	(11546.977,160.002)
	(11548.011,112.691)
	(11549.018,64.360)
	(11550.000,15.000)

\put(7335,10065){\ellipse{150}{150}}
\put(7335,10215){\ellipse{150}{150}}
\put(7335,10365){\ellipse{150}{150}}
\path(7335,9990)	(7334.401,9928.543)
	(7333.957,9869.649)
	(7333.672,9813.276)
	(7333.550,9759.383)
	(7333.593,9707.929)
	(7333.806,9658.873)
	(7334.191,9612.172)
	(7334.753,9567.787)
	(7335.494,9525.677)
	(7336.418,9485.799)
	(7338.829,9412.576)
	(7342.014,9347.791)
	(7345.999,9291.112)
	(7350.812,9242.212)
	(7356.481,9200.760)
	(7370.495,9138.881)
	(7410.000,9090.000)

\path(7410,9090)	(7448.407,9135.771)
	(7466.068,9190.672)
	(7482.997,9262.434)
	(7491.257,9303.588)
	(7499.417,9347.696)
	(7507.505,9394.339)
	(7515.550,9443.095)
	(7523.578,9493.546)
	(7531.617,9545.270)
	(7539.696,9597.847)
	(7547.842,9650.858)
	(7556.084,9703.881)
	(7564.448,9756.497)
	(7572.962,9808.285)
	(7581.655,9858.825)
	(7590.554,9907.697)
	(7599.687,9954.481)
	(7609.082,9998.756)
	(7618.766,10040.102)
	(7639.115,10112.326)
	(7660.955,10167.792)
	(7710.000,10215.000)

\path(7710,10215)	(7758.026,10161.706)
	(7776.996,10101.915)
	(7785.416,10065.267)
	(7793.201,10024.702)
	(7800.406,9980.657)
	(7807.088,9933.571)
	(7813.303,9883.881)
	(7819.107,9832.027)
	(7824.555,9778.446)
	(7829.703,9723.576)
	(7834.608,9667.856)
	(7839.326,9611.722)
	(7843.913,9555.615)
	(7848.423,9499.971)
	(7852.915,9445.229)
	(7857.443,9391.828)
	(7862.063,9340.204)
	(7866.832,9290.796)
	(7871.805,9244.043)
	(7877.039,9200.382)
	(7882.589,9160.252)
	(7888.511,9124.091)
	(7901.697,9065.426)
	(7935.000,9015.000)

\path(7935,9015)	(7985.879,9077.259)
	(7996.702,9113.561)
	(8006.875,9159.035)
	(8016.455,9214.123)
	(8025.495,9279.263)
	(8029.830,9315.739)
	(8034.051,9354.894)
	(8038.165,9396.782)
	(8042.178,9441.457)
	(8046.098,9488.975)
	(8049.930,9539.391)
	(8053.683,9592.759)
	(8057.363,9649.134)
	(8060.977,9708.572)
	(8064.532,9771.128)
	(8068.034,9836.855)
	(8071.491,9905.810)
	(8074.909,9978.047)
	(8076.606,10015.414)
	(8078.295,10053.621)
	(8079.979,10092.677)
	(8081.657,10132.587)
	(8083.330,10173.359)
	(8085.000,10215.000)

\path(7200,9918)(6900,9318)(8100,9318)
	(8400,9918)(7200,9918)
\put(9728,9600){\blacken\ellipse{74}{74}}
\put(9728,9600){\ellipse{74}{74}}
\put(10638,7816){\blacken\ellipse{74}{74}}
\put(10638,7816){\ellipse{74}{74}}
\put(8220,6007){\blacken\ellipse{74}{74}}
\put(8220,6007){\ellipse{74}{74}}
\put(10635,6067){\blacken\ellipse{74}{74}}
\put(10635,6067){\ellipse{74}{74}}
\put(10845,5955){\blacken\ellipse{74}{74}}
\put(10845,5955){\ellipse{74}{74}}
\put(13778,6037){\blacken\ellipse{74}{74}}
\put(13778,6037){\ellipse{74}{74}}
\put(14678,4282){\blacken\ellipse{74}{74}}
\put(14678,4282){\ellipse{74}{74}}
\put(14543,4050){\blacken\ellipse{74}{74}}
\put(14543,4050){\ellipse{74}{74}}
\put(11543,4312){\blacken\ellipse{74}{74}}
\put(11543,4312){\ellipse{74}{74}}
\put(11730,4125){\blacken\ellipse{74}{74}}
\put(11730,4125){\ellipse{74}{74}}
\put(9359,4060){\blacken\ellipse{74}{74}}
\put(9359,4060){\ellipse{74}{74}}
\put(6728,4185){\blacken\ellipse{74}{74}}
\put(6728,4185){\ellipse{74}{74}}
\put(5460,2257){\blacken\ellipse{74}{74}}
\put(5460,2257){\ellipse{74}{74}}
\put(8100,2422){\blacken\ellipse{74}{74}}
\put(8100,2422){\ellipse{74}{74}}
\put(10020,2512){\blacken\ellipse{74}{74}}
\put(10020,2512){\ellipse{74}{74}}
\put(10230,2332){\blacken\ellipse{74}{74}}
\put(10230,2332){\ellipse{74}{74}}
\put(12638,2250){\blacken\ellipse{74}{74}}
\put(12638,2250){\ellipse{74}{74}}
\put(12728,2392){\blacken\ellipse{74}{74}}
\put(12728,2392){\ellipse{74}{74}}
\put(8753,435){\blacken\ellipse{74}{74}}
\put(8753,435){\ellipse{74}{74}}
\put(8839,657){\blacken\ellipse{74}{74}}
\put(8839,657){\ellipse{74}{74}}
\put(11408,622){\blacken\ellipse{74}{74}}
\put(11408,622){\ellipse{74}{74}}
\put(11663,622){\blacken\ellipse{74}{74}}
\put(11663,622){\ellipse{74}{74}}
\dottedline{135}(8955,11610)(8603,11235)
\dottedline{135}(7425,9870)(7095,9420)
\dottedline{135}(7433,9630)(7778,9442)
\dottedline{135}(8333,8040)(7995,7687)
\dottedline{135}(8303,7815)(8670,7620)
\dottedline{135}(8588,7882)(8895,8047)
\dottedline{135}(10740,8025)(10995,7650)
\dottedline{135}(6360,7965)(6750,7972)
\dottedline{135}(6368,7860)(6705,7867)
\dottedline{135}(3968,6165)(4358,6172)
\dottedline{135}(3975,6052)(4305,6060)
\dottedline{135}(3578,6142)(3833,5805)
\dottedline{135}(5955,6255)(5655,6022)
\dottedline{135}(5918,6022)(6248,5797)
\dottedline{135}(6180,6052)(6465,6262)
\dottedline{135}(6443,5977)(6743,6225)
\dottedline{135}(8385,6232)(8573,5850)
\dottedline{135}(8573,6120)(8903,6240)
\dottedline{135}(13455,5835)(13185,6187)
\dottedline{135}(12158,4455)(11865,4012)
\dottedline{135}(10050,4387)(9668,4380)
\dottedline{135}(9675,4275)(9990,4200)
\dottedline{135}(2880,4425)(2483,4350)
\dottedline{135}(2505,4297)(2753,4162)
\dottedline{135}(2153,4425)(2378,4012)
\dottedline{135}(2250,4005)(2063,4245)
\dottedline{135}(7343,4425)(7673,4432)
\dottedline{135}(7350,4305)(7058,4020)
\dottedline{135}(6810,4432)(7140,4297)
\dottedline{135}(5813,2595)(6165,2602)
\dottedline{135}(5813,2475)(5993,2250)
\dottedline{135}(5633,2377)(5303,2617)
\dottedline{135}(10635,2617)(10380,2445)
\dottedline{135}(10635,2347)(10950,2632)
\dottedline{135}(13020,2557)(13350,2535)
\dottedline{135}(8925,495)(9105,817)
\dottedline{135}(9128,615)(9443,757)
\path(9150,12840)(9150,12165)
\path(9150,12840)(9150,12165)
\path(9120.000,12285.000)(9150.000,12165.000)(9180.000,12285.000)
\path(8775,10740)(8325,10065)
\path(8775,10740)(8325,10065)
\path(8366.603,10181.487)(8325.000,10065.000)(8416.526,10148.205)
\path(9000,10740)(9450,9915)
\path(9000,10740)(9450,9915)
\path(9366.201,10005.982)(9450.000,9915.000)(9418.875,10034.713)
\path(7515,9015)(6900,8250)
\path(7515,9015)(6900,8250)
\path(6951.806,8362.322)(6900.000,8250.000)(6998.568,8324.728)
\path(7650,9015)(7860,8040)
\path(7650,9015)(7860,8040)
\path(7805.406,8150.993)(7860.000,8040.000)(7864.061,8163.626)
\path(8265,9180)(10050,7920)
\path(8265,9180)(10050,7920)
\path(9934.663,7964.693)(10050.000,7920.000)(9969.264,8013.711)
\path(10575,9240)(11085,8640)
\path(10575,9240)(11085,8640)
\path(10984.424,8712.003)(11085.000,8640.000)(11030.140,8750.862)
\path(7335,7335)(4230,6600)
\path(7335,7335)(4230,6600)
\path(4339.862,6656.835)(4230.000,6600.000)(4353.683,6598.449)
\path(7530,7110)(6720,6510)
\path(7530,7110)(6720,6510)
\path(6798.570,6605.534)(6720.000,6510.000)(6834.284,6557.321)
\path(7875,7110)(7845,6165)
\path(7875,7110)(7845,6165)
\path(7818.823,6285.891)(7845.000,6165.000)(7878.792,6283.988)
\path(9885,7110)(9225,6555)
\path(9885,7110)(9225,6555)
\path(9297.535,6655.193)(9225.000,6555.000)(9336.151,6609.271)
\path(10140,7140)(10425,6480)
\path(10140,7140)(10425,6480)
\path(10349.886,6578.274)(10425.000,6480.000)(10404.970,6602.061)
\path(11235,7155)(12810,6390)
\path(11235,7155)(12810,6390)
\path(12688.952,6415.443)(12810.000,6390.000)(12715.166,6469.414)
\path(12120,3510)(12495,2850)
\path(12120,3510)(12495,2850)
\path(12409.635,2939.515)(12495.000,2850.000)(12461.803,2969.155)
\path(11745,3555)(11070,2895)
\path(11745,3555)(11070,2895)
\path(11134.827,3000.344)(11070.000,2895.000)(11176.774,2957.444)
\path(7665,3885)(9765,2805)
\path(7665,3885)(9765,2805)
\path(9644.565,2833.203)(9765.000,2805.000)(9672.006,2886.560)
\path(6780,3630)(7200,2685)
\path(6780,3630)(7200,2685)
\path(7123.849,2782.473)(7200.000,2685.000)(7178.678,2806.842)
\path(6450,3540)(6195,2910)
\path(6450,3540)(6195,2910)
\path(6212.215,3032.489)(6195.000,2910.000)(6267.832,3009.978)
\path(9855,1785)(9435,1080)
\path(9855,1785)(9435,1080)
\path(9470.644,1198.446)(9435.000,1080.000)(9522.190,1167.738)
\path(10245,1740)(10560,900)
\path(10245,1740)(10560,900)
\path(10489.775,1001.826)(10560.000,900.000)(10545.955,1022.893)
\path(11520,5685)(13650,4515)
\path(11520,5685)(13650,4515)
\path(13530.379,4546.479)(13650.000,4515.000)(13559.266,4599.068)
\path(11085,5490)(11295,4725)
\path(11085,5490)(11295,4725)
\path(11234.304,4832.778)(11295.000,4725.000)(11292.164,4848.661)
\path(9105,5490)(11115,4290)
\path(9105,5490)(11115,4290)
\path(10996.587,4325.755)(11115.000,4290.000)(11027.344,4377.272)
\path(8505,5430)(8790,4515)
\path(8505,5430)(8790,4515)
\path(8725.671,4620.649)(8790.000,4515.000)(8782.957,4638.492)
\path(8235,5430)(7830,4680)
\path(8235,5430)(7830,4680)
\path(7860.621,4799.843)(7830.000,4680.000)(7913.415,4771.334)
\path(5880,5400)(6345,4410)
\path(5880,5400)(6345,4410)
\path(6266.830,4505.861)(6345.000,4410.000)(6321.138,4531.370)
\path(5550,5340)(5100,4590)
\path(5550,5340)(5100,4590)
\path(5136.015,4708.334)(5100.000,4590.000)(5187.464,4677.464)
\path(4950,5490)(3075,4665)
\path(4950,5490)(3075,4665)
\path(3172.756,4740.788)(3075.000,4665.000)(3196.920,4685.869)
\put(300,13140){\makebox(0,0)[lb]{\smash{{{\SetFigFont{12}{14.4}{rm}17}}}}}
\put(300,9540){\makebox(0,0)[lb]{\smash{{{\SetFigFont{12}{14.4}{rm}15}}}}}
\put(300,7740){\makebox(0,0)[lb]{\smash{{{\SetFigFont{12}{14.4}{rm}14}}}}}
\put(300,5865){\makebox(0,0)[lb]{\smash{{{\SetFigFont{12}{14.4}{rm}13}}}}}
\put(300,4140){\makebox(0,0)[lb]{\smash{{{\SetFigFont{12}{14.4}{rm}12}}}}}
\put(0,14115){\makebox(0,0)[lb]{\smash{{{\SetFigFont{12}{14.4}{rm}\# general}}}}}
\put(225,13740){\makebox(0,0)[lb]{\smash{{{\SetFigFont{12}{14.4}{rm}lines}}}}}
\put(300,2340){\makebox(0,0)[lb]{\smash{{{\SetFigFont{12}{14.4}{rm}11}}}}}
\put(300,540){\makebox(0,0)[lb]{\smash{{{\SetFigFont{12}{14.4}{rm}10}}}}}
\put(5100,7215){\makebox(0,0)[lb]{\smash{{{\SetFigFont{8}{9.6}{rm}1}}}}}
\put(7500,7215){\makebox(0,0)[lb]{\smash{{{\SetFigFont{8}{9.6}{rm}9706}}}}}
\put(2700,5415){\makebox(0,0)[lb]{\smash{{{\SetFigFont{8}{9.6}{rm}36}}}}}
\put(7500,5415){\makebox(0,0)[lb]{\smash{{{\SetFigFont{8}{9.6}{rm}622}}}}}
\put(13200,3615){\makebox(0,0)[lb]{\smash{{{\SetFigFont{8}{9.6}{rm}1}}}}}
\put(10800,3690){\makebox(0,0)[lb]{\smash{{{\SetFigFont{8}{9.6}{rm}28}}}}}
\put(3600,3615){\makebox(0,0)[lb]{\smash{{{\SetFigFont{8}{9.6}{rm}5340}}}}}
\put(6900,1815){\makebox(0,0)[lb]{\smash{{{\SetFigFont{8}{9.6}{rm}430}}}}}
\put(9300,1815){\makebox(0,0)[lb]{\smash{{{\SetFigFont{8}{9.6}{rm}24}}}}}
\put(11700,1815){\makebox(0,0)[lb]{\smash{{{\SetFigFont{8}{9.6}{rm}1}}}}}
\put(10200,15){\makebox(0,0)[lb]{\smash{{{\SetFigFont{8}{9.6}{rm}20}}}}}
\put(7800,15){\makebox(0,0)[lb]{\smash{{{\SetFigFont{8}{9.6}{rm}1}}}}}
\put(9450,8415){\makebox(0,0)[lb]{\smash{{{\SetFigFont{8}{9.6}{rm}$\times 2$}}}}}
\put(4050,6765){\makebox(0,0)[lb]{\smash{{{\SetFigFont{8}{9.6}{rm}$\times 4$}}}}}
\put(7425,6390){\makebox(0,0)[lb]{\smash{{{\SetFigFont{8}{9.6}{rm}$\times 3$}}}}}
\put(12300,6765){\makebox(0,0)[lb]{\smash{{{\SetFigFont{8}{9.6}{rm}$\times 4$}}}}}
\put(12855,5445){\makebox(0,0)[lb]{\smash{{{\SetFigFont{8}{9.6}{rm}4}}}}}
\put(3000,4965){\makebox(0,0)[lb]{\smash{{{\SetFigFont{8}{9.6}{rm}$\times 4$}}}}}
\put(5700,4515){\makebox(0,0)[lb]{\smash{{{\SetFigFont{8}{9.6}{rm}$\times 4$}}}}}
\put(8775,4815){\makebox(0,0)[lb]{\smash{{{\SetFigFont{8}{9.6}{rm}$\times 4$}}}}}
\put(13275,4815){\makebox(0,0)[lb]{\smash{{{\SetFigFont{8}{9.6}{rm}$\times 4$}}}}}
\put(12450,3090){\makebox(0,0)[lb]{\smash{{{\SetFigFont{8}{9.6}{rm}$\times 4$}}}}}
\put(5700,3090){\makebox(0,0)[lb]{\smash{{{\SetFigFont{8}{9.6}{rm}$\times 4$}}}}}
\put(9150,1290){\makebox(0,0)[lb]{\smash{{{\SetFigFont{8}{9.6}{rm}$\times 4$}}}}}
\put(8025,12615){\makebox(0,0)[lb]{\smash{{{\SetFigFont{8}{9.6}{rm}11780}}}}}
\put(8025,10815){\makebox(0,0)[lb]{\smash{{{\SetFigFont{8}{9.6}{rm}11780}}}}}
\put(6525,9015){\makebox(0,0)[lb]{\smash{{{\SetFigFont{8}{9.6}{rm}11110}}}}}
\put(9075,9015){\makebox(0,0)[lb]{\smash{{{\SetFigFont{8}{9.6}{rm}670}}}}}
\put(6450,8415){\makebox(0,0)[lb]{\smash{{{\SetFigFont{8}{9.6}{rm}$\times 64$}}}}}
\put(9975,7215){\makebox(0,0)[lb]{\smash{{{\SetFigFont{8}{9.6}{rm}670}}}}}
\put(10050,5415){\makebox(0,0)[lb]{\smash{{{\SetFigFont{8}{9.6}{rm}32}}}}}
\put(5100,5415){\makebox(0,0)[lb]{\smash{{{\SetFigFont{8}{9.6}{rm}7696}}}}}
\put(10350,4815){\makebox(0,0)[lb]{\smash{{{\SetFigFont{8}{9.6}{rm}$\times 2$}}}}}
\put(2025,3615){\makebox(0,0)[lb]{\smash{{{\SetFigFont{8}{9.6}{rm}51}}}}}
\put(6075,3615){\makebox(0,0)[lb]{\smash{{{\SetFigFont{8}{9.6}{rm}538}}}}}
\put(9150,3615){\makebox(0,0)[lb]{\smash{{{\SetFigFont{8}{9.6}{rm}7}}}}}
\put(9225,3090){\makebox(0,0)[lb]{\smash{{{\SetFigFont{8}{9.6}{rm}$\times 3$}}}}}
\put(5250,1815){\makebox(0,0)[lb]{\smash{{{\SetFigFont{8}{9.6}{rm}9}}}}}
\put(300,11340){\makebox(0,0)[lb]{\smash{{{\SetFigFont{12}{14.4}{rm}16}}}}}
\end{picture}